\documentclass[11pt]{article}

\usepackage[pdftex]{graphicx}
\usepackage{amsmath}
\usepackage{amssymb}
\usepackage{caption}
\usepackage[subrefformat=parens]{subcaption}
\usepackage{here}
\newtheorem{theorem}{Theorem}

\usepackage{version}
\excludeversion{en-text}
\usepackage{comment}
\includecomment{jp-text}
\excludecomment{en-text}
\def\be{\begin{equation}}
\def\ee{\end{equation}}
\def\bea{\begin{eqnarray}}
\def\eea{\end{eqnarray}}
\def\beas{\begin{eqnarray*}}
\def\eeas{\end{eqnarray*}}

\def\dlnt{\bar{\delta}}

\def\be{\begin{equation}}
\def\ee{\end{equation}}
\def\bea{\begin{eqnarray}}
\def\eea{\end{eqnarray}}
\def\beas{\begin{eqnarray*}}
\def\eeas{\end{eqnarray*}}

\newif\ifcol
\colfalse

\newcommand{\argsup}{\mathop{\rm arg~sup}\limits}
\newcommand{\arginf}{\mathop{\rm arg~inf}\limits}

\begin{document}
\title{Adaptive estimator for a parabolic linear SPDE 
with a small noise
}
\author{
${}^{1}${Yusuke Kaino}
 and
 ${}^{1,2,3}${Masayuki Uchida}
 \\
        ${}^1${\small {Graduate School of Engineering Science, Osaka University}}\\
        ${}^2${\small {Center for Mathematical Modeling and Data Science (MMDS),  
         Osaka University} and }\\
          ${}^3${\small {JST CREST}}, 
        {\small {Toyonaka, Osaka 560-8531, Japan}}
}
\maketitle
\noindent
{\bf Abstract.}
{
We deal with parametric estimation for a parabolic linear second order stochastic partial differential equation (SPDE) 
with a small dispersion parameter based on high frequency data which are  observed in time and space. 
By using the thinned data with respect to space obtained from the high frequency data,
the minimum contrast estimators of two coefficient parameters of the SPDE  are proposed. 
With these estimators and the thinned data with respect to time obtained from the high frequency data, 
we construct an approximation of the coordinate process of the SPDE.
Using the approximate coordinate process, 
we {obtain} the adaptive estimator of 
{a coefficient parameter}
of the SPDE.
Moreover, we give simulation results of the proposed estimators of the SPDE.
}

\begin{en-text} 

We consider parametric estimation for a parabolic linear second order stochastic partial differential equation (SPDE) with small perturbations
from high frequency data which are  observed in time and space. 
By using thinned data for space obtained from the high frequency data,
estimators of the part of drift parameters of a parabolic linear SPDE model are proposed. 
With these estimators and thinned data for time, we construct an approximate of coordinate process.
Using this approximate of coordinate process, we calculate an estimator of the rest of drift parameters of a parabolic linear SPDE model.
Moreover, we give simulation results of the adaptive estimators of the SPDE model based on the high frequency data.

We consider parameter estimation of a parabolic second-order linear stochastic partial differential equation model using discrete observations.
There are four parameters to be estimated, and three of them are used to construct normalized volatility parameter and curvature parameter.
We propose two methods to estimate normalized volatility parameter and curvature parameter.
The first method uses data from time 0 to time 1, and the second method uses data from time 0 to time $T$.
A coordinate process is generated from discrete observations and the estimator of the curvature parameter, and we estimate its drift and diffusion coefficient.
The parameters of the original model are estimated from this estimator and the estimator of the normalized volatility parameter.
We verify the asymptotic behavior of the estimators by simulation.
\end{en-text}

\vspace{0.5cm} 


%




\section{Introduction}

We treat parametric estimation of a linear parabolic stochastic partial differential equation (SPDE) with 
one space dimension and a small dispersion parameter $\epsilon \in (0,1]$.
\begin{eqnarray}
& &dX_{t}(y) = \left( \theta_2 \frac{\partial^2 X_{t}(y)}{\partial y^2} + \theta_1 \frac{\partial X_{t}(y)}{\partial y} + \theta_0 X_{t}(y) \right)dt 
+ \epsilon dB_t(y), \quad  (t,y) \in [0, T]\times[0, 1],  \quad
\label{spde0} \\
& & X_t(0) = X_t(1) = 0, \quad t \in [0,T], \qquad X_0(y) = \xi, \quad y \in [0,1],
\nonumber
\end{eqnarray}
where 
{$\epsilon$ 
is known},
$T>0$, 
$B_t$ is defined as a cylindrical Brownian motion in the Sobolev space on $[0,1]$,
an unknown parameter $\theta= (\theta_0, \theta_1, \theta_2)$  
and $\theta_0,\theta_1 \in \mathbb{R}$, $\theta_2 >0$,
and the parameter space $\Theta$ is a compact convex subset of $\mathbb{R}^2 \times (0, \infty)$.
Moreover, the true value of parameter $\theta^*= (\theta_0^*, \theta_1^*, \theta_2^*)$
and we assume that $\theta^* \in \mbox{Int}(\Theta)$.
The data are discrete observations ${\bf X}_{N,M}= \left\{ X_{t_{i:N}}(y_{j:M}) \right\}_{i = 1, ..., N, j = 1, ..., M}$
with $t_{i:N} = i \frac{T}{N}$ and 
${ y_{j:M} = \frac{j}{M}}$.

Statistical inference for SPDE models based on high frequency data has been developed
by some researchers,
see, for example, 
Cont (2005),
Cialenco and Huang (2020), 
Bibinger and Trabs (2020)
and
Hildebrandt (2019).
Recently, Kaino and Uchida {(2020)} studied adaptive maximum likelihood (ML) type estimators 
of the coefficient parameters of the parabolic linear second order SPDE model.
Hildebrandt and Trabs (2019) proposed a contrast function with double increments 
for the parabolic linear second order SPDE model. 
{They obtained the minimum contrast estimators of the coefficient parameters of the SPDE model}
and showed that the estimators have asymptotic normality.

In this paper, we propose adaptive maximum likelihood (ML) type estimator of the coefficient parameter
$\theta= (\theta_0, \theta_1, \theta_2)$ 
of the parabolic linear second order SPDE  (\ref{spde0}) with a small dispersion parameter $\epsilon$.
For $k \in \mathbb{N}$, 
the coordinate process  $x_k(t)$ of the SPDE model (\ref{spde0}) is 
\begin{eqnarray}
x_k(t) 
&=& \int^1_0 X_t(y) \sqrt{2}\sin (\pi k y) \exp \left( \frac{\eta y}{2} \right) dy, \label{coordinate}
\end{eqnarray}
where $\eta = \frac{\theta_1}{\theta_2}$.
Note that the coordinate process (\ref{coordinate}) 
is  the Ornstein-Uhlenbeck process with a small dispersion parameter $\epsilon$ as follows.
{
$$
dx_k(t) = -\lambda_k x_k(t) dt + \epsilon dw_k(t), 
$$
where for $k \in \mathbb{N}$, 
$\lambda_k = - \theta_0 + \frac{\theta_1^2}{4 \theta_2} + \pi^2 k^2 \theta_2$
and 
$(w_k(t))_{t \geq 0}$ is independent real-valued Brownian motions. 
Moreover, the initial value $x_k(0) $ is defined in (\ref{coordinate2}) below. 
For details of 
{the} 
coordinate process, see Bibinger and Trabs (2020).
In fact, $\lambda_k$ is a very important parameter.
Figure $1$ below is the sample paths
{with the initial condition $\xi(y)=4y(1-y)$},
where
$\theta_1$, $\theta_2$ and $\epsilon$ are fixed and only $\theta_0$ is changed.
The rough shape of the sample path depends on the value of $\lambda_1$.
For the case that $\lambda_1$ is positive, when $y$ is fixed 
and $t$ tends to $1$, the value of $X_t(y)$ approaches $0$.
In case that $\lambda_1$ is close to 0, 
for $y$ being fixed and any $t \in [0,1]$, the value of $X_t(y)$ does not change.
When $\lambda_1$ is negative, the value of $X_t(y)$ increases.
Figure $2$ shows the sample path viewed from the $t$-axis side.
Figure $3$ is 
a cross section of 
the sample path at $y=0.5$, which means the sample path $X_t(0.5)$.

\begin{figure}[H]
\begin{minipage}{0.32\hsize}
\begin{center}
\includegraphics[width=4.5cm]{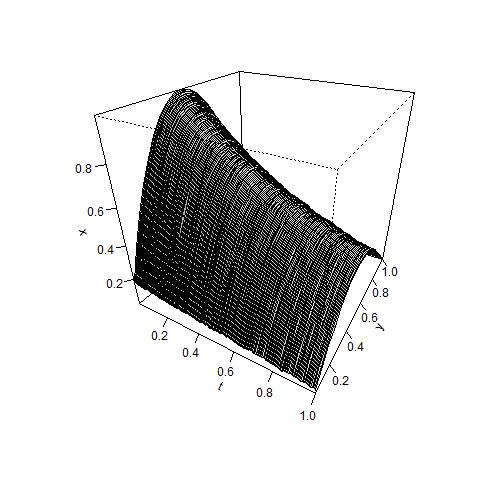}
\captionsetup{labelformat=empty,labelsep=none}
\subcaption{$\theta=(1,0.2,0.2)$, $\lambda_1 \approx 1$, \\ \hspace{0.5cm} $\epsilon=0.01$}
\end{center}
\end{minipage}
\begin{minipage}{0.32\hsize}
\begin{center}
\includegraphics[width=4.5cm]{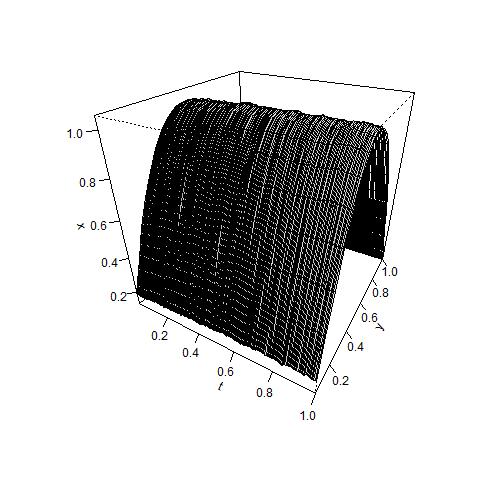}
\captionsetup{labelformat=empty,labelsep=none}
\subcaption{$\theta=(2,0.2,0.2)$, $\lambda_1\approx 0.02$, \\ \hspace{0.5cm}  $\epsilon=0.01$}
\end{center}
\end{minipage}
\begin{minipage}{0.32\hsize}
\begin{center}
\includegraphics[width=4.5cm]{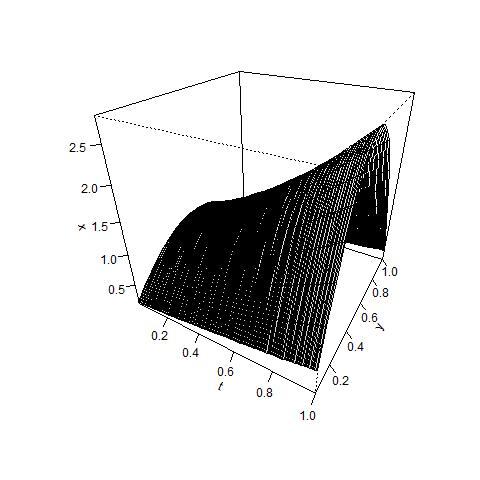}
\captionsetup{labelformat=empty,labelsep=none}
\subcaption{$\theta=(3,0.2,0.2)$, $\lambda_1\approx -1$, \\  \hspace{0.5cm} $\epsilon=0.01$}
\end{center}
\end{minipage}
\caption{Sample paths with $\lambda_1 \approx 1, 0.02, -1$}
\label{t1-11}
\end{figure}

\begin{figure}[H]
\begin{minipage}{0.32\hsize}
\begin{center}
\includegraphics[width=4.5cm]{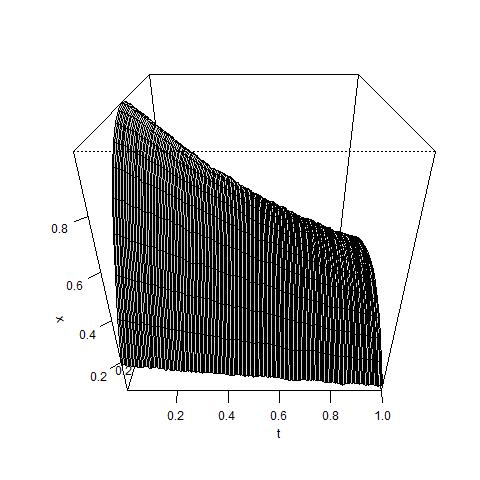}
\captionsetup{labelformat=empty,labelsep=none}
\subcaption{$\theta=(1,0.2,0.2)$, $\lambda_1 \approx 1$, \\ \hspace{0.5cm} $\epsilon=0.01$}
\end{center}
\end{minipage}
\begin{minipage}{0.32\hsize}
\begin{center}
\includegraphics[width=4.5cm]{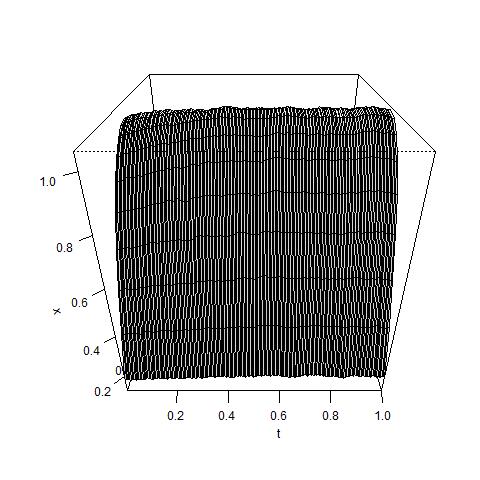}
\captionsetup{labelformat=empty,labelsep=none}
\subcaption{$\theta=(2,0.2,0.2)$, $\lambda_1\approx 0.02$, \\ \hspace{0.5cm}  $\epsilon=0.01$}
\end{center}
\end{minipage}
\begin{minipage}{0.32\hsize}
\begin{center}
\includegraphics[width=4.5cm]{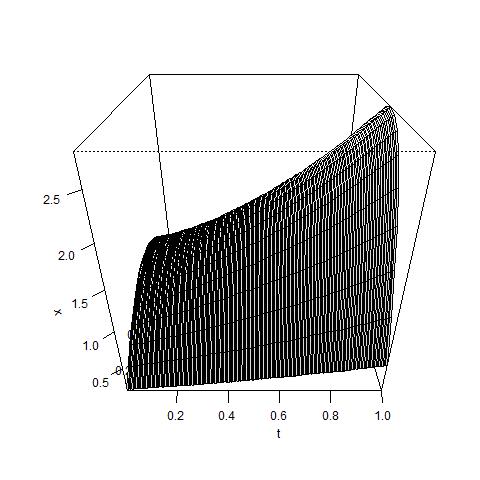}
\captionsetup{labelformat=empty,labelsep=none}
\subcaption{$\theta=(3,0.2,0.2)$, $\lambda_1\approx -1$, \\  \hspace{0.5cm} $\epsilon=0.01$}
\end{center}
\end{minipage}
\caption{Sample paths with $\lambda_1 \approx 1, 0.02, -1$ (t-axis side)}
\end{figure}

\begin{figure}[H]
\begin{minipage}{0.32\hsize}
\begin{center}
\includegraphics[width=4.5cm]{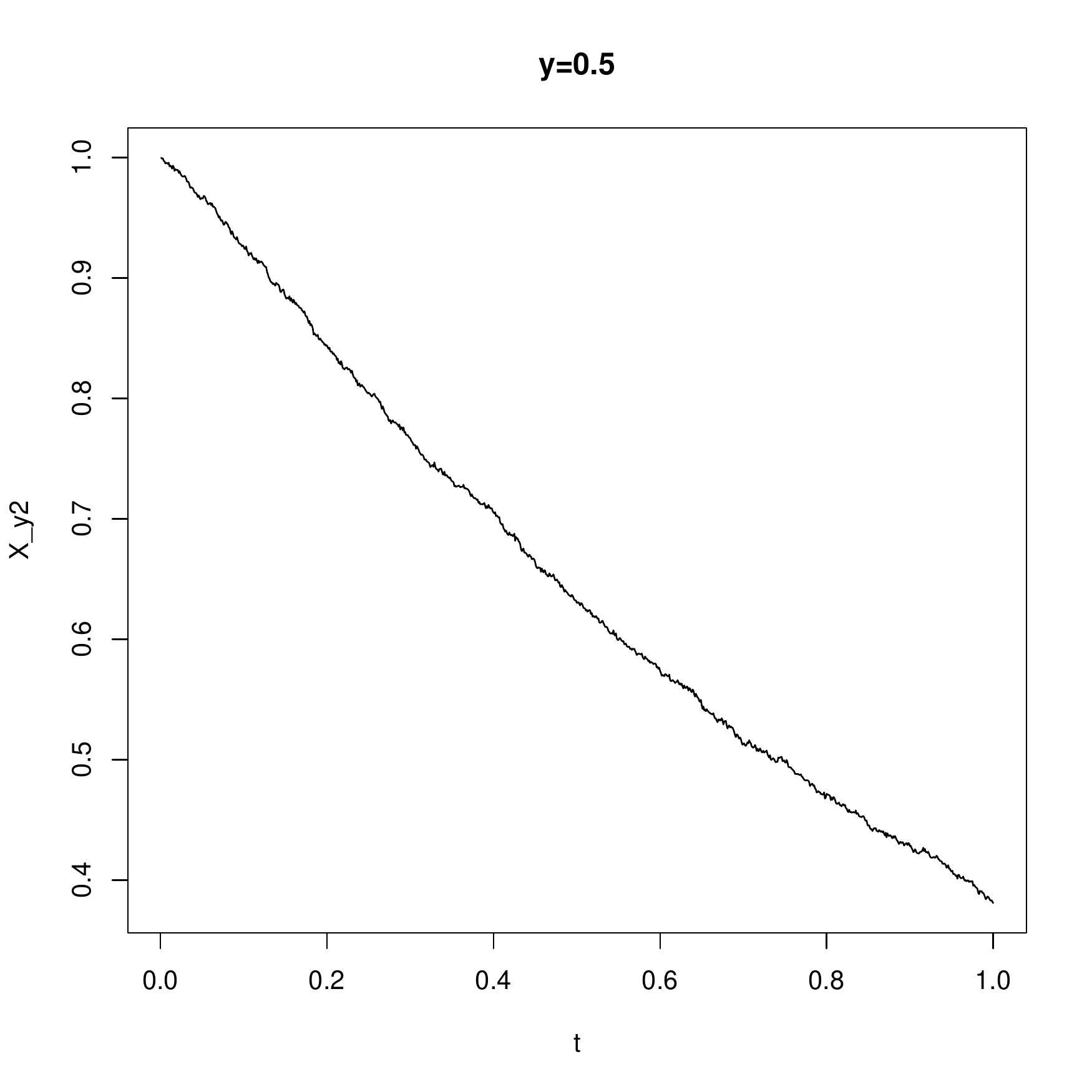}
\captionsetup{labelformat=empty,labelsep=none}
\subcaption{$\theta=(1,0.2,0.2)$, $\lambda_1 \approx 1$, \\ \hspace{0.5cm} $\epsilon=0.01$}
\end{center}
\end{minipage}
\begin{minipage}{0.32\hsize}
\begin{center}
\includegraphics[width=4.5cm]{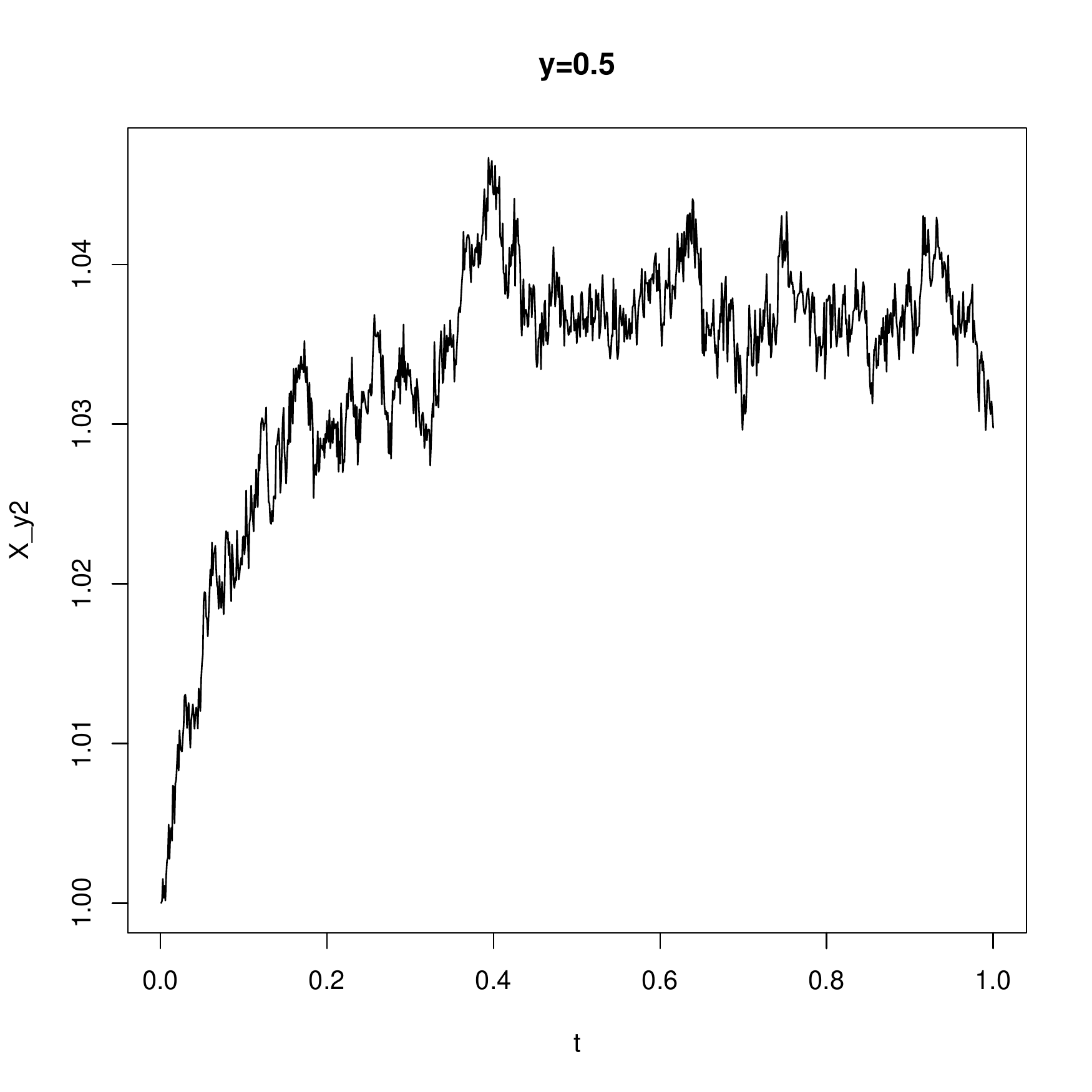}
\captionsetup{labelformat=empty,labelsep=none}
\subcaption{$\theta=(2,0.2,0.2)$, $\lambda_1\approx 0.02$, \\ \hspace{0.5cm}  $\epsilon=0.01$}
\end{center}
\end{minipage}
\begin{minipage}{0.32\hsize}
\begin{center}
\includegraphics[width=4.5cm]{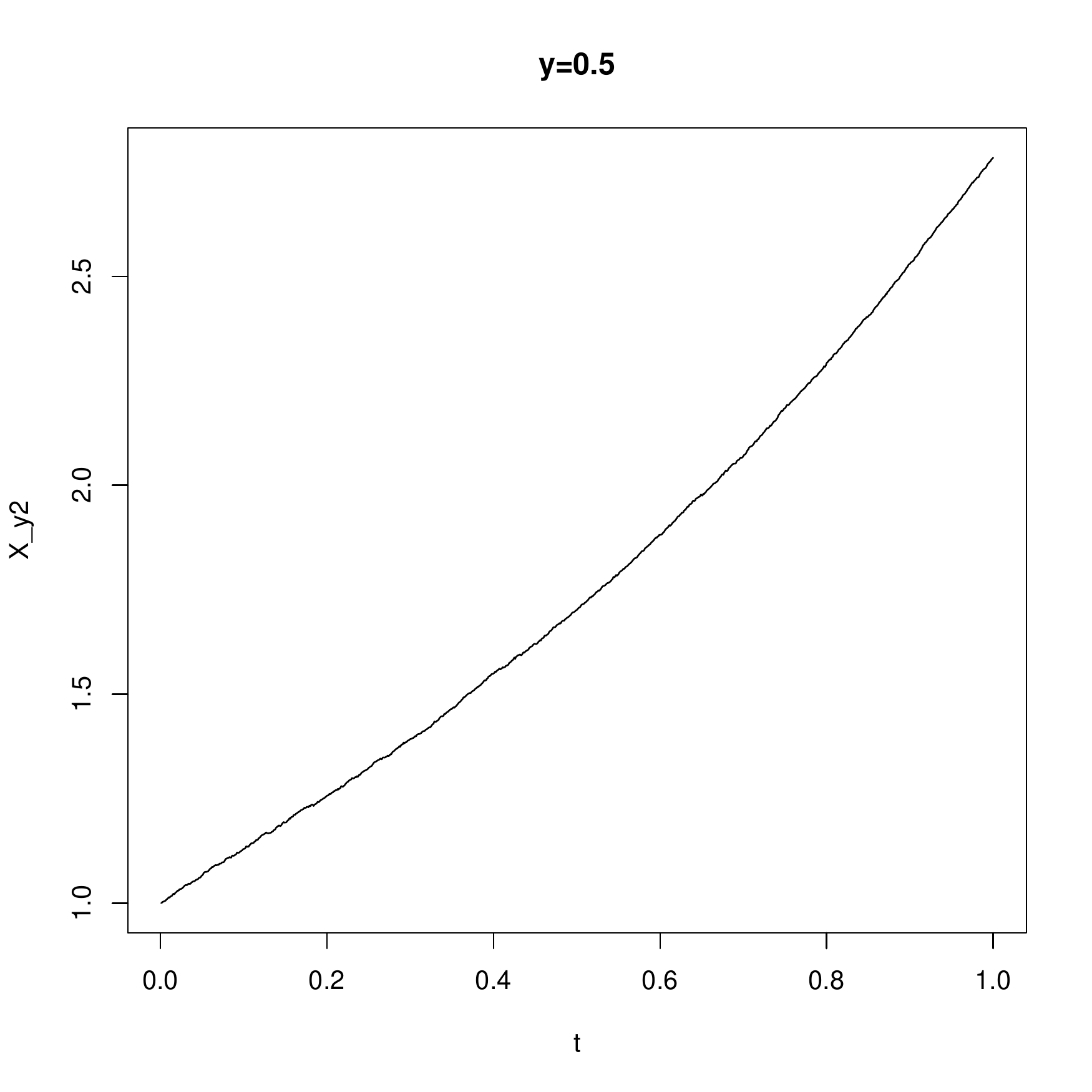}
\captionsetup{labelformat=empty,labelsep=none}
\subcaption{$\theta=(3,0.2,0.2)$, $\lambda_1\approx -1$, \\  \hspace{0.5cm} $\epsilon=0.01$}
\end{center}
\end{minipage}
\caption{Sample paths with $\lambda_1 \approx 1, 0.02, -1$ ($X_t(y)$ at $y = 0.5$)}
\end{figure}

Figures $4$-$6$ are the sample paths
{with the initial condition $\xi(y)=4y(1-y)$}, 
where $\theta_0$, $\theta_1$ and $\theta_2$ are fixed and only $\epsilon$ is changed.
For three kinds of $\lambda_1$, which are positive, near $0$ and negative,
it can be seen that in all cases, the noise of SPDE (\ref{spde0}) increases as $\epsilon$ increases.

\begin{figure}[H]
\begin{minipage}{0.32\hsize}
\begin{center}
\includegraphics[width=4.5cm]{1_1_all.jpeg}
\captionsetup{labelformat=empty,labelsep=none}
\subcaption{$\theta=(1,0.2,0.2)$, $\lambda_1 \approx 1$, \\ \hspace{0.5cm} $\epsilon=0.01$}
\end{center}
\end{minipage}
\begin{minipage}{0.32\hsize}
\begin{center}
\includegraphics[width=4.5cm]{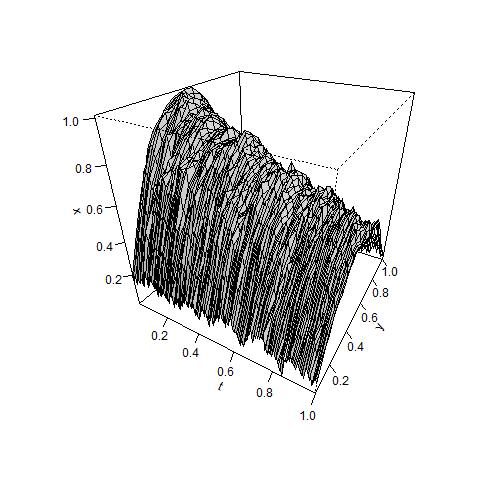}
\captionsetup{labelformat=empty,labelsep=none}
\subcaption{$\theta=(1,0.2,0.2)$, $\lambda_1 \approx 1$, \\ \hspace{0.5cm} $\epsilon=0.1$}
\end{center}
\end{minipage}
\begin{minipage}{0.32\hsize}
\begin{center}
\includegraphics[width=4.5cm]{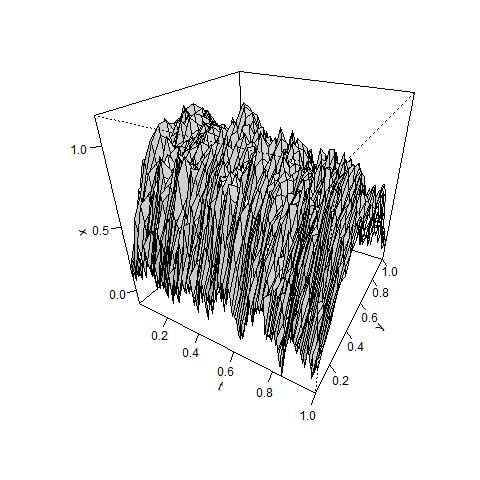}
\captionsetup{labelformat=empty,labelsep=none}
\subcaption{$\theta=(1,0.2,0.2)$, $\lambda_1 \approx 1$, \\ \hspace{0.5cm} $\epsilon=0.25$}
\end{center}
\end{minipage}
\caption{Sample paths with $\epsilon=0.01, 0.1, 0.25$, $\lambda_1\approx$1}
\label{t1-11}
\end{figure}

\begin{figure}[H]
\begin{minipage}{0.32\hsize}
\begin{center}
\includegraphics[width=4.5cm]{2_1_all.jpeg}
\captionsetup{labelformat=empty,labelsep=none}
\subcaption{$\theta=(2,0.2,0.2)$, $\lambda_1\approx 0.02$, \\ \hspace{0.5cm}  $\epsilon=0.01$}
\end{center}
\end{minipage}
\begin{minipage}{0.32\hsize}
\begin{center}
\includegraphics[width=4.5cm]{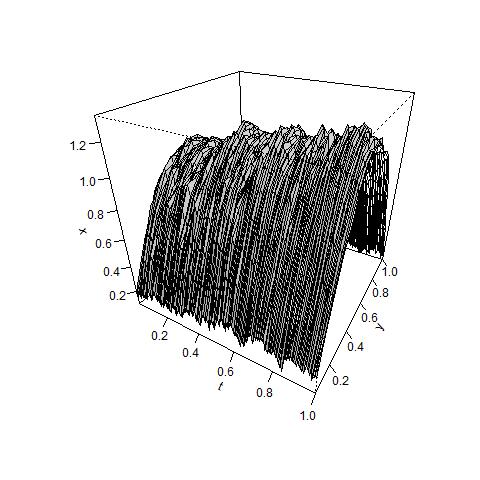}
\captionsetup{labelformat=empty,labelsep=none}
\subcaption{$\theta=(2,0.2,0.2)$, $\lambda_1\approx 0.02$, \\ \hspace{0.5cm}  $\epsilon=0.1$}
\end{center}
\end{minipage}
\begin{minipage}{0.32\hsize}
\begin{center}
\includegraphics[width=4.5cm]{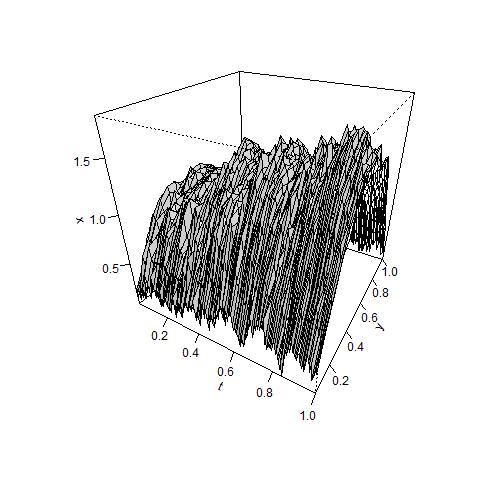}
\captionsetup{labelformat=empty,labelsep=none}
\subcaption{$\theta=(2,0.2,0.2)$, $\lambda_1\approx 0.02$, \\ \hspace{0.5cm}  $\epsilon=0.25$}
\end{center}
\end{minipage}
\caption{Sample paths with $\epsilon=0.01, 0.1, 0.25$, $\lambda_1\approx 0.02$}
\label{t1-11}
\end{figure}

\begin{figure}[H]
\begin{minipage}{0.32\hsize}
\begin{center}
\includegraphics[width=4.5cm]{3_1_all.jpeg}
\captionsetup{labelformat=empty,labelsep=none}
\subcaption{$\theta=(3,0.2,0.2)$, $\lambda_1\approx -1$, \\  \hspace{0.5cm} $\epsilon=0.01$}
\end{center}
\end{minipage}
\begin{minipage}{0.32\hsize}
\begin{center}
\includegraphics[width=4.5cm]{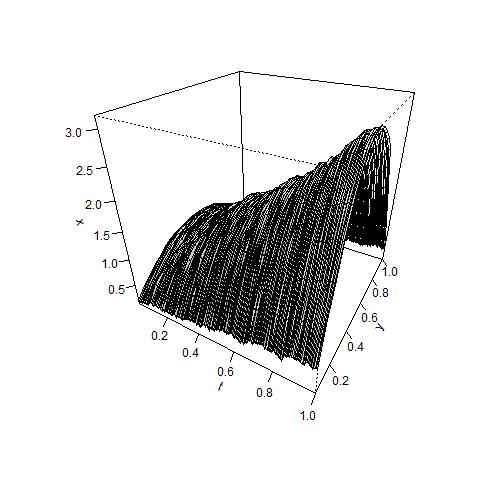}
\captionsetup{labelformat=empty,labelsep=none}
\subcaption{$\theta=(3,0.2,0.2)$, $\lambda_1\approx -1$, \\  \hspace{0.5cm} $\epsilon=0.1$}
\end{center}
\end{minipage}
\begin{minipage}{0.32\hsize}
\begin{center}
\includegraphics[width=4.5cm]{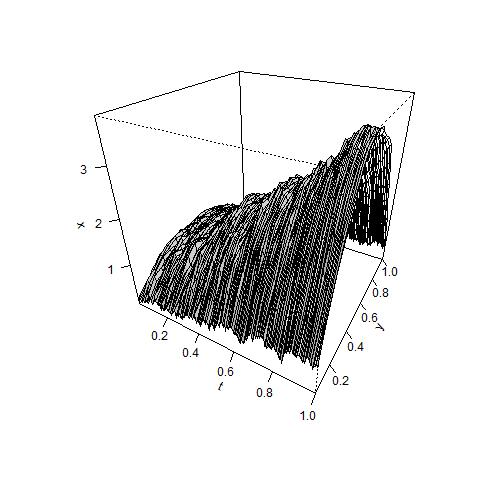}
\captionsetup{labelformat=empty,labelsep=none}
\subcaption{$\theta=(3,0.2,0.2)$, $\lambda_1\approx -1$, \\  \hspace{0.5cm} $\epsilon=0.25$}
\end{center}
\end{minipage}
\caption{Sample paths with $\epsilon=0.01, 0.1, 0.25$, $\lambda_1\approx -1$}
\label{t1-11}
\end{figure}

In order to estimate the unknown parameter $\lambda_1$,
we construct the estimators of  three parameters $\theta_0$, $\theta_1$ and $\theta_2$.
For the properties of $\theta_0$, $\theta_1$ and $\theta_2$,
we can refer Kaino and Uchida {(2020)}. 
First of all, we can get the minimum contrast estimator $(\hat{\theta}_1, \hat{\theta}_2)$
of $(\theta_1, \theta_2)$ 
in the same way as in Bibinger and Trabs (2020).
Next, using the minimum contrast estimator $(\hat{\theta}_1, \hat{\theta}_2)$,
we obtain the following approximate coordinate process
}
\begin{eqnarray*}
\hat{x}_k(t) = \frac{1}{M}\sum^M_{j = 1} X_{t}(y_{j:M}) \sqrt{2}\sin (\pi ky_{j:M}) 
\exp \left( \frac{\hat{\theta}_1 y_{j:M}}{2\hat{\theta}_2} \right)
\end{eqnarray*}
{
and the adaptive ML type estimator $\hat{\theta}_0$ of $\theta_0$ 
is constructed by using  the property that 
the coordinate process (\ref{coordinate}) is a small diffusion process.  
Using statistical inference for small diffusion process
and the thinned data of the approximate coordinate process $\hat{x}_k(t)$,
we can show that the adaptive ML type estimator $\hat{\theta}_0$ of $\theta_0$ has
asymptotic normality under some regularity conditions}. 
For details of statistical inference for small diffusion processes,
see
Kutoyants (1984, 1994), 
Genon-Catalot (1990),  
Laredo (1990), 
S{\o}rensen and Uchida (2003),
Gloter and S{\o}rensen (2009), 
Guy et.\ al.\ (2014)
and 
Kaino and Uchida (2018a). 
For  adaptive ML type estimators and thinned data
for diffusion type processes, see for example, Uchida and Yoshida (2012) 
and Kaino and Uchida (2018b). 
Furthermore, in order to illustrate the asymptotic performance of the estimators of the coefficient parameters 
of the parabolic linear second order SPDE model based on high-frequency data, 
some examples and simulation results of the estimators 
$\hat{\theta}_0$, $\hat{\theta}_1$ and $\hat{\theta}_2$ are given.

{
This paper is organized as follows. 
In Section 2, we first introduce
the minimum contrast estimator of $(\theta_1, \theta_2)$ 
in the SPDE model based on 
the high frequency data 
in the fixed region
$[0,T] \times [0,1]$.
Next, the adaptive ML type estimator of $\theta_0$ 
is constructed by using the minimum contrast estimator $(\hat{\theta}_1, \hat{\theta}_2)$
and the approximate coordinate process. 
It is shown that the adaptive ML type estimator $\hat{\theta}_0$ has asymptotic normality.
In Section 3, we give some examples and simulation results of 
the asymptotic behavior of the estimator $(\hat{\theta}_0, \hat{\theta}_1, \hat{\theta}_2)$ 
proposed in {Section 2}.
The theoretical results in {Section 2} are proved in Section 4.
}

\begin{en-text}
We consider a linear parabolic stochastic partial differential equation (SPDE) with one space dimension.
\begin{eqnarray}
& &dX_{t}(y) = \left( \theta_2 \frac{\partial^2 X_{t}(y)}{\partial y^2} + \theta_1 \frac{\partial X_{t}(y)}{\partial y} + \theta_0 X_{t}(y) \right)dt 
+ \epsilon dB_t(y), \quad  (t,y) \in [0, T]\times[0, 1],  \quad
\label{spde0} \\
& & X_t(0) = X_t(1) = 0, \quad t \in [0,T], \qquad X_0(y) = \xi, \quad y \in [0,1],
\nonumber
\end{eqnarray}
where
$T>0$ is fixed, $B_t$ is defined as a cylindrical Brownian motion in the Sobolev space on $[0,1]$,
an unknown parameter $\theta= (\theta_0, \theta_1, \theta_2)$  
and $\theta_0,\theta_1 \in \mathbb{R}, \theta_2$,
and the parameter space $\Theta$ is a compact convex subset of $\mathbb{R}^2 \times [0, \infty]$.
Moreover, the true value of parameter $\theta^*= (\theta_0^*, \theta_1^*, \theta_2^*)$
and we assume that $\theta^* \in \mbox{Int}(\Theta)$.
The data are discrete observations ${\bf X}_{N,M}= \left\{ X_{t_{i:N}}(y_{j:M}) \right\}_{i = 1, ..., N, j = 1, ..., M}$, 
$t_{i:N} = i \frac{T}{N}$,
${ y_{j:M} = \frac{j}{M}}$.
{
}

Statistical inference for SPDE models based on sampled data has been developed
by some researchers.
For example, 
Cont (2005),
Cialenco and Huang (2017), 
Bibinger and Trabs (2017)
and
Hildebrandt (2019).
Kaino and Uchida {(2020)} studied adaptive maximum likelihood (ML) type estimator of the coefficient parameter
of the parabolic linear second order SPDE model.
Hildebrandt and Trabs (2019) proposed a contrast function with double increments. 
They obtained the estimator of the coefficient parameter of the SPDE model based on the contrast function and showed that the estimator satisfied asymptotic normality.

In this paper, we propose adaptive maximum likelihood (ML) type estimator of the coefficient parameter
$\theta= (\theta_0, \theta_1, \theta_2)$ 
of the parabolic linear second order SPDE model (\ref{spde0}).
For $k \in \mathbb{N}$, 
the coordinate process  $x_k(t)$ of the SPDE model (\ref{spde0}) is 
\begin{eqnarray}
x_k(t) 
&=& \int^1_0 X_t(y) \sqrt{2}\sin (\pi k y) \exp \left( \frac{\eta y}{2} \right) dy, \label{coordinate}
\end{eqnarray}
where $\eta = \frac{\theta_1}{\theta_2}$.
Note that the coordinate process (\ref{coordinate}) 
is  the Ornstein-Uhlenbeck process with small perturbations.
Using the minimum contrast estimator $\hat{\theta}_1$ and $\hat{\theta}_2$,
we obtain the approximate coordinate process 
\begin{eqnarray*}
\hat{x}_k(t) = \frac{1}{M}\sum^M_{j = 1} X_{t}(y_{j:M}) \sqrt{2}\sin (\pi ky_{j:M}) 
\exp \left( \frac{\hat{\theta}_1 y_{j:M}}{2\hat{\theta}_2} \right)
\end{eqnarray*}
and the adaptive estimator is constructed by using  the property that 
the coordinate process (\ref{coordinate}) is a small diffusion process. 
{
It is also shown that the adaptive ML type estimators satisfy asymptotic normality under some regularity conditions. 
Furthermore, in order to illustrate the asymptotic performance of the adaptive ML type estimators of the coefficient parameters of the parabolic linear second order SPDE model based on high-frequency data, some examples and simulation results of the adaptive ML type estimators are given. 
}
For details of statistical inference for small diffusion processes,
see
Gloter and S{\o}rensen (2009), 
Guy et.\ al.\ (2014), 
Kaino and Uchida (2018a),  
{
For  adaptive ML type estimators and thinned data
for diffusion type processes, see for example, Uchida and Yoshida (2012) 
and Kaino and Uchida (2018b). 
}

This paper is organized as follows. 
In Section 2, we consider 
{the adaptive estimator} 
of the SPDE model based on 
{the sampled data 
in the fixed region
$[0,1] \times [0,1]$.}
The adaptive estimator is constructed by using the minimum contrast estimators of $\theta_1$ and $\theta_2$. 
It is shown that the adaptive estimator has asymptotic normality.
In Section 3, concrete examples are given and 
the asymptotic behavior of the estimators proposed in Sections 2 is verified by simulations.
Section 4 is devoted to the proofs of the results presented in Sections 2.

\end{en-text}

\section{Estimation of $\theta_1$, $\theta_2$ and $\theta_0$}

{
For real-valued functions $f$ and $g$ defined on $[0,1]$,
let 
$\langle f,g \rangle_\theta = \int^1_0 e^{y\theta_1/\theta_2}f(y)g(y) dy$ 
and $ \| f \|_\theta = \sqrt{\langle f,f \rangle_\theta}$.
Moreover, set 
$H_\theta = \{ f:[0,1] \to \mathbb{R} : \| f \|_\theta < \infty, f(0) = f(1) = 0  \}$. 
}

\begin{en-text}
In this section, we treat the linear parabolic SPDE (\ref{spde0}) with $T=1$, which is defined as
\begin{eqnarray}
& &dX_{t}(y) = \left( \theta_2 \frac{\partial^2 X_{t}(y)}{\partial y^2} + \theta_1 \frac{\partial X_{t}(y)}{\partial y} + \theta_0 X_{t}(y) \right)dt 
+ \sigma dB_t(y), \quad  (t,y) \in [0, 1]\times[0, 1],  \quad
\label{spde1} \\
& & X_t(0) = X_t(1) = 0, \quad t \in [0,1], \qquad X_0(y) = \xi = 0, \quad y \in [0,1].
\nonumber
\end{eqnarray}
The data are discrete observations ${{\bf \hat{x}}_{N, m}}= \left\{ X_{t_{i:N}}({\bar{y}_{j:M}}) \right\}_{i = 1, \ldots, N, j = 1, \ldots, {m}}$ with
{$t_{i:N} =  \frac{i}{N}$},
{$\bar{y}_{j:M} = \delta+ \frac{j-1}{M}$ and
$m=[(1-2\delta)M]$ for $\delta \in (0,1/2) $.
Note that $\delta \leq \bar{y}_{j:M} \leq 1-\delta$ for $ j = 1, \ldots, m$,
and that ${\bf \hat{x}}_{N, m}$
is generated by the full data  ${\bf X}_{N,M}= \left\{ X_{t_{i:N}}(y_{j:M}) \right\}_{i = 1, ..., N, j = 1, ..., M}$
with $t_{i:N} = \frac{i}{N}$ and $y_{j:M} = \frac{j}{M}$}.
\end{en-text}

{
The differential operator is given by
\begin{eqnarray*}
A_\theta = \theta_0  + \theta_1\frac{\partial}{\partial y} + \theta_2\frac{\partial^2}{\partial y^2},
\end{eqnarray*}
and the eigenfunctions $e_k$ of $A_\theta$ and the corresponding eigenvalues $-\lambda_k$ 
are defined as 
\begin{eqnarray*}
e_k(y) &=& \sqrt{2} \sin(\pi k y)\exp\left( -\frac{\theta_1}{2\theta_2} y \right), \;\;\; y\in [0,1], \\
\lambda_k &=& - \theta_0 + \frac{\theta_1^2}{4 \theta_2} + \pi^2 k^2 \theta_2. 
\end{eqnarray*}
We then obtain that  for $k \in \mathbb{N}$,
$$
A_\theta  e_k =  -\lambda_k e_k.
$$

\noindent
The coordinate process is defined as  
\begin{eqnarray*}
x_k(t) &=& \langle X_t, e_k \rangle_\theta 
= \int^1_0 \exp \left( \frac{\theta_1}{\theta_2} y \right)  X_t(y) e_k(y) dy \\
&=& \int^1_0 X_t(y) \sqrt{2}\sin (\pi k y) \exp \left( \frac{\theta_1 y}{2\theta_2} \right) dy .
\end{eqnarray*}
Here we note that the random field $X_t(y)$ is 
\begin{eqnarray*}
X_t(y) = \sum_{k=1}^\infty x_k(t) e_k(y).
\end{eqnarray*}
Moreover, as we stated in Introduction section, 
we notice that $x_k(t)$ is  the Ornstein-Uhlenbeck process as follows.
\begin{equation}
dx_k(t) = -\lambda_k x_k(t) dt + \epsilon dw_k(t), \quad x_k(0) = \langle \xi, e_k \rangle_\theta, 
\label{coordinate2}
\end{equation}
where $(w_k(t))_{t \geq 0}, k \in \mathbb{N}$ is independent real-valued Brownian motions. 

We assume that $\lambda^*_1 >0$. 
Furthermore,  we make the following assumption. \\
\noindent
$[A1]$ $\xi$ is non-random, $\| A_\theta^{1/2} \xi \|^2_\theta < \infty$
and  $
\langle \xi, e_1 \rangle_\theta 
\ne 0$.
}

\begin{en-text}
The cylindrical Brownian motion $(B_t)_{t \ge 0}$ in (\ref{spde1}) can be defined as 
\begin{eqnarray*}
\langle B_t,f \rangle_\theta = \sum_{k \ge 1} \langle f, e_k \rangle_\theta W^k_t, \quad 
f \in H_\theta, \quad t \ge 0
\end{eqnarray*}
for independent real-valued Brownian motions $(W^k_t)_{t\ge0}$, $k\ge1$. 
$X_t(y)$ is called a mild solution of (\ref{spde1}) on $[0,1]$ if it satisfies that for any $t \in [0,1]$,  
\begin{eqnarray*}
X_t = e^{tA_\theta}\xi + \int^t_0 e^{(t-s)A_\theta}\sigma dB_s \;\;\; a.s.
\end{eqnarray*}
\end{en-text}

{The data are discrete observations ${{\bf \bar{X}}_{N, \bar{M}}}= \left\{X_{t_{i:N}}({\bar{y}_{j:M}}) \right\}_{i = 1, \ldots, N, j = 1, \ldots, {\bar{M}}}$ with
{$t_{i:N} =  i\frac{T}{N}$},
$\bar{y}_{j:M} = \delta+ \frac{j-1}{M}$ and
$\bar{M}=1+[(1-2\delta)M]$ for $\delta \in (0,1/2) $.
Let $m \leq {\bar{M}}$
and
$\tilde{y}_{j:m} = {\delta+} \left[ \frac{{\bar{M}}}{m} \right] {\frac{j-1}{M}} $ 
for $j = 1,...,m$.
}

Assume that 
$m \leq N^\rho$ for $\rho \in (0, 1/2)$.
Let
\begin{eqnarray*}
Z_{j:m} = {\frac{1}{N \sqrt{{t}_{1:N}}} }
\sum_{i=1}^{N} (X_{t_{i:N}}(\tilde{y}_{j:m} )- X_{t_{i-1:N}}(\tilde{y}_{j:m}))^2.
\end{eqnarray*}
{
Setting that the contrast function is 
\begin{eqnarray*}
U_{N,m}(\theta_1,\theta_2) = \frac{1}{m} \sum^{m}_{j = 1}  \left( \frac{1}{\epsilon^2}Z_{j:m} - \frac{1}{\sqrt{\pi \theta_2}} \exp(- \frac{\theta_1}{\theta_2} \tilde{y}_{j:m}) \right)^2,
\end{eqnarray*}
the minimum contrast estimator of $\theta_1$ and $\theta_2$ are defined as 
\begin{eqnarray*}
(\hat{\theta}_1,\hat{\theta}_2) = \arginf_{\theta_1,\theta_2}U_{N,m}(\theta_1,\theta_2).
\end{eqnarray*}
}

\begin{en-text}
Let
\begin{eqnarray*}
U(\zeta^*)  &=& 
\begin{pmatrix}
\int^{1-\delta}_\delta e^{-4 \eta^* y}dy & -(\sigma^*_0)^2\int^{1-\delta}_\delta ye^{-4 \eta^* y}dy\\
 -(\sigma_0^*)^2\int^{1-\delta}_\delta ye^{-4 \eta^* y}dy& (\sigma^*_0)^4\int^{1-\delta}_\delta y^2e^{-4 \eta^* y}dy
\end{pmatrix}, \\
V(\zeta^*) &=& 
\begin{pmatrix}
\int^{1-\delta}_\delta e^{-2 \eta^* y}dy & -(\sigma^*_0)^2\int^{1-\delta}_\delta ye^{-2 \eta^2 y}dy\\
 -(\sigma^*_0)^2\int^{1-\delta}_\delta ye^{-2 \eta^* y}dy& (\sigma^*_0)^4\int^{1-\delta}_\delta y^2e^{-2 \eta^* y}dy
\end{pmatrix}, \\
\Gamma &=& \frac{1}{\pi} \sum^{\infty}_{r = 0} I(r)^2 + \frac{2}{\pi}
\ \mbox{with} \ I(r) = 2\sqrt{r+1} - \sqrt{r+2} - \sqrt{r}.
\end{eqnarray*}

\begin{theorem}[Bibinger and Trabs(2017)] \label{thm1}
Let $\zeta^* = ((\sigma_0^*)^2, \eta^*) \in \Xi$ for {a} compact subset $\Xi \subseteq (0,\infty) \times [0,\infty)$.
Assume that $m \to \infty$ and 
$m =  {O}(N^\rho)$ for some $\rho \in (0,1/2)$.
Then, as $N \to \infty$, 
\begin{eqnarray*}
\sqrt{m N}(\hat{\sigma}_0^2 - (\sigma_0^*)^2,  \hat{\eta}-\eta^* )
\stackrel{d}{\longrightarrow } N(0,(\sigma_0^*)^4 \Gamma \pi V(\zeta^*)^{-1}U(\zeta^*)V(\zeta^*)^{-1}).
\end{eqnarray*}
\end{theorem}

\end{en-text}

\begin{en-text}

Let $k \in \mathbb{N}$.
In order to estimate $\theta_0$,  
we use the minimum contrast estimator $(\hat{\theta}_1,\hat{\theta}_2)$
and  the coordinate process defined that  
\begin{eqnarray*}
x_k(t) &=& \langle X_t, e_k \rangle_\theta 
= \int^1_0 \exp \left( \frac{\theta_1}{\theta_2} y \right)  X_t(y) e_k(y) dy \\
&=& \int^1_0 X_t(y) \sqrt{2}\sin (\pi k y) \exp \left( \frac{\theta_1 y}{2\theta_2} \right) dy .
\end{eqnarray*}
Here we note that the random field $X_t(y)$ is 
\begin{eqnarray*}
X_t(y) = \sum_{k=1}^\infty x_k(t) e_k(y).
\end{eqnarray*}
Moreover, as we stated in Introduction section, 
we notice that $x_k(t)$ is  the Ornstein-Uhlenbeck process as follows.
\begin{equation}
dx_k(t) = -\lambda_k x_k(t) dt + \epsilon dw_k(t), \quad x_k(0) = \langle \xi, e_k \rangle_\theta, 
\label{coordinate2}
\end{equation}
where $(w_k(t))_{t \geq 0}, k \in \mathbb{N}$ is independent real-valued Brownian motions. 

We assume that $\lambda^*_1 >0$. 
Furthermore,  we make the following assumption. \\
\noindent
$[A1]$ $\xi$ is non-random, $\| A_\theta^{1/2} \xi \|^2_\theta < \infty$
and  $
\langle \xi, e_1 \rangle_\theta 
\ne 0$.

\end{en-text}

Let $N_2 \leq N$, ${s}_{i:N_2}=\left[\frac{N}{N_2}\right] t_{i:N}=i \left[\frac{N}{N_2}\right] \frac{T}{N}$ for $i = 1,...,N_2$ 
and $\delta_{N_2} = {s}_{i:N_2} - {s}_{i-1:N_2} = \left[\frac{N}{N_2}\right] \frac{T}{N}$. 
{
The approximate coordinate process $\hat{x}_k(t)$ is given  by
\begin{eqnarray*}
\hat{x}_k(t) = \frac{1}{M}\sum^M_{j = 1} X_{t}(y_{j:M}) \sqrt{2}\sin (\pi ky_{j:M}) 
\exp \left( \frac{\hat{\theta}_1 y_{j:M}}{2\hat{\theta}_2} \right)
\end{eqnarray*}
and 
{
${\bf \hat{x}}_k =\{\hat{x}_k(s_{i:N_2})\}_{i=1,\ldots, N_2}$
are  the thinned data of the approximate coordinate process.} 
The quasi log-likelihood function based on
the thinned data 
${\bf \hat{x}}_k$
are given by }
\begin{eqnarray*}
l_{N_2}(\lambda_k \ | \ {\bf \hat{x}}_k  ) 
&=&
 -\sum^{N_2}_{i=1}\left\{ \frac{1}{2} \log 
 \left( \frac{\epsilon^2(1 - \exp (-2\lambda_k \delta_{N_2} ))}{2\lambda_k} \right) \right. \\
& & \left. 
+ \frac{(\hat{x}_k( s_{i:N_2} ) - \exp (-\lambda_k \delta_{N_2})
\hat{x}_k(s_{i-1:N_2}))^2}{\frac{2\epsilon^2(1-\exp (-2\lambda_k \delta_{N_2}))}{2\lambda_k}} \right\}.
\end{eqnarray*}
{
The adaptive ML type estimator of $\lambda_k$ is defined as 
\begin{eqnarray*}
\hat{\lambda}_k = \argsup_{\lambda_k} l_{N_2}(\lambda_k \ | \ {\bf \hat{x}}_k  ). 
\end{eqnarray*}
If we set $k = 1$,
then $\lambda_1^* = - \theta_0^* + \frac{(\theta_1^*)^2}{4 \theta_2^*} + \pi^2 \theta_2^*$.
The adaptive ML type estimator of $\theta_0$ is given by 
\begin{eqnarray*}
\hat{\theta}_0 = -\hat{\lambda}_1 + \frac{\hat{\theta}_1^2}{4\hat{\theta}_2} + \pi^2\hat{\theta}_2. 
\end{eqnarray*}
}
{
Let
$\eta^* = \frac{\theta_1^*}{\theta_2^*}$, 
$I(r) = 2\sqrt{r+1} - \sqrt{r+2} - \sqrt{r}$,
\begin{eqnarray*}
\Gamma &=& \frac{1}{\pi} \sum^{\infty}_{r = 0} I(r)^2 + \frac{2}{\pi},  \\ 
G(\lambda^*) &=& \frac{x_1(0)^2}{2\lambda^*}(1-e^{-2\lambda^*}), \\
\left(
\begin{array}{cc}
J_{1,1} & J_{1,2} \\
J_{1,2} & J_{2,2} 
\end{array}
\right)
&=& \frac{1}{\theta_2^*} \pi \Gamma W V^{-1} U V^{-1} W^{T},
\end{eqnarray*}
where \begin{eqnarray*}
U &{=}& 
\begin{pmatrix}
\int^{1-\delta}_\delta e^{-4 \eta^* y}dy & -\frac{1}{\sqrt{\theta_2^*}}\int^{1-\delta}_\delta ye^{-4 \eta^* y}dy\\
 -\frac{1}{\sqrt{\theta_2^*}}\int^{1-\delta}_\delta ye^{-4 \eta^* y}dy& 
 \frac{1}{\theta_2^*}\int^{1-\delta}_\delta y^2e^{-4 \eta^* y}dy
\end{pmatrix}
,\\
V &{=}& 
\begin{pmatrix}
\int^{1-\delta}_\delta e^{-2 \eta^* y}dy & -\frac{1}{\sqrt{\theta_2^*}}\int^{1-\delta}_\delta ye^{-2 \eta^* y}dy\\
 -\frac{1}{\sqrt{\theta_2^*}}\int^{1-\delta}_\delta ye^{-2 \eta^* y}dy& 
 \frac{1}{\theta_2^*}\int^{1-\delta}_\delta y^2e^{-2 \eta^* y}dy
\end{pmatrix}
,\\
W &{=}& 
\begin{pmatrix}
-2(\theta_2^*)^{\frac{3}{2}} & 0\\
-2\theta_1^*(\theta_2^*)^{\frac{1}{2}} & \theta_2^*
\end{pmatrix}
\end{eqnarray*}
and $W^T$ is the transpose of $W$.
}

\begin{theorem} \label{thm4}
Assume [A1] , $m \to \infty$ and 
$m = {O}(N^\rho)$ for some $\rho \in (0,1/2)$.
Moreover assume that {$\frac{1}{\epsilon \sqrt{N_2}} =  O(1)$,
$\epsilon \sqrt{N_2} =  O(1)$},
$\frac{Nm}{M^2} =  {O}(1)$,
$\frac{N_2}{M^{1-\rho_1}} \to 0$ for $\rho_1 \in (0,1)$ and
{
$\frac{N_2}{\epsilon^2 N m} \to 0$}.
{As $N_2,m \to \infty$ and $\epsilon \to 0$,} \\
$$
\begin{pmatrix}
\sqrt{Nm}(\hat{\theta}_2 - \theta_{2}^*) \\
\sqrt{Nm}(\hat{\theta}_1 - \theta_{1}^*) \\
\epsilon^{-1}(\hat{\theta}_0 - \theta_{0}^*)
\end{pmatrix}
\stackrel{d}{\longrightarrow } N
\left(
\begin{pmatrix}
0 \\
0 \\
0 
\end{pmatrix}
,
\begin{pmatrix}
J_{1,1} &J_{1,2} & 0\\
{J_{1,2}} &J_{2,2} & 0\\
0 & 0 & {G(\lambda_1^{*})^{-1}}
\end{pmatrix}
\right). 
$$
\end{theorem}



\section{Simulation results}
In the same way as Bibinger and Trabs (2020), the numerical solution of the SPDE (\ref{spde0}) is generated by
\begin{eqnarray}
\tilde{X}_{t_{i:N}}(y_{j:M}) = \sum^{K}_{k = 1}x_k(t_{i:N})e_k(y_{j:M}), \quad i = 1, ..., N, j = 1, ..., M,
\label{appro-SPDE}
\end{eqnarray}
where
\begin{eqnarray*}
x_k(t_{i:N}) = \exp \left( -\lambda_k \frac{T}{N} \right) x_k(t_{i-1:N}) + \sqrt{\frac{\epsilon^2(1-\exp(-2 \lambda_k \frac{T}{N}))}{2 \lambda_k}}N(0,1) ,\quad i = 1, ..., N.
\end{eqnarray*}
The number of iteration is $300$.


\subsection{Example 1}
{The true value of parameter  $\theta^* =(\theta_0^*, \theta_1^*, \theta_2^*) = (0,1,0.2)$, 
$\lambda_1^* = 3.22$}.
We set that $N = 10^4$, $M = 10^4$, $K = 10^5$, $T=1$, $x_1(0) = 3$, 
$\xi(y) = 4.2y(1-y)$.
{When $N=M=10^4$, the size of data ${\bf X}_{N,M}$ is about {1 GB}.
We used R language to compute {the estimators of Theorems 1}.
The personal computer with Intel Gold 6128 (3.40GHz) was used for this simulation.
Figure 7 is {a} sample path of $X_t(y)$ for $(t,y) \in [0,1]\times [0,1]$
when {$(\theta_0^*, \theta_1^*, \theta_2^*, \epsilon) = (0,1,0.2,0)$}.

\begin{figure}[h] 
\begin{center}
\includegraphics[width=9cm]{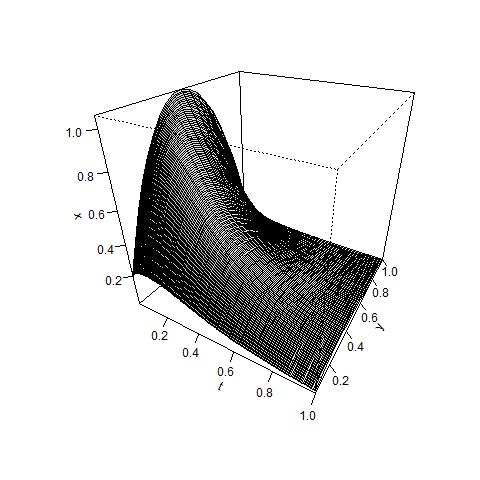} 
\caption{Sample path with $(\theta_0^*, \theta_1^*, \theta_2^*, \epsilon)  = (0,1,0.2,0)$, $\xi(y) = 4.2y(1-y)$ \label{fig5}}
\end{center}
\end{figure}

\subsubsection{$\epsilon=0.1$}

Figure 8 is {a} sample path of $X_t(y)$ for $(t,y) \in [0,1]\times [0,1]$
when $(\theta_0^*, \theta_1^*, \theta_2^*, \epsilon) = (0,1,0.2,0.1)$.
Table 1 is the simulation results of 
{the means and the standard deviations (s.d.s) of} $\hat{\theta}_1$, $\hat{\theta}_2$  and $\hat{\theta}_0$ 
with $(N, m, N_2) = (10^4, 99, 500)$.
Figures 9-11 are the simulation results of 
{the asymptotic distributions of} $\hat{\theta}_1$, $\hat{\theta}_2$  and $\hat{\theta}_0$
with $(N, m, N_2) = (10^4, 99, 500)$.
The left side of Figure 9 is the plot of the empirical distribution function of 
$\sqrt{Nm}(\hat{\theta}_1-\theta_1^*)$ (black line) and the distribution function of $N(0,J_{1,1})$ (red line).
The center of Figure 9 is the Q-Q plot of $\sqrt{Nm}(\hat{\theta}_1-\theta_1^*)$ and $N(0,J_{1,1})$.
The right side of Figure 9 is the plot of the histogram of $\sqrt{Nm}(\hat{\theta}_1-\theta_1^*)$ 
and the density function of $N(0,J_{1,1})$ (red line).
Figures 10 and 11 are the plots of the empirical distribution functions, the Q-Q plots 
and the histograms of $\sqrt{Nm}(\hat{\theta}_2 - \theta_{2}^*)$ 
and $\sqrt{\epsilon^{-1}}(\hat{\theta}_0 - \theta_{0}^*)$, respectively.
From Figures 9-11, 
we can see that the proposed estimators have
the asymptotic distribution in Theorem $1$ and these estimators have good performance.

\begin{figure}[h] 
\begin{center}
\includegraphics[width=9cm]{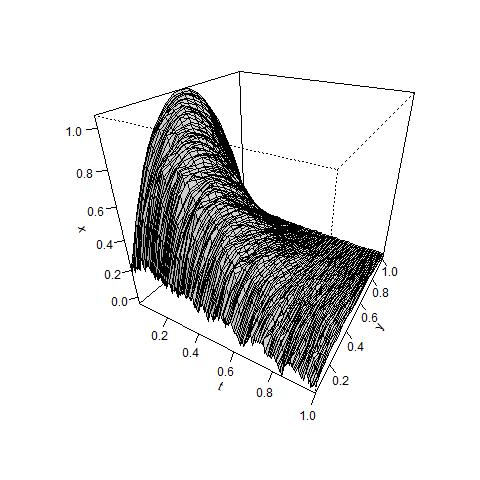} 
\caption{Sample path with 
$(\theta_0^*, \theta_1^*, \theta_2^*, \epsilon) =  (0,1,0.2,0.1)$, 
$\xi(y) = 4.2y(1-y)$ \label{fig5}}
\end{center}
\end{figure}



\begin{table}[h]
\caption{Simulation results of $\hat{\theta}_1$, $\hat{\theta}_2$  and $\hat{\theta}_0$ with $(N, m, N_2) = (10^4, 99, 500)$ \label{table2}}
\begin{center}
\begin{tabular}{c|ccc} \hline
		&$\hat{\theta}_{1}$&$\hat{\theta}_{2}$&$\hat{\theta}_{0}$
\\ \hline
{true value} &1 & 0.2 & 0
\\ \hline
mean &1.001&  0.200&  0.010
 \\
{s.d.} & (0.007)& (0.001)& (0.084)
 \\   \hline
\end{tabular}
\end{center}
\end{table}


\begin{figure}[h] 
\begin{center}
\includegraphics[width=5cm,pagebox=cropbox,clip]{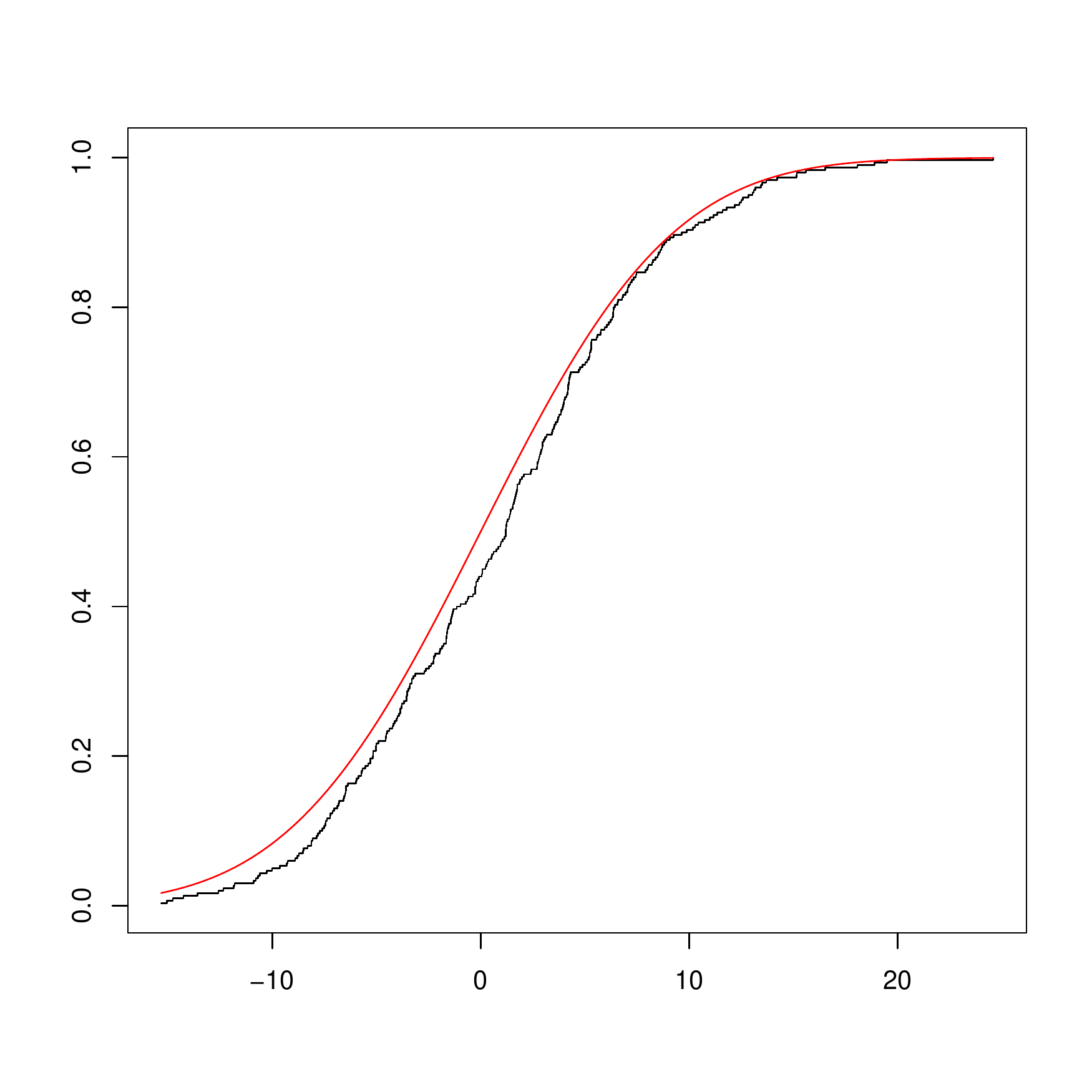}
\includegraphics[width=5cm,pagebox=cropbox,clip]{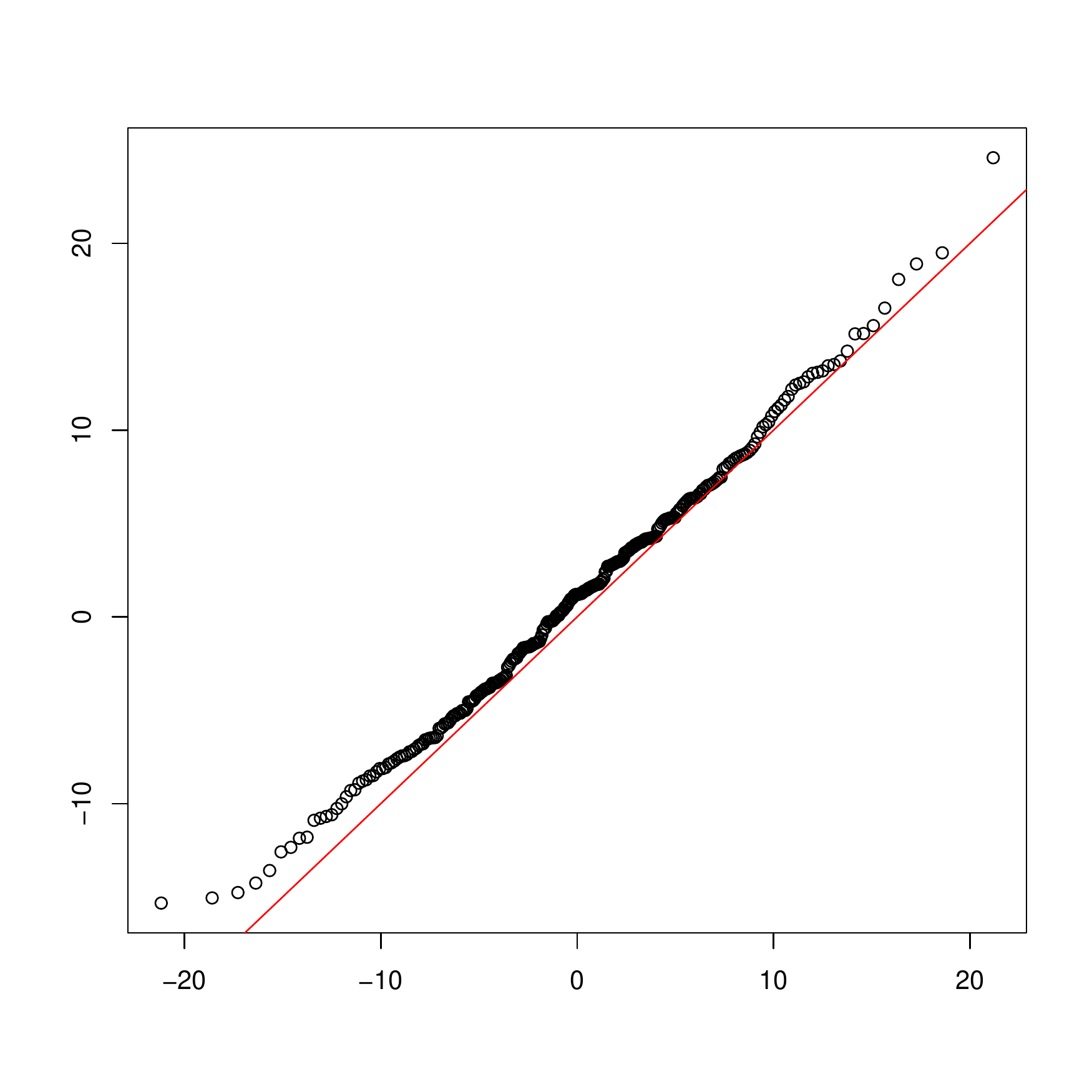}
\includegraphics[width=5cm,pagebox=cropbox,clip]{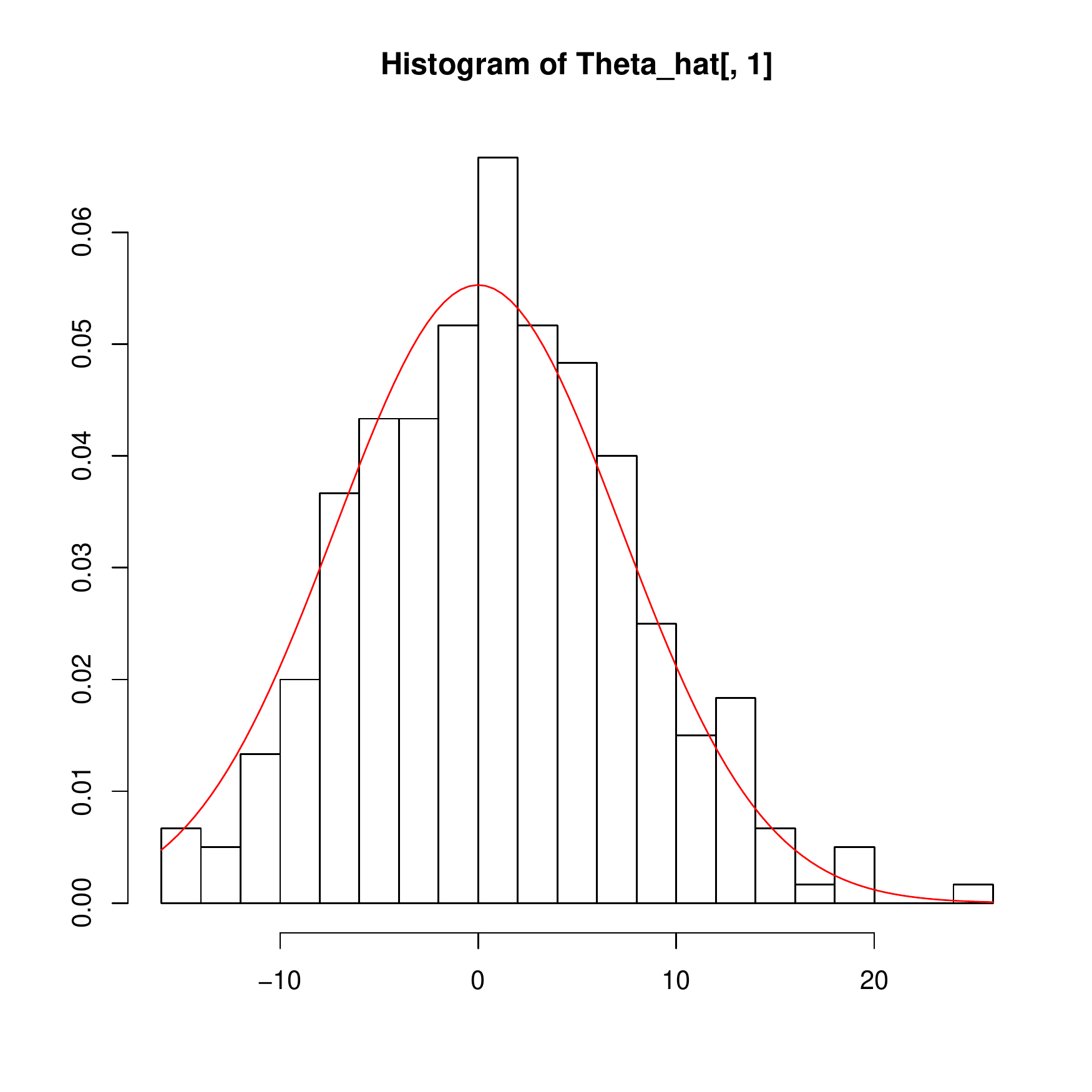}\\
\caption{Simulation results of $\hat{\theta}_1$ \label{fig6}}
\end{center}
\end{figure}

\begin{figure}[h]
\begin{center}
\includegraphics[width=5cm,pagebox=cropbox,clip]{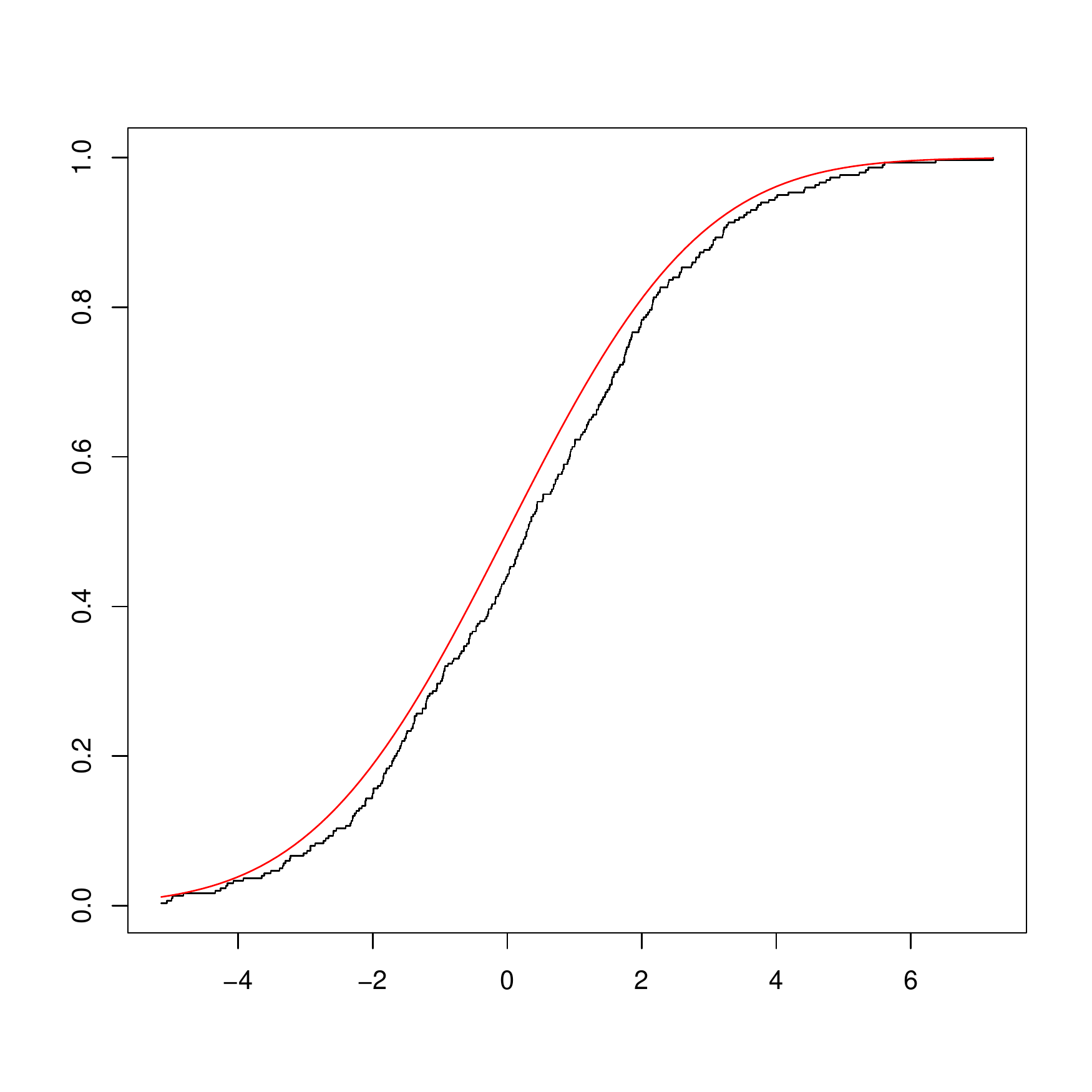}
\includegraphics[width=5cm,pagebox=cropbox,clip]{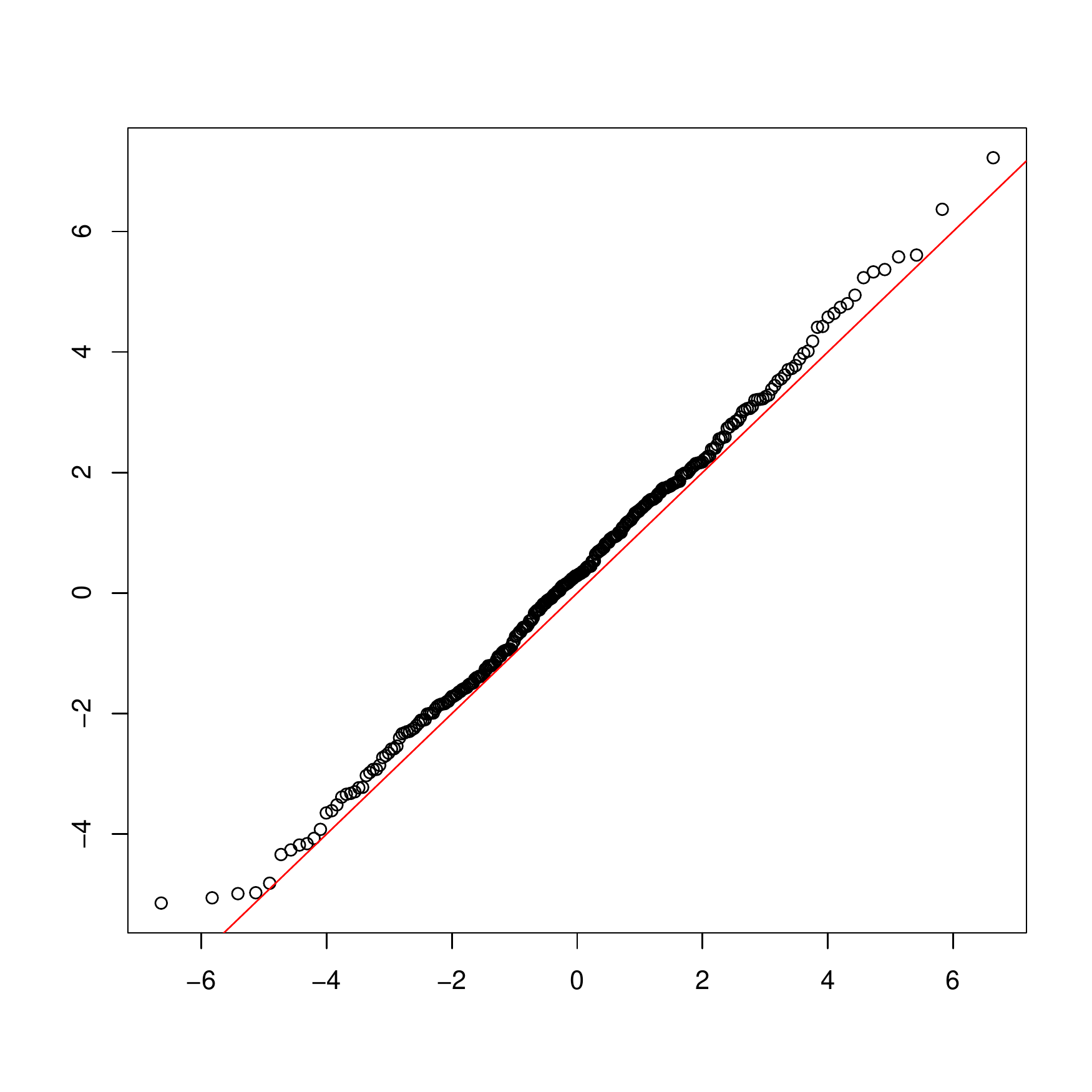}
\includegraphics[width=5cm,pagebox=cropbox,clip]{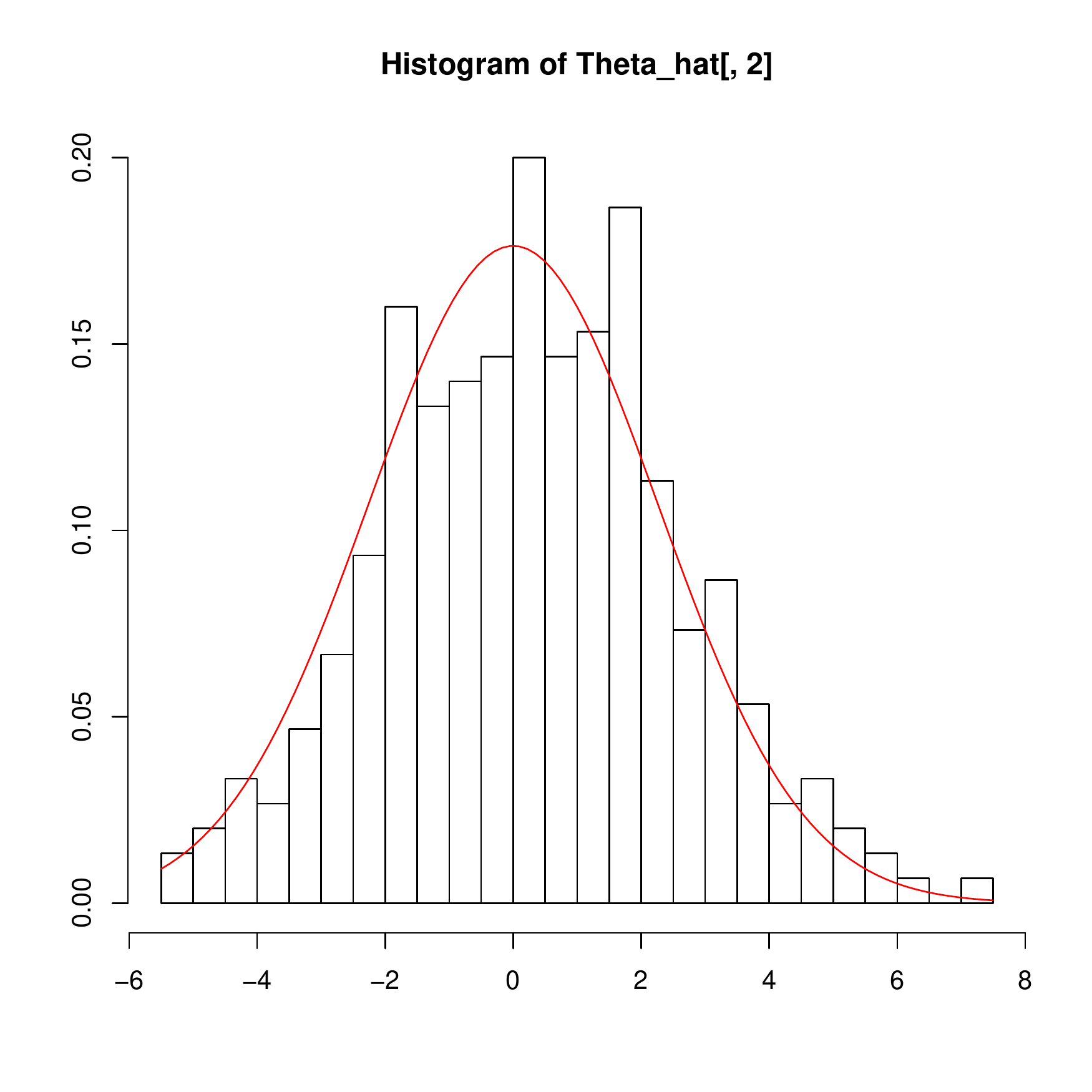}\\
\caption{Simulation results of $\hat{\theta}_2$  \label{fig7}}
\end{center}
\end{figure}


\begin{figure}[t] 
\begin{center}
\includegraphics[width=5cm,pagebox=cropbox,clip]{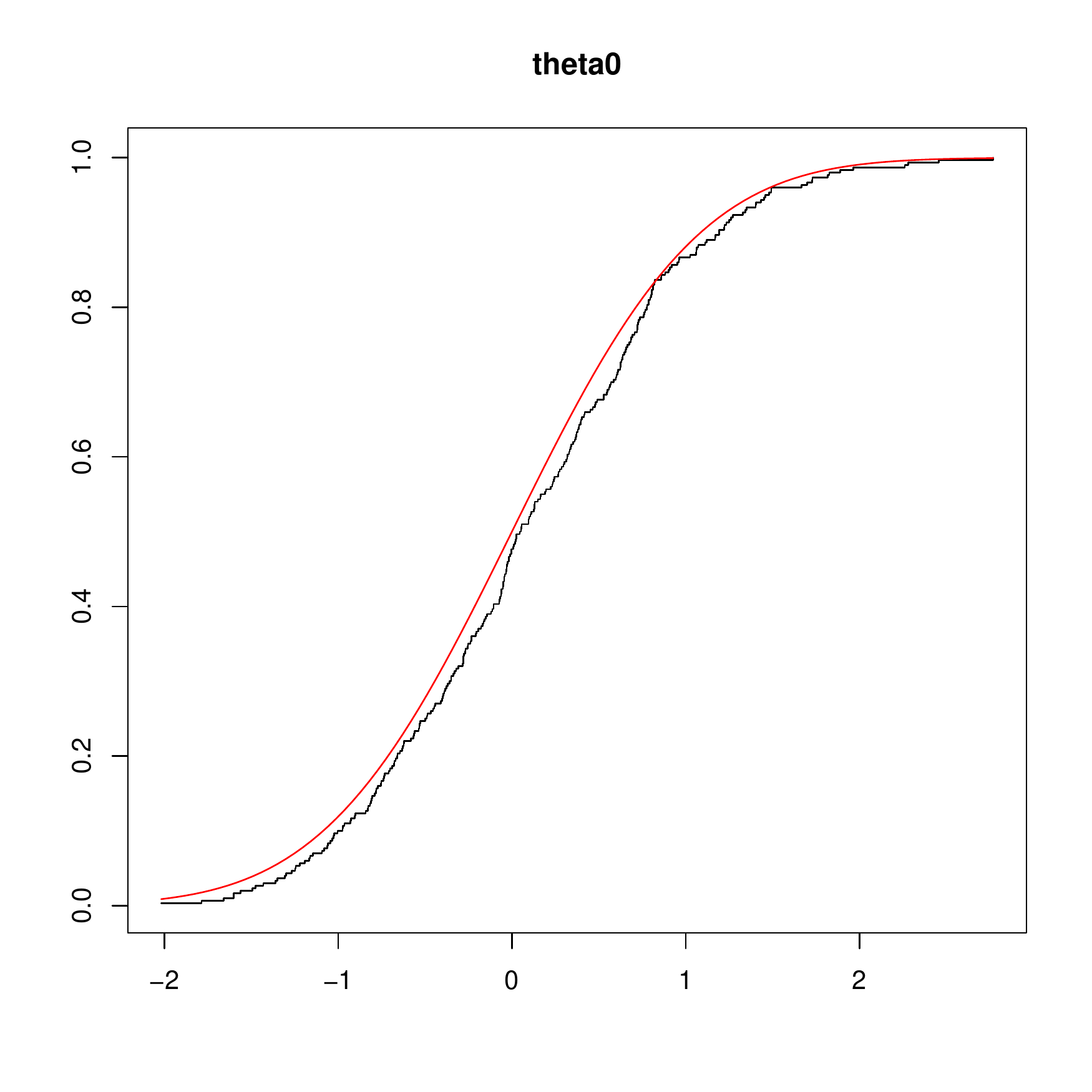}
\includegraphics[width=5cm,pagebox=cropbox,clip]{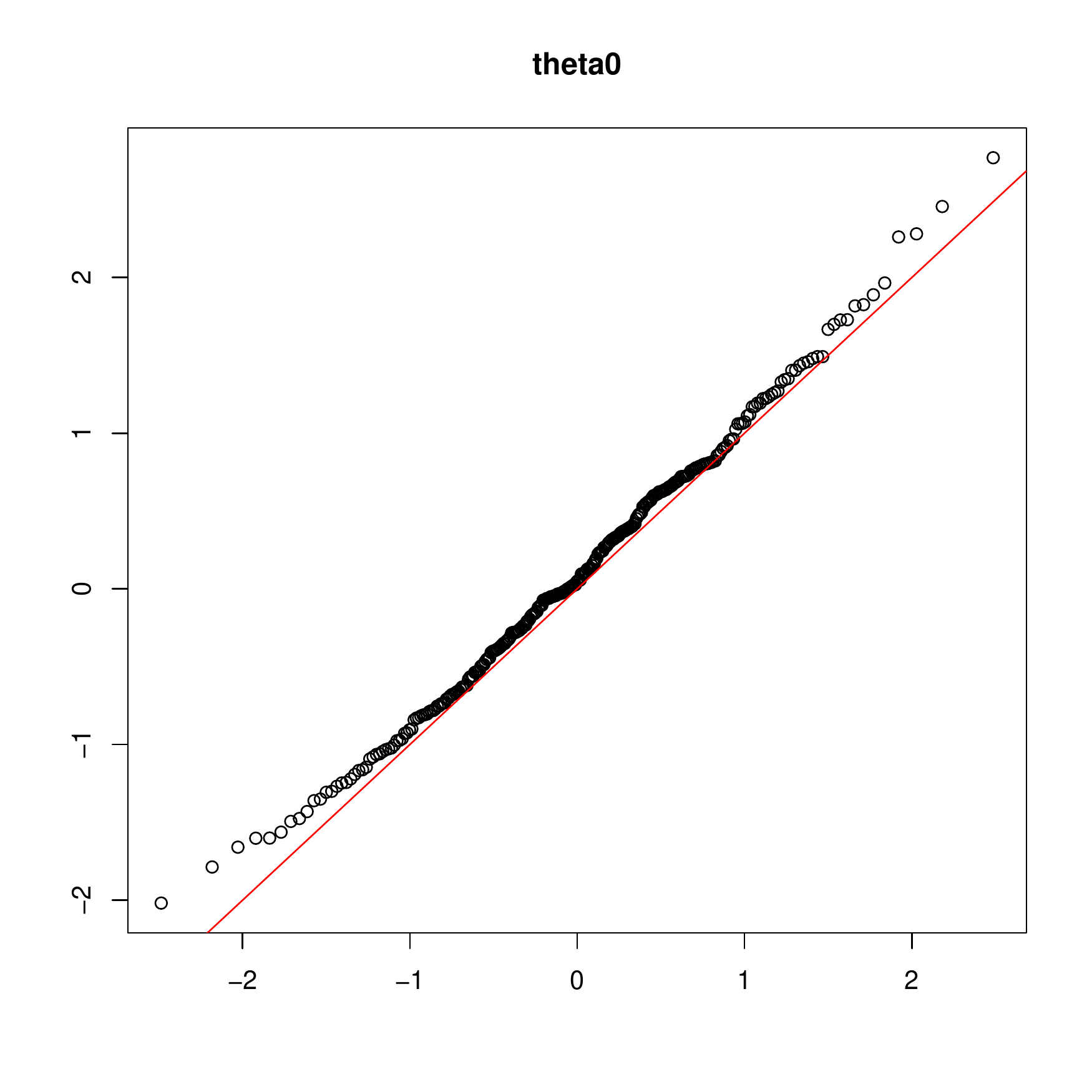}
\includegraphics[width=5cm,pagebox=cropbox,clip]{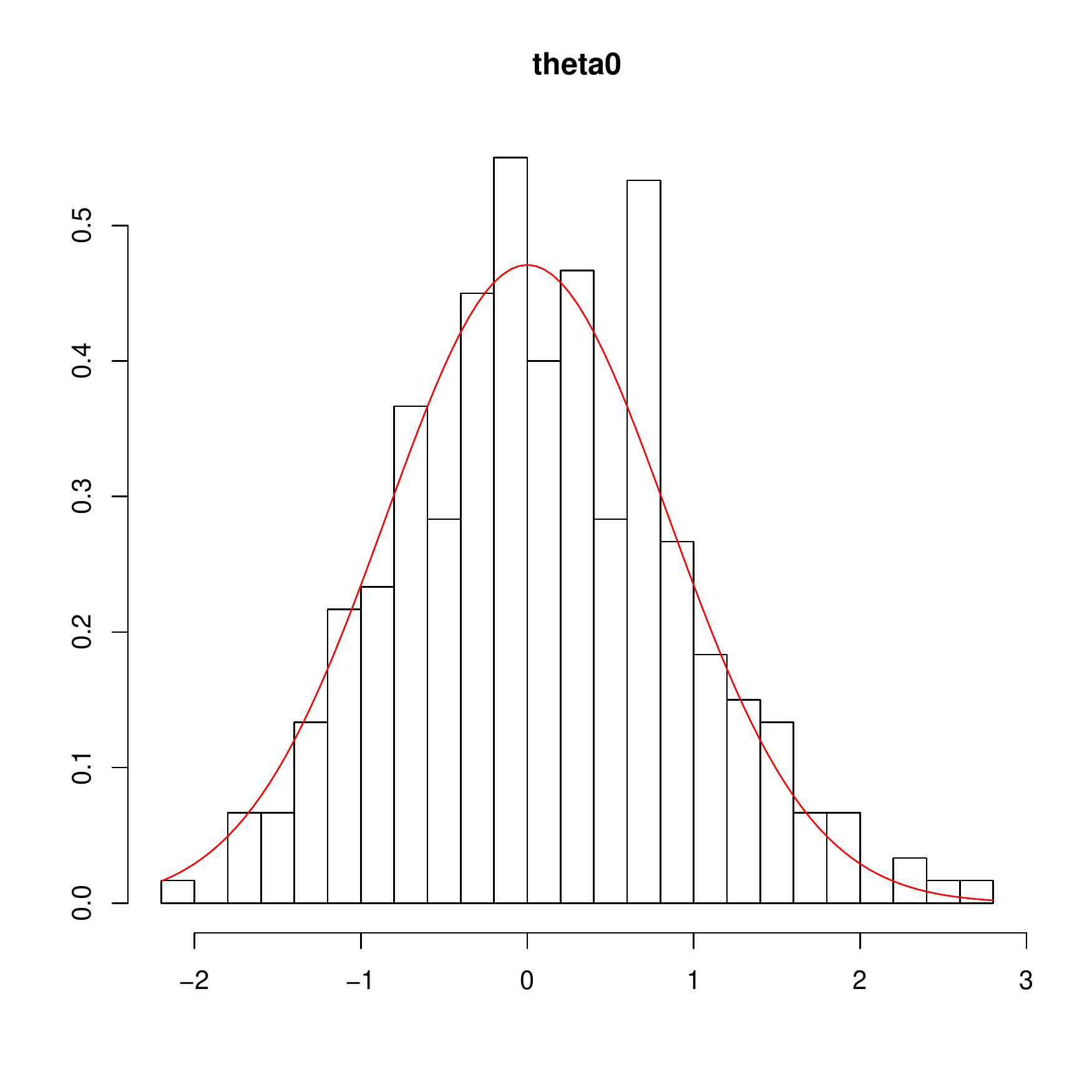}\\
\caption{Simulation results of $\hat{\theta}_0$ \label{fig8}}
\end{center}
\end{figure}

\clearpage

\subsubsection{$\epsilon=0.25$}
Figure 12 is {a} sample path of $X_t(y)$ for $(t,y) \in [0,1]\times [0,1]$
when $(\theta_0^*, \theta_1^*, \theta_2^*, \epsilon) = (0,1,0.2,0.25)$.
Table 2 is the simulation results of {the means and the standard s.d.s of} $\hat{\theta}_1$, $\hat{\theta}_2$  and $\hat{\theta}_0$ with $(N, m, N_2) = (10^4, 99, 500)$.
Figures 13-15 are the simulation results of {the asymptotic distributions of} $\hat{\theta}_1$, $\hat{\theta}_2$  and $\hat{\theta}_0$
with $(N, m, N_2) = (10^4, 99, 500)$.
It seems from Figures 13-15 that 
these estimators have good behaviour.

\begin{figure}[h] 
\begin{center}
\includegraphics[width=9cm]{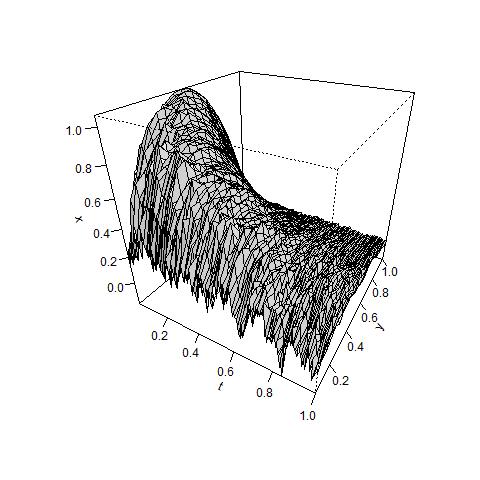} 
\caption{Sample path with $(\theta_0^*, \theta_1^*, \theta_2^*, \epsilon)  = (0,1,0.2, 0.25)$, $\xi(y) = 4.2y(1-y)$}
\end{center}
\end{figure}



\begin{table}[h]
\caption{Simulation results of $\hat{\theta}_1$, $\hat{\theta}_2$  and $\hat{\theta}_0$ with $(N, m, N_2) = (10^4, 99, 500)$ \label{table2}}
\begin{center}
\begin{tabular}{c|ccc} \hline
		&$\hat{\theta}_{1}$&$\hat{\theta}_{2}$&$\hat{\theta}_{0}$
\\ \hline
{true value} &1 & 0.2 & 0
\\ \hline
mean &1.002&  0.201&  0.009
 \\
{s.d.} & (0.007)& (0.002)& (0.188)
 \\   \hline
\end{tabular}
\end{center}
\end{table}


\begin{figure}[h] 
\begin{center}
\includegraphics[width=5cm,pagebox=cropbox,clip]{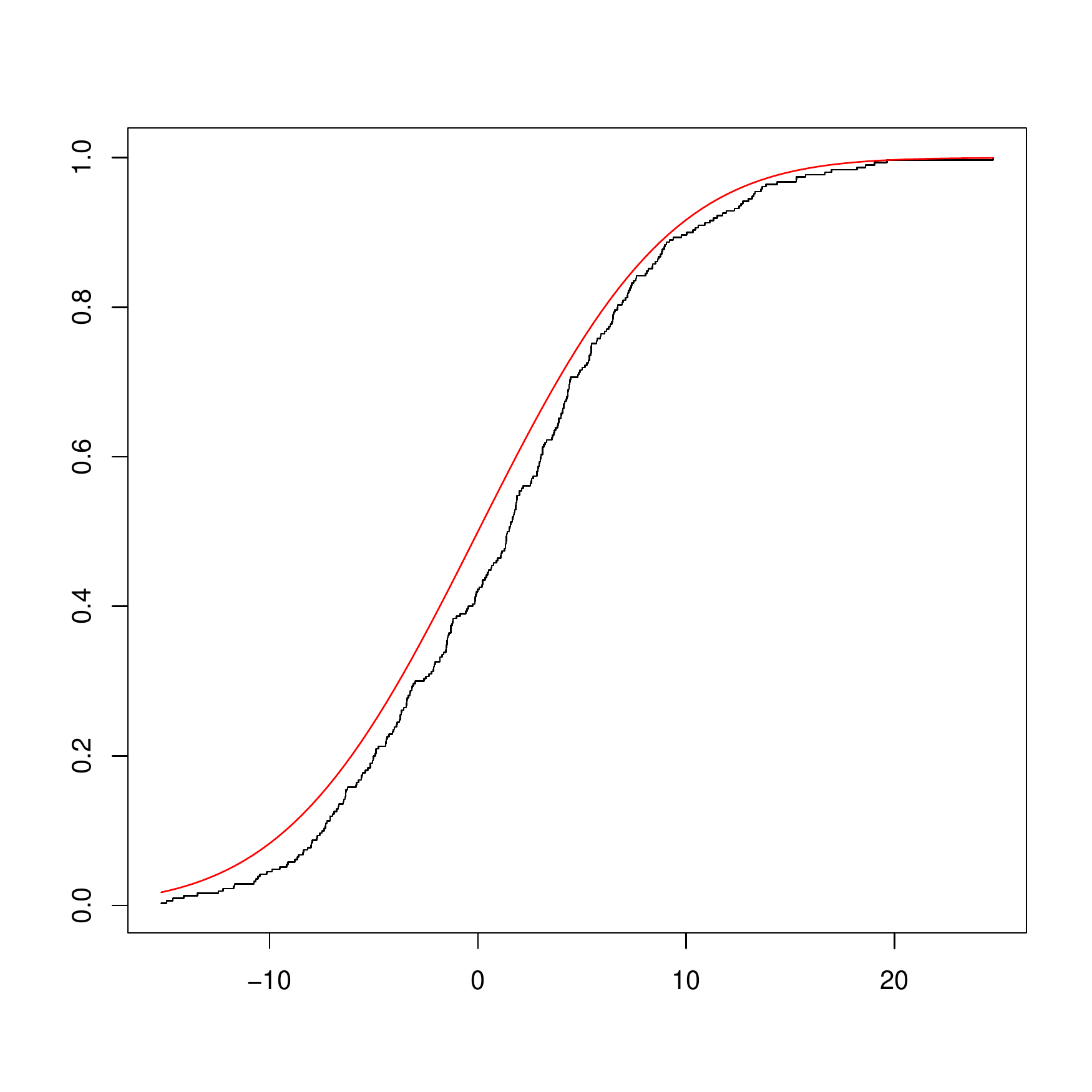}
\includegraphics[width=5cm,pagebox=cropbox,clip]{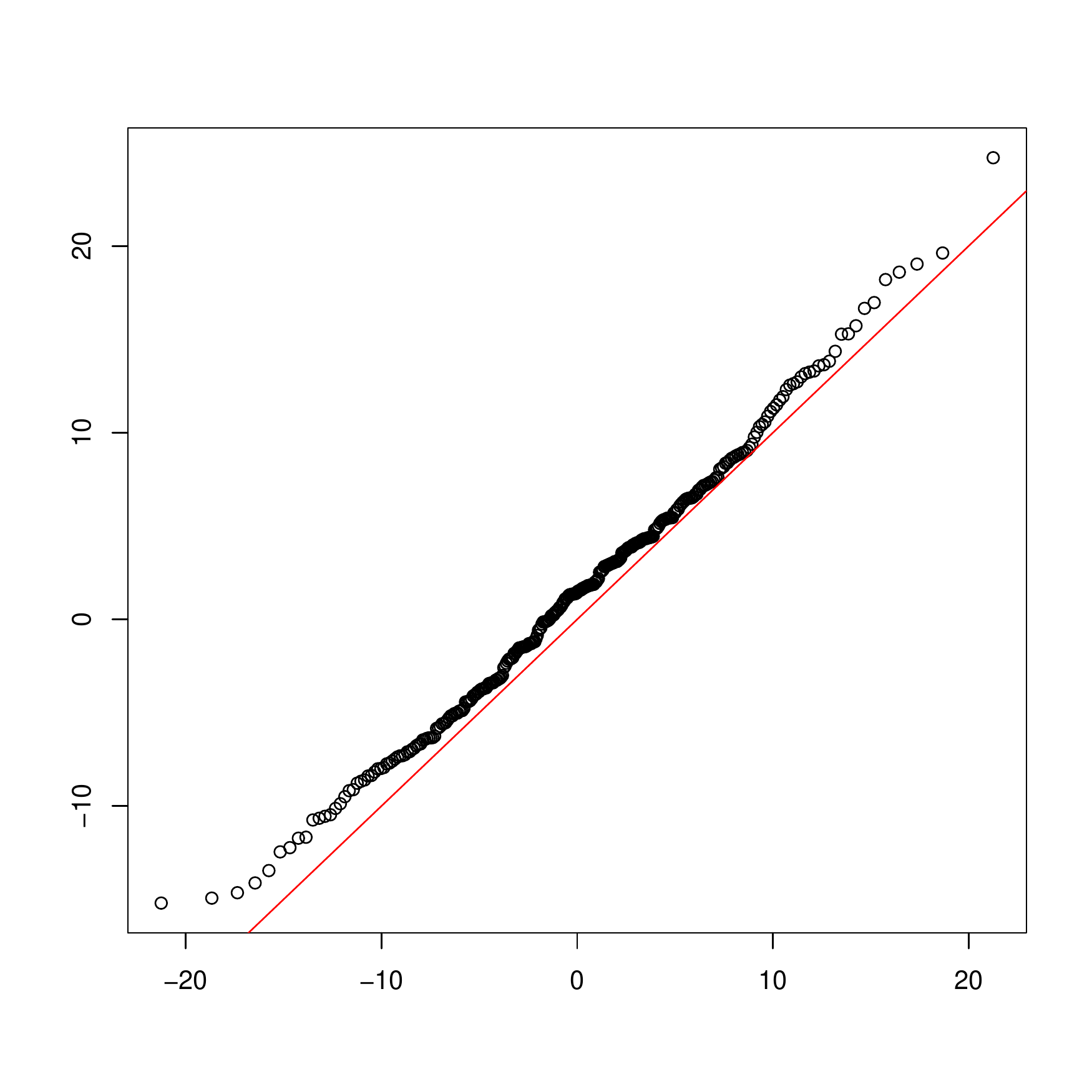}
\includegraphics[width=5cm,pagebox=cropbox,clip]{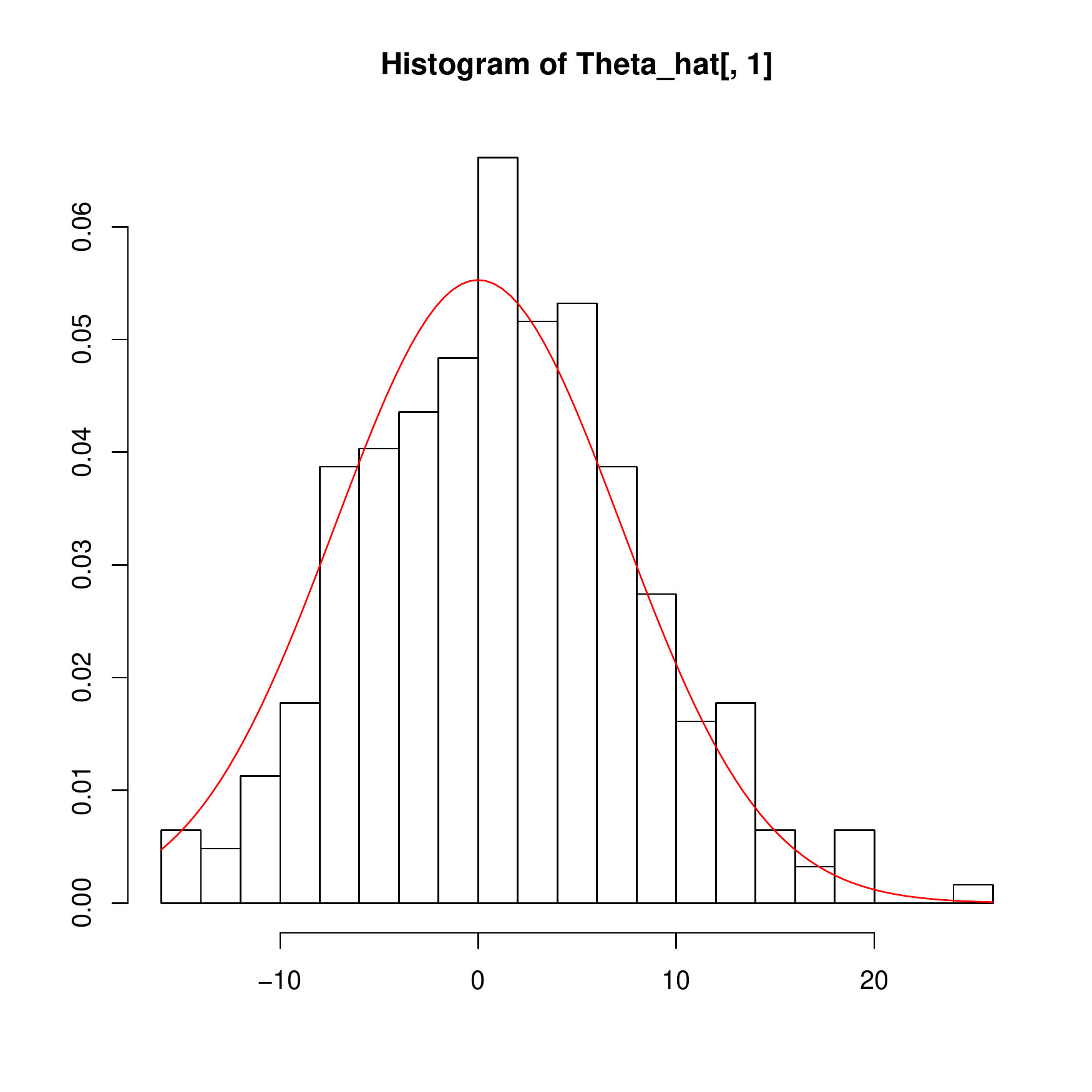}\\
\caption{Simulation results of $\hat{\theta}_1$ \label{fig6}}
\end{center}
\end{figure}

\begin{figure}[h]
\begin{center}
\includegraphics[width=5cm,pagebox=cropbox,clip]{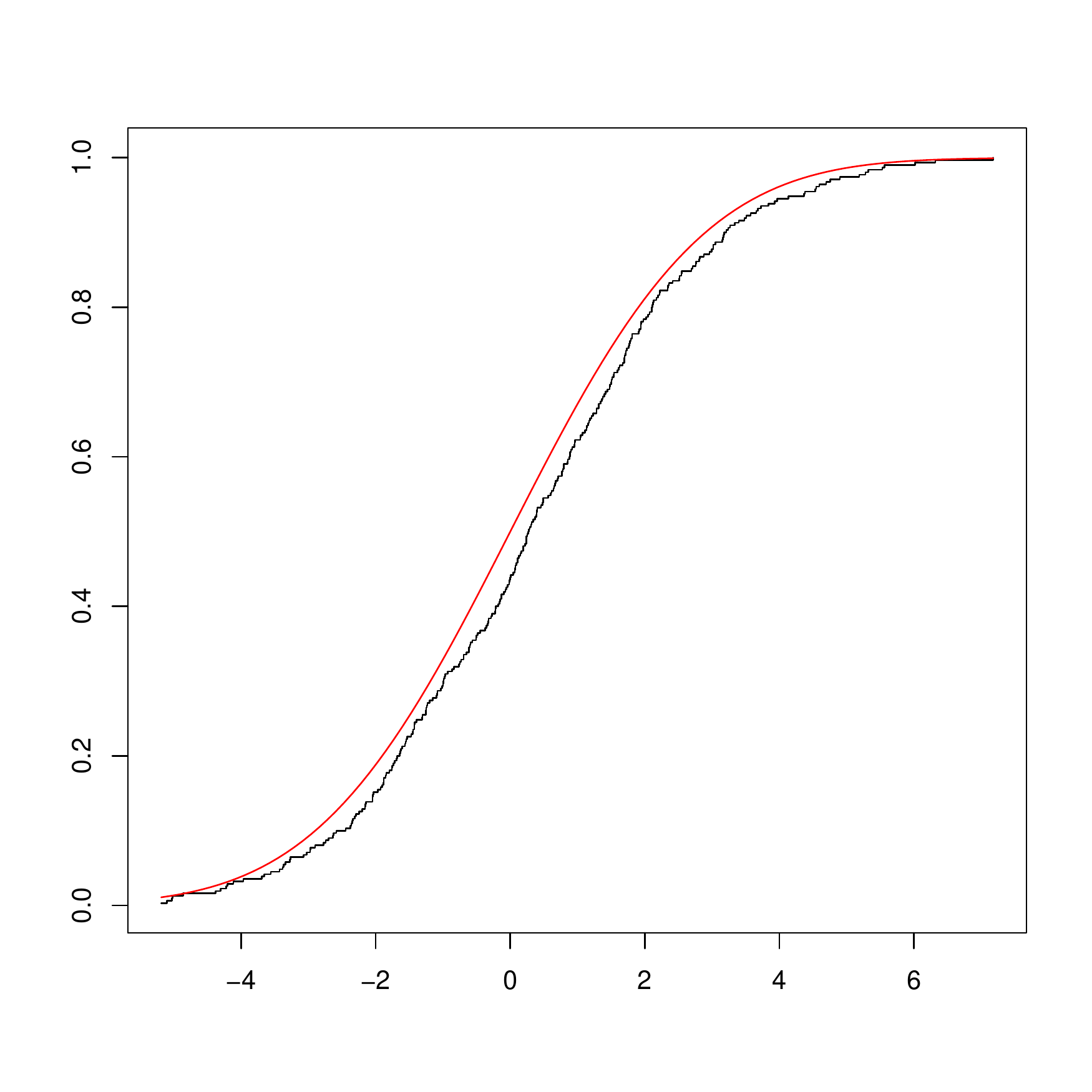}
\includegraphics[width=5cm,pagebox=cropbox,clip]{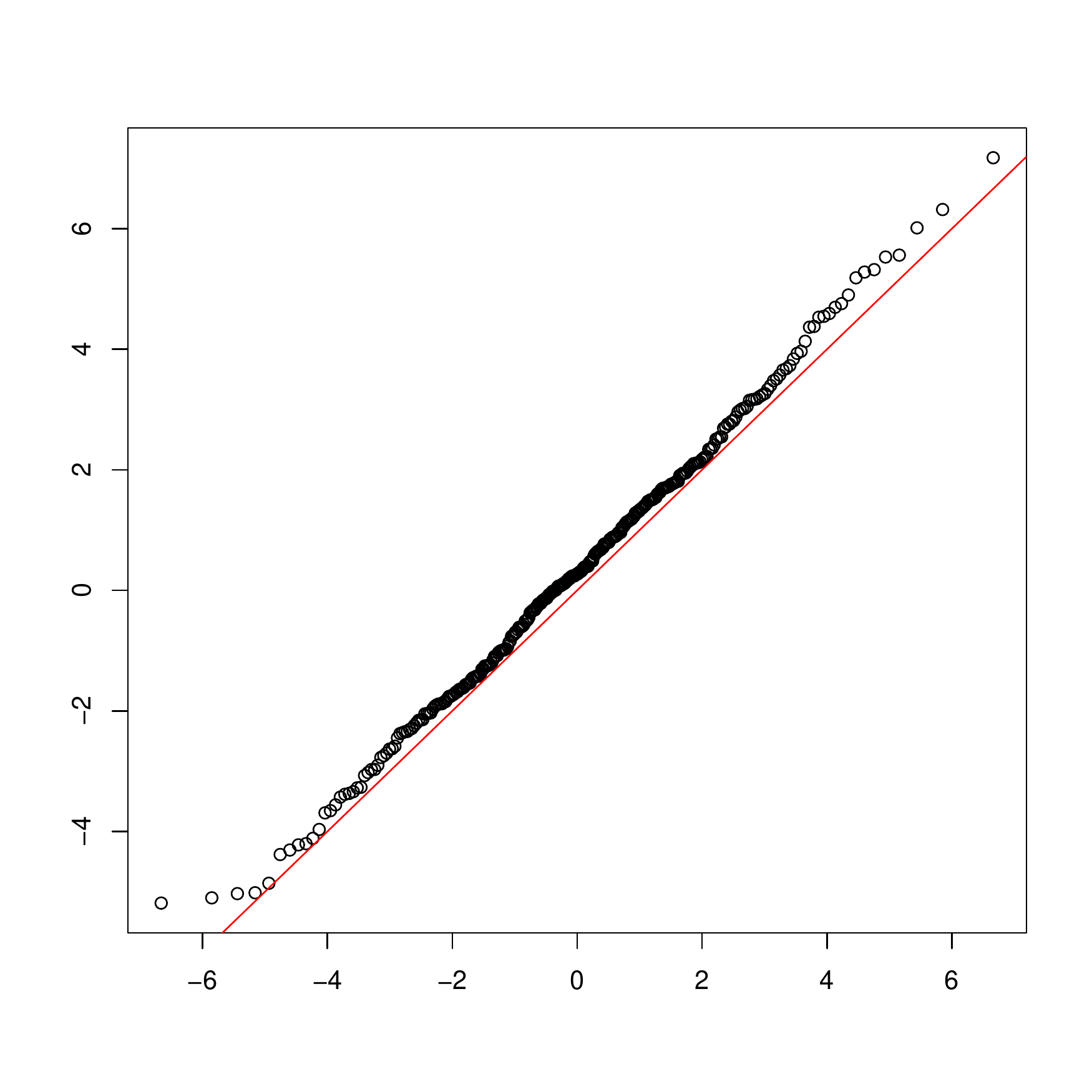}
\includegraphics[width=5cm,pagebox=cropbox,clip]{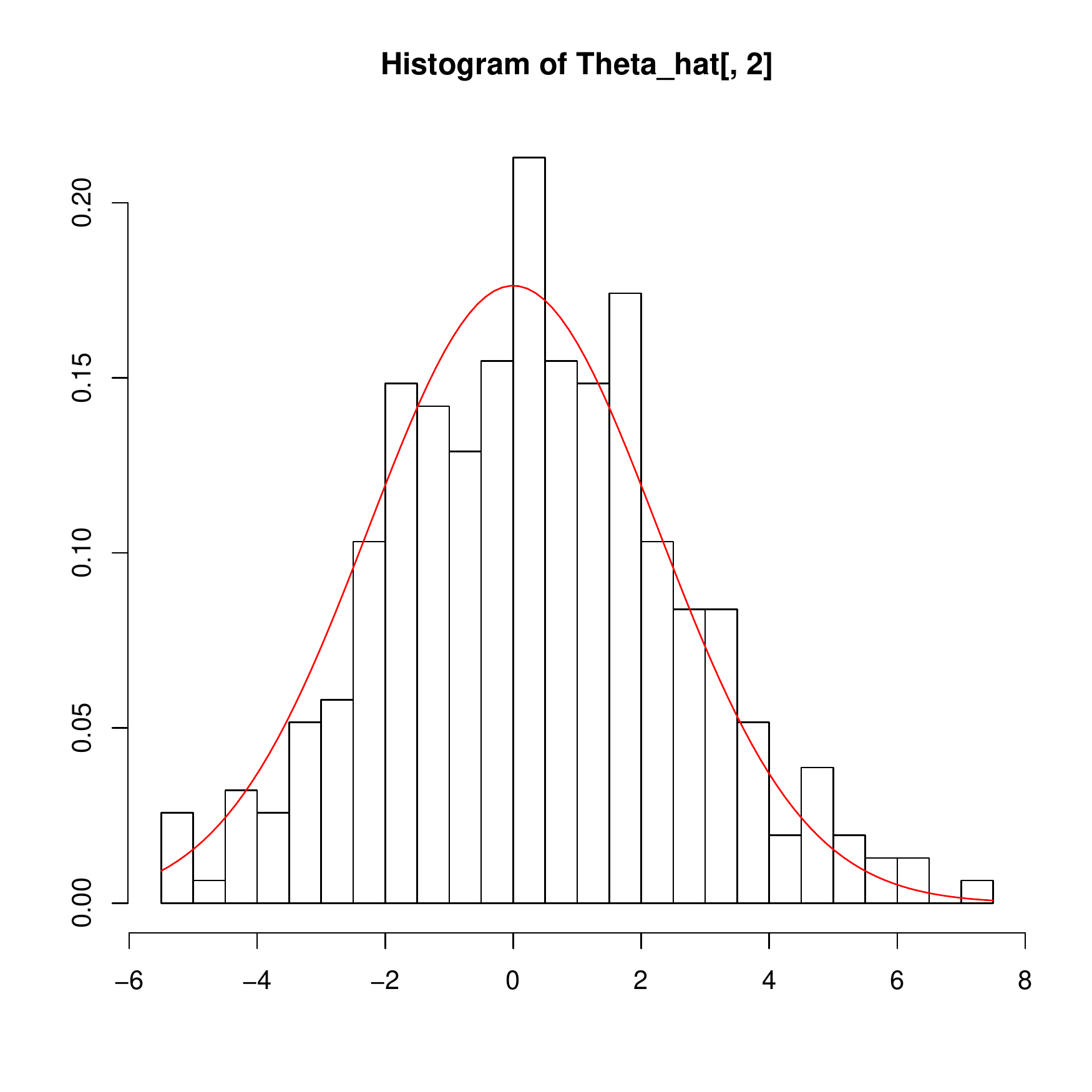}\\
\caption{Simulation results of $\hat{\theta}_2$  \label{fig7}}
\end{center}
\end{figure}


\begin{figure}[t] 
\begin{center}
\includegraphics[width=5cm,pagebox=cropbox,clip]{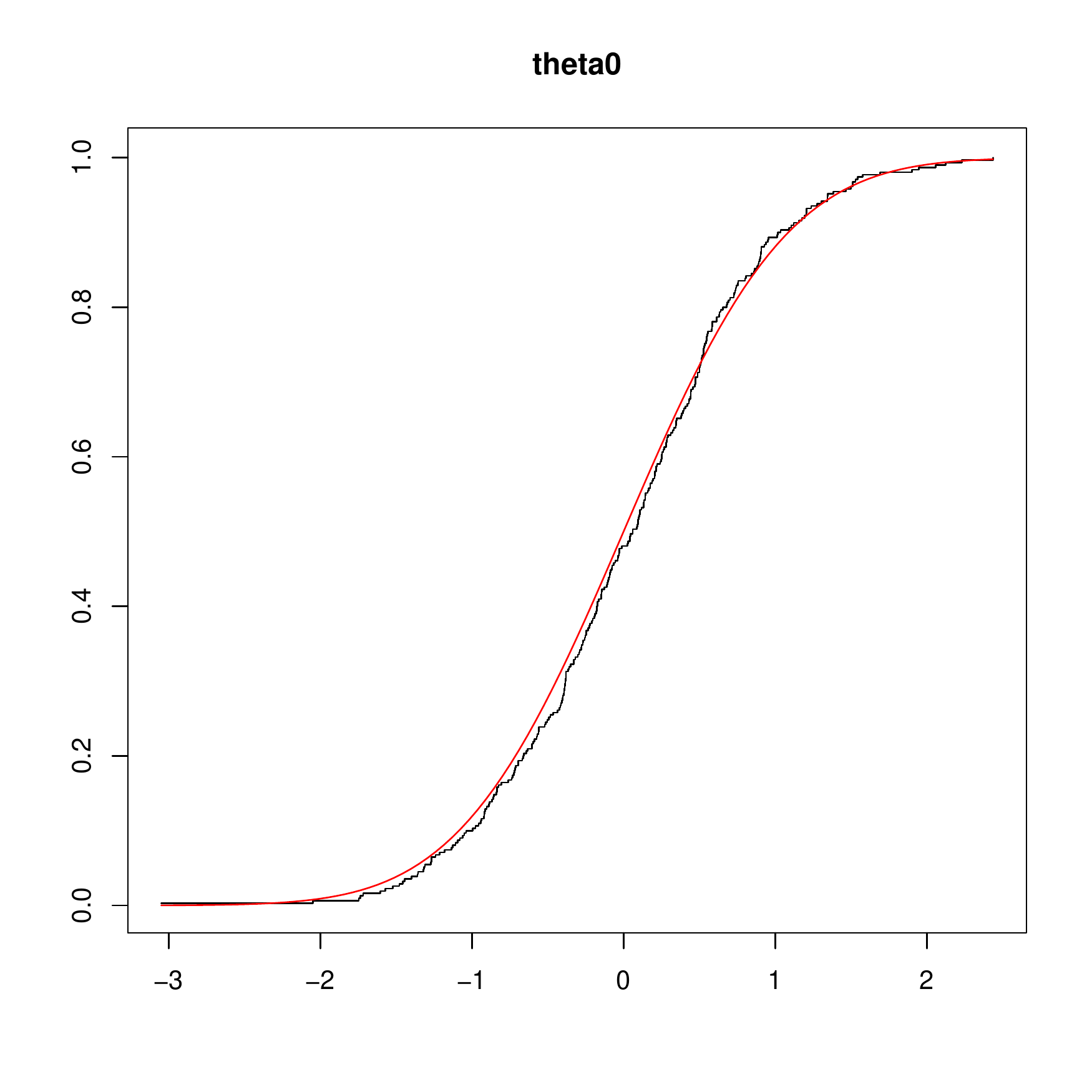}
\includegraphics[width=5cm,pagebox=cropbox,clip]{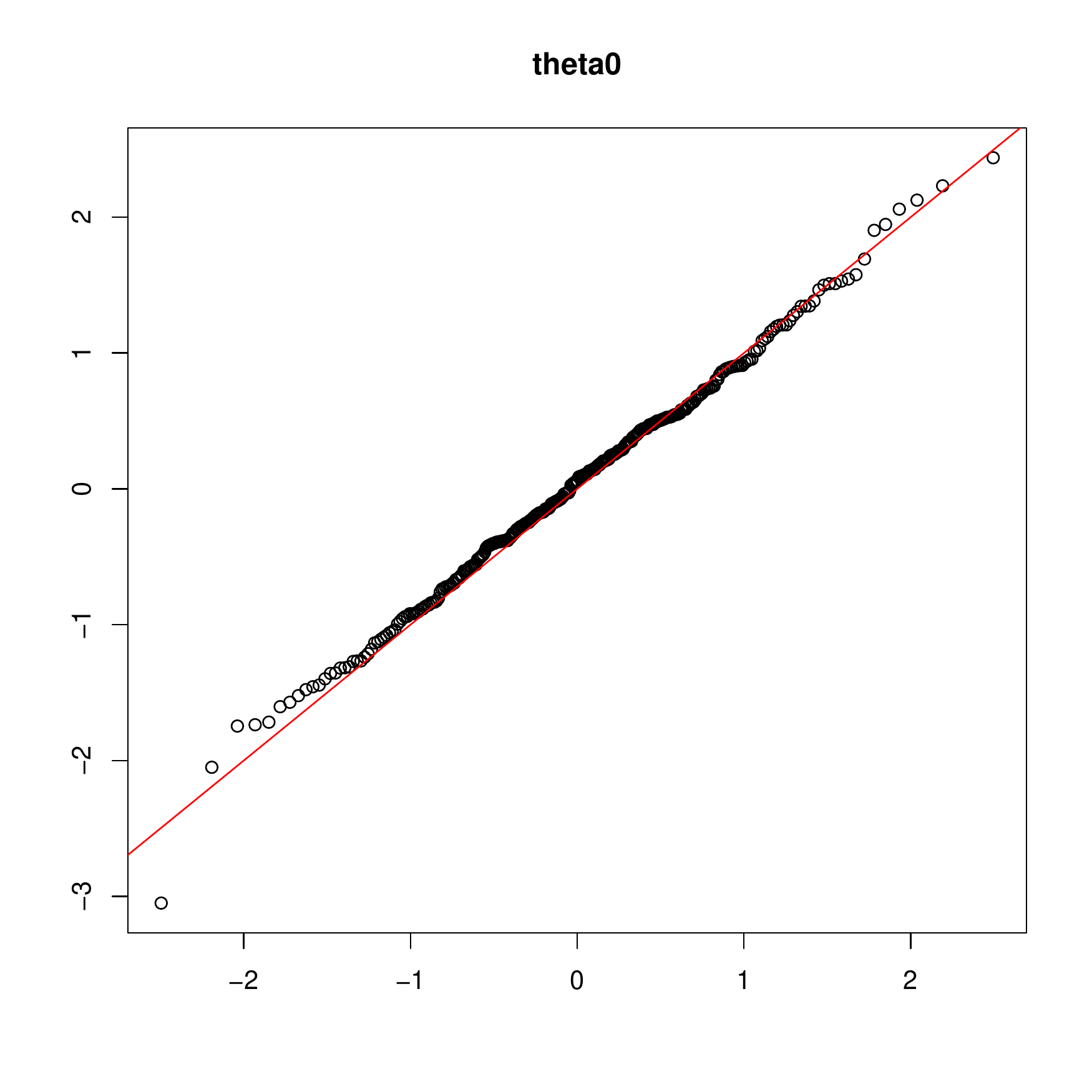}
\includegraphics[width=5cm,pagebox=cropbox,clip]{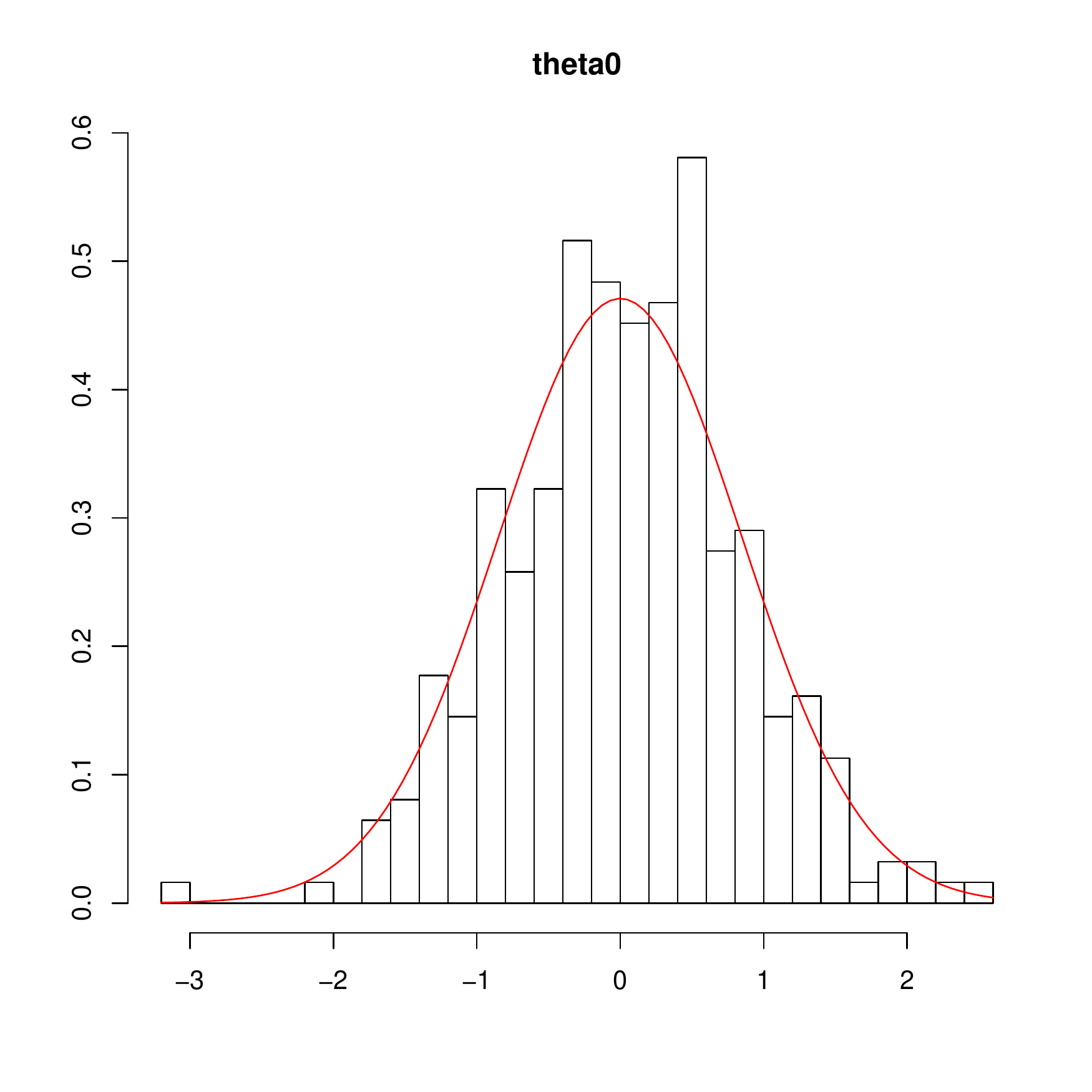}\\
\caption{Simulation results of $\hat{\theta}_0$ \label{fig8}}
\end{center}
\end{figure}


\clearpage

\subsubsection{$\epsilon=0.5$}
Figure 16 is {a} sample path of $X_t(y)$ for $(t,y) \in [0,1]\times [0,1]$
when $(\theta_0^*, \theta_1^*, \theta_2^*, \epsilon) = (0,1,0.2,0.5)$.
Table 3 is the simulation results of {the means and the standard s.d.s of} $\hat{\theta}_1$, $\hat{\theta}_2$  and $\hat{\theta}_0$ with $(N, m, N_2) = (10^4, 99, 500)$.
Figures 17-19 are the simulation results of {the asymptotic distributions of} $\hat{\theta}_1$, $\hat{\theta}_2$  and $\hat{\theta}_0$
with $(N, m, N_2) = (10^4, 99, 500)$.
{
Even if $\epsilon =0.5$, 
we see from Figures 17-19
that 
the estimators stated in Theorem 1 
{have} the asymptotic distribution and 
{they} have good performance.
}

\begin{figure}[h] 
\begin{center}
\includegraphics[width=9cm]{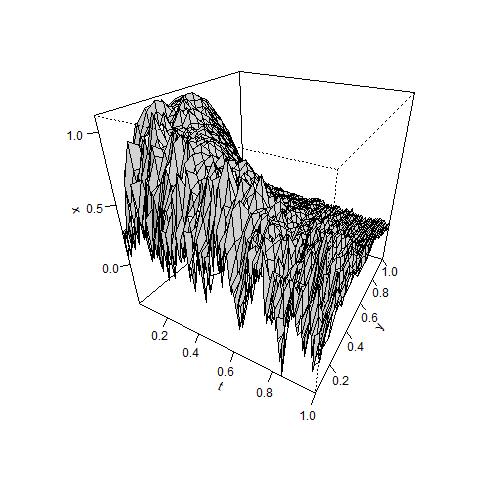} 
\caption{Sample path with $(\theta_0^*, \theta_1^*, \theta_2^*, \epsilon)  = (0,1,0.2,0.5)$, $\xi(y) = 4.2y(1-y)$}
\end{center}
\end{figure}



\begin{table}[h]
\caption{Simulation results of $\hat{\theta}_1$, $\hat{\theta}_2$  and $\hat{\theta}_0$ with $(N, m, N_2) = (10^4, 99, 500)$ \label{table2}}
\begin{center}
\begin{tabular}{c|ccc} \hline
		&$\hat{\theta}_{1}$&$\hat{\theta}_{2}$&$\hat{\theta}_{0}$
\\ \hline
{true value} &1 & 0.2 & 0
\\ \hline
mean &1.002&  0.201&  -0.013
 \\
{s.d.} & (0.007)& (0.002)& (0.367)
 \\   \hline
\end{tabular}
\end{center}
\end{table}


\begin{figure}[h] 
\begin{center}
\includegraphics[width=5cm,pagebox=cropbox,clip]{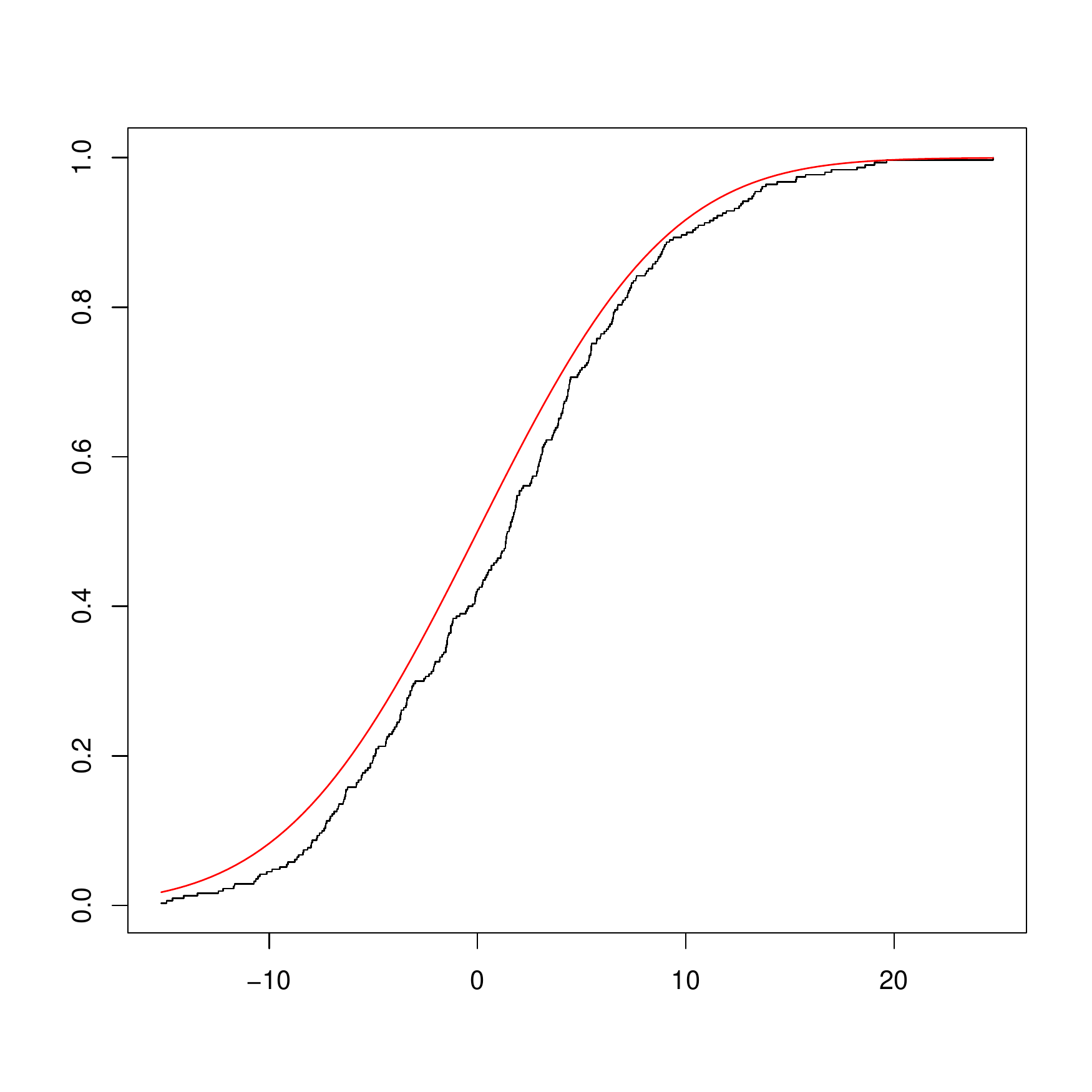}
\includegraphics[width=5cm,pagebox=cropbox,clip]{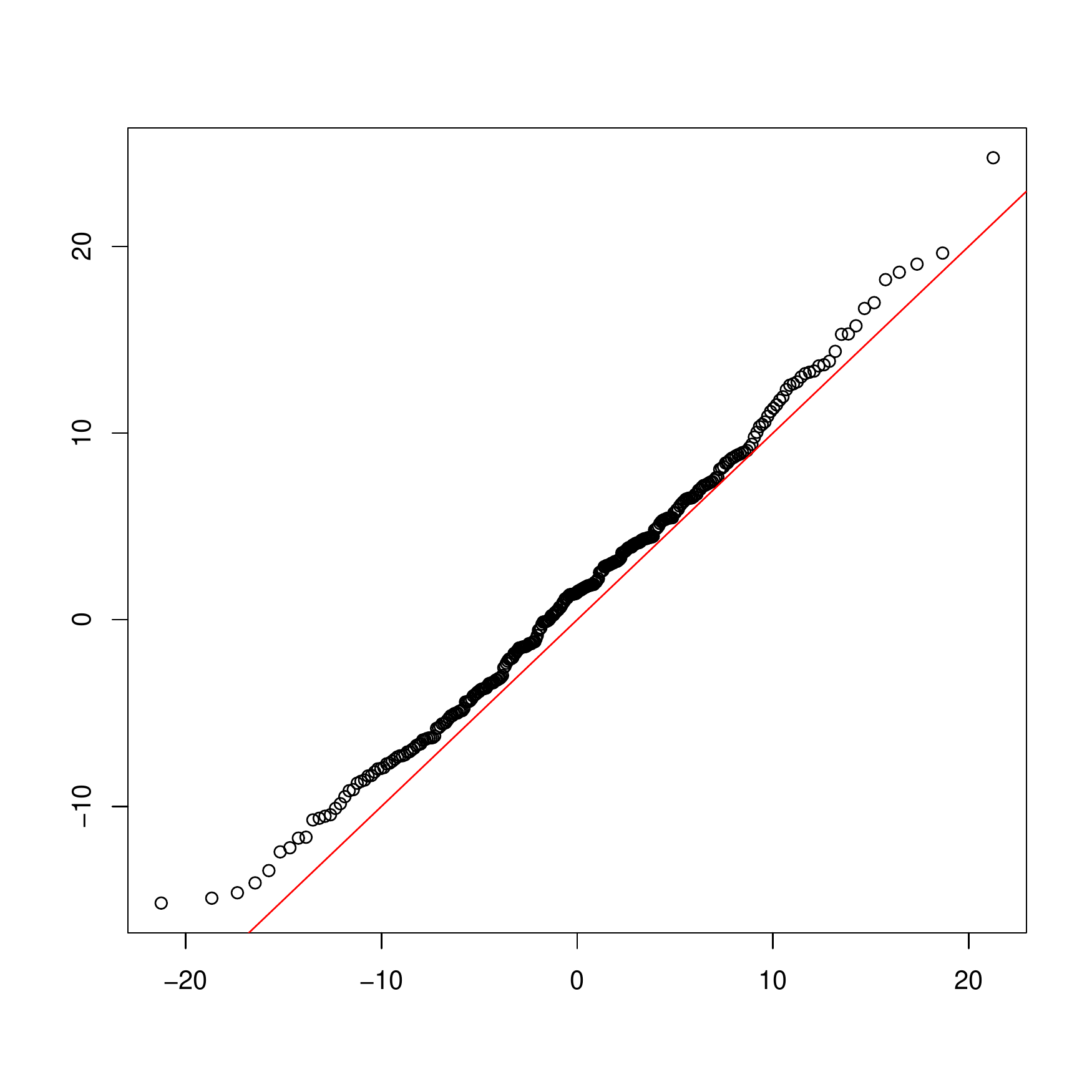}
\includegraphics[width=5cm,pagebox=cropbox,clip]{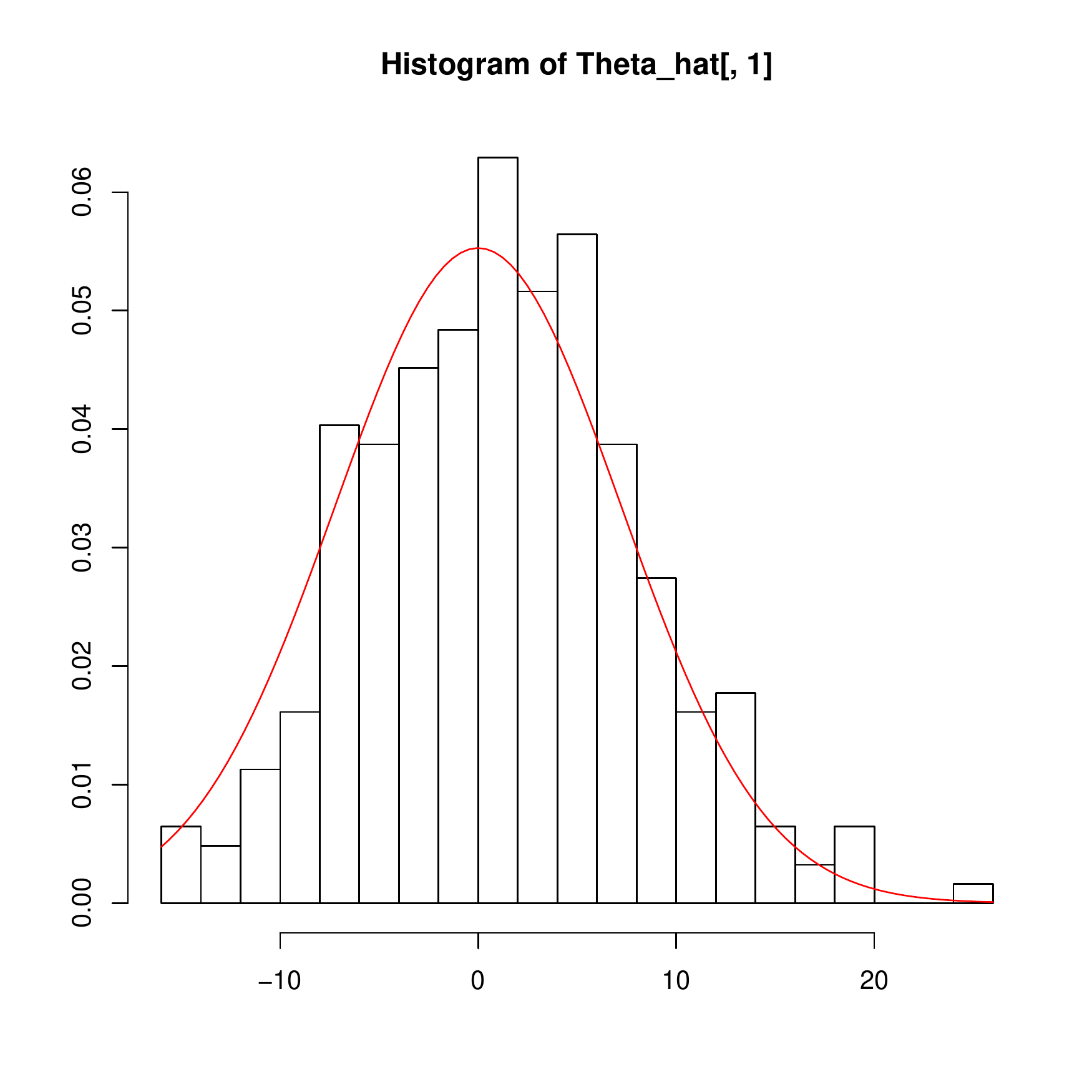}\\
\caption{Simulation results of $\hat{\theta}_1$ \label{fig6}}
\end{center}
\end{figure}

\begin{figure}[h]
\begin{center}
\includegraphics[width=5cm,pagebox=cropbox,clip]{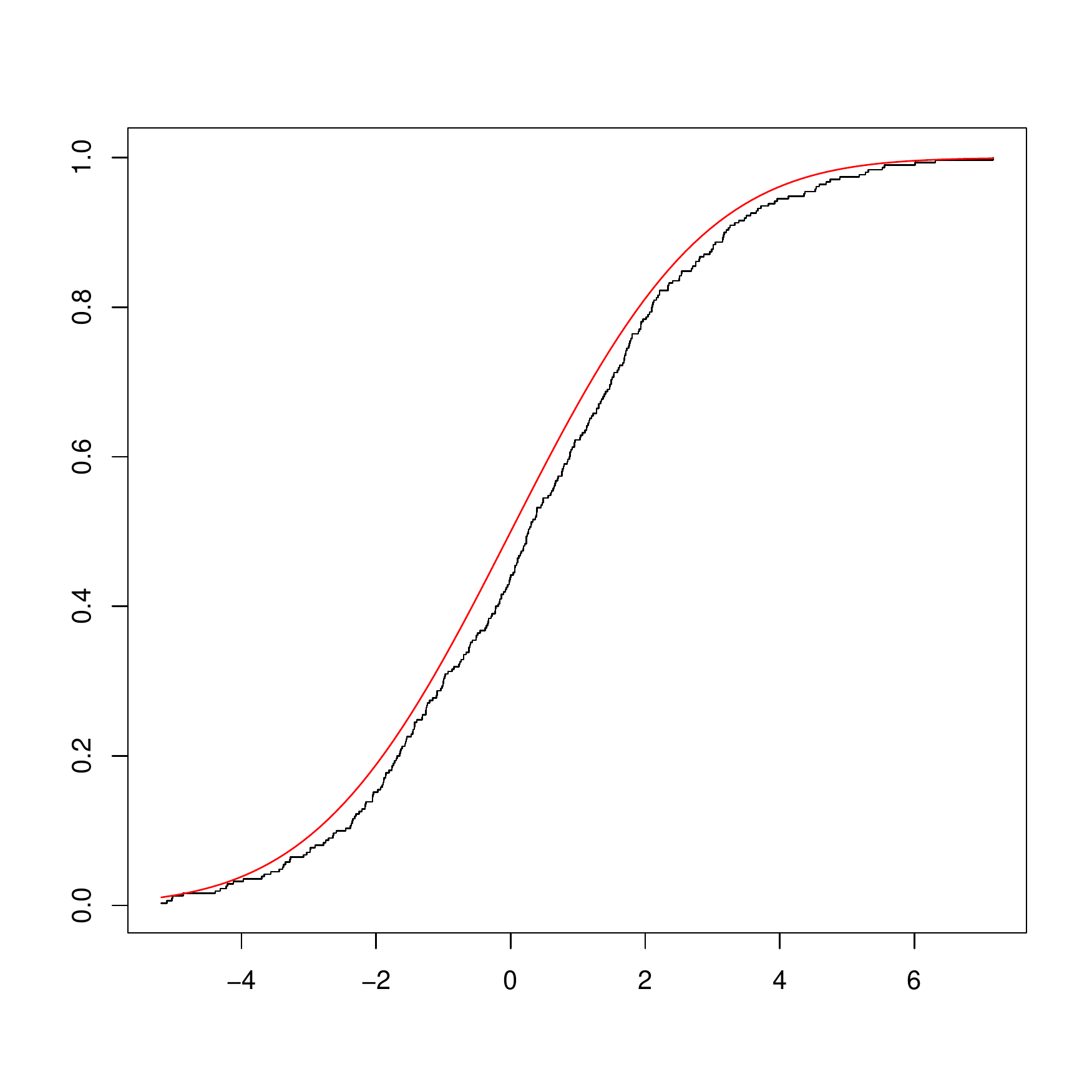}
\includegraphics[width=5cm,pagebox=cropbox,clip]{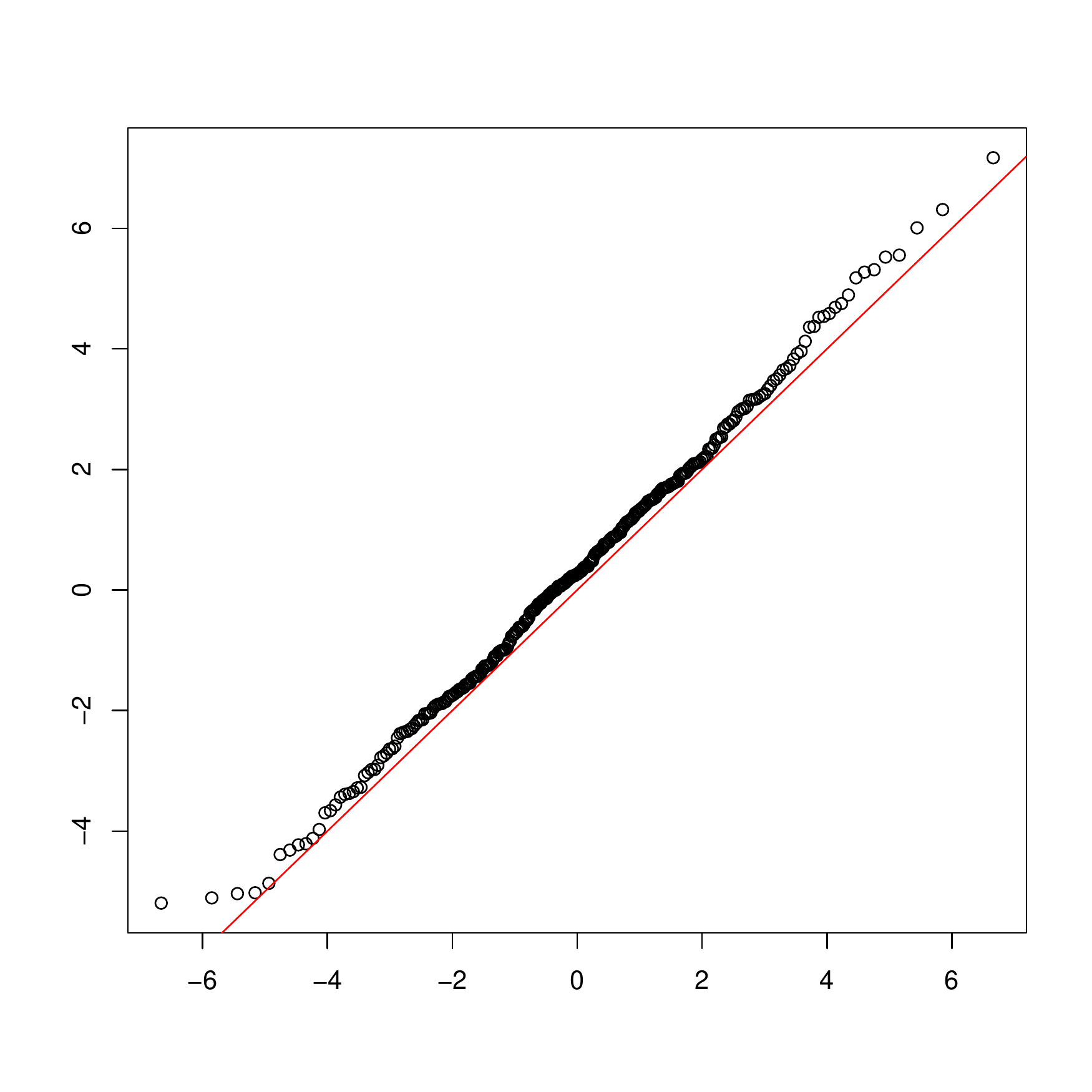}
\includegraphics[width=5cm,pagebox=cropbox,clip]{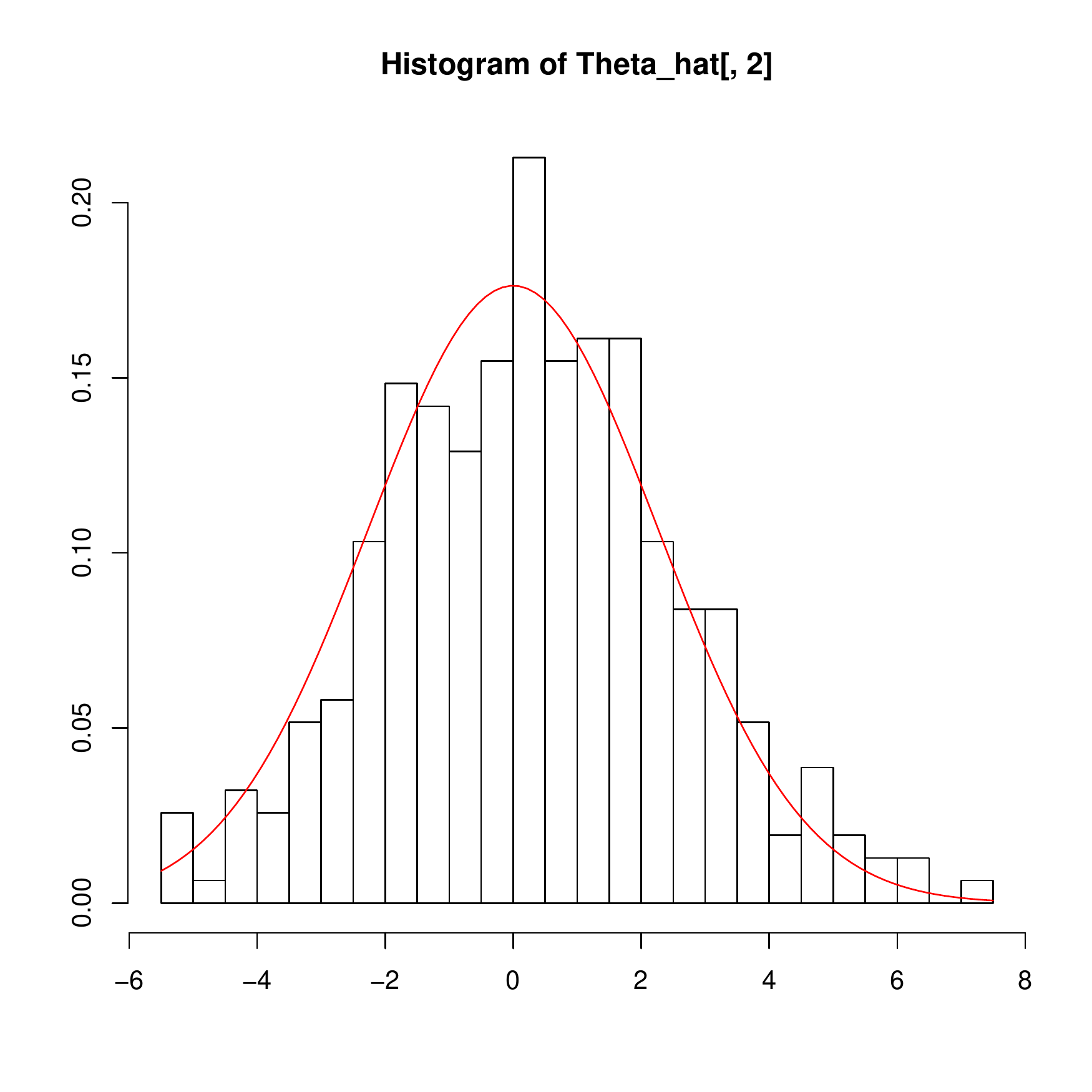}\\
\caption{Simulation results of $\hat{\theta}_2$  \label{fig7}}
\end{center}
\end{figure}


\begin{figure}[t] 
\begin{center}
\includegraphics[width=5cm,pagebox=cropbox,clip]{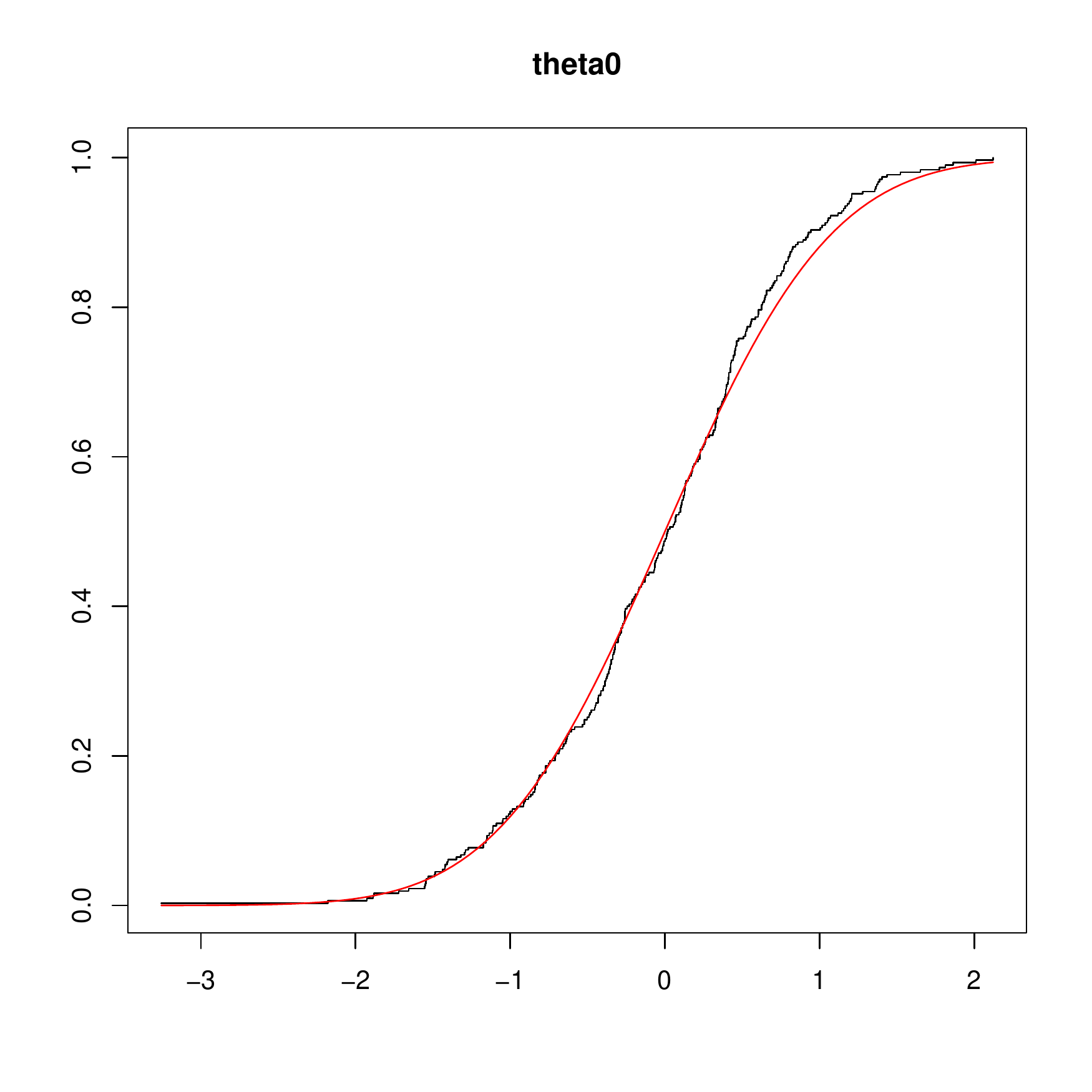}
\includegraphics[width=5cm,pagebox=cropbox,clip]{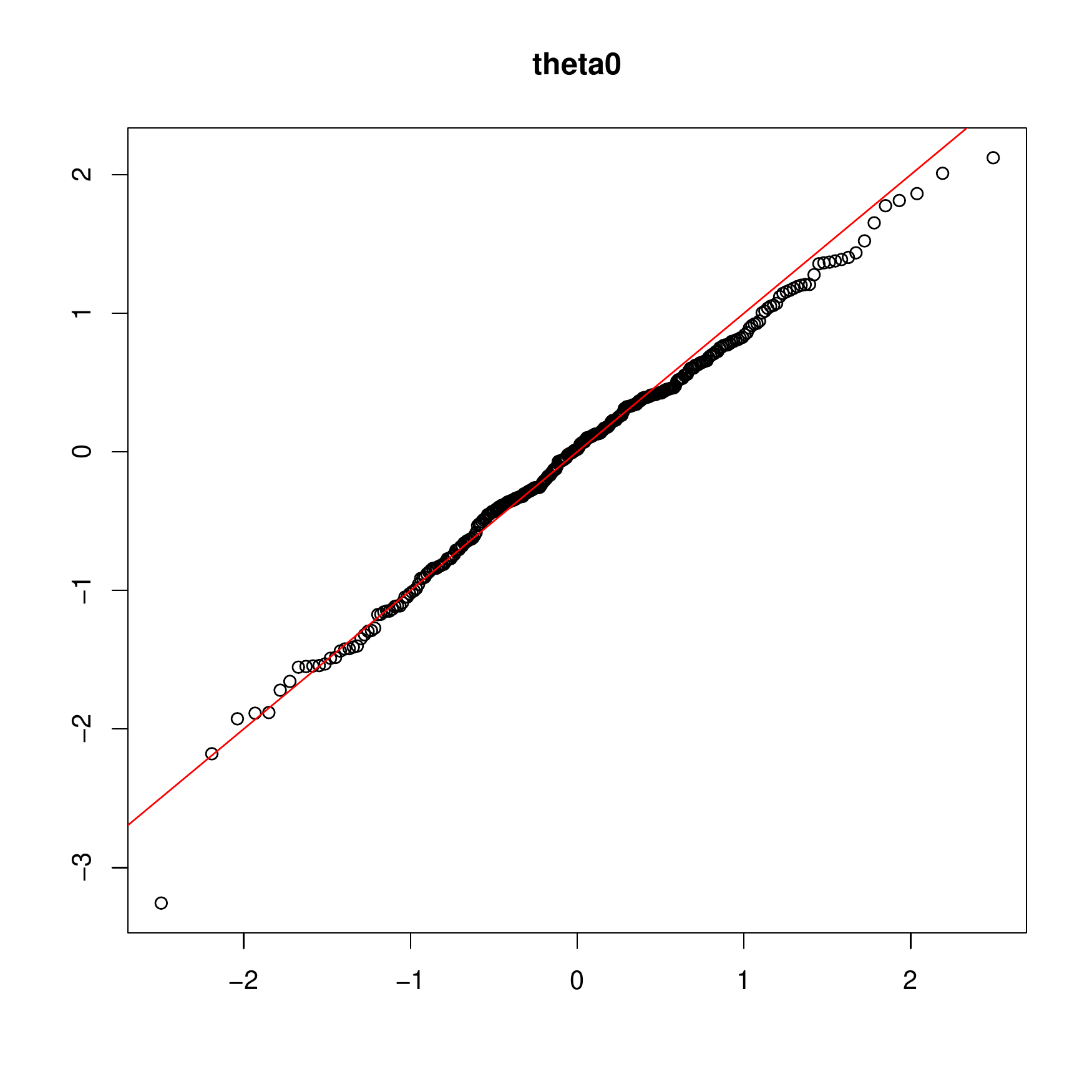}
\includegraphics[width=5cm,pagebox=cropbox,clip]{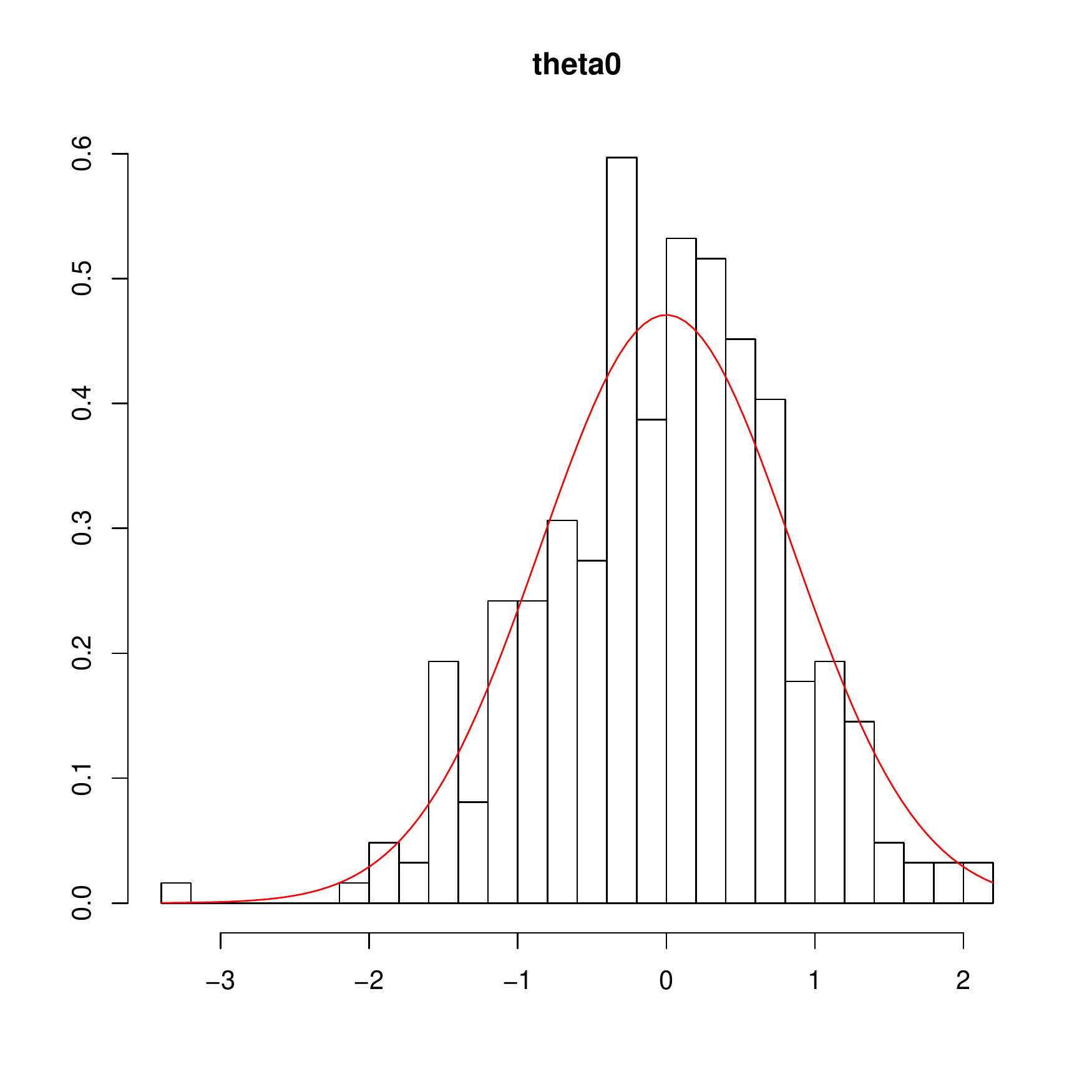}\\
\caption{Simulation results of $\hat{\theta}_0$ \label{fig8}}
\end{center}
\end{figure}


\clearpage

\subsubsection{$\epsilon=0.75$}
Figure 20 is {a} sample path of $X_t(y)$ for $(t,y) \in [0,1]\times [0,1]$
when $(\theta_0^*, \theta_1^*, \theta_2^*, \epsilon) = (0,1,0.2,0.75)$.
Table 4 is the simulation results of {the means and the standard s.d.s of} $\hat{\theta}_1$, $\hat{\theta}_2$  and $\hat{\theta}_0$ with $(N, m, N_2) = (10^4, 99, 500)$.
{Figures 21-23 are the simulation results of {the asymptotic distributions of} $\hat{\theta}_1$, $\hat{\theta}_2$  and $\hat{\theta}_0$
with $(N, m, N_2) = (10^4, 99, 500)$.}
{
From Figures 21-22, 
we can see that 
the estimator of $(\theta_1, \theta_2)$ 
has  the asymptotic distribution in Theorem 1.}
{
Although it seems from Figure 23 that
the estimator of $\theta_0$ has good performance,
the deviation from the red line is larger than that in Figure 19.}

\begin{figure}[h] 
\begin{center}
\includegraphics[width=9cm]{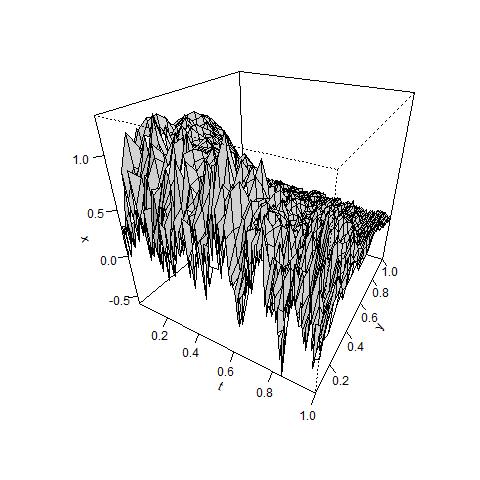} 
\caption{Sample path with $(\theta_0^*, \theta_1^*, \theta_2^*, \epsilon)  = (0,1,0.2,0.75)$, $\xi(y) = 4.2y(1-y)$}
\end{center}
\end{figure}



\begin{table}[h]
\caption{Simulation results of $\hat{\theta}_1$, $\hat{\theta}_2$  and $\hat{\theta}_0$ with $(N, m, N_2) = (10^4, 99, 500)$ \label{table2}}
\begin{center}
\begin{tabular}{c|ccc} \hline
		&$\hat{\theta}_{1}$&$\hat{\theta}_{2}$&$\hat{\theta}_{0}$
\\ \hline
{true value} &1 & 0.2 & 0
\\ \hline
mean &1.002&  0.201&  -0.058
 \\
{s.d.} & (0.007)& (0.002)& (0.548)
 \\   \hline
\end{tabular}
\end{center}
\end{table}


\begin{figure}[h] 
\begin{center}
\includegraphics[width=5cm,pagebox=cropbox,clip]{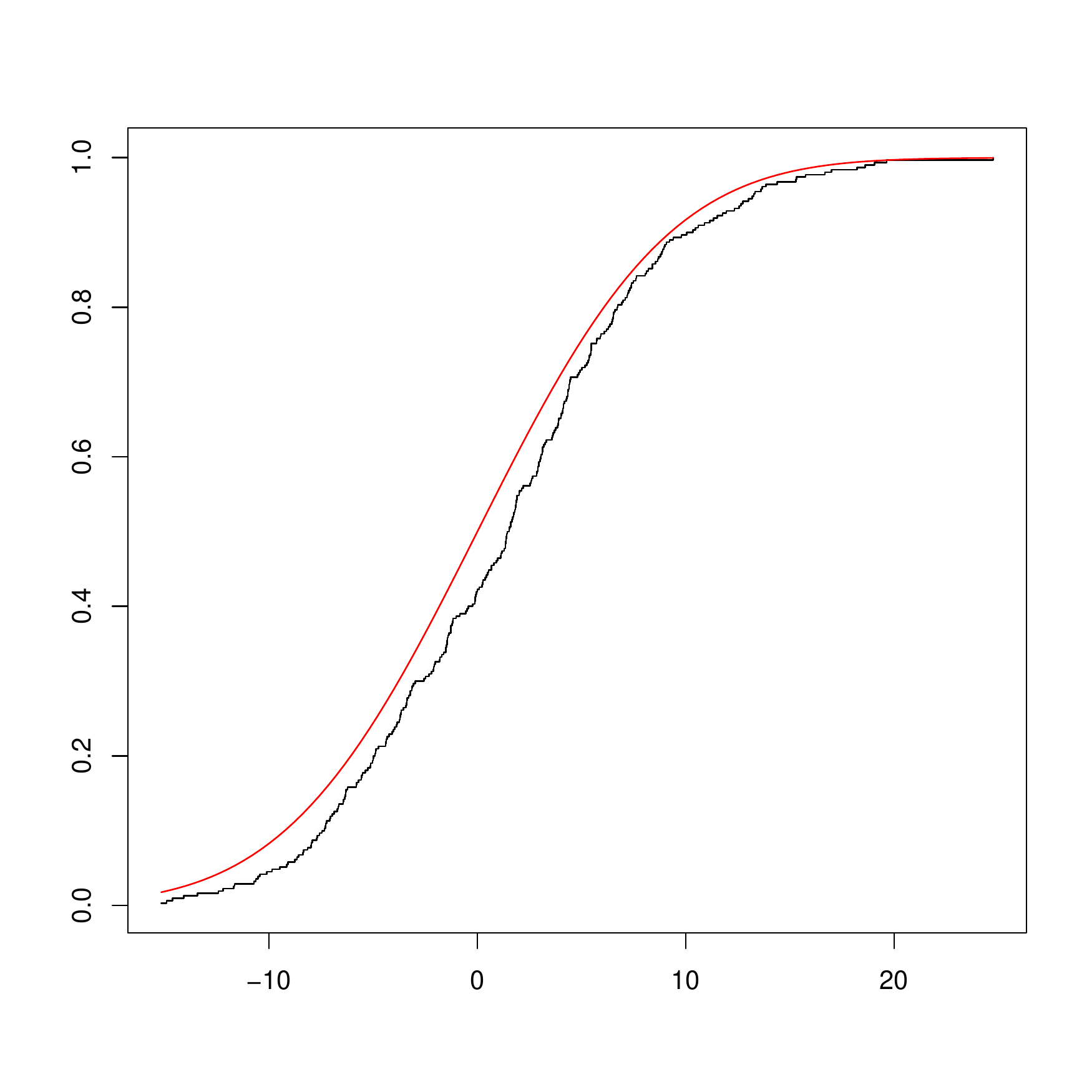}
\includegraphics[width=5cm,pagebox=cropbox,clip]{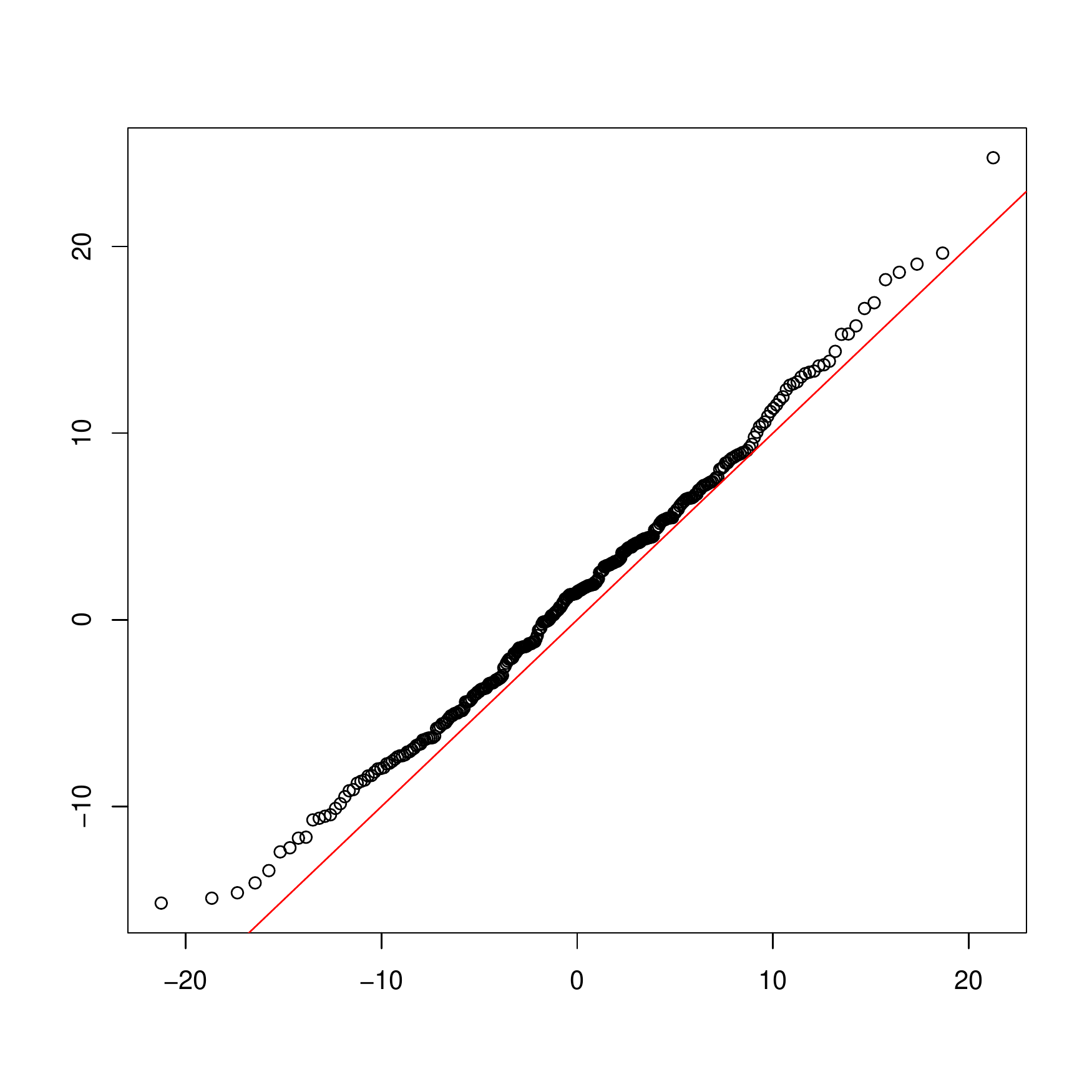}
\includegraphics[width=5cm,pagebox=cropbox,clip]{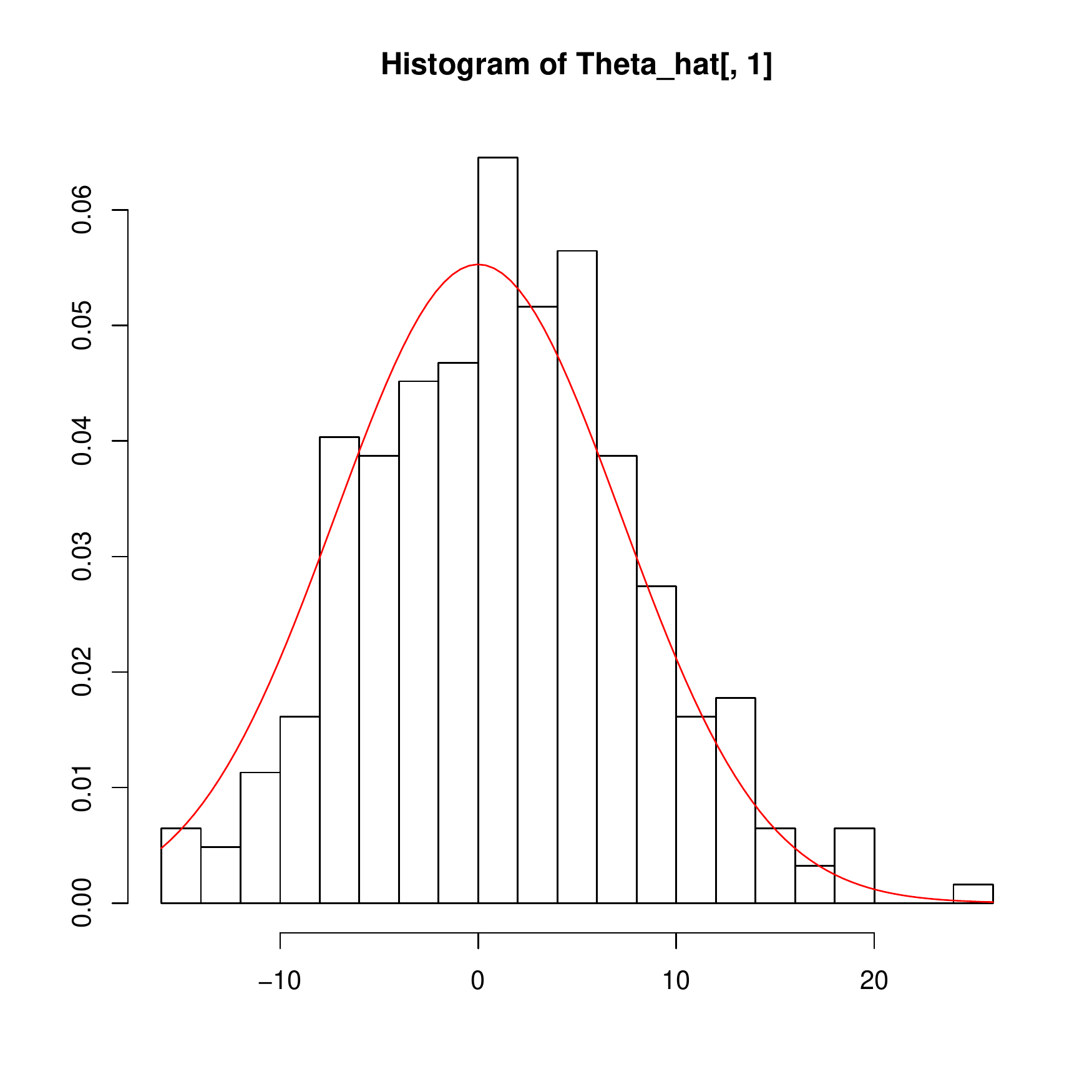}\\
\caption{Simulation results of $\hat{\theta}_1$ \label{fig6}}
\end{center}
\end{figure}

\begin{figure}[h]
\begin{center}
\includegraphics[width=5cm,pagebox=cropbox,clip]{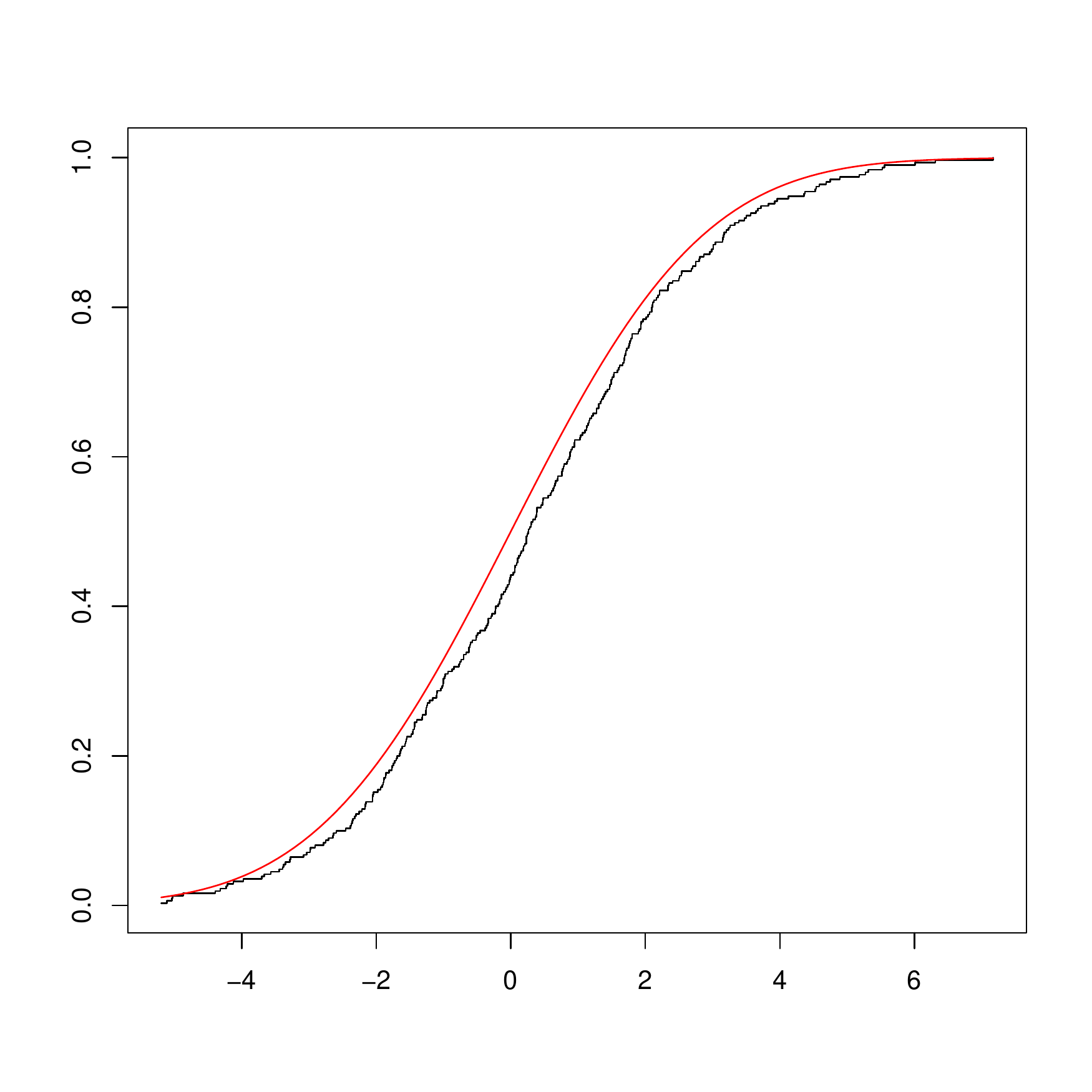}
\includegraphics[width=5cm,pagebox=cropbox,clip]{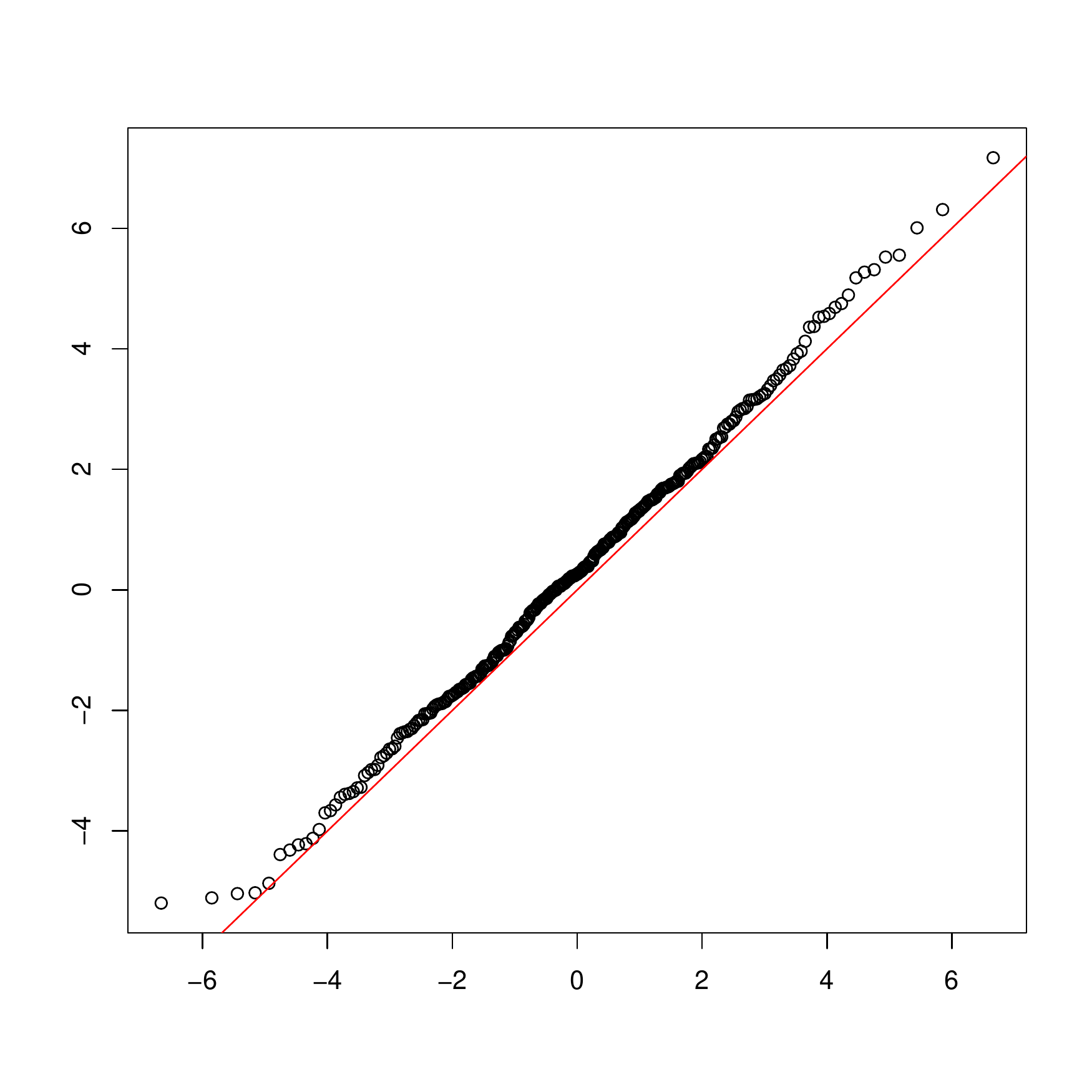}
\includegraphics[width=5cm,pagebox=cropbox,clip]{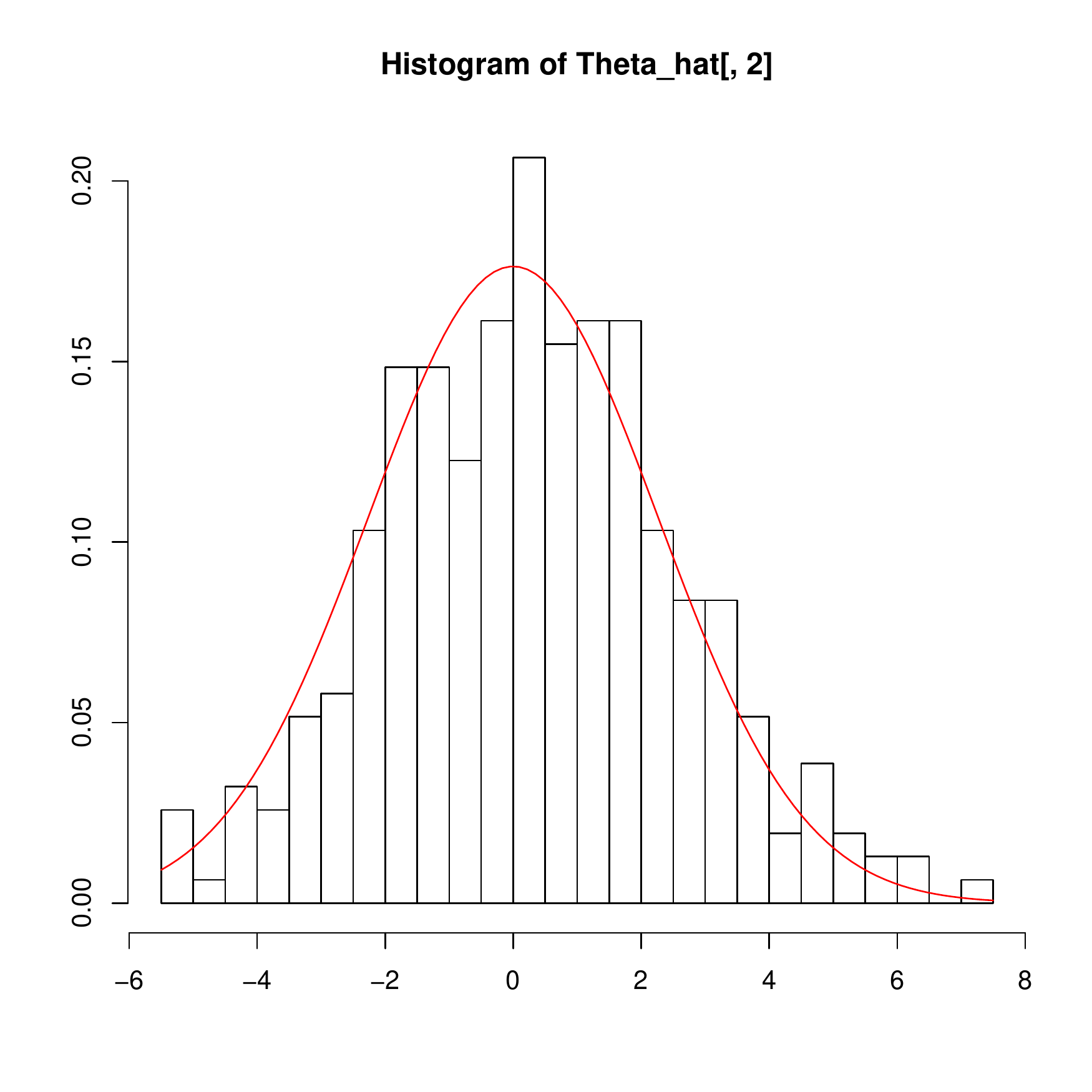}\\
\caption{Simulation results of $\hat{\theta}_2$  \label{fig7}}
\end{center}
\end{figure}


\begin{figure}[t] 
\begin{center}
\includegraphics[width=5cm,pagebox=cropbox,clip]{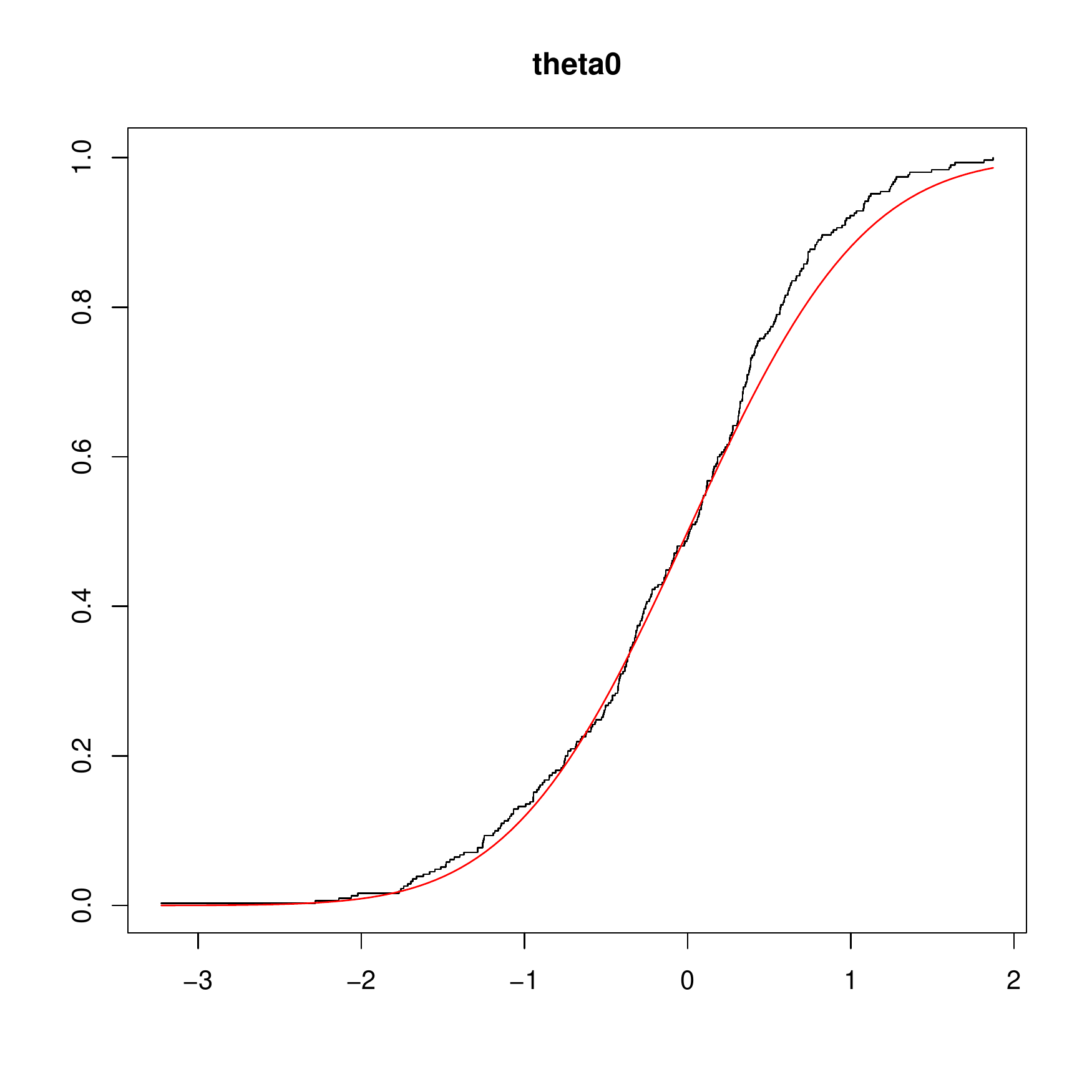}
\includegraphics[width=5cm,pagebox=cropbox,clip]{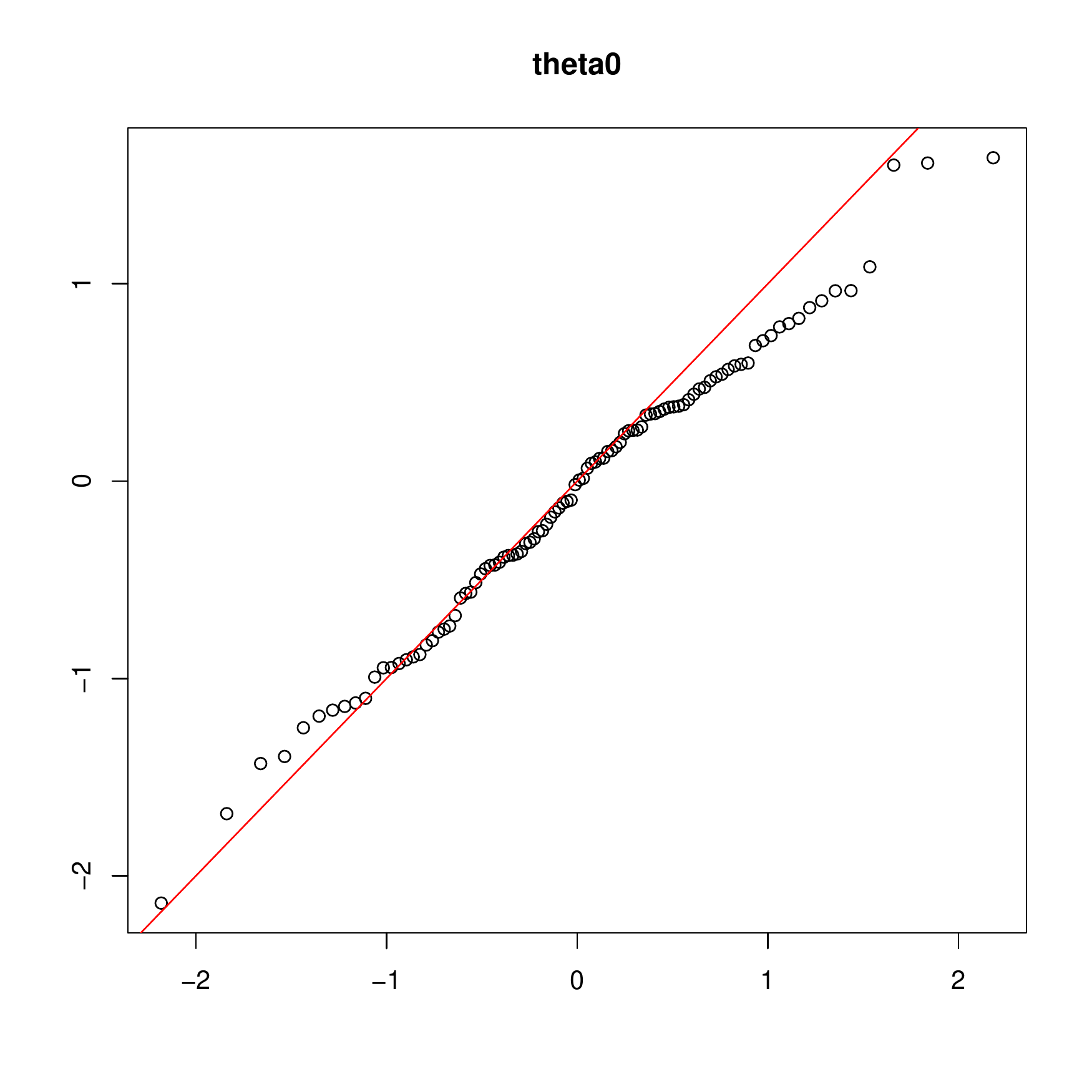}
\includegraphics[width=5cm,pagebox=cropbox,clip]{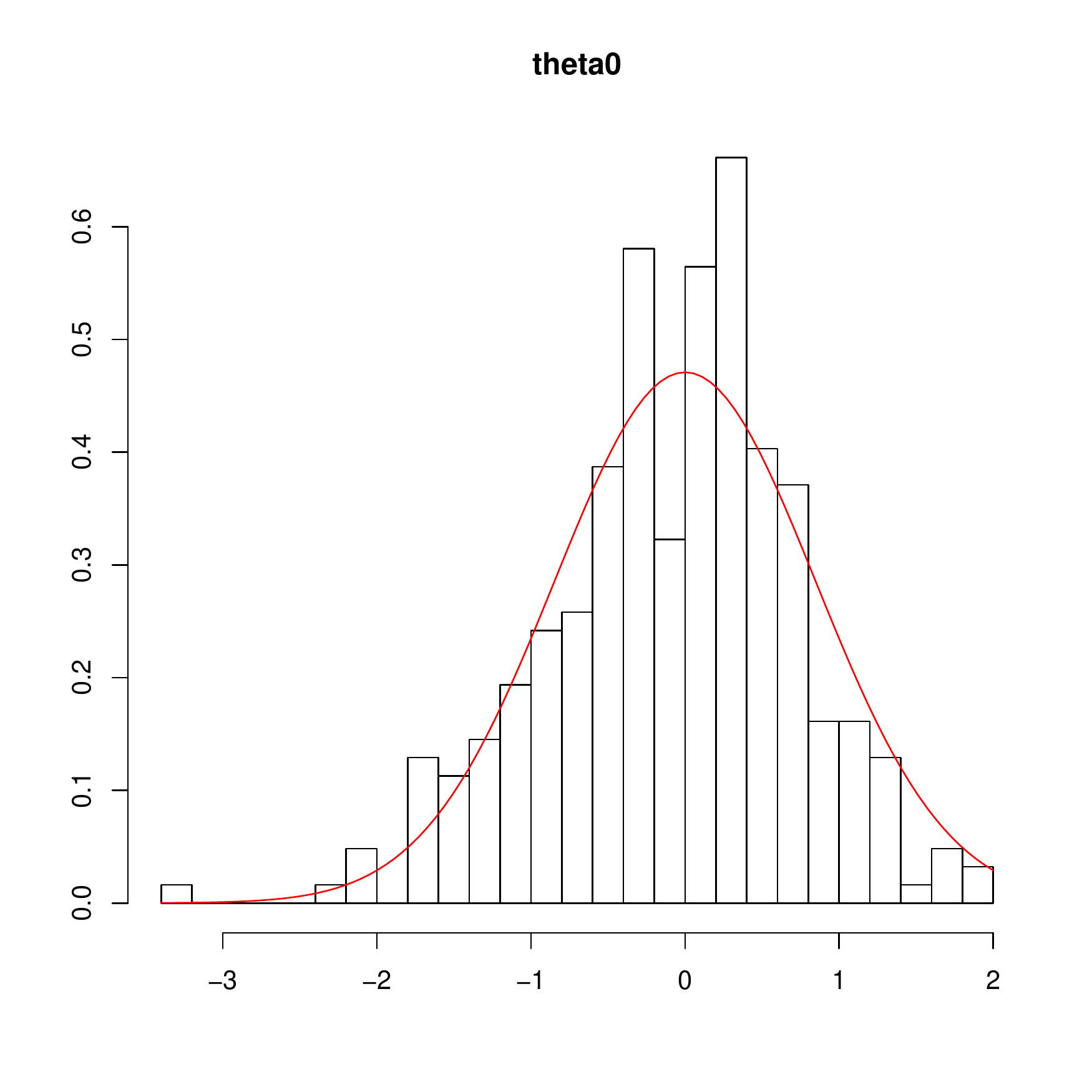}\\
\caption{Simulation results of $\hat{\theta}_0$ \label{fig8}}
\end{center}
\end{figure}


\clearpage

\subsubsection{Summary of example 1}
Table 5 is the simulation results of {the means and the standard s.d.s of} $\hat{\theta}_1$, $\hat{\theta}_2$  and $\hat{\theta}_0$ with $(N, m, N_2) = (10^4, 99, 500)$ from $\epsilon =0.1$ to $0.75$.
{It seems from Table 5 that 
for all $\epsilon$, the estimator of $\theta_0$ has good performance. 
However, it can be seen from Figure 23 
that the asymptotic theory does not work when $\epsilon$ is 0.75.
For this setting, $\epsilon$ 
{should}
be less than 0.5.
}

\begin{table}[h]
\caption{Simulation results of $\hat{\theta}_1$, $\hat{\theta}_2$  and $\hat{\theta}_0$ with $(N, m, N_2) = (10^4, 99, 500)$ \label{table2}}
\begin{center}
\begin{tabular}{cc|ccc} \hline
	&	&$\hat{\theta}_{1}$&$\hat{\theta}_{2}$&$\hat{\theta}_{0}$
\\ \hline
& {true value} &1 & 0.2 & 0
\\ \hline
&mean &1.001&  0.200&  0.010
 \\
$\epsilon=0.1$&{s.d.} & (0.007)& (0.001)& (0.084)
 \\   \hline
&mean &1.002&  0.201&  0.009
 \\
$\epsilon=0.25$&{s.d.} & (0.007)& (0.002)& (0.188)
 \\   \hline
&mean &1.002&  0.201&  -0.013
 \\
$\epsilon=0.5$&{s.d.} & (0.007)& (0.002)& (0.367)
 \\   \hline
&mean &1.002&  0.201&  -0.058
 \\
$\epsilon=0.75$&{s.d.} & (0.007)& (0.002)& (0.548)
 \\   \hline
\end{tabular}
\end{center}
\end{table}


\subsection{Example 2}
The true value of parameter $\theta^*=(\theta_0^*, \theta_1^*, \theta_2^*) = (3.1,1,0.2)$
and 
$\lambda_1^* = 0.12$.
We set that $N = 10^4$, $M = 10^4$, $K = 10^5$, $T=1$, $x_1(0) = 2$, 
$\xi(y) = 2.8y(1-y)$.
Figure 24 is {a} sample path of $X_t(y)$ for $(t,y) \in [0,1]\times [0,1]$
when $(\theta_0^*, \theta_1^*, \theta_2^*, \epsilon) = (3.1,1,0.2,0)$.

\begin{figure}[h] 
\begin{center}
\includegraphics[width=9cm]{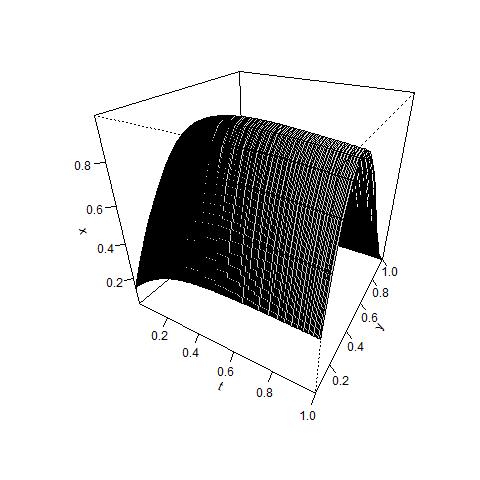} 
\caption{Sample path with 
$(\theta_0^*, \theta_1^*, \theta_2^*, \epsilon) =  (3.1,1,0.2,0)$, 
$\xi(y) = 2.8y(1-y)$ \label{fig5}}
\end{center}
\end{figure}

\subsubsection{$\epsilon=0.1$}
Figure 25 is {a} sample path of $X_t(y)$ for $(t,y) \in [0,1]\times [0,1]$
when $(\theta_0^*, \theta_1^*, \theta_2^*, \epsilon) = (3.1,1,0.2,0.1)$.
Table 6 is the simulation results of {the means and the standard s.d.s of} $\hat{\theta}_1$, $\hat{\theta}_2$  and $\hat{\theta}_0$ with $(N, m, N_2) = (10^4, 99, 500)$.
Figures 26-28 are the simulation results of {the asymptotic distributions of} $\hat{\theta}_1$, $\hat{\theta}_2$  and $\hat{\theta}_0$
with $(N, m, N_2) = (10^4, 99, 500)$.
{
From Figures 26-28, 
we can see that the distributions of the estimators almost 
correspond with the asymptotic distribution in Theorem 1
and the estimators have good performance.
}

\begin{figure}[h] 
\begin{center}
\includegraphics[width=9cm]{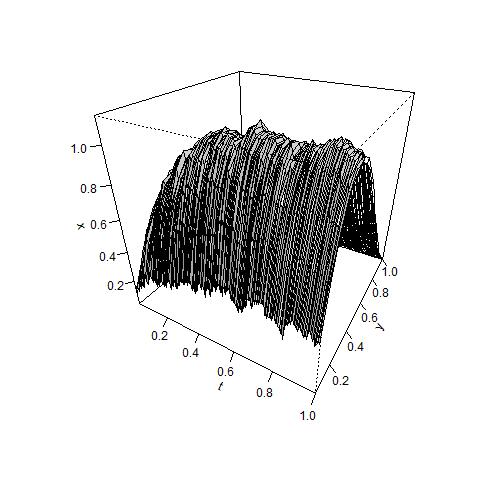} 
\caption{Sample path with $(\theta_0^*, \theta_1^*, \theta_2^*, \epsilon) = (3.1,1,0.2,0.1)$, $\xi(y) = 2.8y(1-y)$ \label{fig5}}
\end{center}
\end{figure}



\begin{table}[h]
\caption{Simulation results of $\hat{\theta}_1$, $\hat{\theta}_2$  and $\hat{\theta}_0$ with $(N, m, N_2) = (10^4, 99, 500)$ \label{table2}}
\begin{center}
\begin{tabular}{c|ccc} \hline
		&$\hat{\theta}_{1}$&$\hat{\theta}_{2}$&$\hat{\theta}_{0}$
\\ \hline
 {true value} &1 & 0.2 & 3.1
\\ \hline
mean &1.001&  0.200&  3.102
 \\
{s.d.} & (0.007)& (0.002)& (0.055)
 \\   \hline
\end{tabular}
\end{center}
\end{table}


\begin{figure}[h] 
\begin{center}
\includegraphics[width=5cm,pagebox=cropbox,clip]{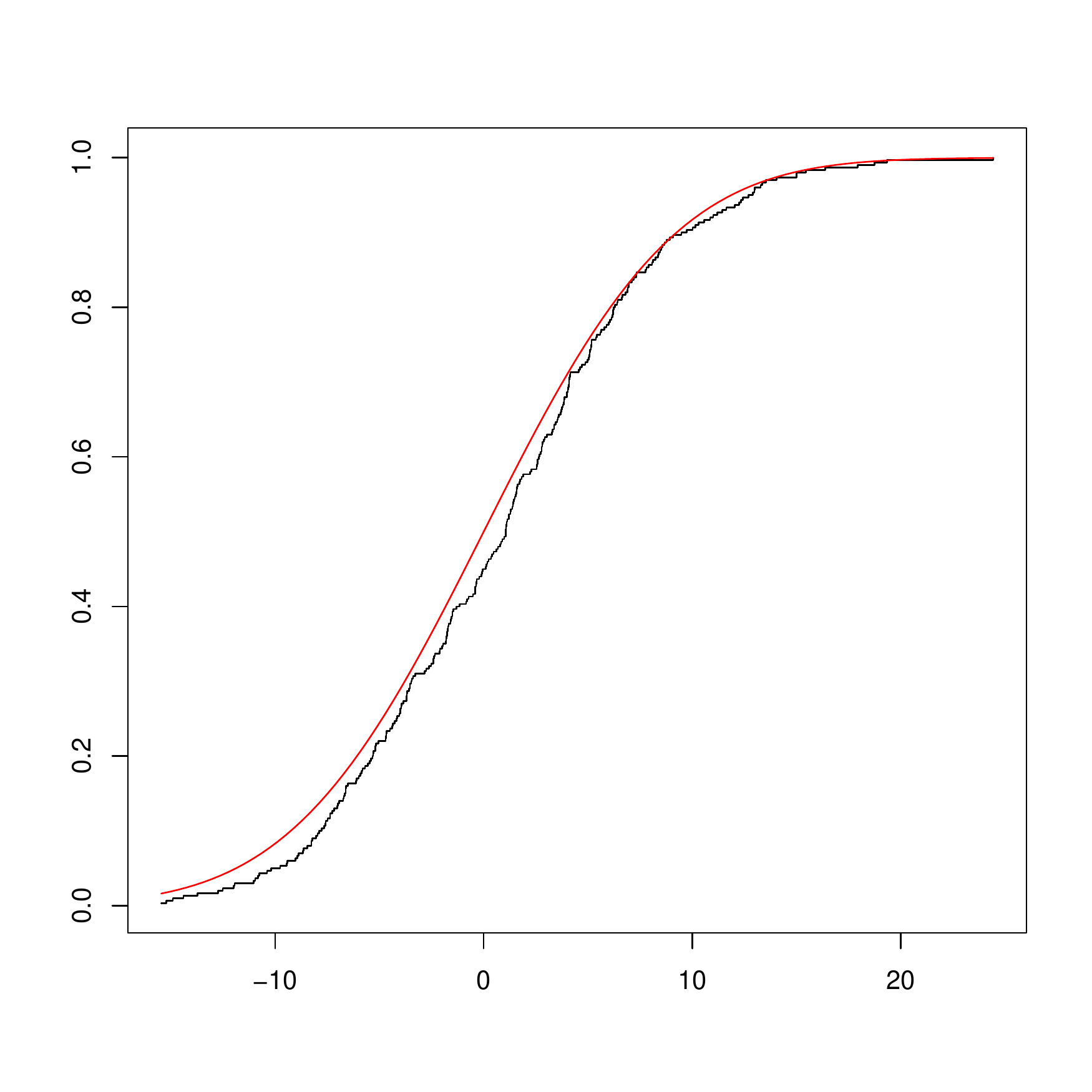}
\includegraphics[width=5cm,pagebox=cropbox,clip]{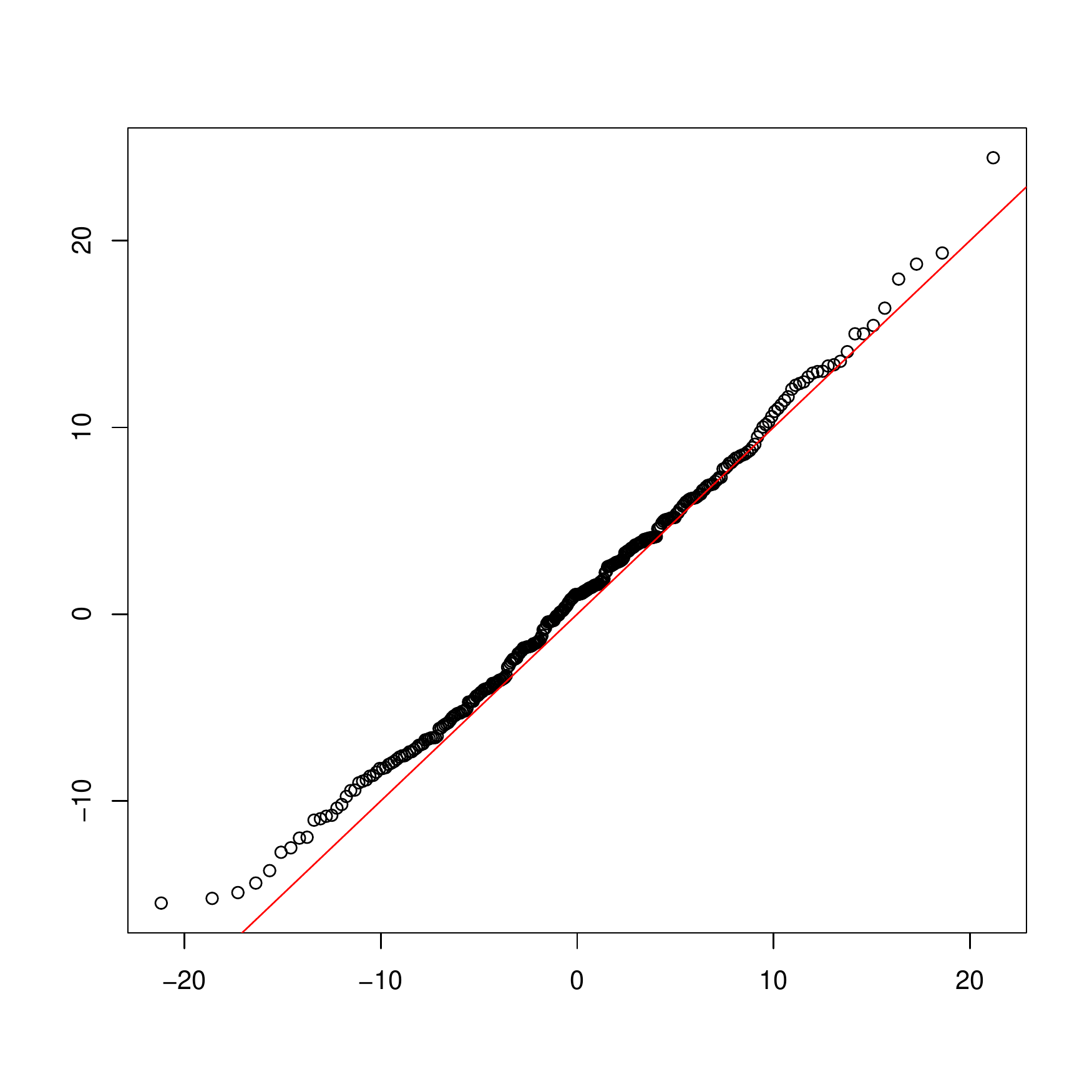}
\includegraphics[width=5cm,pagebox=cropbox,clip]{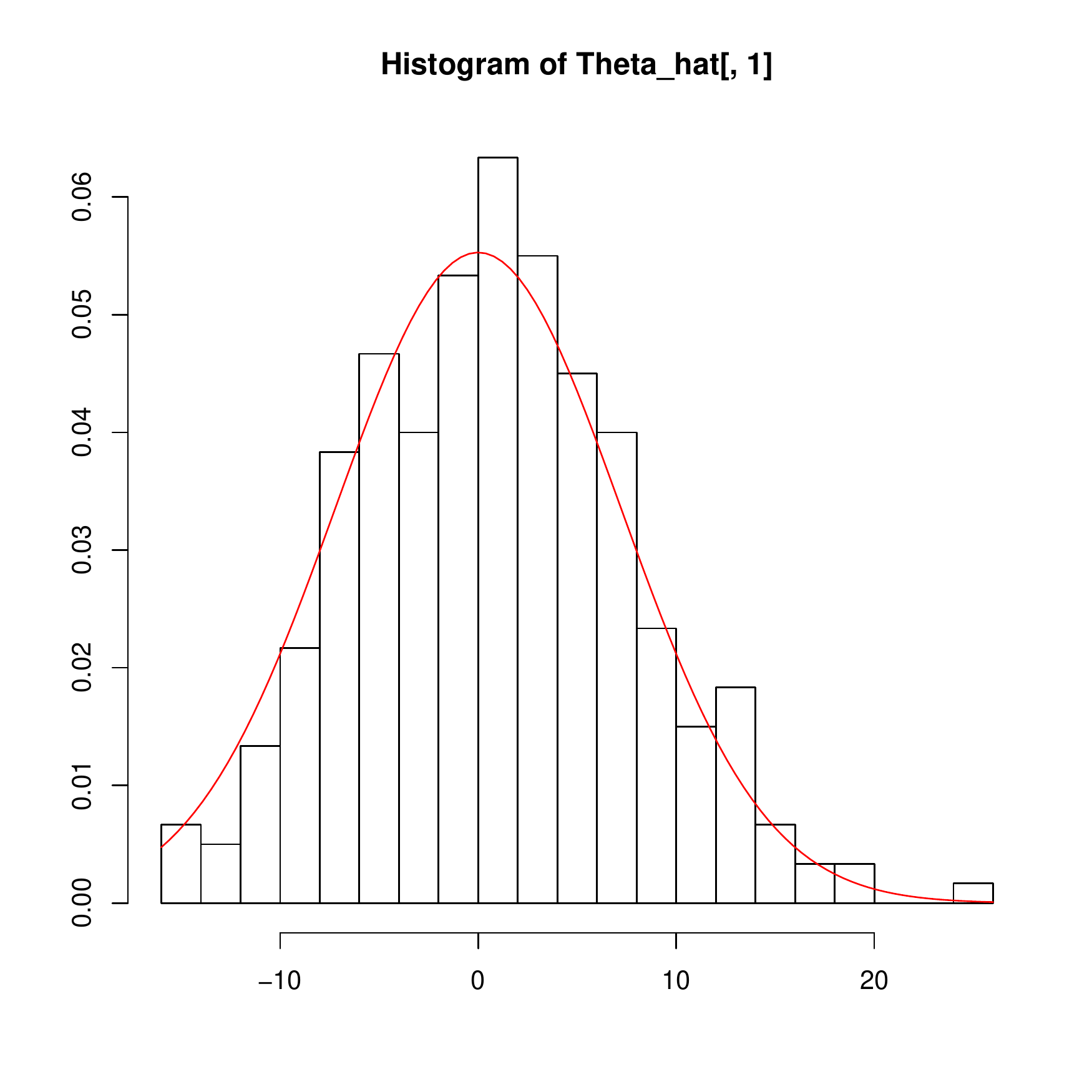}\\
\caption{Simulation results of $\hat{\theta}_1$ \label{fig6}}
\end{center}
\end{figure}

\begin{figure}[h]
\begin{center}
\includegraphics[width=5cm,pagebox=cropbox,clip]{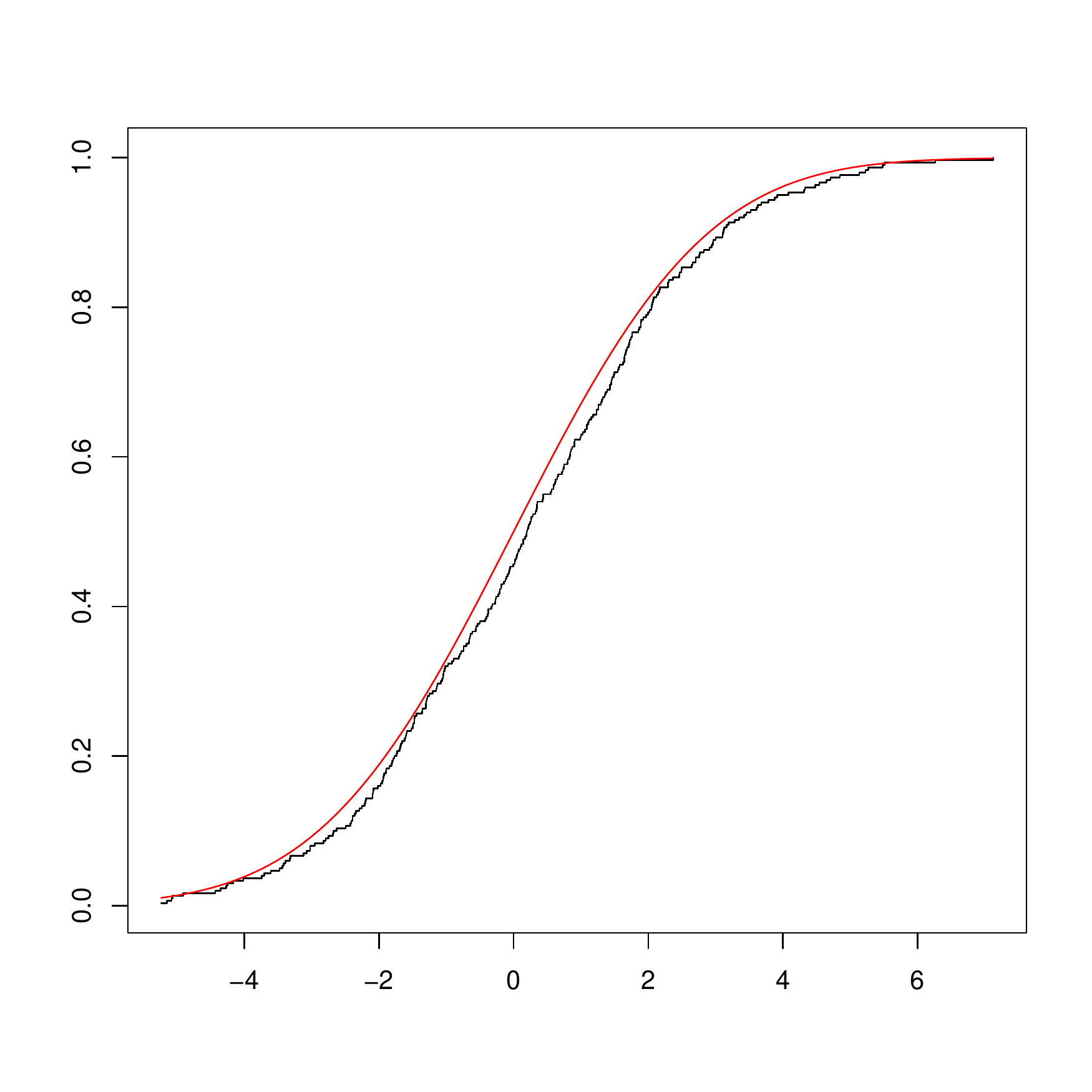}
\includegraphics[width=5cm,pagebox=cropbox,clip]{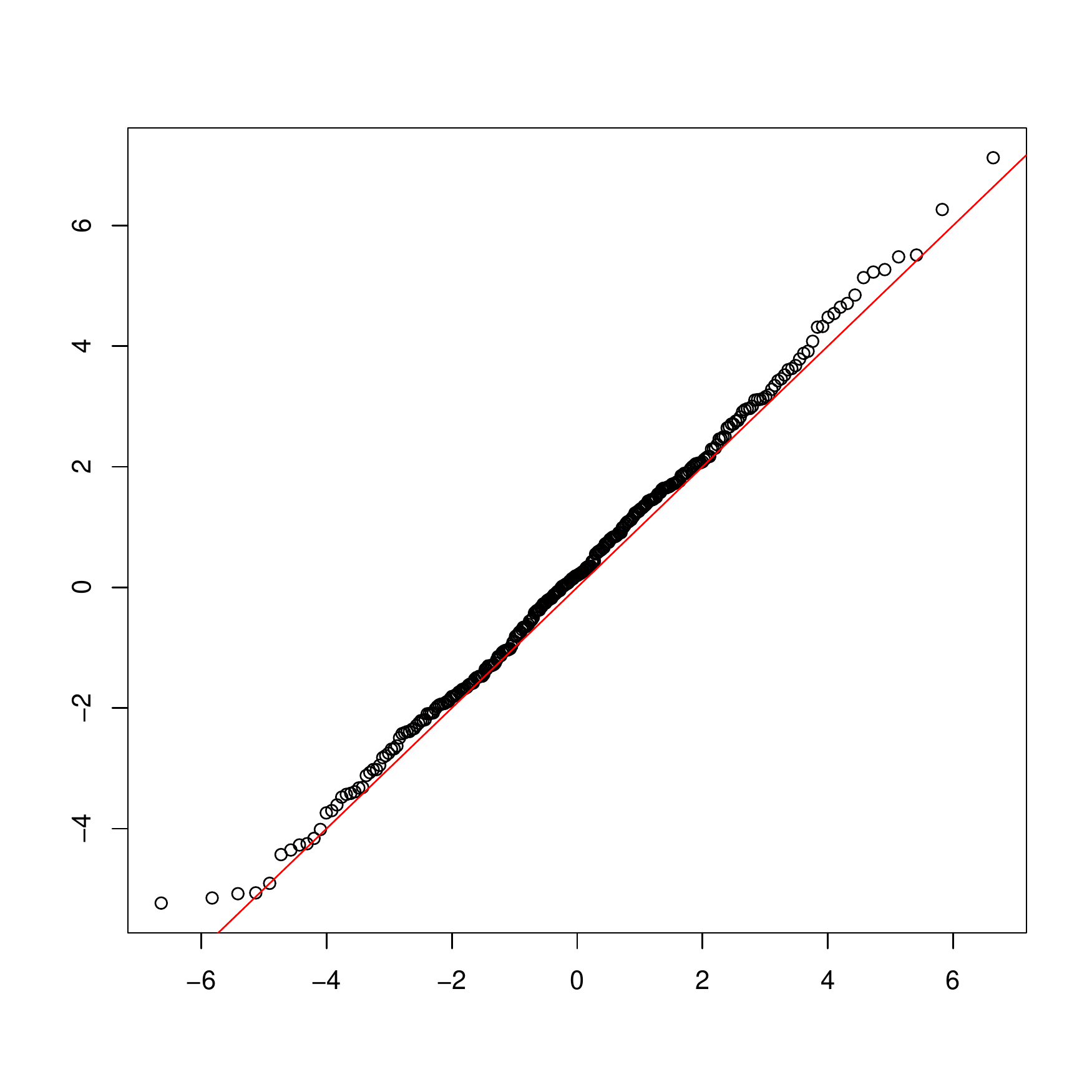}
\includegraphics[width=5cm,pagebox=cropbox,clip]{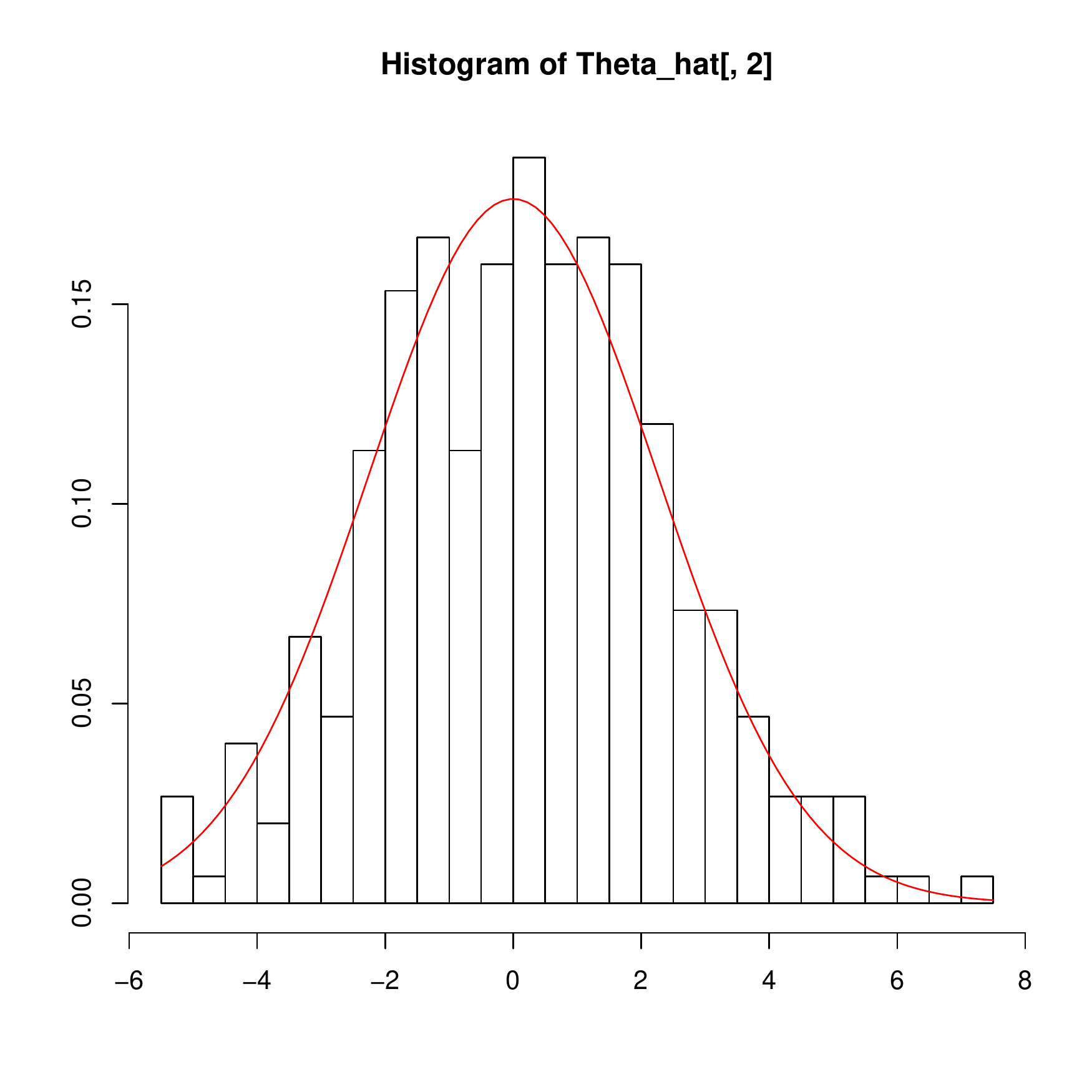}\\
\caption{Simulation results of $\hat{\theta}_2$  \label{fig7}}
\end{center}
\end{figure}


\begin{figure}[t] 
\begin{center}
\includegraphics[width=5cm,pagebox=cropbox,clip]{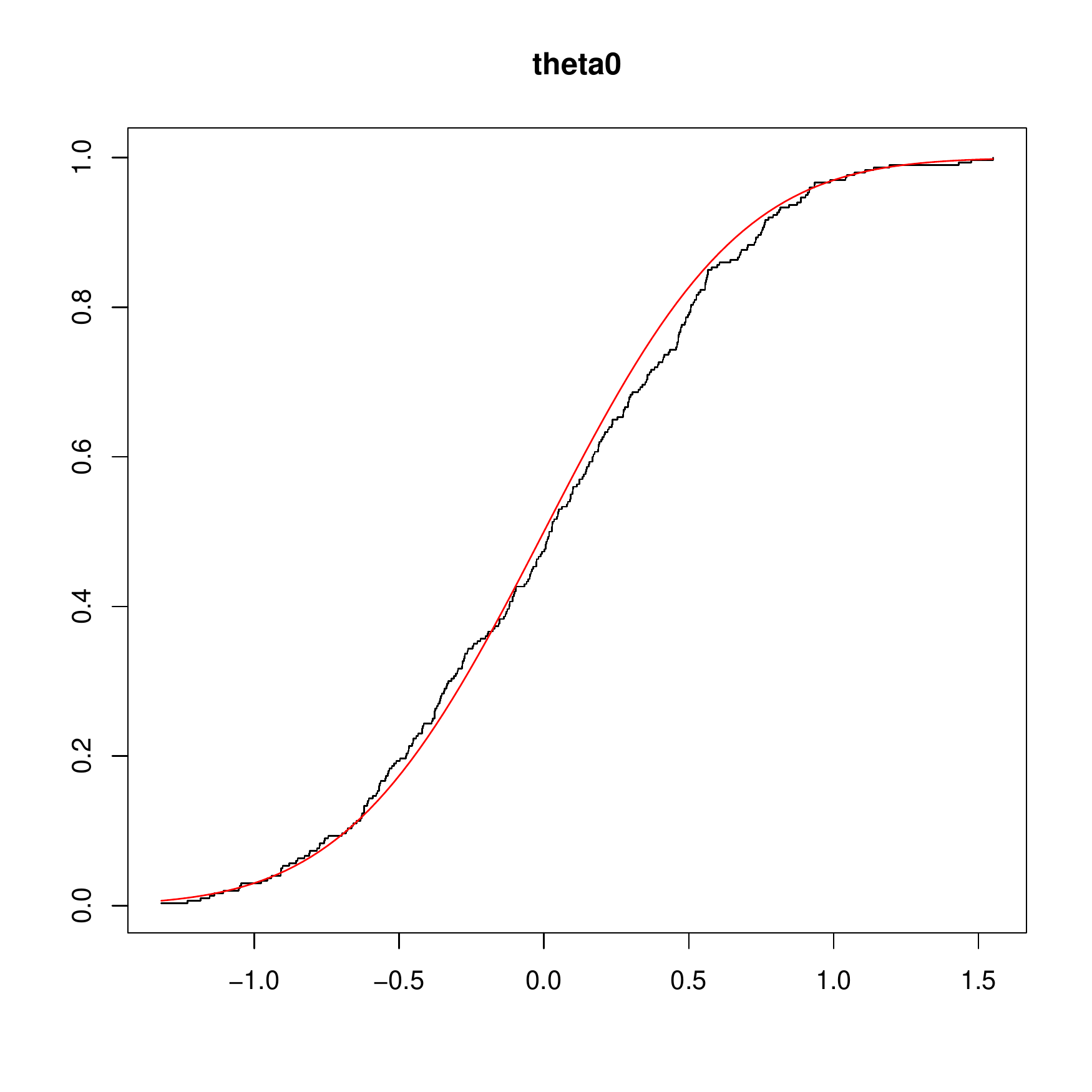}
\includegraphics[width=5cm,pagebox=cropbox,clip]{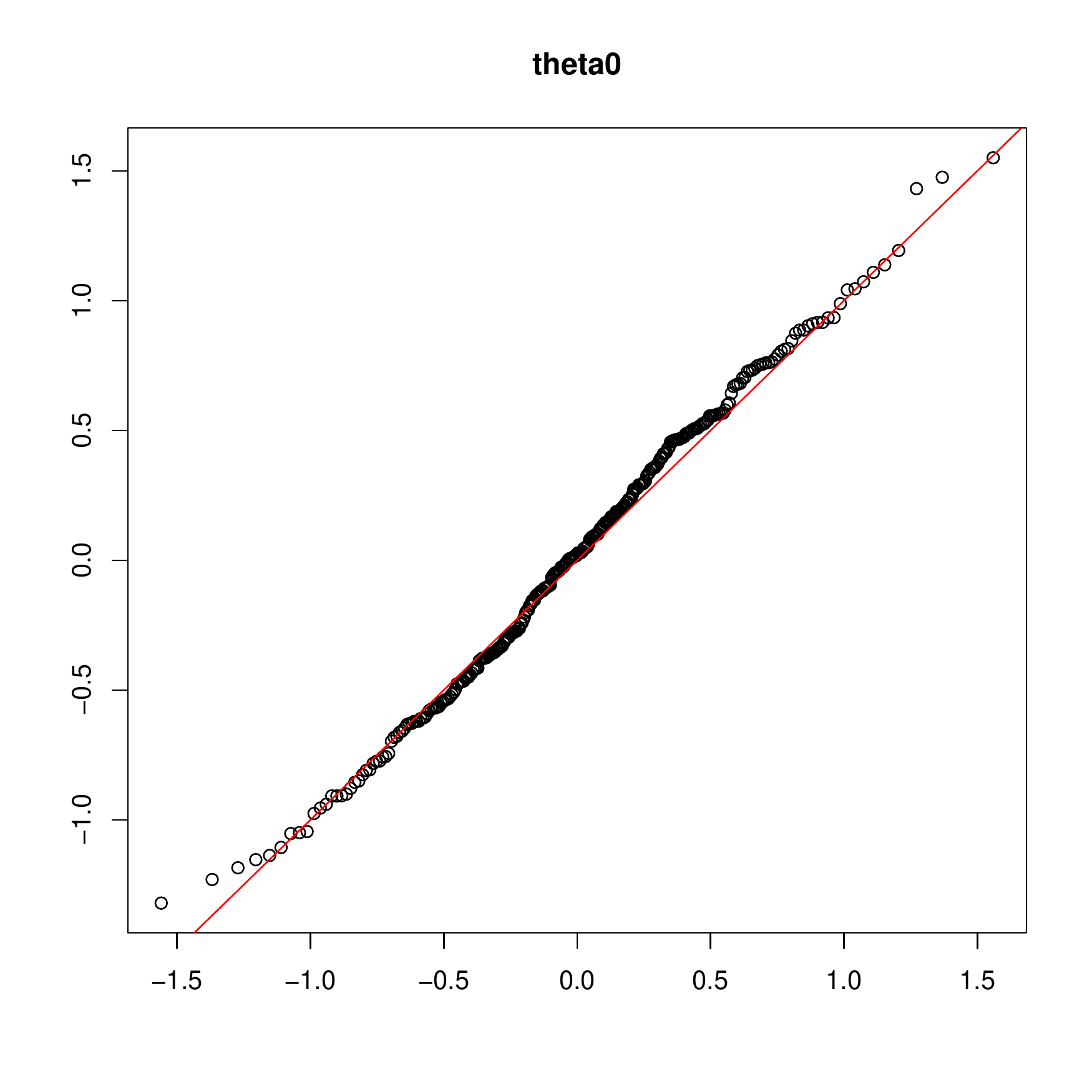}
\includegraphics[width=5cm,pagebox=cropbox,clip]{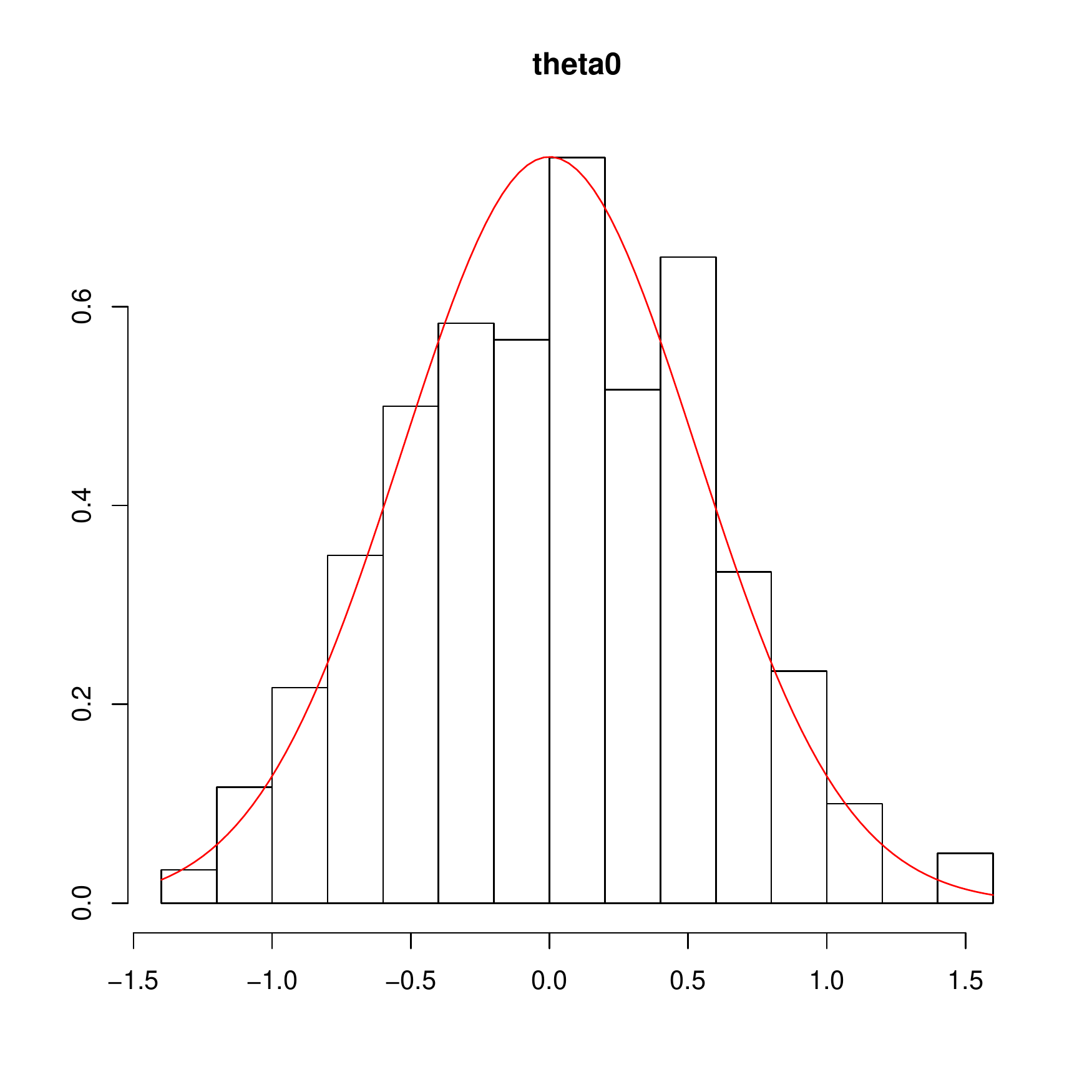}\\
\caption{Simulation results of $\hat{\theta}_0$ \label{fig8}}
\end{center}
\end{figure}


\clearpage

\subsubsection{$\epsilon=0.25$}
Figure 29 is {a} sample path of $X_t(y)$ for $(t,y) \in [0,1]\times [0,1]$
when $(\theta_0^*, \theta_1^*, \theta_2^*, \epsilon) = (3.1,1,0.2,0.25)$.
Table 7 is the simulation results of {the means and the standard s.d.s of} $\hat{\theta}_1$, $\hat{\theta}_2$  and $\hat{\theta}_0$ with $(N, m, N_2) = (10^4, 99, 500)$.
Figures 30-32 are the simulation results of {the asymptotic distributions of} $\hat{\theta}_1$, $\hat{\theta}_2$  and $\hat{\theta}_0$
with $(N, m, N_2) = (10^4, 99, 500)$.
{
It seems from  Figures 30-32 that the estimators 
have the asymptotic distributions in Theorem 1  and 
{their performance is good.}
}

\begin{figure}[h] 
\begin{center}
\includegraphics[width=9cm]{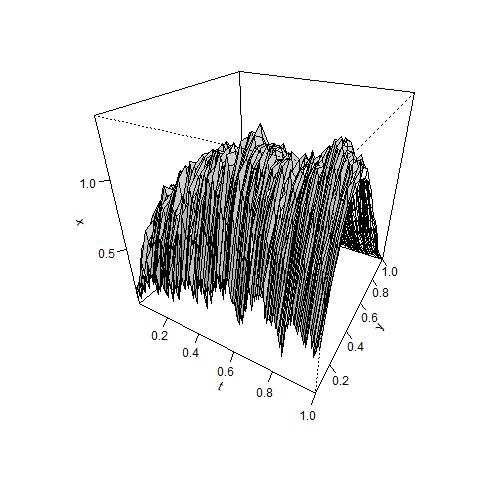} 
\caption{Sample path with $(\theta_0^*, \theta_1^*, \theta_2^*, \epsilon) =  (3.1,1,0.2,0.25)$, $\xi(y) = 2.8y(1-y)$ \label{fig5}}
\end{center}
\end{figure}



\begin{table}[h]
\caption{Simulation results of $\hat{\theta}_1$, $\hat{\theta}_2$  and $\hat{\theta}_0$ with $(N, m, N_2) = (10^4, 99, 500)$ \label{table2}}
\begin{center}
\begin{tabular}{c|ccc} \hline
		&$\hat{\theta}_{1}$&$\hat{\theta}_{2}$&$\hat{\theta}_{0}$
\\ \hline
 {true value} &1 & 0.2 & 3.1

\\ \hline
mean &1.001&  0.200&  3.092
 \\
{s.d.} & (0.007)& (0.002)& (0.126)
 \\   \hline
\end{tabular}
\end{center}
\end{table}


\begin{figure}[h] 
\begin{center}
\includegraphics[width=5cm,pagebox=cropbox,clip]{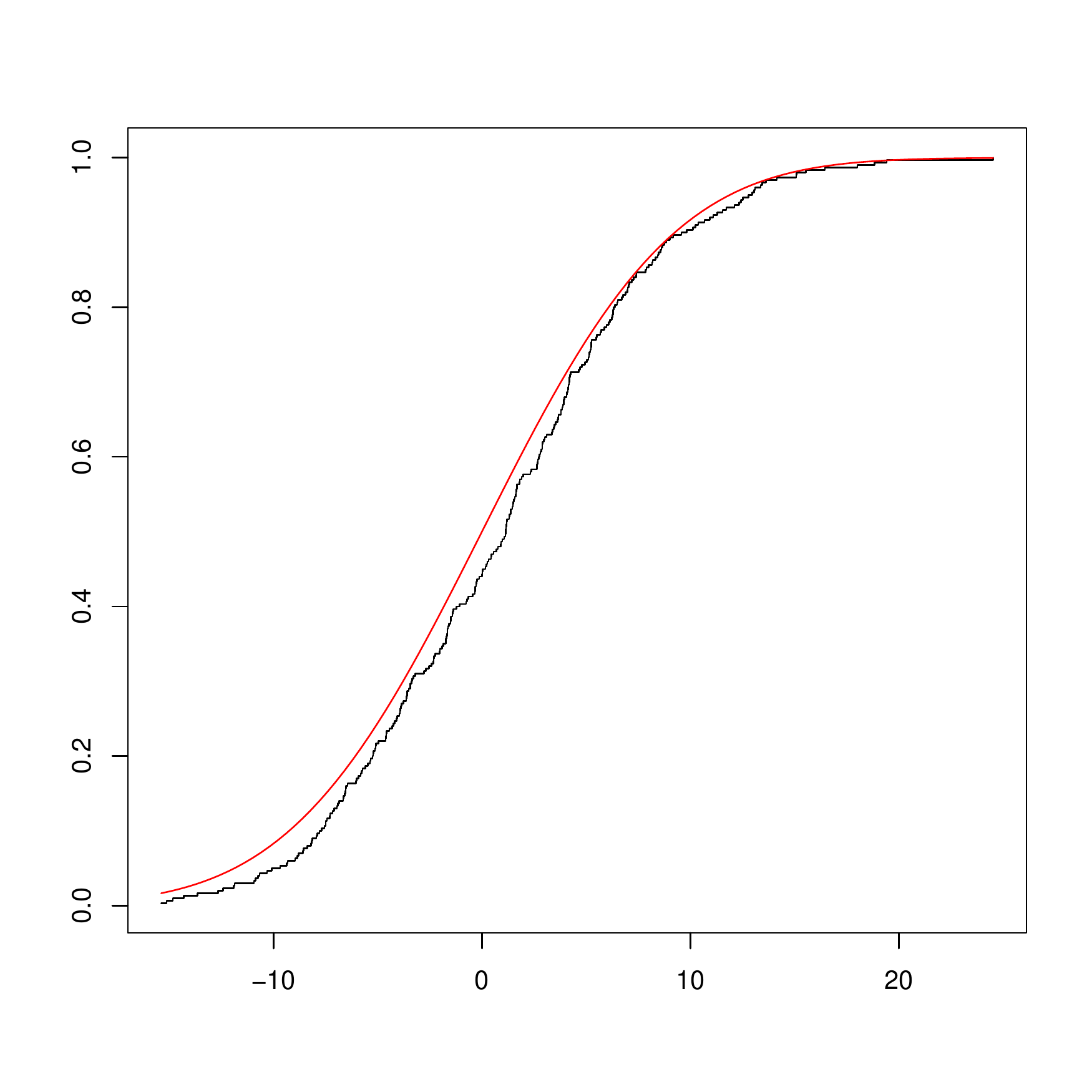}
\includegraphics[width=5cm,pagebox=cropbox,clip]{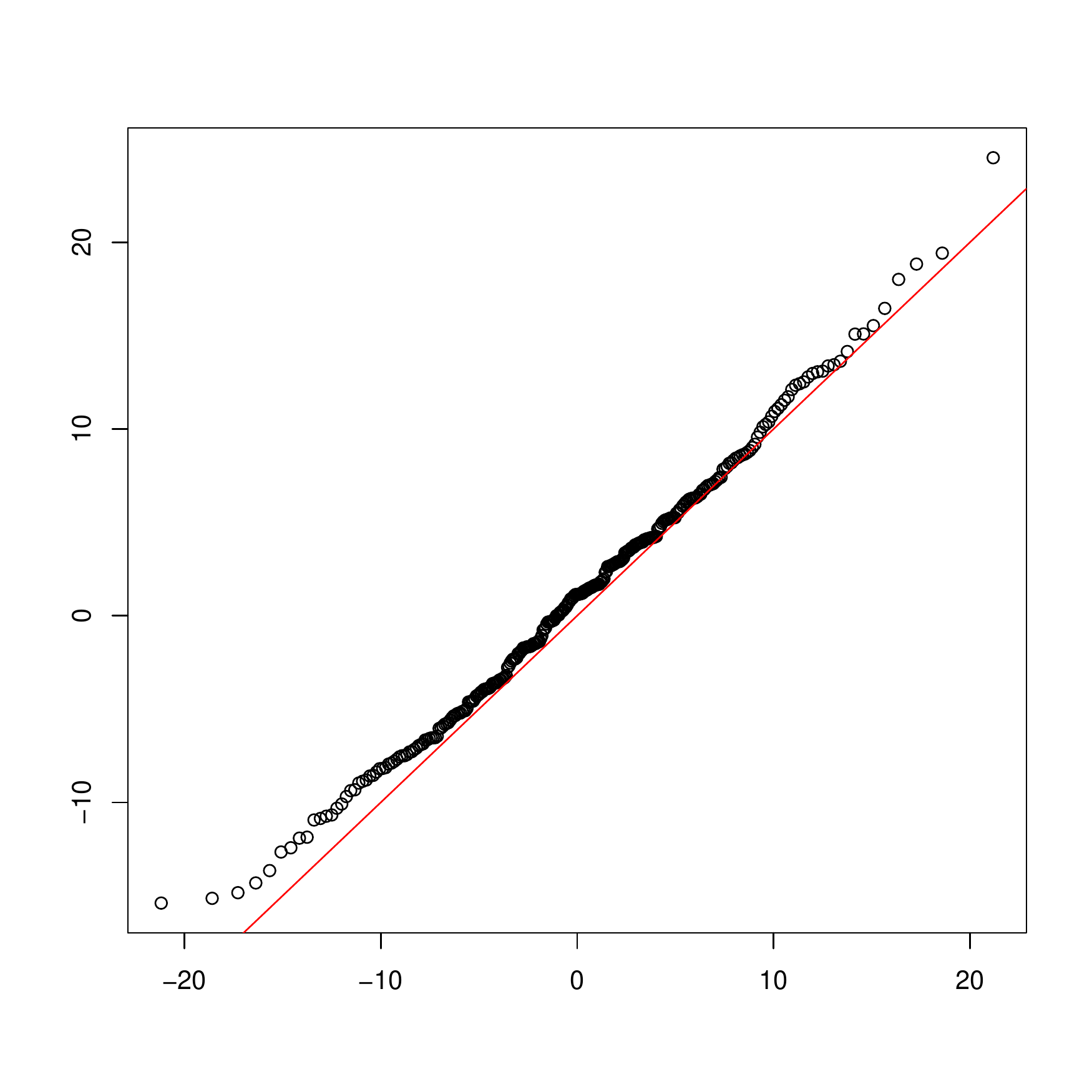}
\includegraphics[width=5cm,pagebox=cropbox,clip]{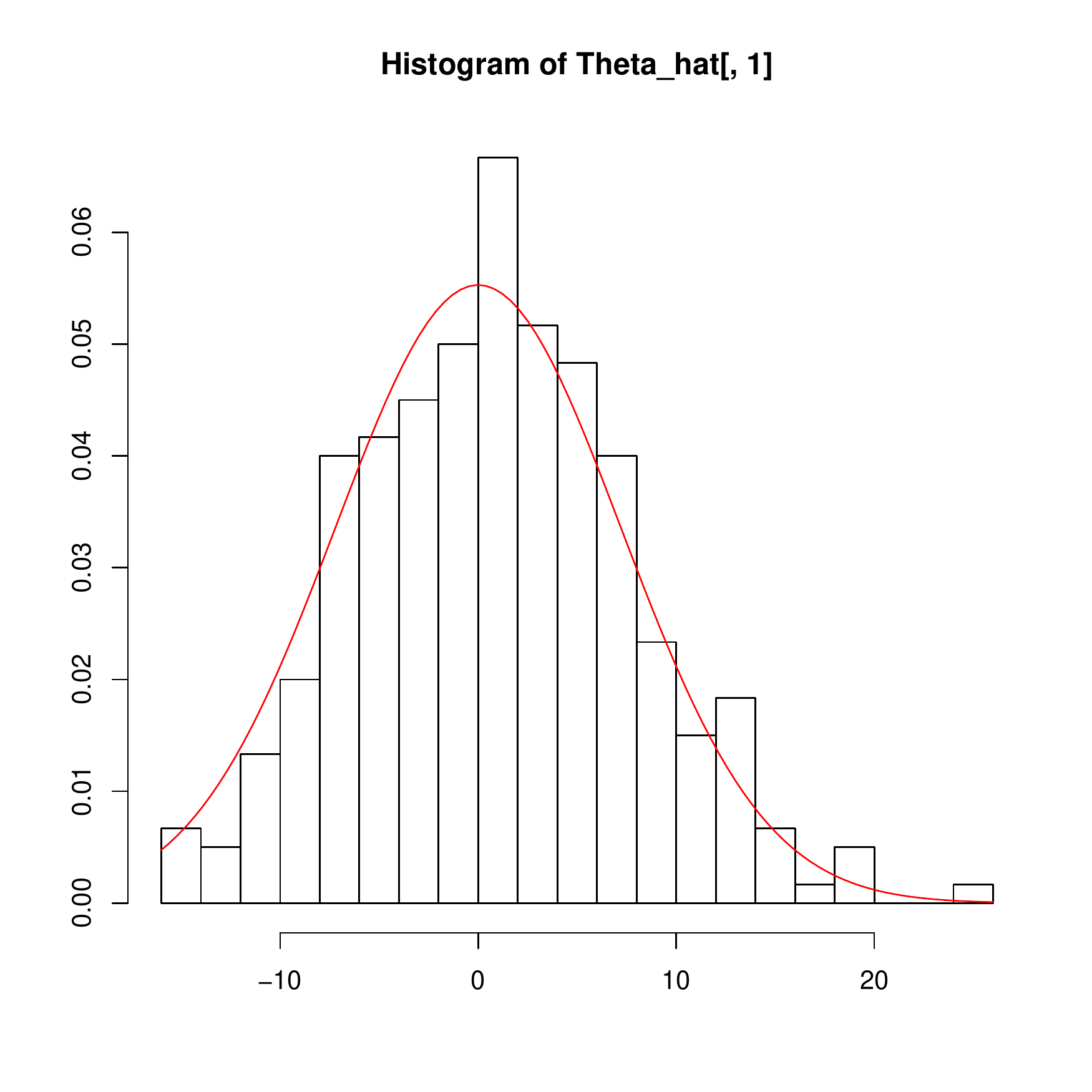}\\
\caption{Simulation results of $\hat{\theta}_1$ \label{fig6}}
\end{center}
\end{figure}

\begin{figure}[h]
\begin{center}
\includegraphics[width=5cm,pagebox=cropbox,clip]{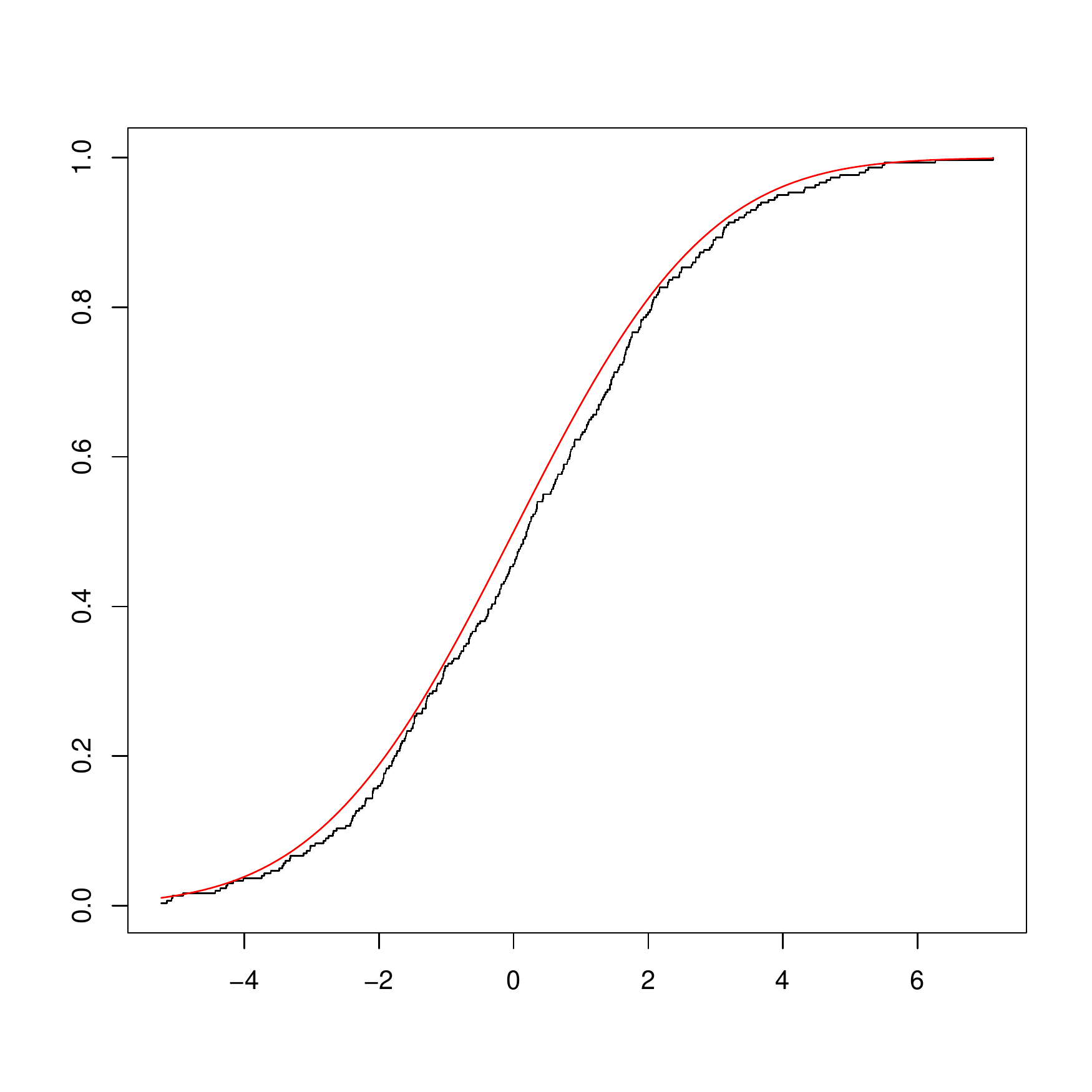}
\includegraphics[width=5cm,pagebox=cropbox,clip]{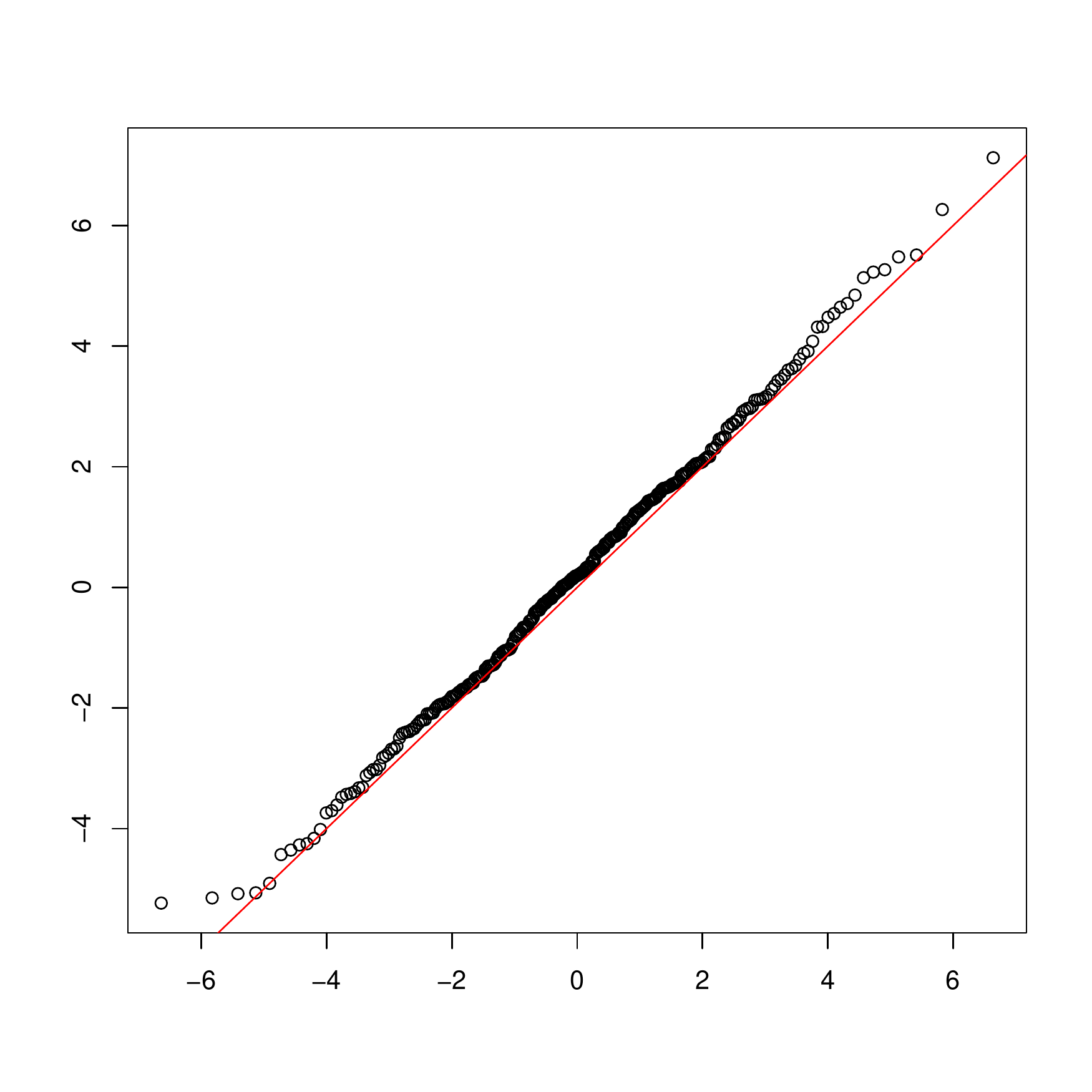}
\includegraphics[width=5cm,pagebox=cropbox,clip]{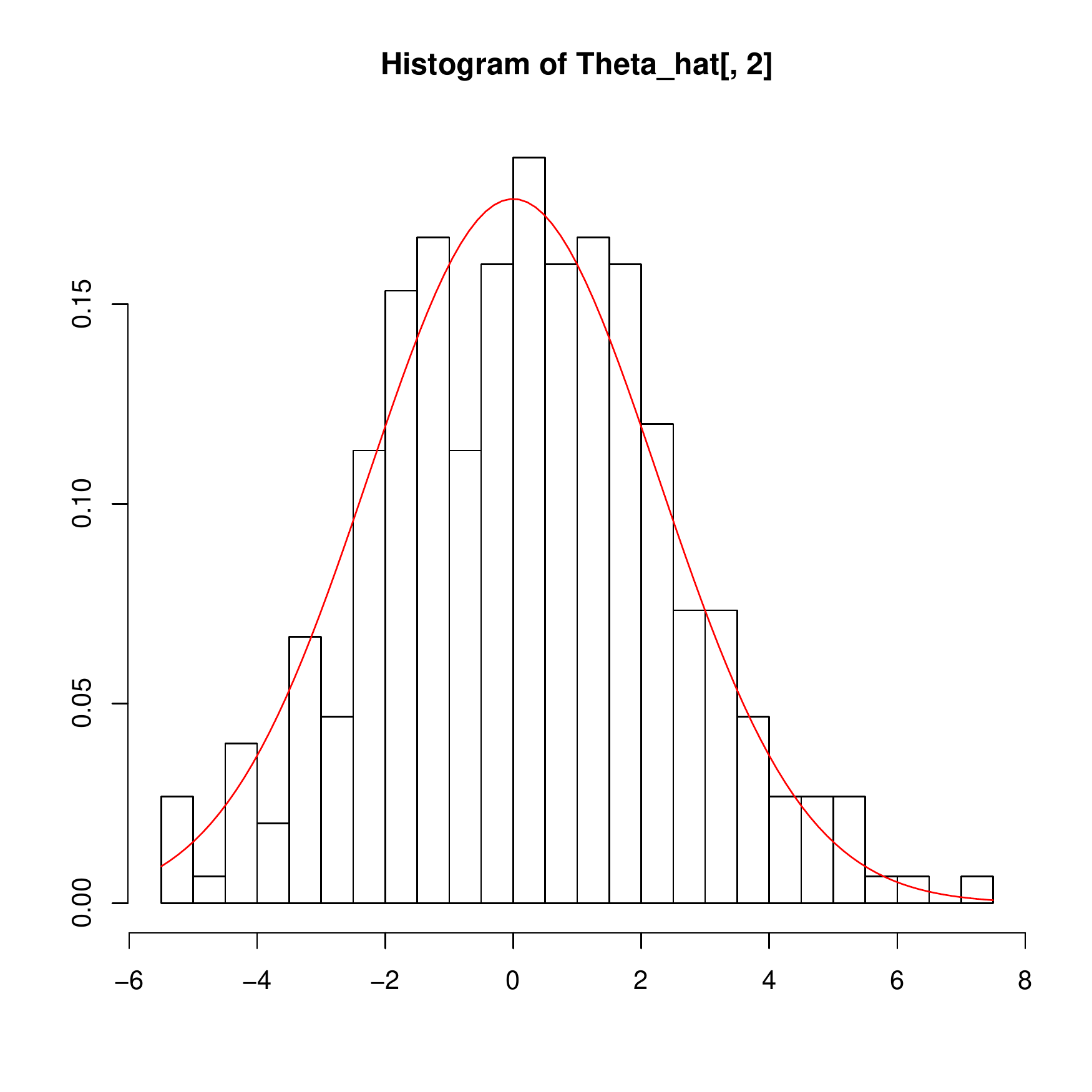}\\
\caption{Simulation results of $\hat{\theta}_2$  \label{fig7}}
\end{center}
\end{figure}


\begin{figure}[t] 
\begin{center}
\includegraphics[width=5cm,pagebox=cropbox,clip]{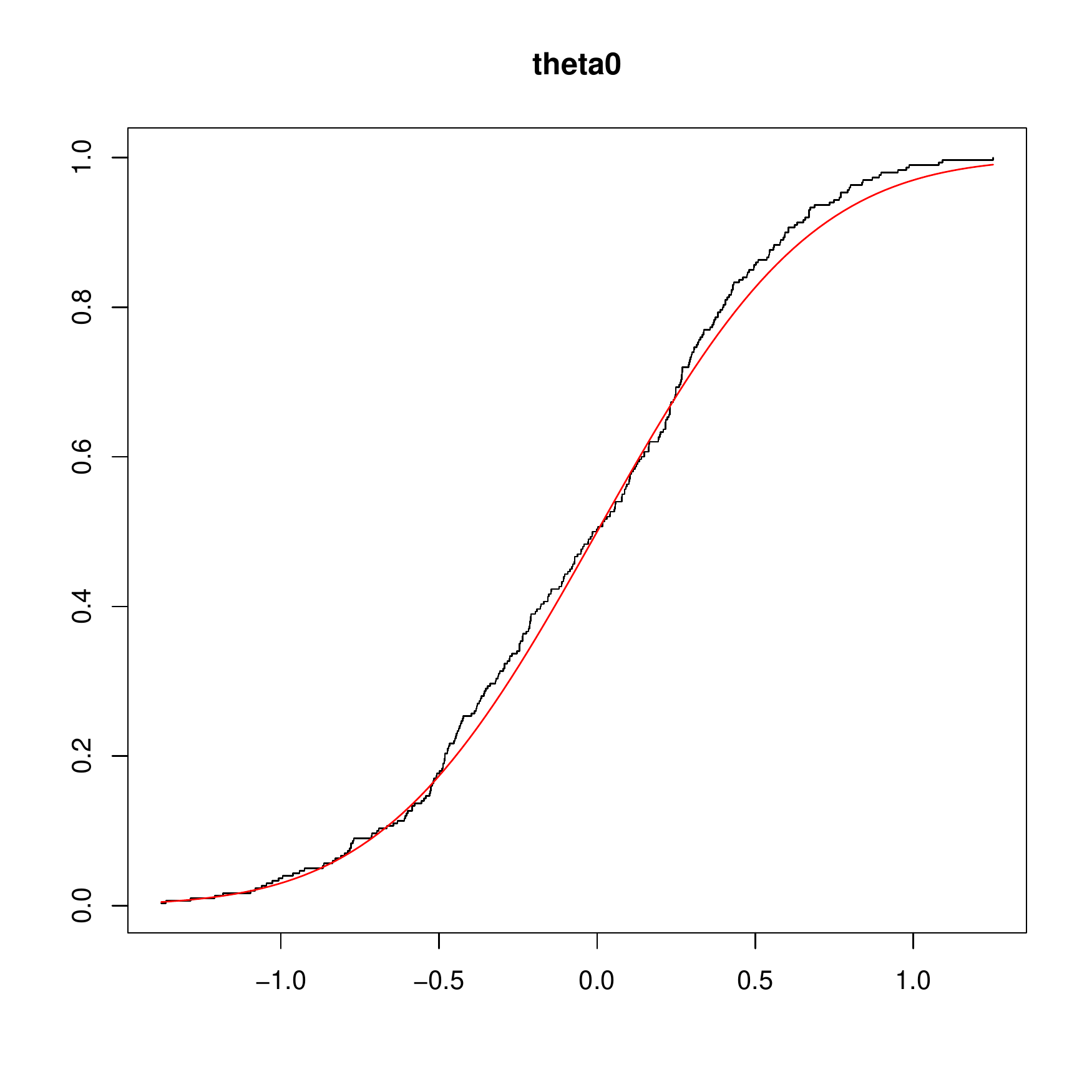}
\includegraphics[width=5cm,pagebox=cropbox,clip]{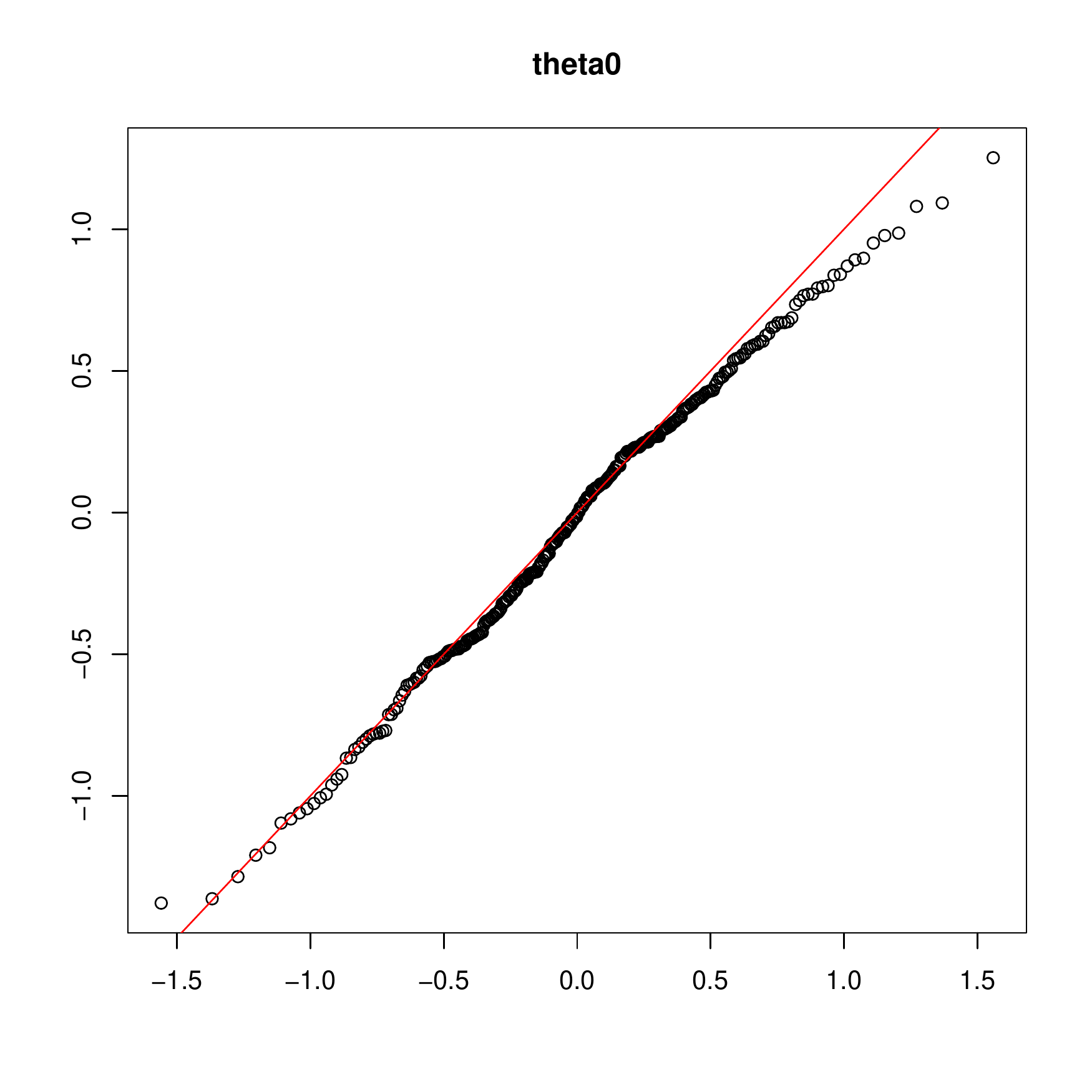}
\includegraphics[width=5cm,pagebox=cropbox,clip]{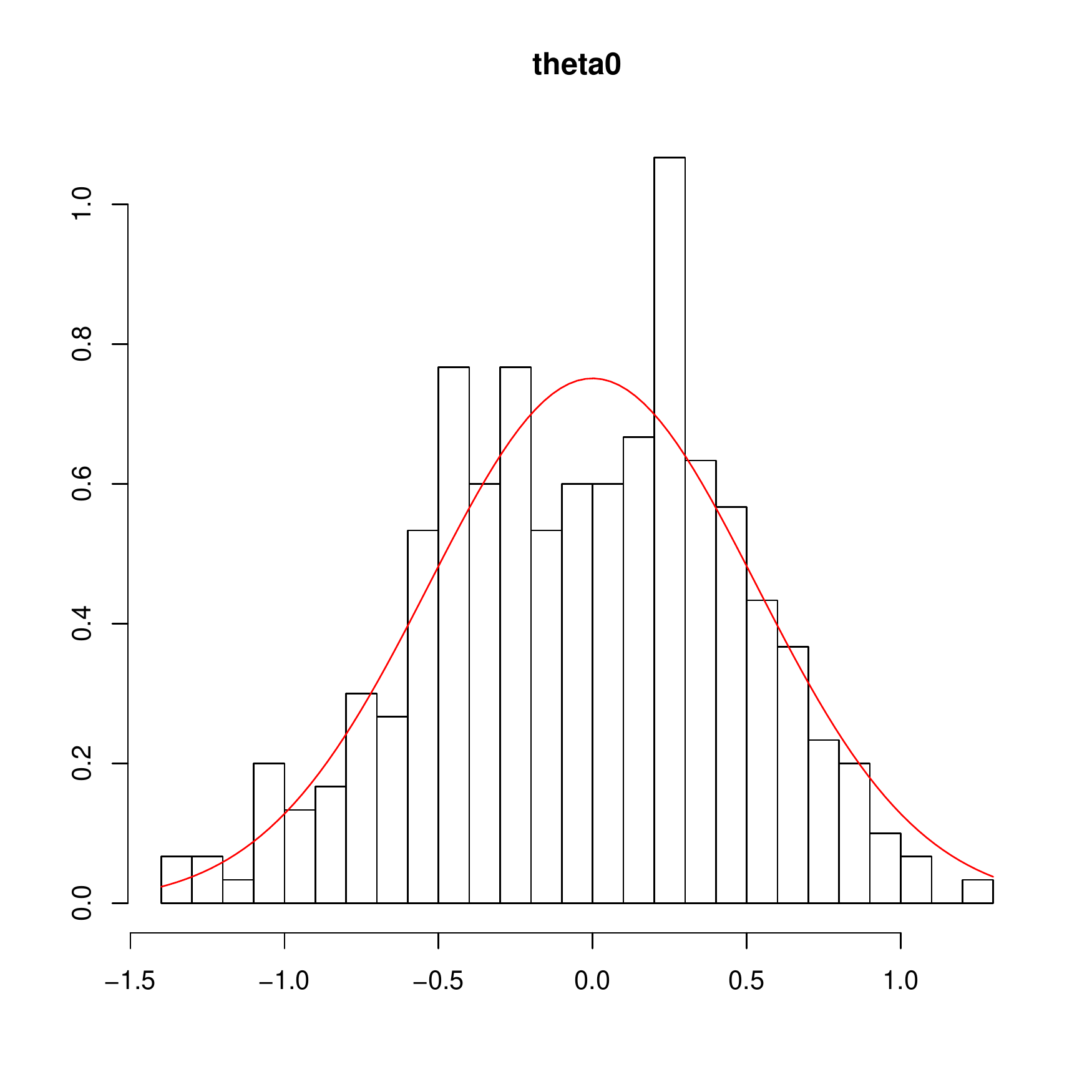}\\
\caption{Simulation results of $\hat{\theta}_0$ \label{fig8}}
\end{center}
\end{figure}


\clearpage

\subsubsection{$\epsilon=0.5$}
Figure 33 is {a} sample path of $X_t(y)$ for $(t,y) \in [0,1]\times [0,1]$
when $(\theta_0^*, \theta_1^*, \theta_2^*, \epsilon) = (3.1,1,0.2,0.5)$.
Table 8 is the simulation results of {the means and the standard s.d.s of} $\hat{\theta}_1$, $\hat{\theta}_2$  and $\hat{\theta}_0$ with $(N, m, N_2) = (10^4, 99, 500)$.
Figures \ref{fig6}-\ref{fig8} are the simulation results of {the asymptotic distributions of} $\hat{\theta}_1$, $\hat{\theta}_2$  and $\hat{\theta}_0$
with $(N, m, N_2) = (10^4, 99, 500)$.
{
{From Figures \ref{fig6}-\ref{fig7}}, 
we can see that the distributions of the estimators of $\theta_1$ and $\theta_2$ 
are almost the same as the asymptotic distributions in Theorem 1 
and the estimators have good performance.
However, it seems from 
{Figure \ref{fig8} that the estimator of $\theta_0$ 
is slightly biased.}
}

\begin{figure}[h] 
\begin{center}
\includegraphics[width=9cm]{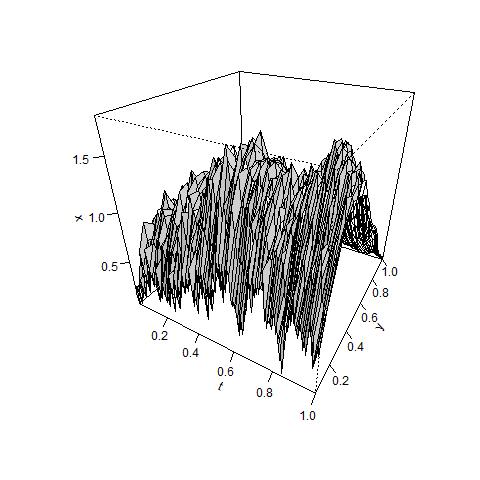} 
\caption{Sample path with $(\theta_0^*, \theta_1^*, \theta_2^*, \epsilon) =  (3.1,1,0.2,0.5)$, $\xi(y) = 2.8y(1-y)$ \label{fig5}}
\end{center}
\end{figure}



\begin{table}[h]
\caption{Simulation results of $\hat{\theta}_1$, $\hat{\theta}_2$  and $\hat{\theta}_0$ with $(N, m, N_2) = (10^4, 99, 500)$ \label{table2}}
\begin{center}
\begin{tabular}{c|ccc} \hline
		&$\hat{\theta}_{1}$&$\hat{\theta}_{2}$&$\hat{\theta}_{0}$
\\ \hline
 {true value} &1 & 0.2 & 3.1

\\ \hline
mean &1.001&  0.200&  3.050
 \\
{s.d.} & (0.007)& (0.002)& (0.258)
 \\   \hline
\end{tabular}
\end{center}
\end{table}


\begin{figure}[h] 
\begin{center}
\includegraphics[width=5cm,pagebox=cropbox,clip]{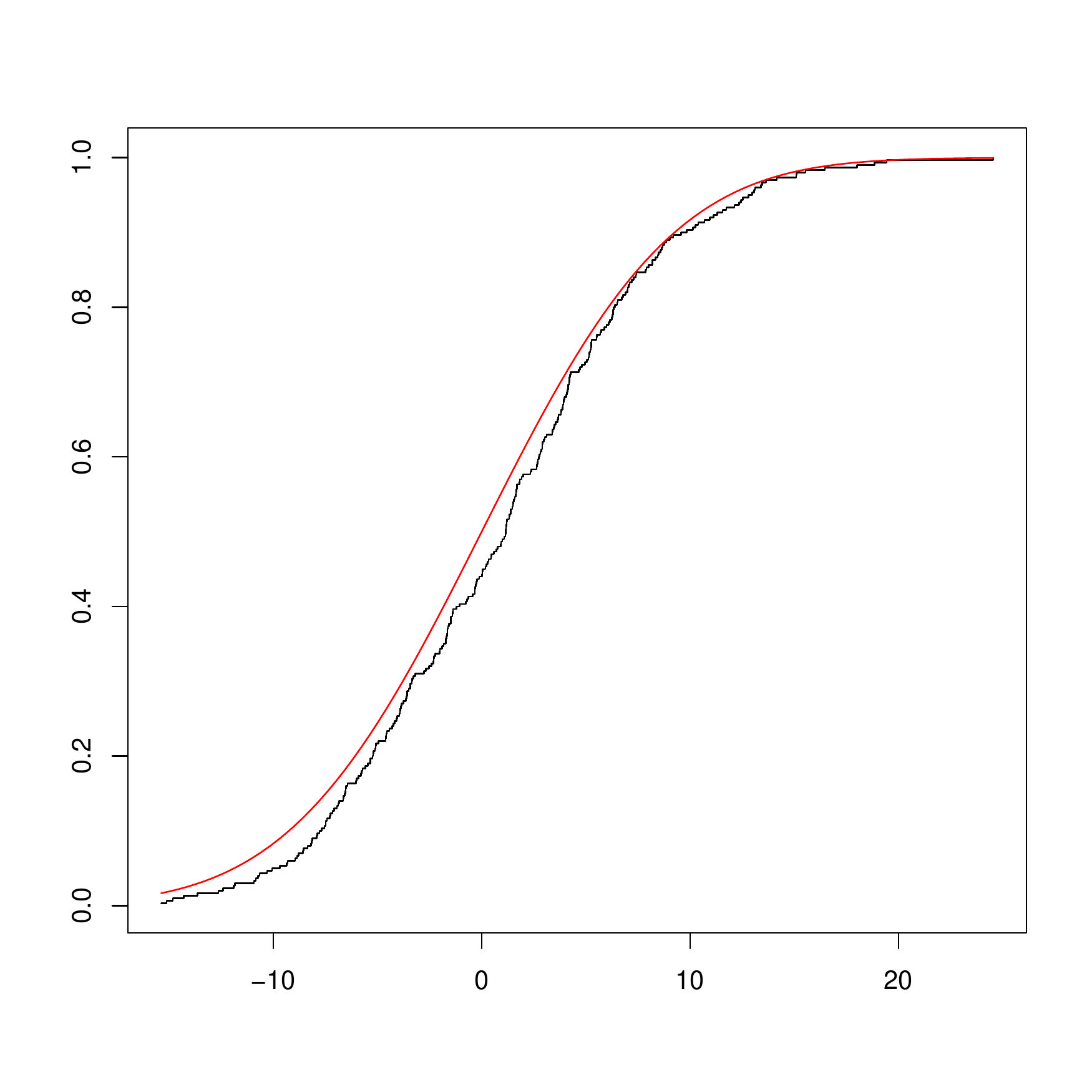}
\includegraphics[width=5cm,pagebox=cropbox,clip]{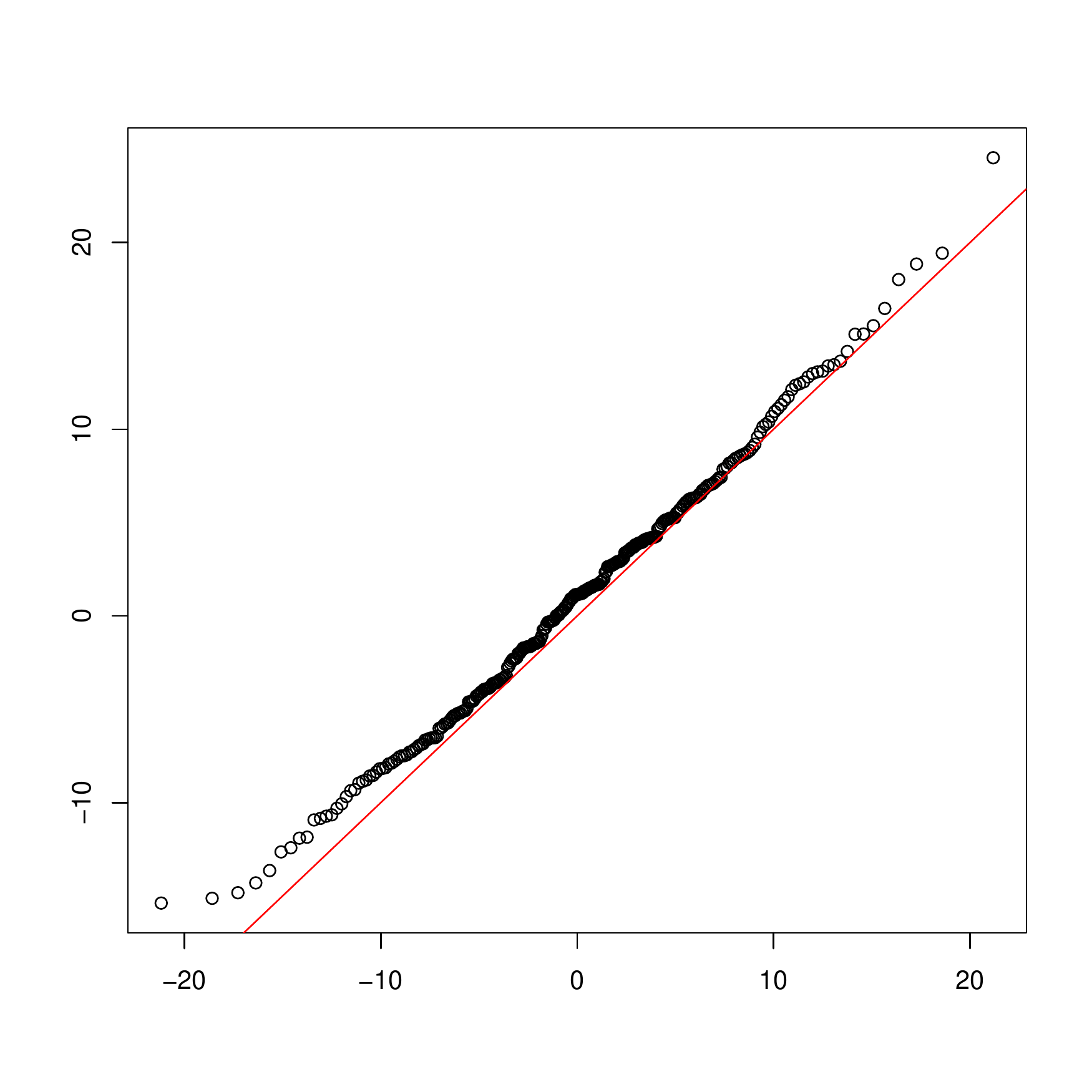}
\includegraphics[width=5cm,pagebox=cropbox,clip]{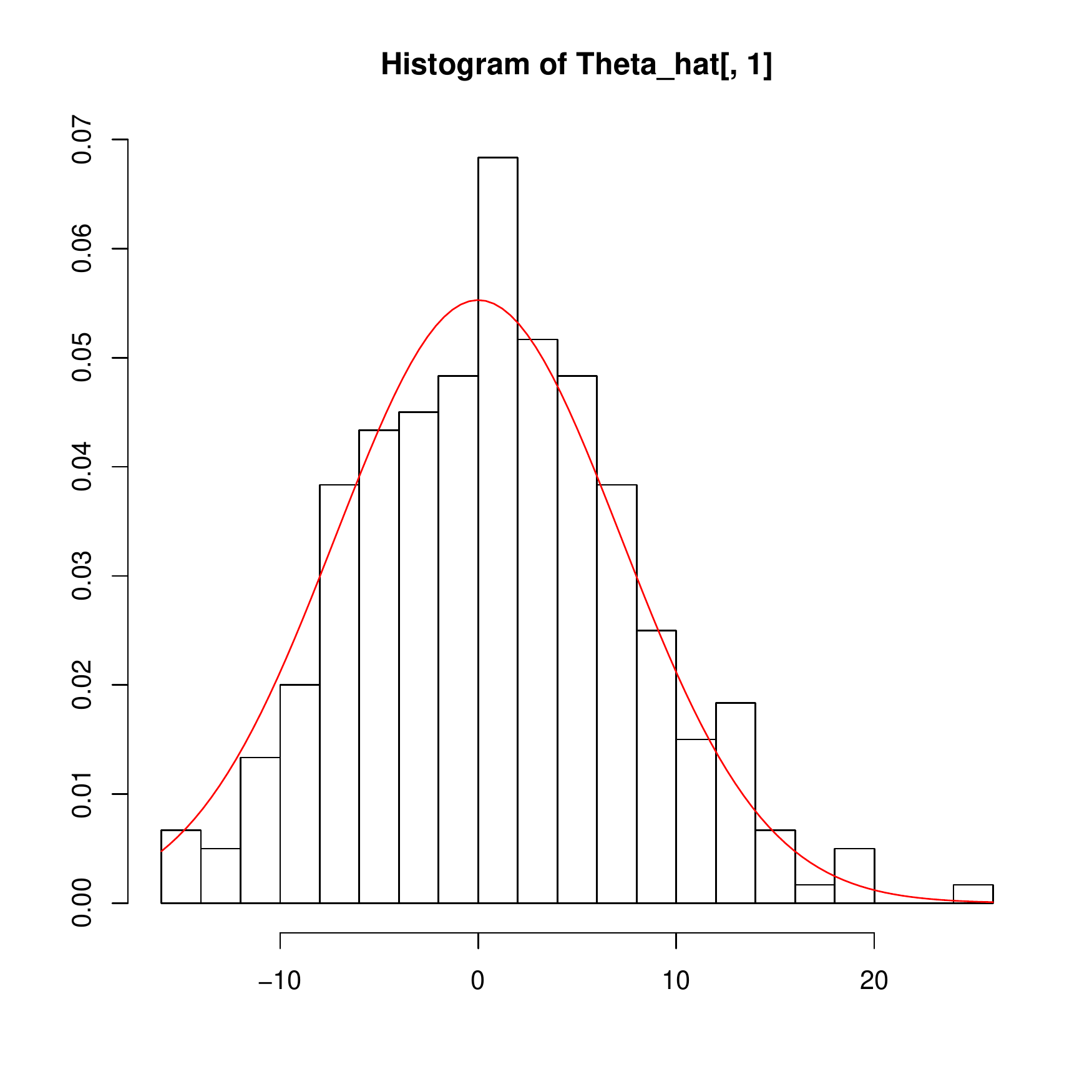}\\
\caption{Simulation results of $\hat{\theta}_1$ \label{fig6}}
\end{center}
\end{figure}

\begin{figure}[h]
\begin{center}
\includegraphics[width=5cm,pagebox=cropbox,clip]{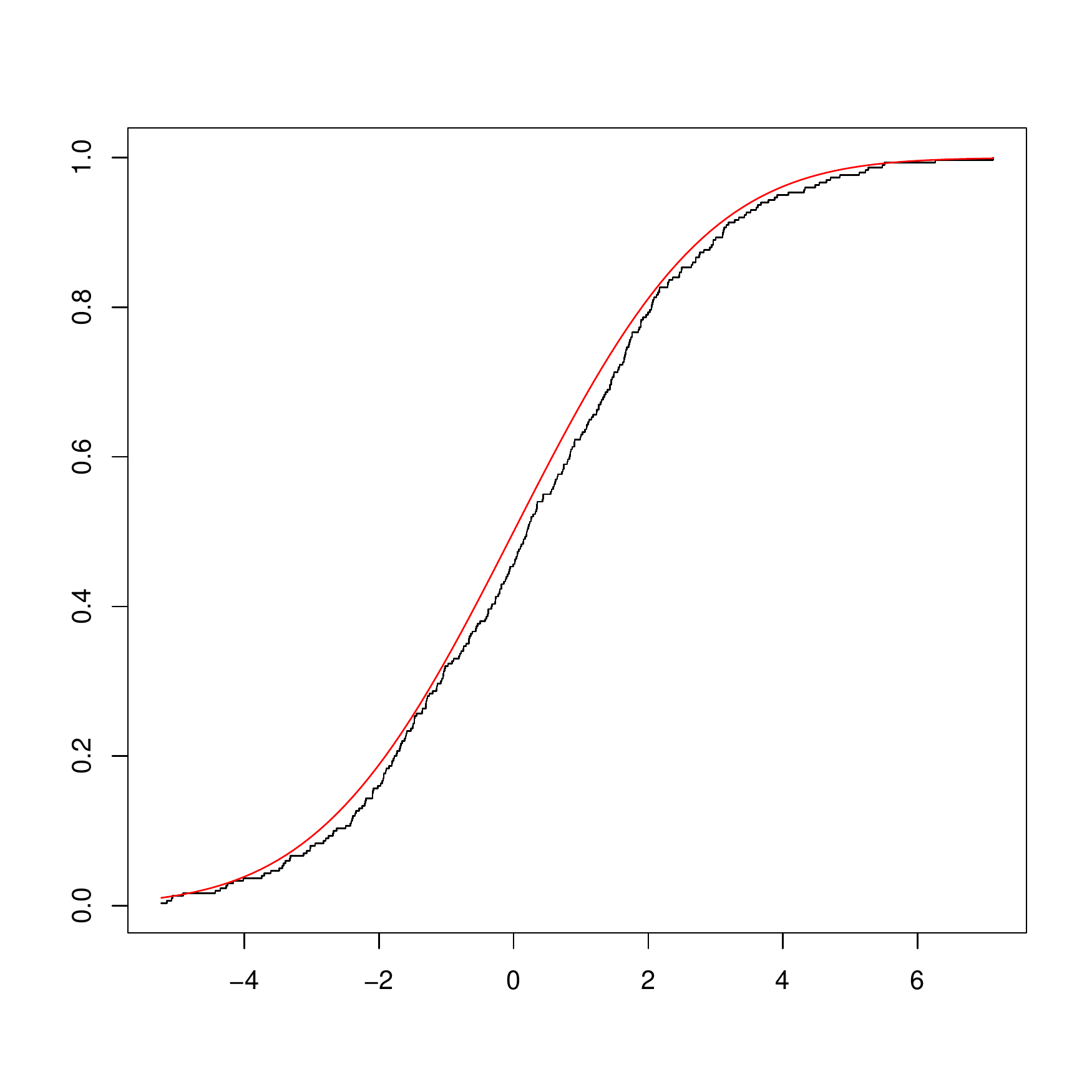}
\includegraphics[width=5cm,pagebox=cropbox,clip]{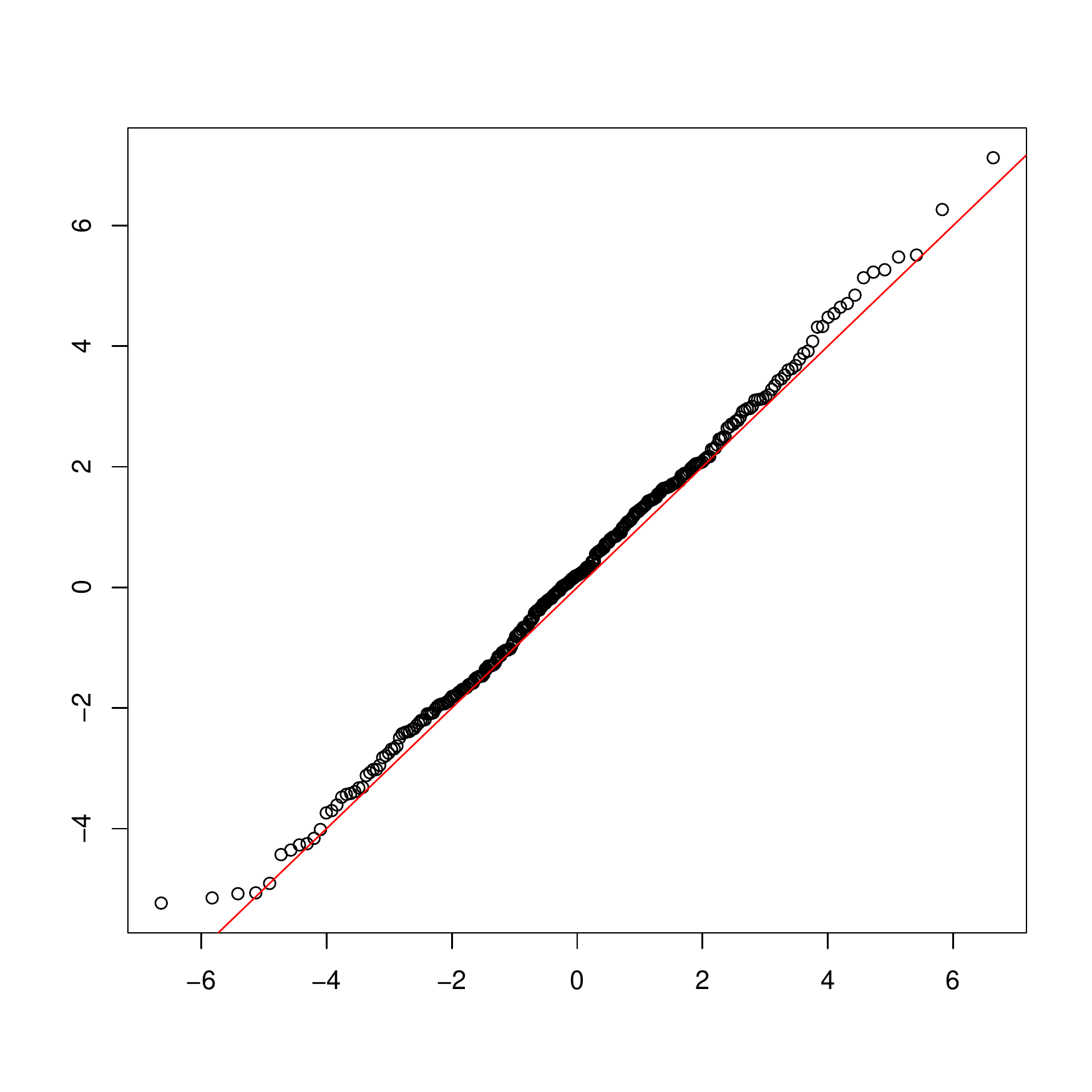}
\includegraphics[width=5cm,pagebox=cropbox,clip]{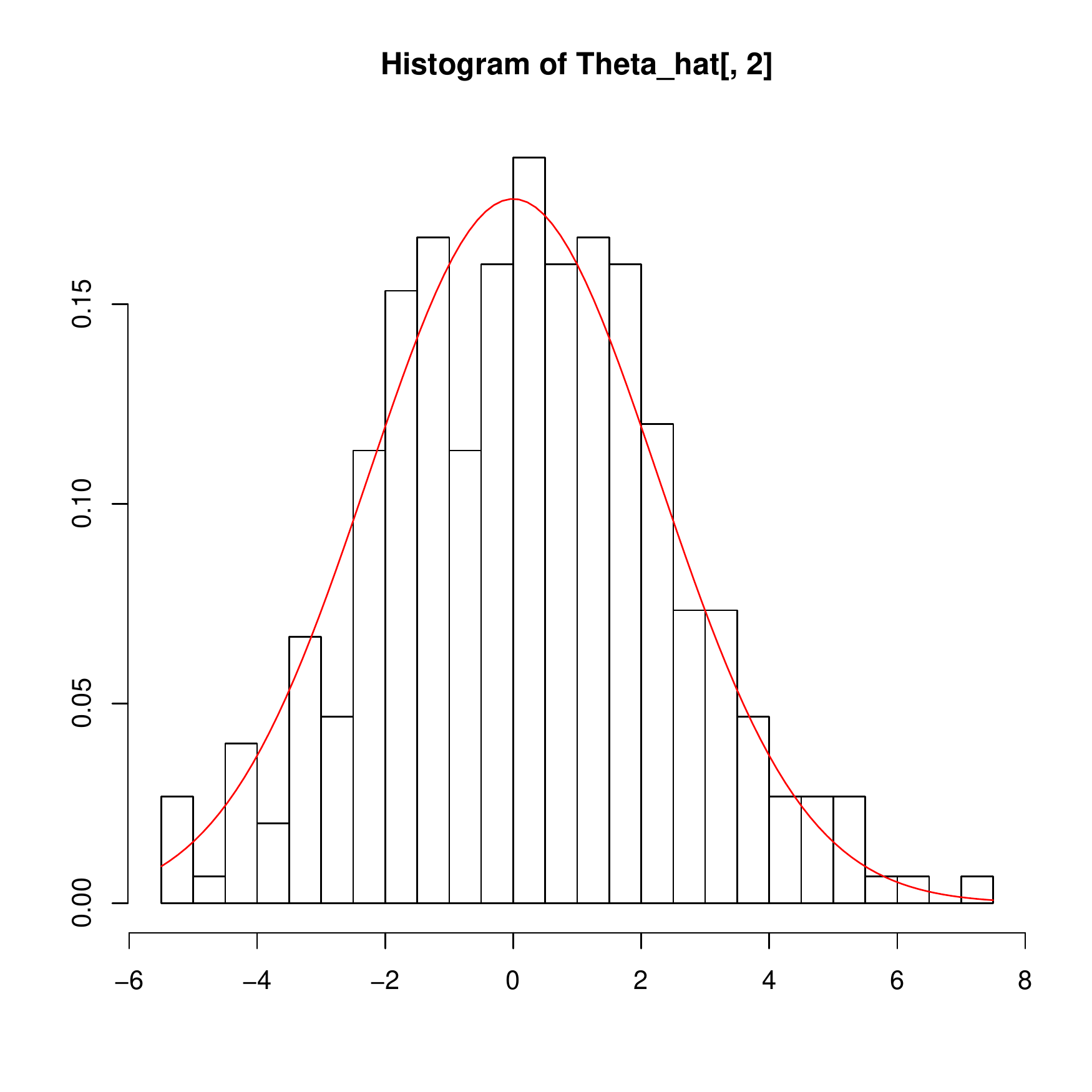}\\
\caption{Simulation results of $\hat{\theta}_2$  \label{fig7}}
\end{center}
\end{figure}


\begin{figure}[t] 
\begin{center}
\includegraphics[width=5cm,pagebox=cropbox,clip]{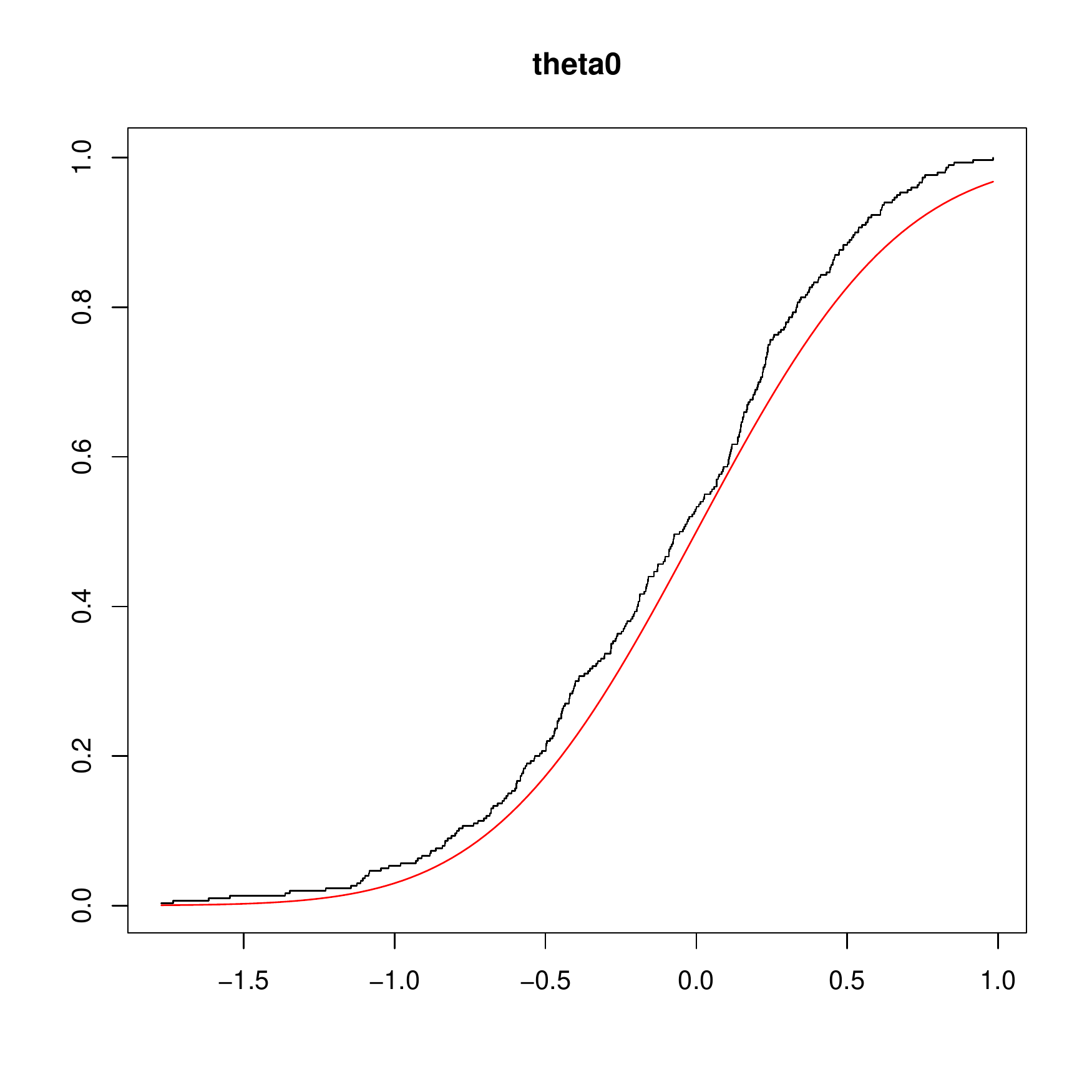}
\includegraphics[width=5cm,pagebox=cropbox,clip]{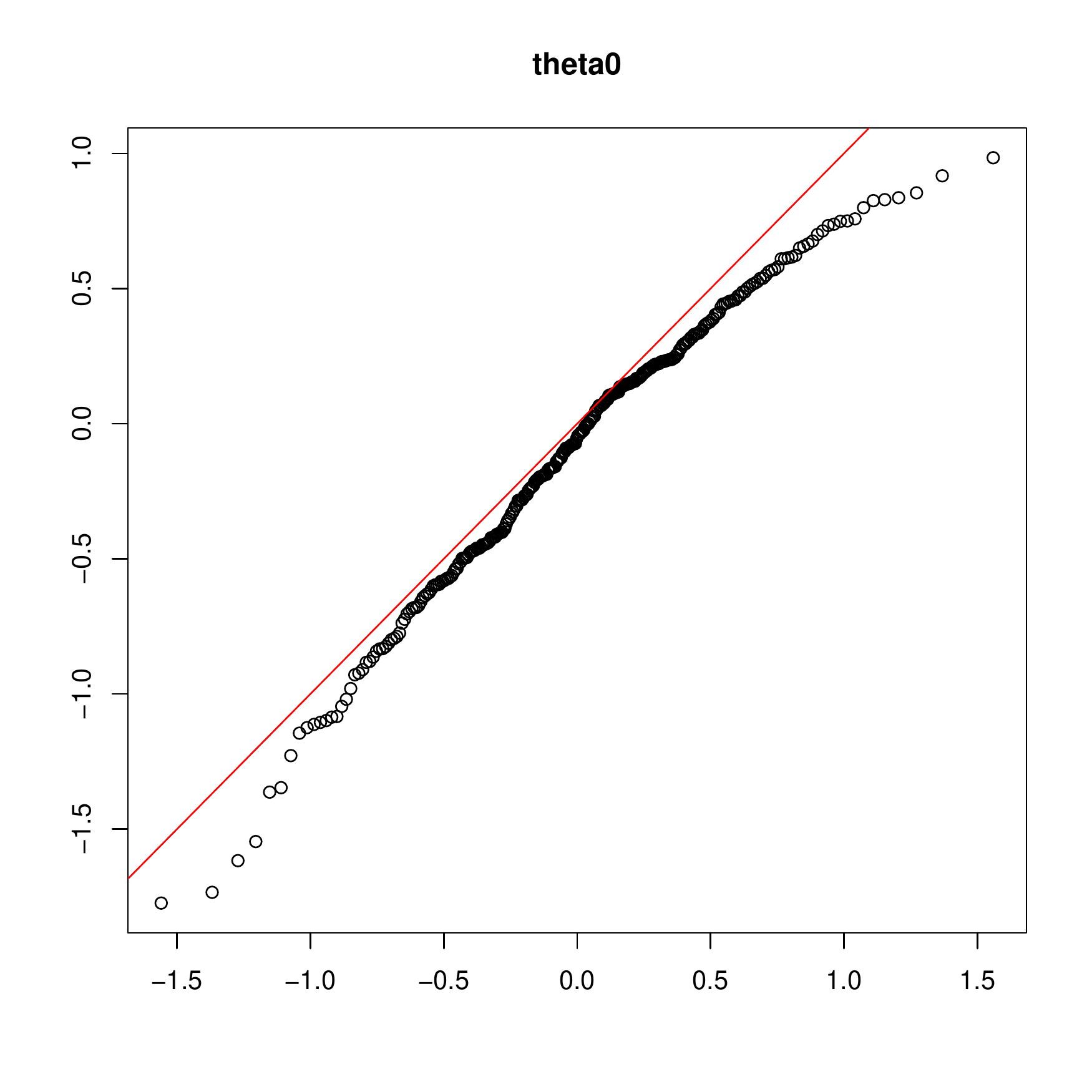}
\includegraphics[width=5cm,pagebox=cropbox,clip]{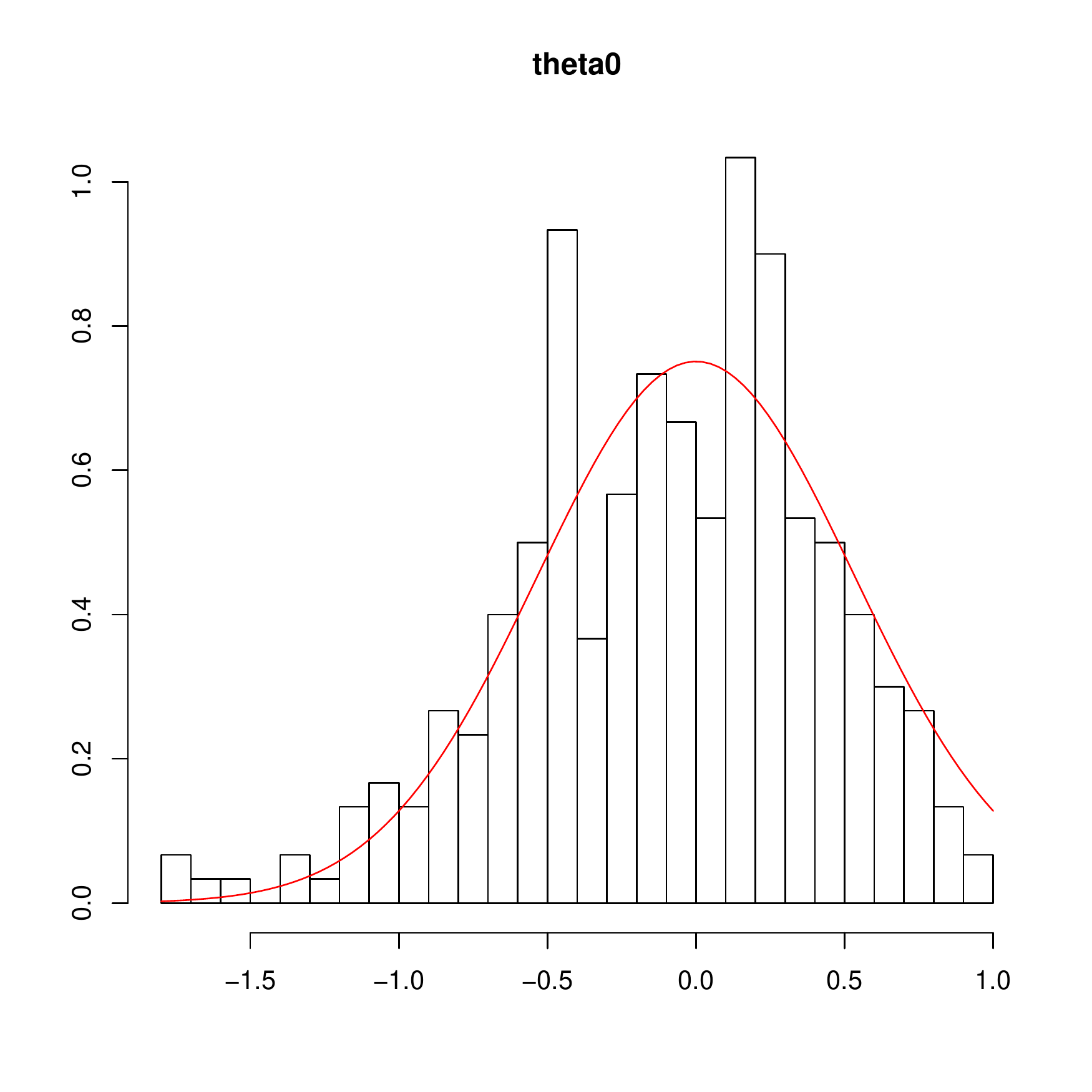}\\
\caption{Simulation results of $\hat{\theta}_0$ \label{fig8}}
\end{center}
\end{figure}


\clearpage

\subsubsection{Summary of example 2}
Table 9 is the simulation results of {the means and the standard s.d.s of} $\hat{\theta}_1$, $\hat{\theta}_2$  and $\hat{\theta}_0$ with $(N, m, N_2) = (10^4, 99, 500)$ from $\epsilon=0.1$ to $0.5$.
{
We can see from Table 9 that
for all $\epsilon$, the estimator of $\theta_0$ has good performance. 
However, it seems from Figure \ref{fig8} that
the asymptotic theory does not work when $\epsilon$ is 0.5.
In this setting,  the estimators work  well if $\epsilon$ is  less than 0.25.
}

\begin{table}[h]
\caption{Simulation results of $\hat{\theta}_1$, $\hat{\theta}_2$  and $\hat{\theta}_0$ with $(N, m, N_2) = (10^4, 99, 500)$ \label{table2}}
\begin{center}
\begin{tabular}{cc|ccc} \hline
	&	&$\hat{\theta}_{1}$&$\hat{\theta}_{2}$&$\hat{\theta}_{0}$
\\ \hline
&{true value} &1 &  0.2 &  3.1
\\ \hline
&mean &1.001&  0.200&  3.102
 \\
$\epsilon=0.1$&{s.d.} & (0.007)& (0.002)& (0.055)
 \\   \hline
&mean &1.001&  0.200&  3.092
 \\
$\epsilon=0.25$&{s.d.} & (0.007)& (0.002)& (0.126)
 \\   \hline
&mean &1.001&  0.200&  3.050
 \\
$\epsilon=0.5$&{s.d.} & (0.007)& (0.002)& (0.258)
 \\   \hline
\end{tabular}
\end{center}
\end{table}





\section{Proofs}
\noindent
{Proof of Theorem 1}. 
{
Let $\sigma_0^2 = \frac{1}{\sqrt{\theta_2}}$, $\eta = \frac{\theta_1}{\theta_2}$,
$(\sigma_0^*)^2 = \frac{1}{\sqrt{\theta_2^*}}$ and $\eta^* = \frac{\theta_1^*}{\theta_2^*}$.
The contrast function is 
\begin{eqnarray*}
U_{N,m}(\sigma_0^2, \eta) = \frac{1}{m} \sum^{m}_{j = 1}  
\left( \frac{1}{\epsilon^2}Z_{j:m} - \frac{\sigma_0^2}{\sqrt{\pi}} 
\exp(- \eta \tilde{y}_{j:m}) \right)^2
\end{eqnarray*}
and the minimum contrast estimator of $\sigma_0^2$ and $\eta$ are defined as 
\begin{eqnarray*}
(\hat{\sigma}_0^2,\hat{\eta}) = \arginf_{\sigma_0^2, \eta}U_{N,m}(\sigma_0^2, \eta).
\end{eqnarray*}
Set
\begin{eqnarray*}
\left(
\begin{array}{cc}
K_{1,1} & K_{1,2} \\
K_{1,2} & K_{2,2} 
\end{array}
\right)
&=& \frac{1}{\theta_2^*} \pi \Gamma V^{-1} U V^{-1}.
\end{eqnarray*}
}
{In an analogous manner to the proof of Theorem 4.2 in Bibinger and Trabs (2020),
it can be shown that} 
under {$\epsilon \rightarrow 0$}, $N \rightarrow \infty$,  
$m \rightarrow \infty$ and 
$m = O(N^\rho)$ 
for
$\rho \in (0, 1/2)$, 
{
$$
\begin{pmatrix}
\sqrt{Nm}(\hat{\sigma}_0^2 - (\sigma_{0}^*)^2) \\
\sqrt{Nm}(\hat{\eta} - \eta^*) \\
\end{pmatrix}
\stackrel{d}{\longrightarrow } N
\left( 
\begin{pmatrix}
0 \\
0 
\end{pmatrix}
,
\begin{pmatrix}
K_{1,1} &K_{1,2} \\
K_{1,2} &K_{2,2}
\end{pmatrix}
\right).
$$
It follows from the delta method that
$$
\begin{pmatrix}
\sqrt{Nm}(\hat{\theta}_2 - \theta_2^*) \\
\sqrt{Nm}(\hat{\theta}_1 - \theta_1^*) \\
\end{pmatrix}
\stackrel{d}{\longrightarrow } N
\left( 
\begin{pmatrix}
0 \\
0 
\end{pmatrix}
,
\begin{pmatrix}
J_{1,1} & J_{1,2} \\
J_{1,2} & J_{2,2}
\end{pmatrix}
\right).
$$

}
{
Next, set $\bar{\delta} =\delta_{N_2}$
and $s_i = s_{i:N_2}$.
Based on the thinned data of the approximate coordinate process
${\bf \hat{x}}_k =\{\hat{x}_k(s_{i:N_2})\}_{i=1,\ldots, N_2}
 =\{\hat{x}_k(s_{i})\}_{i=1,\ldots, N_2}$, 
we have the following quasi log-likelihood function} 
\[
	l_{N_2}(\lambda_k \ | \ {\bf \hat{x}}_k  ) = -\frac{1}{2} \sum_{i=1}^{N_2} 
		\left\{\log \frac{\epsilon^2(1-e^{-2\lambda_k \dlnt})}{2\lambda_k}
		+\frac{\left( \hat{x}_k(s_i)-e^{-\lambda_k \dlnt} \hat{x}_k(s_{i-1})\right)^2}{
		\frac{\epsilon^2(1-e^{-2\lambda_k \dlnt})}{2\lambda_k}}
		\right\}.
\]
{
Let $k=1$, $\lambda=\lambda_1$,
$\hat{\lambda} =\hat{\lambda}_1$, 
${\bf \hat{x}} = {\bf \hat{x}}_1 =\{\hat{x}_1(s_{i:N_2})\}_{i=1,\ldots, N_2}
=\{\hat{x}_1(s_{i})\}_{i=1,\ldots, N_2}$,
${\bf x} = {\bf x}_1 =\{x_1(s_{i:N_2})\}_{i=1,\ldots, N_2}
=\{ x_1(s_{i})\}_{i=1,\ldots, N_2}$, 
}
and
\[
         \Xi(\lambda) = \frac{(1-e^{-2\lambda \dlnt})}{2\lambda \dlnt}.  
\]

{
In order to show the consistency of $\hat{\lambda}$,
it is sufficient to prove that} 
{
\begin{eqnarray}
{\cal Z}:= 
\epsilon^2
\left\{
l_{N_2}(\lambda \ | \ {\bf \hat{x}}  ) 
-
l_{N_2}(\lambda \ | \ {\bf {x}}  ) 
\right\} = o_p(1) \label{consistency-1}
\end{eqnarray}
uniformly in $\lambda$.
Note that 
{
${\cal Z}$ is the difference between the quasi log-likelihood function based on 
the thinned data of the approximate coordinate process
${\bf \hat{x}}$ and  that {based on} the thinned data of the coordinate process
${\bf {x}}$.} 

For the proof of (\ref{consistency-1}), since 
\begin{align*}
{\cal Z}=& \frac{1}{2  \dlnt \Xi(\lambda)}
		\sum_{i=1}^{N_2} \left\{\left(
			\hat{x}_1(s_i)-e^{-\lambda \dlnt} \hat{x}_1(s_{i-1})
		\right)^2 \right. 
	\left. -\left(x_1(s_i)-e^{-\lambda \dlnt} x_1(s_{i-1})\right)^2 \right\}
\end{align*}
and 
\begin{align*}
	& \left(	
		\hat{x}_1(s_i)-e^{-\lambda \dlnt} \hat{x}_1(s_{i-1})
	\right)^2 \\
	=& \left\{
		(\hat{x}_1(s_i)-x_1(s_i))
		-e^{-\lambda \dlnt}(\hat{x}_1(s_{i-1})-x_1(s_{i-1})) \right\}^2 \\
	&+2\{ \hat{x}_1(s_i)-x_1(s_i)
		-e^{-\lambda \dlnt}(\hat{x}_1(s_{i-1})-x_1(s_{i-1})) \}
		(x_1(s_i)-e^{-\lambda \dlnt} x_1(s_{i-1})) \\
	&+(x_1(s_i)-e^{-\lambda \dlnt} x_1(s_{i-1}))^2,
\end{align*}
one has that 
\begin{align}
	{\cal Z} 
	=&{  \frac{1}{2 \dlnt \Xi(\lambda)}
		\sum_{i=1}^{N_2} \left\{ \left(
		\hat{x}_1(s_i)-x_1(s_i) \right)
		-e^{-\lambda \dlnt} \left(\hat{x}_1(s_{i-1})-x_1(s_{i-1}) \right)
		\right\}^2}  \nonumber \\
	&{ +2 \frac{1}{2 \dlnt \Xi(\lambda)}
		\sum_{i=1}^{N_2} \left\{\hat{x}_1(s_i)-\hat{x}_1(s_{i-1})
		-\left( x_1(s_i)-x_1(s_{i-1})\right)\right\}
		\left( x_1(s_i)-e^{-\lambda \dlnt} x_1(s_{i-1})\right) }
		\nonumber \\
	&{ +2\frac{1}{2 \dlnt \Xi(\lambda)}
		\sum_{i=1}^{N_2} (1-e^{-\lambda \dlnt})
		\left( \hat{x}_1(s_{i-1})-x_1(s_{i-1})\right)
		\left( x_1(s_i)-e^{-\lambda \dlnt} x_1(s_{i-1}) \right)  }\nonumber \\
	=:& {\cal W}_1+{\cal W}_2+{\cal W}_3. \nonumber 
\end{align}
First, we will estimate ${\cal W}_1$.
By setting $g_1(t, y, \eta) = X_t(y)\sqrt{2} \sin (\pi y) \exp\left\{\frac{\eta}{2}y \right\}$, 
$y_{j} = y_{j:M}$ and $\hat{\eta} = \frac{\hat{\theta}_1}{\hat{\theta}_2}$,
we obtain that 
\begin{eqnarray*}	
x_1(t) &=& \int_0^1 X_t(y)\sqrt{2} \sin (\pi y) 
	\exp\left\{\frac{\eta}{2}y \right\} dy =  \int_0^1 g_1(t, y, \eta) dy,
	\\
\hat{x}_1(t) &=& \frac{1}{M} \sum_{j=1}^M X_t(y_j) \sqrt{2} \sin (\pi y_j) \exp \left\{ \frac{\hat{\eta}}{2} y_j \right\} 
=  \frac{1}{M} \sum_{j=1}^M g_1(t, y_j, \hat{\eta}).
\end{eqnarray*}
Moreover,  it follows that
\begin{eqnarray*}
	{\cal Z}_1 &:=& \frac{1}{2 \dlnt \Xi(\lambda)} \sum_{i=1}^{N_2} 
		\left( x_1(s_i)-\hat{x}_1(s_i) \right)^2  \\
	&=& \frac{1}{2 \dlnt \Xi(\lambda)} \sum_{i=1}^{N_2}
		\left\{
			M \frac{1}{M} \sum_{j=1}^M \int_{\frac{j-1}{M}}^{\frac{j}{M}}
			\{g_1(s_i,y,\eta)-g_1(s_i,y_j,\hat{\eta})\}dy
		\right\}^2  \\
	&\leq& \frac{1}{2 \dlnt \Xi(\lambda)} \sum_{i=1}^{N_2} M^2 \frac{1}{M}
		\sum_{j=1}^{M} \frac{1}{M} \int_{\frac{j-1}{M}}^{\frac{j}{M}}
		\{g_1(s_i,y,\eta)-g_1(s_i,y_j,\hat{\eta})\}^2 dy  \\
	&=& \frac{1}{2 \dlnt \Xi(\lambda)} \sum_{i=1}^{N_2} \sum_{j=1}^M
		\int_{\frac{j-1}{M}}^{\frac{j}{M}} \{g_1(s_i,y,\eta)-g_1(s_i,y_j,\hat{\eta})\}^2 dy.
\end{eqnarray*}
Noting that 
\begin{eqnarray}
	& & g_1(s_i,y,\eta)-g_1(s_i,y_j,\hat{\eta}) \nonumber\\
	&=& X_{s_i}(y)\sqrt{2}\sin(\pi y) \exp \left\{ \frac{\eta}{2} y\right\}
		-X_{s_i}(y_j)\sqrt{2}\sin(\pi y_j) \exp \left\{ \frac{\hat{\eta}}{2} y_j
		\right\} \nonumber\\
	&=&{  \left( X_{s_i}(y)-X_{s_i}(y_j) \right) \sqrt{2} \sin(\pi y) 
		\exp \left\{ \frac{\eta}{2} y\right\}  } \nonumber \\
	& &{ + X_{s_i}(y_j) \left( \sqrt{2}\sin(\pi y) 
		\exp\left\{ \frac{\eta}{2} y \right\} -\sqrt{2}\sin(\pi y_j) 
		\exp\left\{ \frac{\eta}{2} y_j\right\} \right) } \nonumber\\
	& &{ + X_{s_i}(y_j) \sqrt{2}\sin(\pi y_j)
		\left( \exp\left\{\frac{\eta}{2} y_j \right\} 
		-\exp\left\{\frac{\hat{\eta}}{2} y_j\right\}\right) } \nonumber \\
	&=&: \bar{g}^{(1)}(s_i,y) + \bar{g}^{(2)}(s_i,y) + \bar{g}^{(3)}(s_i,y_j,\hat{\eta}),  \nonumber
\end{eqnarray}
one has that 
\begin{eqnarray*}
	{\cal Z}_1 &\leq& \frac{C}{2 \dlnt \Xi(\lambda)} \sum_{i=1}^{N_2} \sum_{j=1}^M
		\int_{\frac{j-1}{M}}^{\frac{j}{M}} \{\bar{g}_1(s_i,y)\}^2 dy  \nonumber \\
		& & + \frac{C}{2 \dlnt \Xi(\lambda)} \sum_{i=1}^{N_2} \sum_{j=1}^M
		\int_{\frac{j-1}{M}}^{\frac{j}{M}} \{\bar{g}_2(s_i,y)\}^2 dy  \nonumber \\
		& & + \frac{C}{2 \dlnt \Xi(\lambda)} \sum_{i=1}^{N_2} \sum_{j=1}^M
		\int_{\frac{j-1}{M}}^{\frac{j}{M}} \{\bar{g}_3(s_i,y_j,\hat{\eta})\}^2 dy \nonumber \\
		&=& : {\cal Z}_{11} + {\cal Z}_{12} + {\cal Z}_{13} . \nonumber \\
\end{eqnarray*}
It follows that 
\begin{eqnarray*} 
  E\left[{\cal Z}_{11} \right] & = & \frac{C}{2 \dlnt \Xi(\lambda)} \sum_{i=1}^{N_2} \sum_{j=1}^M
		\int_{\frac{j-1}{M}}^{\frac{j}{M}} E\left[ \{\bar{g}^{(1)}(s_i,y)\}^2 \right] dy \\
		& = & \frac{C}{2 \dlnt \Xi(\lambda)} \sum_{i=1}^{N_2} \sum_{j=1}^M
		\int_{\frac{j-1}{M}}^{\frac{j}{M}} E\left[ \left( X_{s_i}(y)-X_{s_i}(y_j)\right)^2 \right] 2 \sin^2(\pi y_j) \exp(\eta y) dy \\
		& = & \frac{C}{2 \dlnt \Xi(\lambda)} \sum_{i=1}^{N_2} \sum_{j=1}^M \frac{1}{M} 
		C_1 \left( \frac{\epsilon^2}{M^{1-\rho_1}} + \frac{1}{M^2}\right) \\
		& \leq & \frac{C}{2 \Xi(\lambda)} \left( \frac{\epsilon^2 N_2^2}{M^{1-\rho_1}} + \frac{N_2^2}{M^2} \right)
\end{eqnarray*}
and that 
\begin{eqnarray*} 
  E\left[{\cal Z}_{12} \right] & = & \frac{C}{2 \dlnt \Xi(\lambda)} \sum_{i=1}^{N_2} \sum_{j=1}^M
		\int_{\frac{j-1}{M}}^{\frac{j}{M}} E\left[ \{\bar{g}^{(2)}(s_i,y)\}^2 \right] dy \\
		& = & \frac{C}{2 \dlnt \Xi(\lambda)} \sum_{i=1}^{N_2} \sum_{j=1}^M
		\int_{\frac{j-1}{M}}^{\frac{j}{M}} E\left[ \left( X_{s_i}(y) \right)^2 \right] 2 (\sin^2(\pi y) \exp(\frac{\eta}{2} y) - \sin^2(\pi y_j) \exp(\frac{\eta}{2} y_j))^2 dy \\
		& = & \frac{C}{2 \dlnt \Xi(\lambda)} \sum_{i=1}^{N_2} \sum_{j=1}^M \frac{1}{M} C_2(y - y_j)^2 \\
		& \leq & \frac{C}{2 \Xi(\lambda)} \frac{N_2^2}{M^2}.
\end{eqnarray*}
{
Let
$
R(y_j, \hat{\eta})
:=\int_0^1 \frac{y_j}{2} \exp \left\{ \frac{y_j}{2} (\eta + u ( \hat{\eta} - \eta)) \right\} du
$
and  
\begin{align*}
&	{\cal Z}_2 :=	\frac{1}{{\dlnt}} \sum_{i=1}^{N_2} \frac{1}{M}
		\sum_{j=1}^M X_{s_i}^2(y_j)2\sin^2(\pi y_j) 
(R(y_j, \hat{\eta}))^2 \frac{1}{N m}.
\end{align*}
{We} obtain that for $\epsilon_1>0$, 
\begin{align}
%
P[|{\cal Z}_2| > \varepsilon_1] 
		& \leq C_3 \frac{N_2^2}{N m} \frac{1}{\epsilon_1}.
		\label{Z-Res1}
\end{align}
}
It follows from  (\ref{Z-Res1}) that 
\begin{eqnarray*}
{\cal Z}_{13} & = & \frac{C}{2 \dlnt \Xi(\lambda)} \sum_{i=1}^{N_2} \sum_{j=1}^M \frac{1}{M}
		\left( X_{s_i}(y_j) \right)^2 2 \sin^2(\pi y_j)( \exp(\frac{\eta}{2} y_j)  - \exp(\frac{\hat{\eta}}{2} y_j))^2 \\
		& = & \frac{C}{2 \Xi(\lambda)} \frac{1}{\dlnt} \sum_{i=1}^{N_2} \frac{1}{M} \sum_{j=1}^M	 \left( X_{s_i}(y_j) \right)^2 2\sin^2(\pi y_j) 
(R(y_j, \hat{\eta}))^2 \frac{1}{N m}  (\sqrt{Nm}(\hat{\eta} - \eta))^2 \\
		& = & \frac{C}{2 \Xi(\lambda)} {\cal Z}_2 (\sqrt{Nm}(\hat{\eta} - \eta))^2  = O_p(\frac{N_2^2}{Nm}).
\end{eqnarray*}
Therefore, 
\begin{eqnarray*}
	{\cal Z}_1 &=& O_p\left(\frac{N_2^2}{M^2}\right)
		+O_p \left(\epsilon^2 \frac{N_2^2}{M^{1-\rho_1}} \right)
		+O_p\left(\frac{N_2^2}{N m}\right), \\
	{\cal W}_1 &=& O_p\left(\frac{N_2^2}{M^2}\right)
		+O_p \left(\epsilon^2 \frac{N_2^2}{M^{1-\rho_1}} \right)
		+O_p\left(\frac{N_2^2}{N m}\right). 
\end{eqnarray*}

For the estimate of ${\cal W}_3$,  we have that
\begin{align*}
	{\cal W}_3^2 \leq& \frac{(1 - e^{-\lambda \dlnt)^2}}{({\dlnt})^2 \Xi(\lambda)^2}
		\sum_{i=1}^{N_2} \left( \hat{x}_1(s_i)-x_1(s_{i-1}) \right)^2 
		\sum_{i=1}^{N_2} \left( x_1(s_i)-e^{-\lambda \dlnt}x_1(s_{i-1}) \right)^2 \\
	=& C \dlnt {\cal Z}_1 O_p(1) \\
	=& O_p \left(\frac{{N_2}}{M^{1-\rho_1}} \right)
		+O_p\left(\frac{{N_2}}{N m}\right). 
\end{align*}

For the estimate of { ${\cal W}_2$}, 
setting $\Delta X_{s_i}(y) = X_{s_i}(y)-X_{s_{i-1}}(y)$, 
we obtain that 
\begin{align*}
     {\cal W}_2
	&= \frac{1}{\dlnt \Xi(\lambda)} \sum_{i=1}^{N_2} \frac{1}{M}
		\sum_{j=1}^M 
                  \Delta X_{s_i}(y_j) 
		\sqrt{2}\sin(\pi y_j)
		\left(\exp\left\{\frac{\hat{\eta}}{2}y_j \right\}
		-\exp\left\{\frac{\eta}{2} y_j\right\}\right) 
		\left(x_1(s_i)-e^{-\lambda \dlnt}x_1(s_{i-1})\right) \nonumber \\
		&+ \frac{1}{\dlnt \Xi(\lambda)} 
		\sum_{i=1}^{N_2} \sum_{j=1}^{M}
		\int_{\frac{j-1}{M}}^{\frac{j}{M}}
                  \Delta X_{s_i}(y_j) 
		\left(\sqrt{2}\sin(\pi y_j)-\sqrt{2}\sin(\pi y)\right)
		\exp\left\{\frac{\eta}{2}y_j\right\} dy 
                 \left(x_1(s_i)-e^{-\lambda \dlnt}x_1(s_{i-1})\right) \nonumber \\
	&+ \frac{1}{\dlnt \Xi(\lambda)}
		\sum_{i=1}^{N_2} \sum_{j=1}^{M}
		\int_{\frac{j-1}{M}}^{\frac{j}{M}}
                  \Delta X_{s_i}(y_j) 
		\sqrt{2}\sin(\pi y) \left(\exp\left\{\frac{\eta}{2}y_j\right\}
		-\exp\left\{\frac{\eta}{2}y\right\} \right)dy 
		\left(x_1(s_i)-e^{-\lambda \dlnt}x_1(s_{i-1})\right) \nonumber \\
	&+ \frac{1}{\dlnt \Xi(\lambda)}
		\sum_{i=1}^{N_2} \sum_{j=1}^{M}
		\int_{\frac{j-1}{M}}^{\frac{j}{M}}
		\left\{ \Delta X_{s_i}(y_j) -\Delta X_{s_i}(y) \right\}		
		\sqrt{2}\sin(\pi y)\exp\left\{\frac{\eta}{2}y\right\}dy 
		\left(x_1(s_i)-e^{-\lambda \dlnt}x_1(s_{i-1})\right) \nonumber \\
	&=: \ ({\rm I})+({\rm II})+({\rm III})+({\rm IV}). \nonumber
\end{align*}

For the estimate of (I), one has that
\begin{align}
	({\rm I})^2 \leq&{  \frac{1}{\dlnt^2 \Xi(\lambda)^2} \sum_{i=1}^{N_2}\frac{1}{M}
		\sum_{j=1}^{M} 
                \left( \Delta X_{s_i}(y_j) \right)^2 
                \left(\sqrt{2}\sin(\pi y_j) R(y_j,\hat{\eta} )  \right)^2 
		\frac{1}{N m}   } \nonumber \\
	&{  \sum_{i=1}^{N_2} \left(x_1(s_i)-e^{-\lambda \dlnt}x_1(s_{i-1})\right)^2 
		\times \left( \sqrt{N m}(\hat{\eta}-\eta)\right)^{2} } \nonumber \\
		=:&{  \mathcal{B}_	1 \times \sum_{i=1}^{N_2} \left(x_1(s_i)-e^{-\lambda \dlnt}x_1(s_{i-1})\right)^2  \times \left( \sqrt{N m}(\hat{\eta}-\eta)\right)^{2}.  } \nonumber
\end{align}
{
Since 
$E \left[ 
\sum^{N_2}_{i=1}
\left( 
\Delta X_{s_i}(y_j)
\right)^2 
\right] = 
O \left( \sqrt{N_2} \epsilon^2 \right) + O \left(\sqrt{{\dlnt}} \right)
$,
one has that
\begin{align*}
	{ \mathcal{B}_1} \leq&  C_1 \frac{({\dlnt})^{-2}}{N m}
	\sum_{i=1}^{N_2}\frac{1}{M}
		\sum_{j=1}^{M} \left( X_{s_i}(y_j)-X_{s_{i-1}}(y_j) \right)^2 
	= O_p \left( {\frac{\epsilon^2 {N_2}^{\frac{5}{2}}}{N m}  } \right)
	 + O_p \left( \frac{N_2^{\frac{3}{2}}}{N m} \right). 
\end{align*}
}

\begin{en-text}
Let  $\eta_1 >0$ and $\epsilon_1>0$.
On $J=\{|\hat{\eta}-\eta| < \eta_1\}$,
\begin{align*}
	{ \mathcal{B}_1} \leq&  C_1 \frac{({\dlnt})^{-2}}{N m}
	\sum_{i=1}^{N_2}\frac{1}{M}
		\sum_{j=1}^{M} \left( X_{s_i}(y_j)-X_{s_{i-1}}(y_j) \right)^2 
	= O_p \left( {\frac{\epsilon^2 {N_2}^{\frac{5}{2}}}{N m}  } \right)
	 + O_p \left( \frac{N_2^{\frac{3}{2}}}{N m} \right)
\end{align*}
because 
$E \left[ 
\sum^{N_2}_{i=1}
\left( 
\Delta X_{s_i}(y_j)
\right)^2 
\right] = 
{O \left( \sqrt{N_2} \epsilon^2 \right)} + O \left(\sqrt{{\dlnt}} \right).
$
\end{en-text}

For the estimate of (II), we obtain that
\begin{align*}
	({\rm II})^2 \leq& C_1 \frac{1}{\dlnt^2 \Xi(\lambda)^2}
		\sum_{i=1}^{N_2} \frac{1}{M}
		\sum_{j=1}^{M}
		 \left( \Delta X_{s_i}(y_j) \right)^2 
		  |y-y_j|^2 
		\sum_{i=1}^{N_2} \left(x_1(s_i)-e^{-\lambda \dlnt}x_1(s_{i-1})\right)^2 \\
	=& O_p \left( {\frac{\epsilon^2 {N_2}^{\frac{5}{2}}}{N m}  } \right)
	 + O_p \left( \frac{N_2^{\frac{3}{2}}}{N m} \right) . 
\end{align*}

For the estimate of (III), we have that 
\begin{align*}
	({\rm III})^2 \leq& C_1 \frac{1}{\dlnt^2 \Xi(\lambda)^2}
		\sum_{i=1}^{N_2} \frac{1}{M} \sum_{j=1}^{M}
		 \left( \Delta X_{s_i}(y_j) \right)^2 
		|y-y_j|^2
		\sum_{i=1}^{N_2} \left(x_1(s_i)-e^{-\lambda \dlnt}x_1(s_{i-1})\right)^2 \\
	=& O_p \left( {\frac{\epsilon^2 {N_2}^{\frac{5}{2}}}{N m}  } \right)
	 + O_p \left( \frac{N_2^{\frac{3}{2}}}{N m} \right) . 
\end{align*}

For the estimate of (IV),  it follows that 
\begin{eqnarray*}
	({\rm IV})^2 &\leq&  \frac{1}{\dlnt^2 \Xi(\lambda)^2} {N_2}
		\sum_{i=1}^{N_2} M^2 \frac{1}{M} \\
		& & \times \sum_{j=1}^{M} \frac{1}{M}
		\int_{\frac{j-1}{M}}^{\frac{j}{M}}
		\left( \left( \sum_{k=1}^\infty (x_k(s_i)-x_k(s_{i-1}))(e_k(y_j) -e_k(y)) \right)^2
		(x_1(s_i)-{e^{-\lambda \delta} }x_1(s_{i-1}))^2 \right) dy.
\end{eqnarray*}
Let
\begin{equation*}
	{\cal Z}_3 := \left( \sum_{k=1}^\infty (x_k(s_i)-x_k(s_{i-1}))(e_k(y_j) -e_k(y)) \right)^2
		(x_1(s_i)-{e^{-\lambda \delta} }x_1(s_{i-1}))^2. 
\end{equation*}
We obtain that 
\begin{align*}
	E[{\cal Z}_3] 
	\leq& C_1   \frac{1}{M^{1-\rho_1}}\left(\frac{\epsilon^2}{N_2^2} + \frac{\epsilon}{N_2^3} + \frac{1}{N_2^4}\right) .
\end{align*}
Moreover, one has that 

\begin{align*}
	E[({\rm IV})^2] \leq& C {\dlnt}^{-2} {N_2^2} \cdot \frac{1}{M^{1-\rho_1}}\left(\frac{\epsilon^2}{N_2^2} + \frac{\epsilon}{N_2^3} + \frac{1}{N_2^4}\right)
		= \frac{{N_2^2\epsilon^2 + N_2\epsilon^2 + 1}}{ M^{1-\rho_1}}. 
\end{align*}
Therefore, 
\begin{equation*}
 \epsilon^2 \left\{ 
		 l_{N_2}(\lambda \ | \ {\bf \hat{x}}  )
		- l_{N_2}(\lambda \ | \ {\bf {x}}  ) \right\} 	
= o_p(1)
\end{equation*}
uniformly in $\lambda$. 
This completes the proof of consistency of 
$\hat{\lambda}$.

Next, in order to prove the  {asymptotic normality of $\hat{\lambda}$,}
we consider the following derivatives of the quasi  log-likelihood function with respect to the parameter $\lambda$.
{Note that}
\begin{align*}
	\partial_\lambda l_{N_2}(\lambda \ | \ {\bf \hat{x}}  ) 
		=& -\frac{1}{2} \sum_{i=1}^{N_{2}}
		\left\{
			\frac{\partial_{\lambda}\Xi(\lambda)}{\Xi(\lambda)}
			-\frac{(\partial_\lambda \Xi(\lambda)) \left(\hat{x}_1(s_i)-e^{-\lambda \dlnt} \hat{x}_1(s_{i-1})\right)^2}{\epsilon^2 \Xi(\lambda)^2 \dlnt} \right.\\
			& \left. +\frac{2\dlnt e^{-\lambda \dlnt} \hat{x}_1(s_{i-1})\left(\hat{x}_1(s_i)-e^{-\lambda \dlnt} \hat{x}_1(s_{i-1})\right)}{\epsilon^2 \Xi(\lambda) \dlnt}
    \right\},
\end{align*}
and 
\begin{align*}
	\partial_{\lambda}^2 l_{N_2}(\lambda, \epsilon^2 \ | \ {\bf \hat{x}}  ) 
		=&-\frac{1}{2}\sum_{i=1}^{N_2} \left\{ \partial_\lambda\left(
		\frac{\partial_\lambda\Xi(\lambda)}{\Xi(\lambda)}\right)
		-\partial_\lambda \left(
			\frac{\partial_\lambda \Xi(\lambda)}{\Xi(\lambda)^2}
		\right)
		\frac{\left(\hat{x}_1(s_i)-e^{-\lambda \dlnt} \hat{x}_1(s_{i-1})\right)^2}{\epsilon^2 \dlnt} 
	\right. \\
		&\left. 
		+\left( \frac{\partial_\lambda \Xi(\lambda)}{\Xi(\lambda)^2} \right)
		\frac{2 \dlnt e^{-\lambda \dlnt} \hat{x}_1(s_{i-1})
		\left(\hat{x}_1(s_i)-e^{-\lambda \dlnt}\hat{x}_1(s_{i-1})\right)} {\epsilon^2 \dlnt} 
	\right.\\
		&\left. -\frac{\partial_\lambda\Xi(\lambda)}{\Xi(\lambda)^2}
		\frac{2 \dlnt e^{-\lambda \dlnt} \hat{x}_1(s_{i-1})
		\left(\hat{x}_1(s_i)-e^{-\lambda \dlnt}\hat{x}_1(s_{i-1})\right)}{\epsilon^2 \dlnt}
	\right.\\
		&\left. -\frac{1}{\Xi(\lambda)}
		\frac{2 {\dlnt}^2 e^{-\lambda \dlnt} \hat{x}_1(s_{i-1}) 
		\left(\hat{x}_1(s_i)-e^{-\lambda \dlnt} \hat{x}_1(s_{i-1}) \right)}{\epsilon^2 \dlnt}
	\right.\\
		&\left. -\frac{1}{\Xi(\lambda)} 
		\frac{2 {\dlnt}^2 e^{-2 \lambda \dlnt}\hat{x}_1^2 (s_{i-1})}{\epsilon^2 \dlnt}	
\right\}.
\end{align*}

\begin{en-text}
Let the discrete data of the coordinate process ${\bf x} =\{ x_1(s_{i:N_2:T})\}_{i=1,\ldots, N_2}$. 
The difference between the score functions of the volatility parameter $\epsilon^2$ 
based on ${\bf \hat{x}}$ and  ${\bf {x}}$ 
is as follows.
\begin{align*}
	(A) :=& \frac{1}{\sqrt{N_2}} \left\{ 
		\partial_{\epsilon^2} l_{N_2}(\lambda,\epsilon^2 \ | \ {\bf \hat{x}}  )
		-\partial_{\epsilon^2} l_{N_2}(\lambda,\epsilon^2 \ | \ {\bf {x}}  ) \right\} \\
	=& \frac{1}{\sqrt{N_2}} \frac{1}{2\epsilon^4 \dlnt \Xi(\lambda)}
		\sum_{i=1}^{N_2} \left\{\left(
			\hat{x}_1(s_i)-e^{-\lambda \dlnt} \hat{x}_1(s_{i-1})
		\right)^2 \right. 
	\left. -\left(x_1(s_i)-e^{-\lambda \dlnt} x_1(s_{i-1})\right)^2 \right\}. 
\end{align*}

One has that 
\begin{align*}
	& \left(	
		\hat{x}_1(s_i)-e^{-\lambda \dlnt} \hat{x}_1(s_{i-1})
	\right)^2 \\
	=& \left\{
		\hat{x}_1(s_i)-x_1(s_i)
		-e^{-\lambda \dlnt}(\hat{x}_1(s_i)-x_1(s_{i-1}))
		+(x_1(s_i)-e^{-\lambda \dlnt} x_1(s_{i-1})) 
	\right\}^2 \\
	=& \left\{
		(\hat{x}_1(s_i)-x_1(s_i))
		-e^{-\lambda \dlnt}(\hat{x}_1(s_{i-1})-x_1(s_{i-1})) \right\}^2 \\
	&+2\{ \hat{x}_1(s_i)-x_1(s_i)
		-e^{-\lambda \dlnt}(\hat{x}_1(s_{i-1})-x_1(s_{i-1})) \}
		(x_1(s_i)-e^{-\lambda \dlnt} x_1(s_{i-1})) \\
	&+(x_1(s_i)-e^{-\lambda \dlnt} x_1(s_{i-1}))^2.
\end{align*}
It follows that 
\begin{align}
	(A) =& \frac{({\dlnt})^{-1}}{\sqrt{N_2}} \frac{1}{2\epsilon^2\Xi(\lambda)}
		\sum_{i=1}^{N_2} \left[\left\{
		\left( \hat{x}_1(s_i)-x_1(s_i) \right)
		-e^{-\lambda \dlnt} \left( \hat{x}_1(s_{i-1})-x_1(s_{i-1}) \right)
		\right\}^2
	\right. \nonumber \\
	&\left. +2 \left( \hat{x}_1(s_i)-x_1(s_i) \right) 
		\left( x_1(s_i)-e^{-\lambda \dlnt} x_1(s_{i-1}) \right) \right. \nonumber\\
	&\left. -2e^{-\lambda \dlnt} 
		\left( \hat{x}_1(s_{i-1})-x_1(s_{i-1})\right)
		\left( x_1(s_i)-e^{-\lambda \dlnt} x_1(s_{i-1}) \right)
	\right] \nonumber\\
	=& \frac{({\dlnt})^{-1}}{\sqrt{N_2}} \frac{1}{2\epsilon^2 \Xi(\lambda)}
		\sum_{i=1}^{N_2} \left[ \left\{ \left(
		\hat{x}_1(s_i)-x_1(s_i) \right)
		-e^{-\lambda \dlnt} \left(\hat{x}_1(s_{i-1})-x_1(s_{i-1}) \right)
		\right\}^2 \right.  \label{Thm4-1}\\
	&\left. +2 \left\{\hat{x}_1(s_i)-\hat{x}_1(s_{i-1})
		-\left( x_1(s_i)-x_1(s_{i-1})\right)\right\}
		\left( x_1(s_i)-e^{-\lambda \dlnt} x_1(s_{i-1})\right) 
		\right. \label{Thm4-2} \\
	&\left. +2(1-e^{-\lambda \dlnt})
		\left( \hat{x}_1(s_{i-1})-x_1(s_{i-1})\right)
		\left( x_1(s_i)-e^{-\lambda \dlnt} x_1(s_{i-1}) \right)
	\right].  \label{Thm4-3}
\end{align}

For the evaluation of \eqref{Thm4-1}, 
we set that
\begin{eqnarray*}
g_1(t, y, \eta) &=& X_t(y)\sqrt{2} \sin (\pi y) 
	\exp\left\{\frac{\eta}{2}y \right\}.
\end{eqnarray*}
Noting that
\begin{eqnarray*}	
x_1(t) &=& \int_0^1 X_t(y)\sqrt{2} \sin (\pi y) 
	\exp\left\{\frac{\eta}{2}y \right\} dy =  \int_0^1 g_1(t, y, \eta) dy,
	\\
\hat{x}_1(t) &=& \frac{1}{M} \sum_{j=1}^M X_t(y_j) \sqrt{2} \sin (\pi y_j) \exp \left\{ \frac{\hat{\eta}}{2} y_j \right\} 
=  \frac{1}{M} \sum_{j=1}^M g_1(t, y_j, \hat{\eta}),
\end{eqnarray*}
we have that 
\begin{eqnarray*}
	(B) &:=& \frac{1}{\sqrt{N_2}} \frac{1}{\dlnt} \sum_{i=1}^{N_2} 
		\left( x_1(s_i)-\hat{x}_1(s_i) \right)^2  \\
	&=& \frac{1}{\sqrt{N_2}} \frac{1}{\dlnt} \sum_{i=1}^{N_2}
		\left\{
			M \frac{1}{M} \sum_{j=1}^M \int_{\frac{j-1}{M}}^{\frac{j}{M}}
			\{g_1(s_i,y,\eta)-g_1(s_i,y_j,\hat{\eta})\}dy
		\right\}^2  \\
	&\leq& \frac{1}{\sqrt{N_2}} \frac{1}{\dlnt} \sum_{i=1}^{N_2} M^2 \frac{1}{M}
		\sum_{j=1}^{M} \frac{1}{M} \int_{\frac{j-1}{M}}^{\frac{j}{M}}
		\{g_1(s_i,y,\eta)-g_1(s_i,y_j,\hat{\eta})\}^2 dy  \\
	&=& \frac{1}{\sqrt{N_2}} \frac{1}{\dlnt} \sum_{i=1}^{N_2} \sum_{j=1}^M
		\int_{\frac{j-1}{M}}^{\frac{j}{M}} \{g_1(s_i,y,\eta)-g_1(s_i,y_j,\hat{\eta})\}^2 dy.
\end{eqnarray*}
Moreover, 
\begin{eqnarray}
	& & g_1(s_i,y,\eta)-g_1(s_i,y_j,\hat{\eta}) \nonumber\\
	&=& X_{s_i}(y)\sqrt{2}\sin(\pi y) \exp \left\{ \frac{\eta}{2} y\right\}
		-X_{s_i}(y_j)\sqrt{2}\sin(\pi y_j) \exp \left\{ \frac{\hat{\eta}}{2} y_j
		\right\} \nonumber\\
	&=& \left( X_{s_i}(y)-X_{s_i}(y_j) \right) \sqrt{2} \sin(\pi y) 
		\exp \left\{ \frac{\eta}{2} y\right\}  \label{thm4-i} \\
	& &+ X_{s_i}(y_j) \left( \sqrt{2}\sin(\pi y) 
		\exp\left\{ \frac{\eta}{2} y \right\} -\sqrt{2}\sin(\pi y_j) 
		\exp\left\{ \frac{\eta}{2} y_j\right\} \right) \label{thm4-ii}\\
	& &+ X_{s_i}(y_j) \sqrt{2}\sin(\pi y_j)
		\left( \exp\left\{\frac{\eta}{2} y_j \right\} 
		-\exp\left\{\frac{\hat{\eta}}{2} y_j\right\}\right). \label{thm4-iii}
\end{eqnarray}

Set 
$
R(y_j, \hat{\eta})
:=\int_0^1 \frac{y_j}{2} \exp \left\{ \frac{y_j}{2} (\eta + u ( \hat{\eta} - \eta)) \right\} du
$.
Let  $\delta_1>0$.
Since
on $J=\{|\hat{\eta}-\eta| < \delta_1 \}$
\begin{align*}
&	(D) :=	\frac{({\dlnt})^{-1}}{\sqrt{N_2}} \sum_{i=1}^{N_2} \frac{1}{M}
		\sum_{j=1}^M X_{s_i}^2(y_j)2\sin^2(\pi y_j) 
(R(y_j, \hat{\eta}))^2
		\leq C_1 \frac{{N_2}^{\frac{3}{2}}}{T},
\end{align*}
we obtain that 
\begin{align}
	P(|({\rm D})| > \varepsilon) &= P(|({\rm D})| > \varepsilon \cap J)+P(|({\rm D})| > 
		\varepsilon \cap J^c) 
		\leq C_1 \frac{{N_2}^{\frac{3}{2}}}{T} \frac{1}{\varepsilon}+o(1). \label{Res1}
\end{align}

\noindent
It follows that
\begin{eqnarray*} 
  E[(\ref{thm4-i})^2 \times \sqrt{N_2} ({\dlnt})^{-1}] &\leq& 
 C_1 E[ \left( X_{s_i}(y)-X_{s_i}(y_j)\right)^2] \sqrt{N_2}  ({\dlnt})^{-1} 
		\leq \frac{C_2}{M^{1-\rho_1}} \frac{N_2^{\frac{3}{2}}}{T},
\\
   E[(\ref{thm4-ii})^2 \times \sqrt{N_2} ({\dlnt})^{-1}] &\leq& C_1 (y-y_j)^2 \sqrt{N_2} ({\dlnt})^{-1} 
		\leq \frac{C_1}{M^2} \frac{N_2^{\frac{3}{2}}}{T},
\\
\frac{({\dlnt})^{-1}}{\sqrt{N_2}} \sum_{i=1}^{N_2} (\ref{thm4-iii})^2 
&=&
		\frac{({\dlnt})^{-1}}{\sqrt{N_2}} \sum_{i=1}^{N_2} \frac{1}{M}
		\sum_{j=1}^M X_{s_i}^2(y_j)2\sin^2(\pi y_j) 
(R(y_j, \hat{\eta}))^2 \frac{1}{N m} 
		\left(\sqrt{N m}(\hat{\eta}-\eta) \right)^{2} 
		\\
	&=& O_p\left( \frac{{N_2}^{\frac{3}{2}}}{T} \frac{1}{N m} \right),
\end{eqnarray*}
where we use  (\ref{Res1}) for the last estimate.  

Hence
\begin{equation*}
	(B) = O_p \left(\frac{{N_2}^{\frac{3}{2}}}{T M^{1-\rho_1}} \right)
		+O_p\left(\frac{{N_2}^{\frac{3}{2}}}{T N m}\right), \quad 
	\eqref{Thm4-1} = O_p \left(\frac{{N_2}^{\frac{3}{2}}}{T M^{1-\rho_1}} \right)
		+O_p\left(\frac{{N_2}^{\frac{3}{2}}}{T N m}\right). 
\end{equation*}

For the evaluation of \eqref{Thm4-3},  one has that
\begin{align*}
	\eqref{Thm4-3}^2 \leq& \frac{({\dlnt})^{-1}}{\sqrt{N_2}} ({\dlnt})^2
		\frac{({\dlnt})^{-1}}{\sqrt{N_2}}\sum_{i=1}^{N_2}
		\left( \hat{x}_1(s_i)-x_1(s_{i-1}) \right)^2 \sum_{i=1}^{N_2} 
		\left( x_1(s_i)-e^{-\lambda \dlnt}x_1(s_{i-1}) \right)^2 \\
	=& \frac{T^2}{\sqrt{N_2} N_2} \times \left(  O_p \left(\frac{{N_2}^{\frac{3}{2}}}{T M^{1-\rho_1}} \right)
		+O_p\left(\frac{{N_2}^{\frac{3}{2}}}{T N m}\right) \right) 
		\\
	=& O_p\left(\frac{T}{M^{1-\rho_1}}\right)+O_p\left(\frac{T}{N m}\right). 
\end{align*}

For the evaluation of \eqref{Thm4-2}, 
setting that $\Delta X_{s_i}(y) = X_{s_i}(y_j)-X_{s_{i-1}}(y_j)$, 
we obtain that 
\begin{align*}
	\eqref{Thm4-2} =& \frac{({\dlnt})^{-1}}{\sqrt{N_2}} \sum_{i=1}^{N_2}
		\left\{\hat{x}_1(s_i)-\hat{x}_1(s_{i-1})
		-\left(x_1(s_i)-x_1(s_{i-1})\right) \right\} 
	\left(x_1(s_i)-e^{-\lambda \dlnt}x_1(s_{i-1})\right) \nonumber \\
	=& \frac{({\dlnt})^{-1}}{\sqrt{N_2}} \sum_{i=1}^{N_2} \frac{1}{M}
		\sum_{j=1}^M 
                  \Delta X_{s_i}(y_j) 
		\sqrt{2}\sin(\pi y_j)
		\left(\exp\left\{\frac{\hat{\eta}}{2}y_j \right\}
		-\exp\left\{\frac{\eta}{2} y_j\right\}\right) 
		\left(x_1(s_i)-e^{-\lambda \dlnt}x_1(s_{i-1})\right) \nonumber \\
		&+ \frac{({\dlnt})^{-1}}{\sqrt{N_2}} 
		\sum_{i=1}^{N_2} \sum_{j=1}^{M}
		\int_{\frac{j-1}{M}}^{\frac{j}{M}}
                  \Delta X_{s_i}(y_j) 
		\left(\sqrt{2}\sin(\pi y_j)-\sqrt{2}\sin(\pi y)\right)
		\exp\left\{\frac{\eta}{2}y_j\right\} dy 
                 \left(x_1(s_i)-e^{-\lambda \dlnt}x_1(s_{i-1})\right) \nonumber \\
	&+ \frac{({\dlnt})^{-1}}{\sqrt{N_2}}
		\sum_{i=1}^{N_2} \sum_{j=1}^{M}
		\int_{\frac{j-1}{M}}^{\frac{j}{M}}
                  \Delta X_{s_i}(y_j) 
		\sqrt{2}\sin(\pi y) \left(\exp\left\{\frac{\eta}{2}y_j\right\}
		-\exp\left\{\frac{\eta}{2}y\right\} \right)dy 
		\left(x_1(s_i)-e^{-\lambda \dlnt}x_1(s_{i-1})\right) \nonumber \\
	&+ \frac{({\dlnt})^{-1}}{\sqrt{N_2}}
		\sum_{i=1}^{N_2} \sum_{j=1}^{M}
		\int_{\frac{j-1}{M}}^{\frac{j}{M}}
		\left\{ \Delta X_{s_i}(y_j) -\Delta X_{s_i}(y) \right\}		
		\sqrt{2}\sin(\pi y)\exp\left\{\frac{\eta}{2}y\right\}dy 
		\left(x_1(s_i)-e^{-\lambda \dlnt}x_1(s_{i-1})\right) \nonumber \\
	&=: \ ({\rm I})+({\rm II})+({\rm III})+({\rm IV}). \nonumber
\end{align*}

For the evaluation of (I), one has that
\begin{align}
	({\rm I})^2 \leq& \left(\frac{({\dlnt})^{-1}}{\sqrt{N_2}}\right)^2 \sum_{i=1}^{N_2}\frac{1}{M}
		\sum_{j=1}^{M} 
                \left( \Delta_{s_i}(y_j) \right)^2 
                \left(\sqrt{2}\sin(\pi y_j) R(y_j,\hat{\eta} )  \right)^2 
		\frac{T}{N m}  \label{C1} \\
	& \times \frac{1}{T} \sum_{i=1}^{N_2} \left(x_1(s_i)-e^{-\lambda \dlnt}x_1(s_{i-1})\right)^2 
		\times \left( \sqrt{N m}(\hat{\eta}-\eta)\right)^{2}.  \nonumber
\end{align}
Let  $\delta_1>0$.
On $J=\{|\hat{\eta}-\eta| < \delta_1\}$,
\begin{align*}
	(\ref{C1}) \leq&  C_1 \frac{({\dlnt})^{-2}}{N m} \frac{T}{N_2 }
	\sum_{i=1}^{N_2}\frac{1}{M}
		\sum_{j=1}^{M} \left( X_{s_i}(y_j)-X_{s_{i-1}}(y_j) \right)^2 
	=& O_p \left( \frac{T}{({\dlnt})^{\frac{3}{2}} N m} \right)
	= O_p \left( \frac{{N_2}^{\frac{3}{2}}}{T^{\frac{1}{2}} N m} \right)
\end{align*}
because 
$E \left[ 
\left( 
\Delta X_{s_i}(y_j)
\right)^2 
\right] = \sqrt{{\dlnt}}=\sqrt{T/N_2}$.
It follows that under 
$\frac{{N_2}^{\frac{3}{2}}}{T^{\frac{1}{2}} N m} \rightarrow 0$, 
\begin{align*}
	P(|(\ref{C1})| > \varepsilon) &= P(|(\ref{C1})| > \varepsilon \cap J)+P(|(\ref{C1})| > 
		\varepsilon \cap J^c) 
		\leq \frac{{N_2}^{\frac{3}{2}}}{T^{\frac{1}{2}}N m} \frac{1}{\varepsilon}+o(1) 
		\rightarrow 0.
\end{align*}
Therefore, $({\rm I})=o_p(1)$.

For the evaluation of (II), we obtain that
\begin{align*}
	({\rm II})^2 \leq& C_1 \left( \frac{{\dlnt}^{-1}}{\sqrt{N_2}}\right)^2
		\sum_{i=1}^{N_2} \frac{1}{M}
		\sum_{j=1}^{M}
		 \left( \Delta X_{s_i}(y_j) \right)^2 
		  |y-y_j|^2 
		\sum_{i=1}^{N_2} \left(x_1(s_i)-e^{-\lambda \dlnt}x_1(s_{i-1})\right)^2 \\
	=& O_p \left( \left(\frac{\sqrt{N_2}}{T}\right)^2
		N_2 \sqrt{\frac{T}{N_2}} \frac{1}{M^2} T \right)
		=O_p \left(\frac{{N_2}^{\frac{3}{2}}}{T^{\frac{1}{2}}M^2} \right). 
\end{align*}

For the evaluation of (III), one has that 
\begin{align*}
	({\rm III})^2 \leq& C_1 \left( \frac{{\dlnt}^{-1}}{\sqrt{N_2}}\right)^2
		\sum_{i=1}^{N_2} \frac{1}{M} \sum_{j=1}^{M}
		 \left( \Delta X_{s_i}(y_j) \right)^2 
		|y-y_j|^2
		\sum_{i=1}^{N_2} \left(x_1(s_i)-e^{-\lambda \dlnt}x_1(s_{i-1})\right)^2 \\
		=& O_p \left(\frac{{N_2}^{\frac{3}{2}}}{T^{\frac{1}{2}}M^2} \right). 
\end{align*}

For the evaluation of (IV),  we obtain that 
\begin{eqnarray*}
	({\rm IV})^2 &\leq&  \left( \frac{{\dlnt}^{-1}}{\sqrt{N_2}}\right)^2
		\sum_{i=1}^{N_2} M^2 \frac{1}{M} \\
		& & \times \sum_{j=1}^{M} \frac{1}{M}
		\int_{\frac{j-1}{M}}^{\frac{j}{M}}
		\left( \left( \sum_{k=1}^\infty (x_k(s_i)-x_k(s_{i-1}))(e_k(y_j) -e_k(y)) \right)^2
		(x_1(s_i)-x_1(s_{i-1}))^2 \right) dy.
\end{eqnarray*}
By setting that
\begin{equation*}
	({\rm E}) := \left( \sum_{k=1}^\infty (x_k(s_i)-x_k(s_{i-1}))(e_k(y_j) -e_k(y)) \right)^2
		(x_1(s_i)-x_1(s_{i-1}))^2, 
\end{equation*}
it follows that 
\begin{align*}
	E[({\rm E})] =& \sum_{k=2}^\infty E[(x_k(s_i)-x_k(s_{i-1}))^2]
		E[(x_1(s_i)-x_1(s_{i-1}))^2](e_k(y_j) -e_k(y))^2 \\
	&+ E[(x_1(s_i)-x_1(s_{i-1}))^4](e_1(y_j) -e_1(y))^2 \\
	\leq& C_1 \left(  \frac{1}{M^{1-\rho_1}} \dlnt+\frac{{\dlnt}^2}{M^2} \right).
\end{align*}
Therefore, 
under 
$\frac{{N_2}^{\frac{3}{2}}} {T M^{1-\rho_1}} \rightarrow 0$, 
\begin{align*}
	E[({\rm IV})^2] \leq& \frac{{\dlnt}^{-2}}{N_2} N_2 \cdot \frac{\dlnt}{M^{1-\rho_1}}
		=\frac{1}{\dlnt M^{1-\rho_1}} = \frac{N_2}{T M^{1-\rho_1}} \rightarrow 0. 
\end{align*}
\end{en-text}

The difference between 
{the score function of the drift parameter $\lambda$ 
based on { 
the thinned data of the approximate coordinate process
${\bf \hat{x}}$ and  that base on the thinned data of the coordinate process
${\bf {x}}$} 
is as follows.
\begin{align*}
	{\cal F} :=& \epsilon \left\{
		\partial_\lambda l_{N_2}(\lambda \ | \ {\bf \hat{x}}  )
		-\partial_\lambda l_{N_2}(\lambda \ | \ {\bf {x}}  ) \right\} \\
	=&\left(-\frac{1}{2}\right)\epsilon\sum_{i=1}^{N_2}
		\left[\left(-\frac{\partial_\lambda \Xi(\lambda)}{\epsilon^2\Xi(\lambda)^2}\right)
		\left\{\left(\hat{x}_1(s_i)-e^{-\lambda \dlnt}\hat{x}_1(s_{i-1})\right)^2
		-\left(x_1(s_i)-e^{-\lambda \dlnt}x_1(s_{i-1}) \right)^2 \right\} {\dlnt}^{-1} \right.\\
	&+2 \left. \frac{e^{-\lambda \dlnt}}{\epsilon^2 \Xi(\lambda)}
		\left\{\hat{x}_1(s_i) \left(\hat{x}_1(s_i)
		-e^{-\lambda \dlnt}  \hat{x}_1(s_{i-1})\right)
		-x_1(s_i)\left(x_1(s_i)-e^{-\lambda \dlnt}x_1(s_{i-1}) \right)\right\}\right] \\
	=:& {\cal F}_1+{\cal F}_2.
\end{align*}
{In an analogous way to (\ref{consistency-1}),
it follows that}
$$
		{\cal F}_1 = o_p(1).
$$
For the estimate of ${\cal F}_2$, we have that  
\begin{align}
	{\cal F}_2 =& \left( -\epsilon \right) 
		\frac{e^{-\lambda \dlnt}}{\epsilon^2\Xi(\lambda)}
		\sum_{i=1}^{N_2} \left[ \left( \hat{x}_1(s_i)-x_1(s_i)  \right)
		\left( \hat{x}_1(s_i)-e^{-\lambda \dlnt}  \hat{x}_1(s_{i-1}) 
		\right) \right. \nonumber\\
	&+ \left. x_1(s_i) \left\{ \hat{x}_1(s_i)-x_1(s_i)
		-e^{-\lambda \dlnt}( \hat{x}_1(s_{i-1})-x_1(s_{i-1})) \right\}\right] \nonumber \\
=&{   \left( -\epsilon \right) 
		\frac{e^{-\lambda \dlnt}}{\epsilon^2\Xi(\lambda)}
		\sum_{i=1}^{N_2} \left( \hat{x}_1(s_i)-x_1(s_i)  \right)^2   }
		 \nonumber \\
	&{ + \epsilon 
		\frac{e^{-2\lambda \dlnt}}{\epsilon^2\Xi(\lambda)}
		\sum_{i=1}^{N_2} \left( \hat{x}_1(s_i)-x_1(s_i)\right)
		\left( \hat{x}_1(s_{i-1})-x_1(s_{i-1})\right) }
		 \nonumber  \\
	&{ + \left( -\epsilon \right) 
		\frac{e^{-\lambda \dlnt}}{\epsilon^2\Xi(\lambda)}
		\sum_{i=1}^{N_2} \left( \hat{x}_1(s_i)-x_1(s_i)  \right)
		\left( x_1(s_i)-e^{-\lambda \dlnt}x_1(s_{i-1})\right) }
		 \nonumber  \\
	&{ + \left( -\epsilon \right) 
		\frac{e^{-\lambda \dlnt}}{\epsilon^2\Xi(\lambda)}
		\sum_{i=1}^{N_2} x_1(s_i) \left\{ \hat{x}_1(s_i)-x_1(s_i)
		-e^{-\lambda \dlnt}( \hat{x}_1(s_{i-1})-x_1(s_{i-1})) \right\}  } \nonumber \\
	=:& {\cal H}_1 +{\cal H}_2+{\cal H}_3+{\cal H}_4.  \nonumber 
\end{align}
For the estimate of  ${\cal H}_1$, it follows that
\begin{align*}
	|{\cal H}_1| &\leq C_1 \frac{\dlnt}{\epsilon} {\cal Z}_1
		= O_p\left(\frac{{N_2}}{\epsilon M^2} \right)
		+O_p\left(\frac{\epsilon{N_2}}{M^{1-\rho_1}} \right)
		+O_p\left( \frac{{N_2}}{\epsilon N m} \right) \stackrel{p}{\rightarrow} 0. 
\end{align*}
For the estimate of ${\cal H}_2$, one has that
\begin{align*}
	{\cal H}_2^2 &\leq C \frac{1}{\epsilon^2} \sum_{i=1}^{N_2}
		\left( \hat{x}_1(s_i)-x_1(s_i) \right)^2
		\sum_{i=1}^{N_2} \left( \hat{x}_1(s_{i-1})-x_1(s_{i-1}) \right)^2 \\
	&= O_p\left(\frac{N_2^2}{\epsilon^2 M^{4}}\right)
		+O_p\left(\frac{\epsilon^2 N_2^2}{M^{2(1-\rho_1)}}\right)
		+O_p\left( \frac{N_2^2}{\epsilon^2 (N m)^2} \right) 
		\stackrel{p}{\rightarrow} 0.
\end{align*}
For the estimate of  ${\cal H}_3$, 
by letting that 
\begin{equation*}
	\Delta x ( s_i, s_{i-1})  = x_1(s_i)-e^{-\lambda \dlnt}x_1(s_{i-1}),
\end{equation*}
it follows that 
\begin{align}
	& -{\cal H}_3 
	=
	 \epsilon\frac{e^{-\lambda \dlnt}}{\epsilon^2 \Xi(\lambda)}
		\sum_{i=1}^{N_2} \left\{\left( \hat{x}_1(s_i)-x_1(s_i)\right)
		\left(x_1(s_i)-e^{-\lambda \dlnt}x_1(s_{i-1})\right)\right\} \nonumber
		\\
	 &=\epsilon\frac{e^{-\lambda \dlnt}}{\epsilon^2 \Xi(\lambda)}
	 \sum_{i=1}^{N_2}\frac{1}{M}
		\sum_{j=1}^M X_{s_i}(y_j)\sqrt{2}\sin(\pi y_j)
		\left(\exp\left\{\frac{\hat{\eta}}{2}y_j\right\}
		-\exp\left\{\frac{\eta}{2}y_j\right\}\right)
		\Delta x ( s_i, s_{i-1})  \nonumber \\
	&{ +\epsilon\frac{e^{-\lambda \dlnt}}{\epsilon^2 \Xi(\lambda)}
	\sum_{i=1}^{N_2}\sum_{j=1}^M\int_{\frac{j-1}{M}}^{\frac{j}{M}}
		X_{s_i}(y_j)\left(\sqrt{2}\sin(\pi y_j)-\sqrt{2}\sin(\pi y)\right)
		\exp\left\{\frac{\eta}{2}y_j\right\}dy  \Delta x ( s_i, s_{i-1})   \nonumber }\\
	&{ +\epsilon\frac{e^{-\lambda \dlnt}}{\epsilon^2 \Xi(\lambda)}
	\sum_{i=1}^{N_2}\sum_{j=1}^M\int_{\frac{j-1}{M}}^{\frac{j}{M}}
		X_{s_i}(y_j)\sqrt{2}\sin(\pi y)\left(\exp\left\{\frac{\eta}{2}y_j\right\}
		-\exp\left\{\frac{\eta}{2}y\right\}\right)dy \Delta x ( s_i, s_{i-1})   \nonumber } \\
	&{ +\epsilon\frac{e^{-\lambda \dlnt}}{\epsilon^2 \Xi(\lambda)}
	\sum_{i=1}^{N_2}\sum_{j=1}^M \int_{\frac{j-1}{M}}^{\frac{j}{M}}
		\left(X_{s_i}(y_j)-X_{s_i}(y)\right)\sqrt{2}\sin(\pi y)
		\exp\left\{\frac{\eta}{2}y\right\}dy \Delta x ( s_i, s_{i-1})   \nonumber } \\
		&{ =: \mathcal{B}_2 +\mathcal{B}_3 + \mathcal{B}_4 + \mathcal{B}_5. \nonumber}
\end{align}
For the estimate of $\mathcal{B}_2$, we obtain that 
\begin{align*}
	{ \mathcal{B}_2^2} \leq& C\frac{1}{\epsilon^2}\sum_{i=1}^{N_2}\frac{1}{M}\sum_{j=1}^M
		\left(X_{s_i}(y_j)\sqrt{2}\sin(\pi y_j) 
                 R(y_j, \hat{\eta}) 
		 (\hat{\eta}-\eta)\right)^2
		\sum_{i=1}^{N_2} (\Delta x ( s_i, s_{i-1}) )^2 \\
	=& O_p\left(\frac{N_2}{N m}\right)\stackrel{p}{\rightarrow} 0.
\end{align*}
For the estimate of  $\mathcal{B}_3$, one has that 
\begin{align*}	
	{ \mathcal{B}_3^2} \leq& \frac{1}{\epsilon^2}\sum_{i=1}^{N_2} \sum_{j=1}^M  \int_{\frac{j-1}{M}}^{\frac{j}{M}} (X_{s_i}(y_j))^2(y-y_j)^2 dy
		\sum_{i=1}^{N_2} (\Delta x ( s_i, s_{i-1}) )^2 \\
	=& O_p\left(\frac{N_2}{M^2} \right) \stackrel{p}{\rightarrow} 0.
\end{align*}
By the same manner as the estimate of $\mathcal{B}_3$, 
it is proved that 
$$ 
\mathcal{B}_4 \stackrel{p}{\rightarrow} 0.
$$ 
For the estimate of  $\mathcal{B}_5$, 
since one has that 
\begin{eqnarray*}
	{ \mathcal{B}_5^2} &\leq& \frac{C_1}{\epsilon^2}N_2^2\frac{1}{N_2}
		\sum_{i=1}^{N_2}\sum_{j=1}^M \int_{\frac{j-1}{M}}^{\frac{j}{M}}
		\left(X_{s_i}(y_j)-X_{s_i}(y)\right) dy
		(\Delta x ( s_i, s_{i-1}) )^2 \\
		&=& O_p\left(\frac{N_2}{M^2} \right) + O_p\left(\frac{\epsilon^2 N_2}{{M^{1-\rho_1}}} \right)  \stackrel{p}{\rightarrow} 0,
\end{eqnarray*}
it follows that
$$
{\cal H}_3 \stackrel{p}{\rightarrow}  0.
$$ 
For the estimate of  ${ {\cal H}_4}$, 
set
\begin{align}
         -{\cal H}_4 =&{   \epsilon\frac{e^{-\lambda \dlnt}}{\epsilon^2 \Xi(\lambda)}
	\sum_{i=1}^{N_2}
		\left\{\hat{x}_1(s_i)-\hat{x}_1(s_{i-1})
		-\left(x_1(s_i)-x_1(s_{i-1})\right) \right\} 
	x_1(s_i) \nonumber } \\
	+&{   \epsilon\frac{e^{-\lambda \dlnt}}{\epsilon^2 \Xi(\lambda)}
	\sum_{i=1}^{N_2}
		\left(1-e^{\lambda \delta} \right) \left(\hat{x}_1(s_{i-1})-x_1(s_{i-1})\right)  
	x_1(s_i)  \nonumber } \\
	=:&{  \mathcal{B}_6 + \mathcal{B}_7}. \nonumber 
\end{align}
We obtain that 
\begin{align*}
{ \mathcal{B}_7^2} \leq& C \frac{\delta^2}{\epsilon^2}  \sum_{i=1}^{N_2} (x_1(s_i))^2 
\delta \frac{1}{\delta} \sum_{i=1}^{N_2} \left(\hat{x}_1(s_{i-1})-x_1(s_{i-1})\right)^2 \\
=& O_p \left( \frac{1}{\epsilon^2 N_2 M^{2}} \right)
+ O_p \left( \frac{1}{N_2 M^{1-\rho_1}} \right)
+ O_p \left( \frac{1}{\epsilon^2 N m} \right)
 \stackrel{p}{\rightarrow}  0.
\end{align*}
Furthermore, let
\begin{align*}
	{ \mathcal{B}_6} =&  \epsilon\frac{e^{-\lambda \dlnt}}{\epsilon^2 \Xi(\lambda)}
	 \sum_{i=1}^{N_2} \frac{1}{M}
		\sum_{j=1}^M 
                  \Delta X_{s_i}(y_j) 
		\sqrt{2}\sin(\pi y_j)
		\left(\exp\left\{\frac{\hat{\eta}}{2}y_j \right\}
		-\exp\left\{\frac{\eta}{2} y_j\right\}\right) 
		x_1(s_i) \nonumber \\
		&+ \epsilon\frac{e^{-\lambda \dlnt}}{\epsilon^2 \Xi(\lambda)}
		\sum_{i=1}^{N_2} \sum_{j=1}^{M}
		\int_{\frac{j-1}{M}}^{\frac{j}{M}}
                  \Delta X_{s_i}(y_j) 
		\left(\sqrt{2}\sin(\pi y_j)-\sqrt{2}\sin(\pi y)\right)
		\exp\left\{\frac{\eta}{2}y_j\right\} dy 
                 x_1(s_i) \nonumber \\
	&+  \epsilon\frac{e^{-\lambda \dlnt}}{\epsilon^2 \Xi(\lambda)}
		\sum_{i=1}^{N_2} \sum_{j=1}^{M}
		\int_{\frac{j-1}{M}}^{\frac{j}{M}}
                  \Delta X_{s_i}(y_j) 
		\sqrt{2}\sin(\pi y) \left(\exp\left\{\frac{\eta}{2}y_j\right\}
		-\exp\left\{\frac{\eta}{2}y\right\} \right)dy 
		x_1(s_i) \nonumber \\
	&+  \epsilon\frac{e^{-\lambda \dlnt}}{\epsilon^2 \Xi(\lambda)}
		\sum_{i=1}^{N_2} \sum_{j=1}^{M}
		\int_{\frac{j-1}{M}}^{\frac{j}{M}}
		\left\{ \Delta X_{s_i}(y_j) -\Delta X_{s_i}(y) \right\}		
		\sqrt{2}\sin(\pi y)\exp\left\{\frac{\eta}{2}y\right\}dy 
		x_1(s_i) \nonumber \\
	&=: \ ({\rm V})+({\rm VI})+({\rm VII})+({\rm VIII}). \nonumber
\end{align*}
For the estimate of (V), it follows that
\begin{align}
	({\rm V})^2 \leq&{  C \frac{1}{\epsilon^2} \sum_{i=1}^{N_2}\frac{1}{M}
		\sum_{j=1}^{M} 
                \left( \Delta X_{s_i}(y_j) \right)^2 
                \left(\sqrt{2}\sin(\pi y_j) R(y_j,\hat{\eta} )  \right)^2 
		\frac{N_2}{N m} \nonumber } \\
	&{  \times  \frac{1}{N_2} \sum_{i=1}^{N_2} \left(x_1(s_i) \right)^2 
		\times \left( \sqrt{N m}(\hat{\eta}-\eta)\right)^{2}  \nonumber } \\
		=:&{  \mathcal{B}_8  \times  \frac{1}{N_2} \sum_{i=1}^{N_2} \left(x_1(s_i) \right)^2 
		\times \left( \sqrt{N m}(\hat{\eta}-\eta)\right)^{2}}.   \nonumber 
\end{align}
Since
\begin{align*}
	{ \mathcal{B}_8} \leq&  C_1 \frac{N_2}{N m} \frac{1}{\epsilon^2}
	\sum_{i=1}^{N_2}\frac{1}{M}
		\sum_{j=1}^{M} \left( X_{s_i}(y_j)-X_{s_{i-1}}(y_j) \right)^2 
	= O_p \left( \frac{{N_2}^{\frac{3}{2}}}{N m} \right)
	+ O_p \left( \frac{{N_2}^{\frac{1}{2}}}{\epsilon^2 N m} \right),
\end{align*}
one has that 
$({\rm V})=o_p(1)$.

\begin{en-text}
It follows that 
under
\begin{align*}
	P(|{ \mathcal{B}_8}| > \varepsilon) &= P(|{ \mathcal{B}_8}| > \varepsilon \cap J)+P(|{ \mathcal{B}_8}| > 
		\varepsilon \cap J^c) 
		\leq \frac{{N_2}^{\frac{3}{2}}}{N m}
	+ \frac{{N_2}^{\frac{1}{2}}}{\epsilon^2 N m}
		+o(1) 
		\rightarrow 0.
\end{align*}
\end{en-text}

For the estimate of (VI), it follows that
\begin{align*}
	({\rm VI})^2 \leq& \frac{C_1}{\epsilon^2}
		\sum_{i=1}^{N_2} M \sum_{j=1}^M \frac{1}{M} \int_{\frac{j-1}{M}}^{\frac{j}{M}}
		\left(\Delta X_{s_i}(y_j)\right)^2 (y_j - y)^2 dy
		\sum_{i=1}^{N_2} (x_i(s_i))^2 \\
	=& O_p \left( \frac{N_2^{\frac{3}{2}}}{M^2}\right) + O_p \left( \frac{N_2^{\frac{1}{2}}}{\epsilon^2 M^2}\right)
		=o_p(1). 
\end{align*}
In the same manner as the estimate of  (VI),  it is proved that ${ ({\rm VII}) \stackrel{p}{\rightarrow} 0}$.
 
For the estimate of (VIII),  one has that 
\begin{eqnarray*}
	({\rm VIII})^2 &\leq&  \frac{1}{\epsilon^2}
		\sum_{i=1}^{N_2} M \sum_{j=1}^{M} \frac{1}{M}
		\int_{\frac{j-1}{M}}^{\frac{j}{M}}
		\left( \Delta X_{s_i}(y_j) - \Delta X_{s_i}(y) \right)^2 dy
		\sum_{i=1}^{N_2} (x_1(s_i))^2.
\end{eqnarray*}
Set
\begin{equation*}
	{\cal Z}_4 := \left( \Delta X_{s_i}(y_j) - \Delta X_{s_i}(y) \right)^2.
\end{equation*}
Since we obtain that
\begin{align*}
	E[{\cal Z}_4] 
	\leq& C_1 \left(  \frac{1}{M^{1-\rho_1}}\left( \frac{\epsilon^2}{N_2^2} 
	+ \frac{\epsilon}{N_2^3} + \frac{1}{N_2^4}\right) \right),
\end{align*}
it follows that 
\begin{align*}
	E[({\rm VIII})^2] \leq& C\left( \frac{1}{M^{1-\rho_1}} + \frac{1}{\epsilon N_2 M^{1-\rho_1}} + 
	\frac{1}{{\epsilon^2 N_2^2 M^{1-\rho_1}}} \right) \rightarrow 0. 
\end{align*}

\noindent
Hence,
$$
{ \mathcal{B}_6}  \stackrel{p}{\rightarrow} 0,  \quad  {\cal H}_4  \stackrel{p}{\rightarrow} 0,
\quad {\cal F}_2 \stackrel{p}{\rightarrow} 0,
$$
which yields that  ${\cal F} \stackrel{p}{\rightarrow} 0$.

{In a similar manner}, one has that
\begin{eqnarray*}
\epsilon^2
\left\{
\partial_{\lambda}^2 l_{N_2}(\lambda \ | \ {\bf \hat{x}}  ) 
-
\partial_{\lambda}^2 l_{N_2}(\lambda \ | \ {\bf {x}}  ) 
\right\} = o_p(1)
\end{eqnarray*}
uniformly in $\lambda$.
{
Therefore, one has that
\begin{eqnarray*}
 \epsilon^{-1} \left( \hat{\lambda} - \lambda^* \right) 
 &=&  G(\lambda^*)^{-1} \epsilon \partial_\lambda  l_{N_2}(\lambda \ | \ {\bf {x}} ) +o_p(1),
 \\
 &\stackrel{d}{\rightarrow} &
N \left( 
0, G(\lambda^*)^{-1}
\right). \label{AN_1}
\end{eqnarray*}
}
{
Moreover, it is shown that
\begin{equation*}
\begin{pmatrix}
\sqrt{Nm}(\hat{\theta}_2 - \theta_{2}^*) \\
\sqrt{Nm}(\hat{\theta}_1 - \theta_{1}^*) \\
\epsilon^{-1}(\hat{\lambda}-\lambda^*)
\end{pmatrix}
\stackrel{d}{\longrightarrow } N
\left(
\begin{pmatrix}
0 \\
0 \\
0 
\end{pmatrix}
,
\begin{pmatrix}
J_{1,1} &J_{1,2} & 0 \\
J_{1,2} &J_{2,2} & 0 \\
0 & 0 & G(\lambda^{*})^{-1}
\end{pmatrix}
\right).
 \label{AN_2}
\end{equation*}
}
Since we obtain that
\begin{eqnarray*}
\epsilon^{-1} (\hat{\theta}_0 -\theta_0^*) 
&=& \epsilon^{-1} \left( \hat{\lambda} - \lambda^* 
+  \frac{\left( \hat{\theta}_1 \right)^2}{  4 \hat{\theta}_2}  - \frac{\left( \theta_1 \right)^2}{  4 \theta_2^*}
+ \pi^2 (\hat{\theta}_2 -\theta_2^*) \right)
\\
&=&
 \epsilon^{-1} \left( \hat{\lambda} - \lambda^* \right)
+ \frac{ \epsilon^{-1} \left( (\hat{\theta}_1)^2 - ({\theta_1^*})^2 \right)}{4 \hat{\theta}_2}  
+ 
\frac{(\theta_1^*)^2}{4} 
\epsilon^{-1}
\left( \frac{1}{\hat{\theta}_2}   - \frac{1}{\theta_2^*} \right) 
+ \pi^2 \epsilon^{-1} \left( \hat{\theta}_2 -\theta_2^* \right)
\\
&=&
\epsilon^{-1} \left( \hat{\lambda} - \lambda^* \right) +o_p(1),
\end{eqnarray*}
one has that
\begin{equation*}
\begin{pmatrix}
\sqrt{N m}(\hat{\theta}_2- \theta_2^*) \\
\sqrt{N m}(\hat{\theta}_1- \theta_1^*) \\
\epsilon^{-1}(\hat{\theta}_0-\theta_0^*)
\end{pmatrix}
=
\begin{pmatrix}
\sqrt{N m}(\hat{\theta}_2- \theta_2^*) \\
\sqrt{N m}(\hat{\theta}_1- \theta_1^*) \\
\epsilon^{-1}(\hat{\lambda} - \lambda^*)
\end{pmatrix}
+o_p(1).
\end{equation*}

\begin{en-text}
\begin{equation}
\epsilon^{-1}(\hat{\lambda}-\lambda^*)
\stackrel{d}{\rightarrow} 
N \left( 
0, G(\lambda^*)
\right). \label{AN_1}
\end{equation}
\end{en-text}

\noindent
Consequenctly, 
\begin{equation*}
\begin{pmatrix}
\sqrt{Nm}(\hat{\theta}_2 - \theta_{2}^*) \\
\sqrt{Nm}(\hat{\theta}_1 - \theta_{1}^*) \\
\epsilon^{-1}(\hat{\theta}_0 - \theta_{0}^*)
\end{pmatrix}
\stackrel{d}{\longrightarrow } N
\left(
\begin{pmatrix}
0 \\
0 \\
0 
\end{pmatrix}
,
\begin{pmatrix}
J_{1,1} &J_{1,2} & 0\\
{J_{1,2}} &J_{2,2} & 0\\
0 & 0 & G(\lambda^{*})^{-1}
\end{pmatrix}
\right).
 \label{AN_2}
\end{equation*}
This completes the proof.

\section*{Acknowledgement}
This work 
was partially supported by 
JST CREST,
JSPS KAKENHI Grant Number 
JP17H01100 
and Cooperative Research Program
of the Institute of Statistical Mathematics.

\section*{References}
\begin{description}{}{
}


\item
Bibinger, M.\ and Trabs, M.\ 
(2020).\ 
Volatility estimation for stochastic pdes using high-frequency observations.\
{\it Stochastic Processes and their Applications}, {\bf 130}, 3005--3052.


\item
Cialenco, I.\ and Huang, Y.\ 
(2020).\ 
A note on parameter estimation for discretely sampled SPDEs.\
{\it Stochastics and Dynamics},  {\bf 20}, 2050016.

\item
Cont, R. (2005). 
Modeling term structure dynamics: an infinite dimensional approach.\ 
{\it International Journal of Theoretical and Applied Finance}, 
{\bf 8}, 357--380.

\item
Genon-Catalot, V.\ 
(1990).\ 
Maximum contrast estimation for diffusion processes
from discrete observations.\
{\it Statistics }{\bf 21}, 99--116.

\item
Gloter A.\ and S{\o}rensen, M.\
(2009).\
Estimation for stochastic differential equations with a small diffusion coefficient.\
{\it Stochastic Processes and their Applications}, {\bf 119}, 679-699. 

\item
Guy, R., Laredo, C.\ and Vergu, E.\
(2014).\
Parametric inference for discretely observed multidimensional diffusions 
with small diffusion coefficient.\
{\it Stochastic Processes and their Applications}, {\bf 124}, 51--80.


\item
Hildebrandt, F.\  
(2020).\
On generating fully discrete samples of the stochastic heat equation on an interval.\ 
{\it Statistics \& Probability Letters}, 
{\bf 162}, 108750.

\item
Hildebrandt, F.\ and Trabs, M.\  (2019).\  
Parameter estimation for SPDEs based on discrete observations in time and space.\ 
arXiv:1910.01004.

\item
Kaino, Y.\ and Uchida, M.\ (2018a).\ 
Hybrid estimators for small diffusion processes based on reduced data.\ 
{\it Metrika}, {\bf 81}, 745--773.

\item
Kaino, Y.\ and Uchida, M.\ (2018b).\ 
Hybrid estimators for stochastic differential equations from reduced data.\ 
{\it Statistical Inference for Stochastic Processes}, {\bf 21}, 435--454.


\item
{
Kaino, Y.\ and Uchida, M.\ (2020).\ 
Parametric estimation for a parabolic linear SPDE model based on discrete observations.\ 
To appear in {\it Journal of Statistical Planning and Inference}.}

\item
Kutoyants, Yu.\ A.\ 
(1984).\ 
Parameter estimation for stochastic processes.\ 
Prakasa Rao, B.L.S.\ (ed.\ )
Heldermann, Berlin.

\item
Kutoyants, Yu.\ A.\ 
(1994).\ 
{\it Identification of dynamical systems with small noise.\ }
Kluwer, Dordrecht.\

\item
Laredo, C.\ F.\ 
(1990).\ 
A sufficient condition for asymptotic sufficiency of incomplete observations 
of a diffusion process.\ 
{\it  Annals of Statistics, }{\bf 18}, 1158--1171.

\item
S{\o}rensen, M.\ and Uchida, M.\ 
(2003).\
Small diffusion asymptotics for discretely sampled 
stochastic differential equations.\
{\it Bernoulli }{\bf 9}, 1051--1069.

 \item
Uchida, M.\ and Yoshida, N.\ 
(2012).\ 
Adaptive estimation of an ergodic diffusion process
based on sampled data.\ 
{\it Stochastic Processes and their Applications}, {\bf 122}, 2885--2924.

\end{description}


\begin{en-text}

\section{Appendix} \label{appendix}
Figure 44 is {the} sample paths,
{where}  $\theta_1$, $\theta_2$ and $\epsilon$ are fixed and only $\theta_0$ is changed.
The rough shape of the sample path depends on the value of $\lambda_1$.
When $\lambda_1$ is positive, when y is fixed and t is moved from 0 to 1, the value of $X_t(y)$ approaches 0.
When $\lambda_1$ is close to 0, when y is fixed and t is moved from 0 to 1, the value of $X_t(y)$ does not change.
When $\lambda_1$ is negative, the value of $X_t(y)$ increases.
Figure 45 shows the sample path viewed from the t-axis side.
Figure 46 is a cross section of the sample path at y=0.5.

\begin{figure}[H]
\begin{minipage}{0.32\hsize}
\begin{center}
\includegraphics[width=4.5cm]{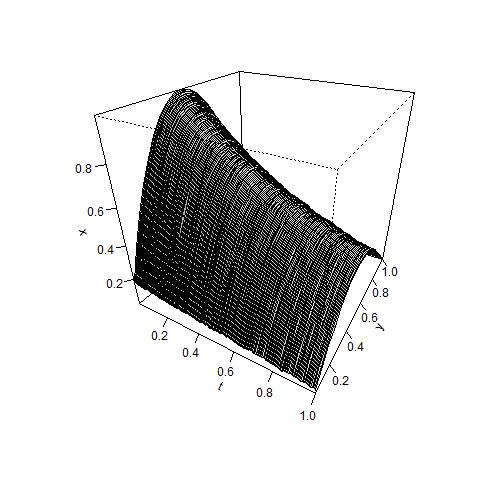}
\captionsetup{labelformat=empty,labelsep=none}
\subcaption{$\theta=$(1,0.2,0.2,0.01), $\lambda_1\approx$1}
\end{center}
\end{minipage}
\begin{minipage}{0.32\hsize}
\begin{center}
\includegraphics[width=4.5cm]{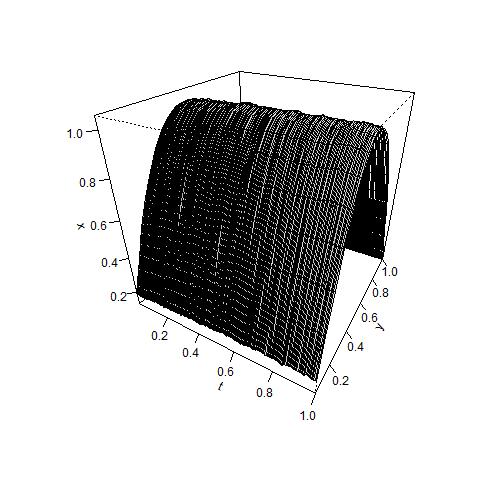}
\captionsetup{labelformat=empty,labelsep=none}
\subcaption{$\theta=$(2,0.2,0.2,0.01), $\lambda_1\approx$0}
\end{center}
\end{minipage}
\begin{minipage}{0.32\hsize}
\begin{center}
\includegraphics[width=4.5cm]{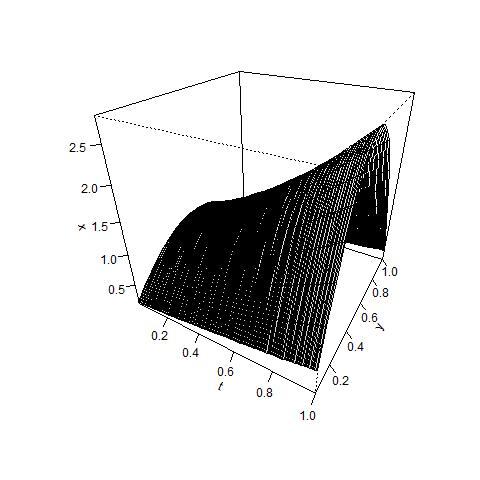}
\captionsetup{labelformat=empty,labelsep=none}
\subcaption{$\theta=$(3,0.2,0.2,0.01), $\lambda_1\approx$-1}
\end{center}
\end{minipage}
\caption{Sample paths with $\lambda_1=1, 0, -1$}
\label{t1-11}
\end{figure}

\begin{figure}[H]
\begin{minipage}{0.32\hsize}
\begin{center}
\includegraphics[width=4.5cm]{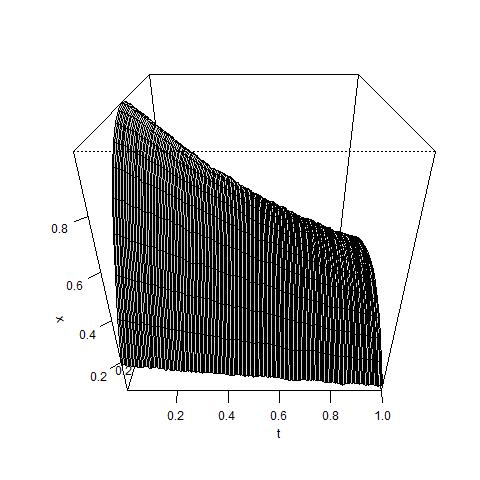}
\captionsetup{labelformat=empty,labelsep=none}
\subcaption{$\theta=$(1,0.2,0.2,0.01), $\lambda_1\approx$1}
\end{center}
\end{minipage}
\begin{minipage}{0.32\hsize}
\begin{center}
\includegraphics[width=4.5cm]{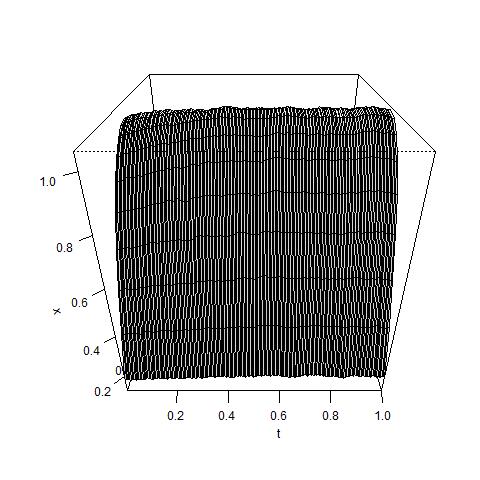}
\captionsetup{labelformat=empty,labelsep=none}
\subcaption{$\theta=$(2,0.2,0.2,0.01), $\lambda_1\approx$0}
\end{center}
\end{minipage}
\begin{minipage}{0.32\hsize}
\begin{center}
\includegraphics[width=4.5cm]{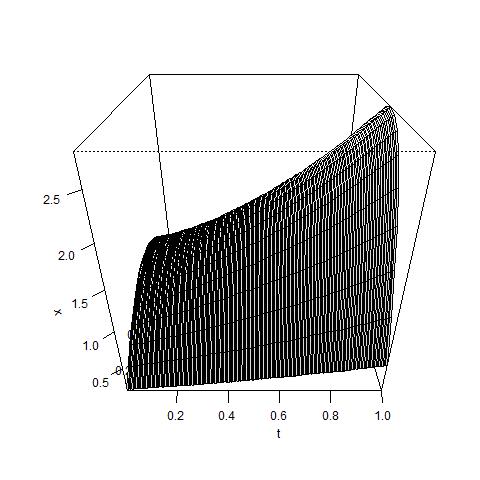}
\captionsetup{labelformat=empty,labelsep=none}
\subcaption{$\theta=$(3,0.2,0.2,0.01), $\lambda_1\approx$-1}
\end{center}
\end{minipage}
\caption{Sample paths with $\lambda_1=1, 0, -1$ (t-axis side)}
\end{figure}

\begin{figure}[H]
\begin{minipage}{0.32\hsize}
\begin{center}
\includegraphics[width=4.5cm]{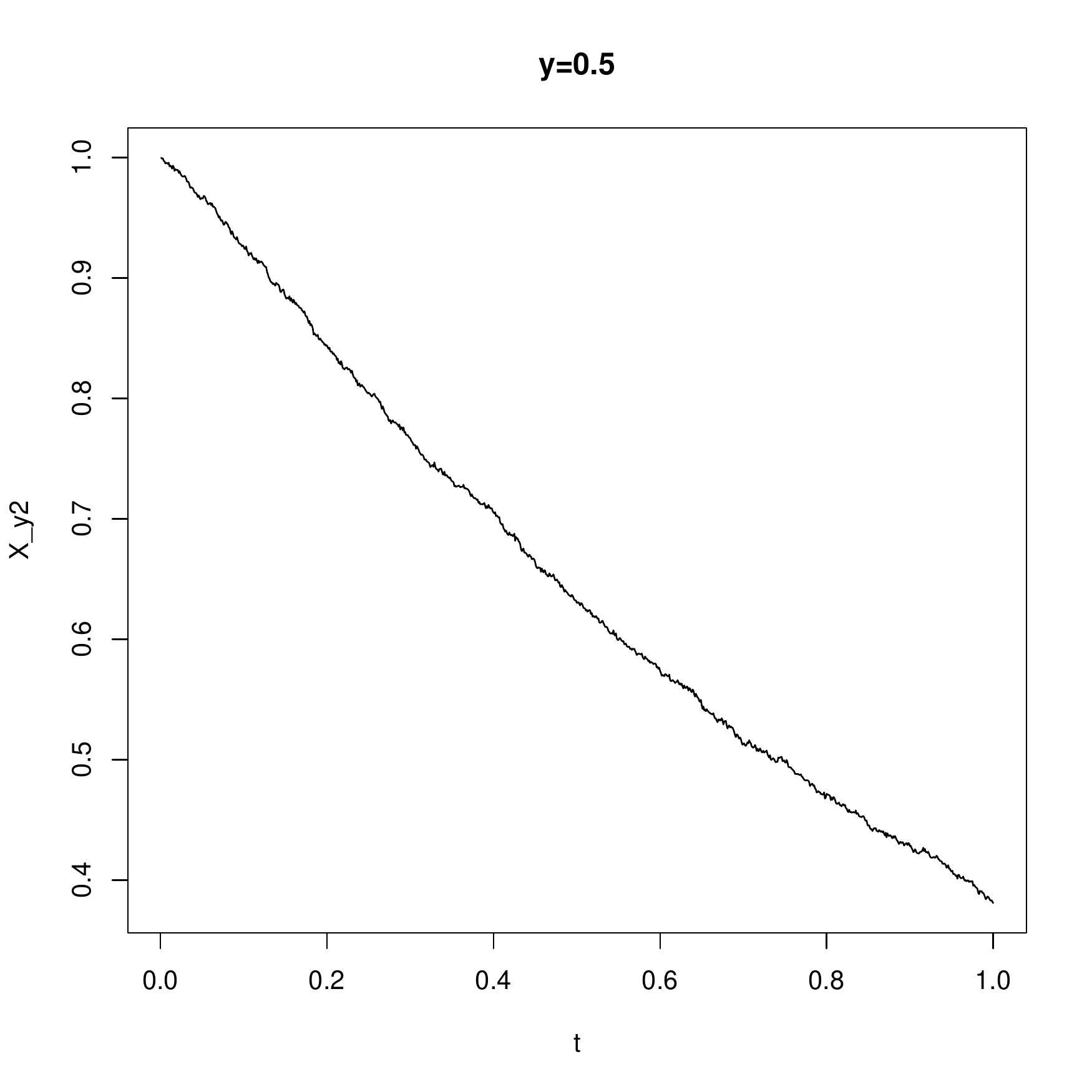}
\captionsetup{labelformat=empty,labelsep=none}
\subcaption{$\theta=$(1,0.2,0.2,0.01), $\lambda_1\approx$1}
\end{center}
\end{minipage}
\begin{minipage}{0.32\hsize}
\begin{center}
\includegraphics[width=4.5cm]{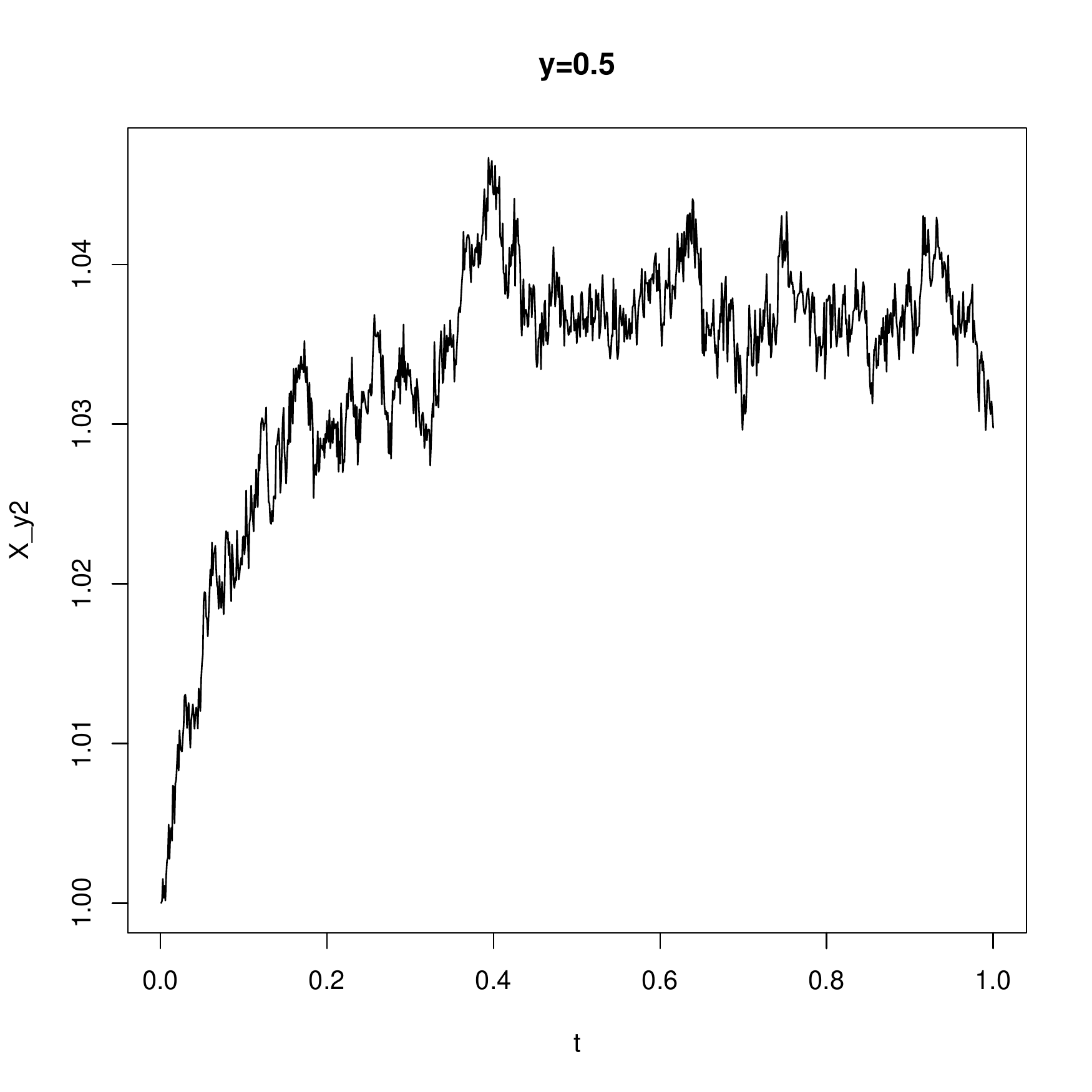}
\captionsetup{labelformat=empty,labelsep=none}
\subcaption{$\theta=$(2,0.2,0.2,0.01), $\lambda_1\approx$0}
\end{center}
\end{minipage}
\begin{minipage}{0.32\hsize}
\begin{center}
\includegraphics[width=4.5cm]{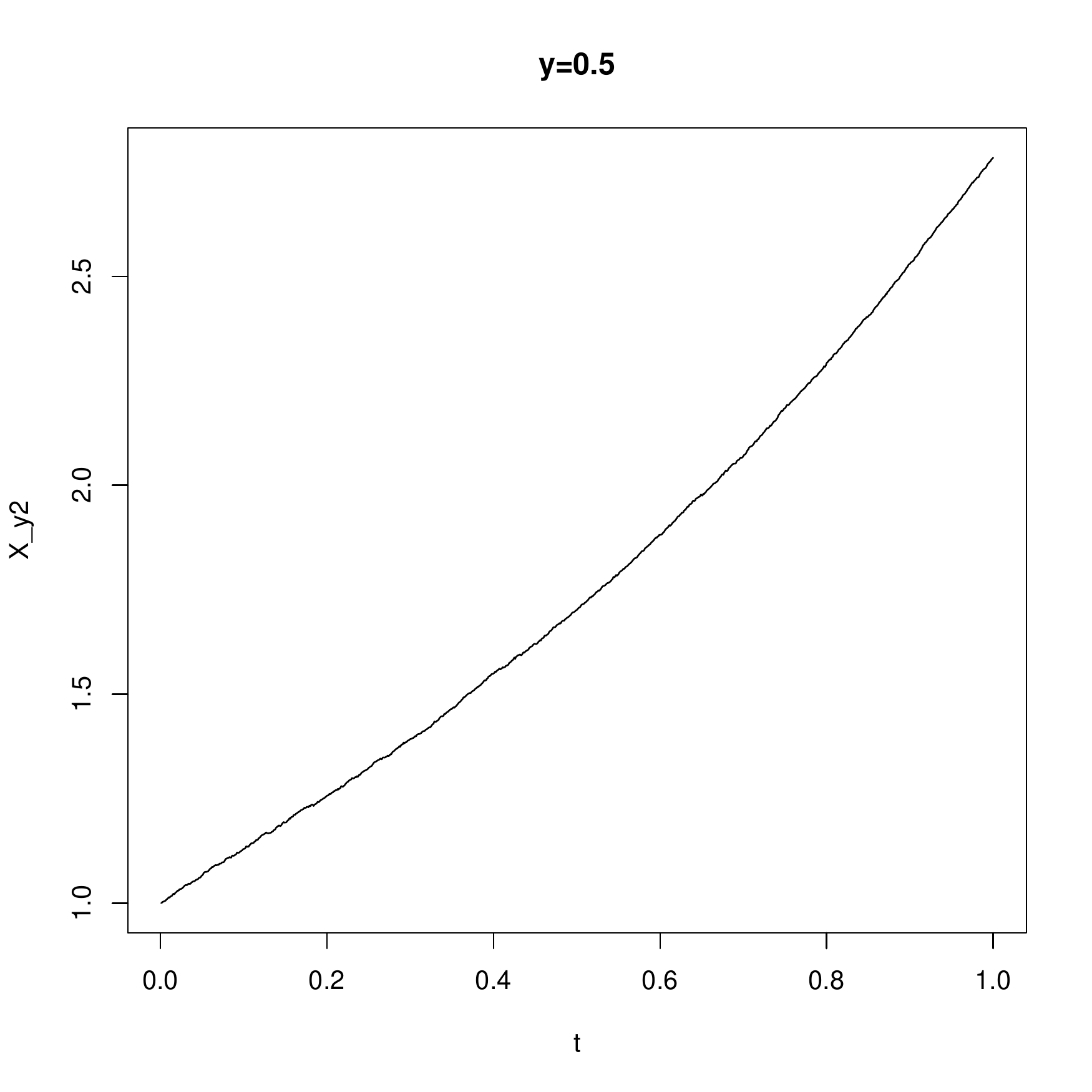}
\captionsetup{labelformat=empty,labelsep=none}
\subcaption{$\theta=$(3,0.2,0.2,0.01), $\lambda_1\approx$-1}
\end{center}
\end{minipage}
\caption{Sample paths with $\lambda_1=1, 0, -1$ (cross section at $y = 0.5$)}
\end{figure}

Figures 47-49 {the} sample paths,
{ where}  $\theta_0$, $\theta_1$ and $\theta_2$ are fixed and only $\epsilon$ is changed.
It can be seen that the noise increases as $\epsilon$ increases.

\begin{figure}[H]
\begin{minipage}{0.32\hsize}
\begin{center}
\includegraphics[width=4.5cm]{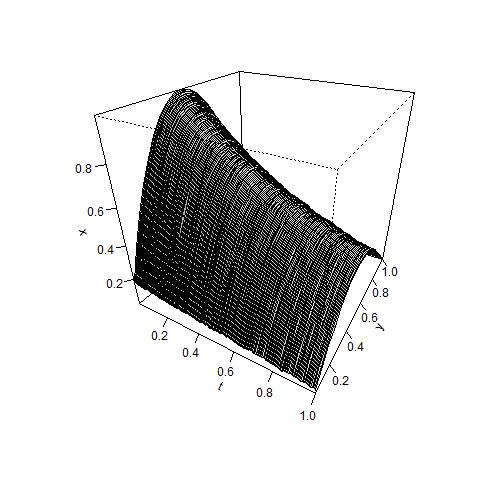}
\captionsetup{labelformat=empty,labelsep=none}
\subcaption{$\theta=$(1,0.2,0.2,0.01)}
\end{center}
\end{minipage}
\begin{minipage}{0.32\hsize}
\begin{center}
\includegraphics[width=4.5cm]{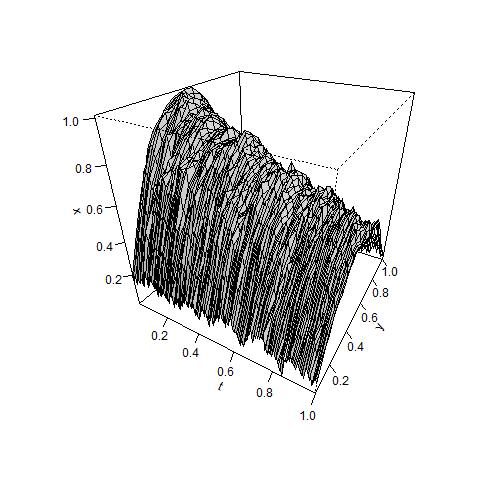}
\captionsetup{labelformat=empty,labelsep=none}
\subcaption{$\theta=$(1,0.2,0.2,0.1)}
\end{center}
\end{minipage}
\begin{minipage}{0.32\hsize}
\begin{center}
\includegraphics[width=4.5cm]{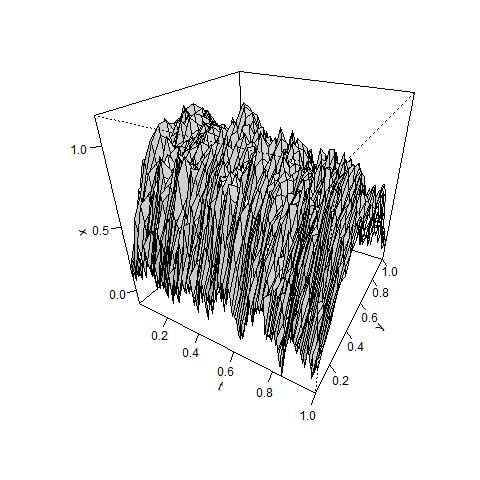}
\captionsetup{labelformat=empty,labelsep=none}
\subcaption{$\theta=$(1,0.2,0.2,0.25)}
\end{center}
\end{minipage}
\caption{Sample paths with $\epsilon=0.01, 0.1, 0.25$, $\lambda_1\approx$1}
\label{t1-11}
\end{figure}

\begin{figure}[H]
\begin{minipage}{0.32\hsize}
\begin{center}
\includegraphics[width=4.5cm]{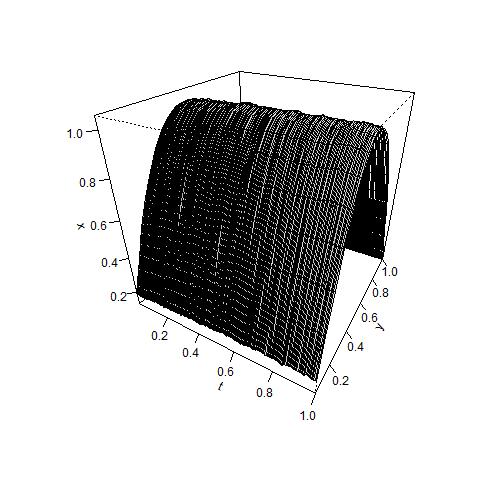}
\captionsetup{labelformat=empty,labelsep=none}
\subcaption{$\theta=$(2,0.2,0.2,0.01)}
\end{center}
\end{minipage}
\begin{minipage}{0.32\hsize}
\begin{center}
\includegraphics[width=4.5cm]{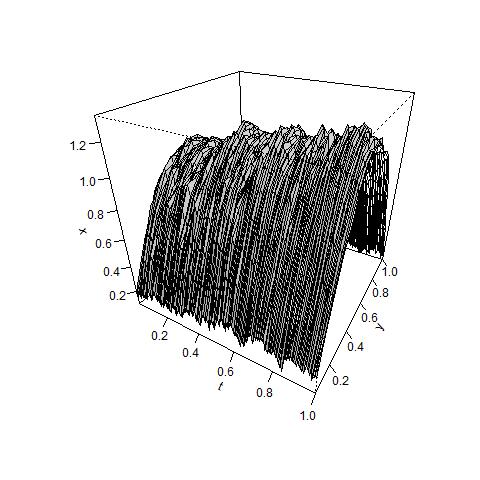}
\captionsetup{labelformat=empty,labelsep=none}
\subcaption{$\theta=$(2,0.2,0.2,0.1)}
\end{center}
\end{minipage}
\begin{minipage}{0.32\hsize}
\begin{center}
\includegraphics[width=4.5cm]{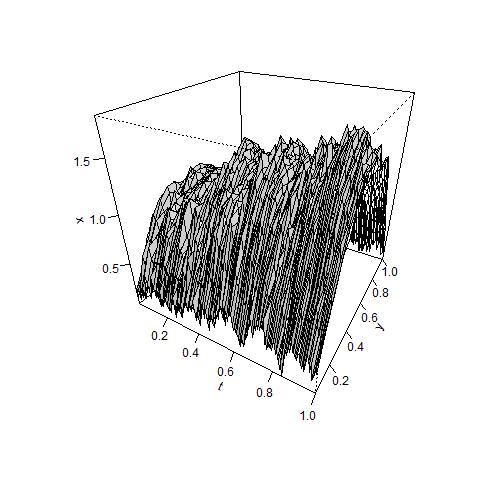}
\captionsetup{labelformat=empty,labelsep=none}
\subcaption{$\theta=$(2,0.2,0.2,0.25)}
\end{center}
\end{minipage}
\caption{Sample paths with $\epsilon=0.01, 0.1, 0.25$, $\lambda_1\approx$0}
\label{t1-11}
\end{figure}

\begin{figure}[H]
\begin{minipage}{0.32\hsize}
\begin{center}
\includegraphics[width=4.5cm]{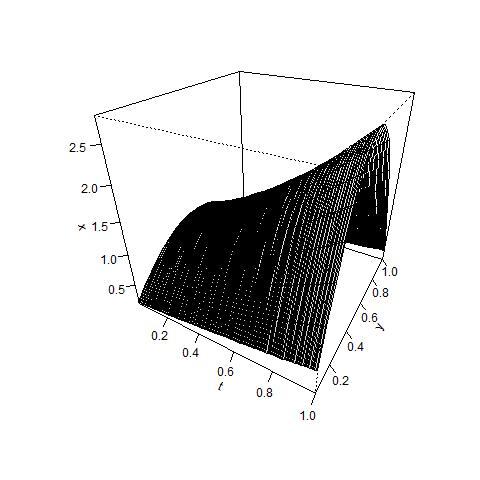}
\captionsetup{labelformat=empty,labelsep=none}
\subcaption{$\theta=$(3,0.2,0.2,0.01)}
\end{center}
\end{minipage}
\begin{minipage}{0.32\hsize}
\begin{center}
\includegraphics[width=4.5cm]{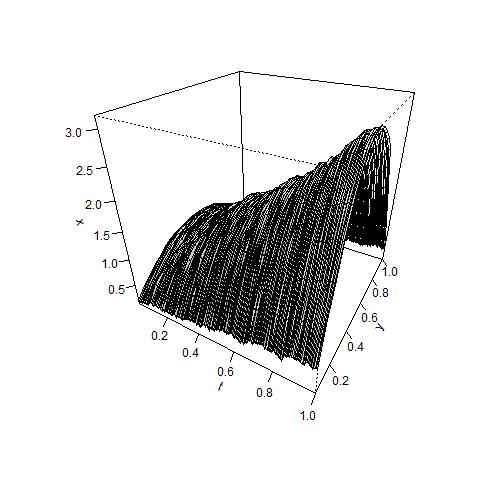}
\captionsetup{labelformat=empty,labelsep=none}
\subcaption{$\theta=$(3,0.2,0.2,0.1)}
\end{center}
\end{minipage}
\begin{minipage}{0.32\hsize}
\begin{center}
\includegraphics[width=4.5cm]{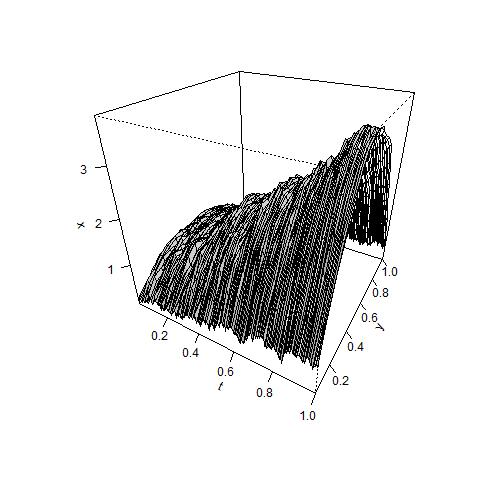}
\captionsetup{labelformat=empty,labelsep=none}
\subcaption{$\theta=$(3,0.2,0.2,0.25)}
\end{center}
\end{minipage}
\caption{Sample paths with $\epsilon=0.01, 0.1, 0.25$, $\lambda_1\approx$-1}
\label{t1-11}
\end{figure}

See Kaino and Uchida {(2020)} for the properties of $\theta_0$, $\theta_1$ and $\theta_2$.

\end{en-text}

\end{document}